\documentclass[11pt]{article}
\usepackage{amssymb}
\usepackage{mathrsfs}
\usepackage{amsthm}
\usepackage{CJK}
\usepackage{amsmath}
\usepackage{graphicx}
\usepackage{needspace}
\usepackage{longtable}
\usepackage{booktabs}
\usepackage{array}

\usepackage[hidelinks]{hyperref}
\hypersetup{pdftitle={Finite-Time Blowup for Navier--Stokes with Smooth Forcing: Part II. Residual Correction via Oscillatory Pulses},pdfauthor={Zhen Lei and Xiao Ren}}

\usepackage{enumerate}
\usepackage{tikz}

\usepackage{pgfplots}
\pgfplotsset{compat=1.18}

\usetikzlibrary{arrows.meta,decorations.pathreplacing}

\usepackage{xcolor}
\definecolor{RevisionPurple}{RGB}{128,0,128}
\definecolor{NoteGreen}{RGB}{0,128,128}

 \newtheorem{thm}{Theorem}[section]
 
 \newtheorem{lem}[thm]{Lemma}
 \newtheorem{prop}[thm]{Proposition}
 \theoremstyle{definition}

 \numberwithin{equation}{section}

\newtheorem{theorem}{Theorem}[section]

\theoremstyle{definition}

\newtheorem{example}[theorem]{Example}

\theoremstyle{remark}

\newcommand{\RR}{\mathbb R}
\newcommand{\ZZ}{\mathbb Z}
\newcommand{\TT}{\mathbb T}
\begin{document}
\begin{CJK*}{GBK}{song}
\title{{\Large Finite-Time Blowup for Navier--Stokes with Smooth Forcing}\\[2pt]
{\large Part II. Residual Correction via Oscillatory Pulses}}

\author{ Zhen Lei\footnote{Center for Applied Mathematics \& School of Mathematics, Fudan University, Shanghai 200433, P.
 R. China. Email: zlei@fudan.edu.cn} \and Xiao Ren\footnote{Center for Applied Mathematics, Fudan University, Shanghai 200433, P. R. China. Email: xren@fudan.edu.cn}}

\date{}

\maketitle

\begin{abstract}
This paper is a companion to Part~I:
\emph{Construction of Self-Similar Solutions with Admissible Stress and Flat Remainder}~\cite{PartI}
(\href{https://arxiv.org/abs/2609.35406}{arXiv:2609.35406}).
Part~I constructs approximately self-similar flow profiles with a singular stress in divergence form.
The associated velocity becomes unbounded as \(t\uparrow1\).
Here, localized oscillatory pulses cancel this stress contribution through their averaged quadratic momentum flux.
Successive corrections preserve the velocity growth and leave a smooth external force \(f\in C^\infty\).
This gives a local solution of the three-dimensional incompressible
Navier--Stokes equations with forcing \(f\).
Its velocity satisfies
\(\limsup_{t\uparrow1}\|u(t)\|_{L^\infty}=\infty\).
Together, the two papers present a readable and accessible account of OpenAI's
\emph{Finite Time Blowup for Navier--Stokes}~\cite{OpenAI}.
We claim no original mathematical contribution.

\end{abstract}

\tableofcontents

\section{Introduction}
\label{sec:introduction}

This paper continues the account of OpenAI's
\emph{Finite Time Blowup for Navier--Stokes}~\cite{OpenAI}
begun in Part~I~\cite{PartI}. Part~I first constructs a leading
self-similar profile and then constructs and sums lower-order corrections
to obtain a collapsing axisymmetric background with an admissible
stress. Here we explain
how oscillatory pulses and successive corrections remove its singular
momentum residual while preserving the velocity growth. Our aim is a
readable and accessible presentation of this part of the construction;
we claim no original mathematical contribution. This exposition is self-contained: the reader can start with Part II directly, assuming the required properties of the background solution constructed in Part I.

\subsection{The background input from Part I}
\label{intro:flux-cone}

We first state the properties of the complete background from Part~I
that will be used below. We use cylindrical coordinates
\((t,r,\theta,z)\), with velocity components in the orthonormal frame
\((e_r,e_\theta,e_z)\). The angular and axial directions are called
\emph{tangential}, and
\(\nabla_{\theta,z}=(r^{-1}\partial_\theta,\partial_z)\).
Fix the parameter \(0<\delta<1/200\) supplied by Part~I and use
\begin{equation}\label{intro:background-variables}
 R=\frac{r^2}{2\lambda^2},\qquad
 Z=\frac{z}{\lambda^{1-\delta}},\qquad
 1-t=\lambda^2(1-Z^2).
\end{equation}
Here \(R\) and \(Z\) are the similarity variables, and \(\lambda=\lambda(t,z)\) is the scale parameter.
We use the following background properties and estimates for \(0<\lambda\le\lambda_*\), with \(\lambda_*>0\) sufficiently small,
\(R\ge0\), and \(|Z|\le1\). Derivatives in \(Z\) at \(Z=\pm1\) are understood as one-sided derivatives.

The background fields \((u_B,p_B,\mathsf T_B,\mathsf E_B)\) are real,
axisymmetric and smooth away from the point $|x|=0, t = 1$. For
\[
\mathsf{NS}(u,p)=\partial_tu+(u\cdot\nabla)u-\Delta u+\nabla p,
\]
they satisfy
\begin{equation}\label{intro:background}
	\nabla\cdot u_B=0,\qquad
	\mathsf{NS}(u_B,p_B)=-\mathcal D\mathsf T_B+\mathsf E_B,
\end{equation}
where \(\mathsf T_B=(\mathsf T_B^\theta,\mathsf T_B^z)\) and
\[
\mathcal D\mathsf T
=\left(\partial_r+\frac2r\right)\mathsf T^\theta e_\theta
+\left(\partial_r+\frac1r\right)\mathsf T^z e_z.
\]
Thus the radial component of
\( \mathsf{NS}(u_B,p_B)-\mathsf E_B\) is exactly zero.

\begin{enumerate}[(i)]
\item \textit{Flat error.}
The error $\mathsf E_B$ is flat in $\lambda$ in the sense that for every  multi-index
\(\nu\), and every \(N>0\),
\begin{equation}\label{intro:flat-background-remainder}
 |\partial_{t,x}^{\nu}\mathsf E_B|
 \le C_{\nu,N}\lambda^N.
\end{equation}

\item \textit{Bounds for the rescaled velocity.}
The rescaled components \(\lambda^{1+\delta}u_B^\theta,\lambda^{1+\delta}u_B^z,\lambda u_B^r\),
together with all their derivatives in \((\log\lambda,\sqrt{2R},Z)\), are uniformly bounded on the domain
\begin{equation}\label{intro:background-velocity-bounds}
 0<\lambda\le\lambda_*,\qquad R\ge0,\qquad |Z|\le1.
\end{equation}

\item \textit{The active annulus and weighted stress bounds.}
The stress $\mathsf T_B$ is supported in \([R_a,R_h]\) with $0<R_a<R_h<+\infty$. On the open
\emph{active annulus} \(R_a<R<R_h\), define
\begin{equation}\label{intro:background-edge-weight}
 \begin{aligned}
 d_{\mathrm{edge}}(R)
   &=\min\{1,\log(R/R_a),\log(R_h/R)\},\\
 \omega(R)
   &=\exp\left(-\frac{c_{\mathrm{collar}, a}}{\log^2(R/R_a)}
               -\frac{c_{\mathrm{collar}, h}}{\log^2(R_h/R)}\right),
 \qquad c_{\mathrm{collar}, a},c_{\mathrm{collar}, h}>0.
 \end{aligned}
\end{equation}
The constants \(c_{\mathrm{collar}, a},c_{\mathrm{collar}, h}\) are fixed by the background. Extend
\(\omega\) by zero outside the annulus. There exists a fixed \(c>0\)
and positive  constants \(C_M>0\)
such that
\begin{equation}\label{intro:weighted-background-stress}
 \begin{aligned}
 \lambda^{2+\delta}|\mathsf T_B|&\ge c\omega,\\
 \left|\left\{\lambda \partial_\lambda, \sqrt{2R} \partial_R, \partial_Z\right\}^{\le M}(\lambda^{2+\delta}\mathsf T_B)\right|
   &\lesssim_M \omega d_{\mathrm{edge}}^{-C_M}.
 \end{aligned}
\end{equation}
These estimates hold uniformly for \(|Z|\le1\) and imply smooth
zero extension of $\mathsf T_B$ at both radial edges. 

\item \textit{A uniform strict cone condition.}
Define the normalized angular velocity and shear of the complete
background by
\begin{equation}\label{intro:full-background-shear}
 \begin{gathered}
 F_B=\frac{\lambda^{1+\delta}u_B^\theta}{\sqrt{2R}},\qquad
 \mathcal S_B=
 \left(2R\partial_RF_B,
       \sqrt{2R}\,\partial_R(\lambda^{1+\delta}u_B^z)\right),\\
 \mathcal S_B^\perp=(-\mathcal S_B^z,\mathcal S_B^\theta),\qquad
 \kappa_B=-\frac{|\mathcal S_B|^2}{F_B\mathcal S_B^\theta}.
 \end{gathered}
\end{equation}
On the closed stress annulus, both \(F_B\) and \(-\mathcal S_B^\theta\)
lie in \([C^{-1},C]\) for some fixed \(C\ge1\). For some fixed \(c_{\mathrm{cone}}>0\),
the following inequalities hold throughout the open active annulus:
\begin{equation}\label{intro:admissible-cone}
 \begin{gathered}
 \kappa_B\ge2+c_{\mathrm{cone}},\qquad
 -\mathsf T_B\cdot\frac{\mathcal S_B}{|\mathcal S_B|}
       \ge c_{\mathrm{cone}} |\mathsf T_B|,\\
 \sqrt{\frac{\kappa_B-2}{2}}
 \left|\frac{\mathsf T_B\cdot\mathcal S_B^\perp}
                 {\mathsf T_B\cdot\mathcal S_B}\right|
       \le1-c_{\mathrm{cone}}.
 \end{gathered}
\end{equation}

\item \textit{A reserved interval for the five moment corrections.}
There is a fixed nonempty open interval
\(I_3\Subset(R_a,R_h)\) such that the complete background satisfies
the exact identities
\begin{equation}\label{intro:background-moment-patch}
 \begin{gathered}
 u_B^z=0,\qquad
 u_B^\theta=\lambda^{-1-\delta}A_*(Z)(2R)^{-1/2-\mu}
       \quad (R\in I_3),\\
 \mu>0,\qquad A_*\in C^\infty([-1,1]),\qquad
 \inf_{[-1,1]}A_*>0.
 \end{gathered}
\end{equation}
In the construction of Part~I, this factor is
\begin{equation}\label{intro:part-i-amplitude}
 A_*(Z)=\frac{a_*}{1+Z^2},\qquad
 a_*>0,
\end{equation}
where \(a_*\) is a fixed constant determined by the
 construction procedure there; see \cite[Interval~O.3]{PartI}.
In Part II, we use only the smoothness and positive lower bound for $A_*$ in
\eqref{intro:background-moment-patch}, not the particular formula \eqref{intro:part-i-amplitude}.

\item \textit{The exact heat exterior.}
For \(R\ge R_h\), the background agrees with the heat swirl
\begin{equation}\label{intro:background-heat}
 \begin{aligned}
 u_B^r=u_B^z&=0,\qquad u_B^\theta=K(t,r),\\
 p_B(t,r)&=-\int_r^\infty K(t,s)^2\,\frac{ds}{s},\\
 K(t,r)&=c_\infty\left(\frac{r^2}{2}\right)^{-(1+\delta)/2}
     H_\delta\left(\frac{4(1-t)}{r^2}\right),\qquad c_\infty>0,
 \end{aligned}
\end{equation}
where
\[
 H_\delta(\xi)
 =\frac1{\Gamma(1+\delta/2)}
  \int_0^\infty e^{-v}v^{\delta/2}
                  (1+\xi v)^{-\delta/2}\,dv>0,
 \qquad \xi\ge0.
\]
This velocity and pressure solve the Navier--Stokes equations
exactly, so \(\mathsf E_B=0\) in the exterior.
\end{enumerate}

Part~I uses global moment conditions to obtain the annular support
of its stress. We do not impose those profile identities as
additional background hypotheses here. The exact equation
\eqref{intro:background} and the smooth annular support of \(\mathsf T_B\)
already give the background force identities
needed below: at each fixed \((t,z)\),
\begin{equation}\label{intro:background-force-moments}
 \begin{aligned}
 \int_0^\infty r^2
   \bigl(\mathsf{NS}(u_B,p_B)-\mathsf E_B\bigr)^\theta\,dr
   &=-[r^2\mathsf T_B^\theta]_0^\infty=0,\\
 \int_0^\infty r
   \bigl(\mathsf{NS}(u_B,p_B)-\mathsf E_B\bigr)^z\,dr
   &=-[r\mathsf T_B^z]_0^\infty=0.
 \end{aligned}
\end{equation}
The later moment corrections in Sections~4--5 will control defects relative to this fixed background.

\subsection{Main result of Part II}
\label{intro:main-result}

The background conditions above suffice for the following construction.

\begin{thm}\label{thm:part-ii-main}
Suppose that \((u_B,p_B,\mathsf T_B,\mathsf E_B)\) satisfy the
conditions in Subsection~\ref{intro:flux-cone}. There exist a compact set
\(K\subset\mathbb R^3\), $t_0 \in (0,1)$, and real fields
\[
 u\in C^\infty(\mathbb R^3\times[0,1);\mathbb R^3),\qquad
 p\in C^\infty(\mathbb R^3\times[0,1)),\qquad
 f\in C_c^\infty(\mathbb R^3\times(0,\infty);\mathbb R^3)
\]
with the following properties.
\begin{enumerate}[(i)]
\item \textit{The equation and compact supports.}
On \(\mathbb R^3\times[0,1)\),
\begin{equation}\label{intro:main-equations}
 \partial_tu+(u\cdot\nabla)u-\Delta u+\nabla p=f,
 \qquad \nabla\cdot u=0.
\end{equation}
Moreover,
\[
 (u,p)=0\quad\text{for }0\le t\le t_0,\qquad
 \operatorname{supp}u(\cdot,t)\cup\operatorname{supp}p(\cdot,t)
 \subset K\quad(0\le t<1).
\]

\item \textit{The correction and remaining flat force.}
For some sufficiently small \(\varepsilon>0\), set
\[
 Q_\varepsilon=B_\varepsilon(0)\times(1-\varepsilon^2,1).
\]
On \(Q_\varepsilon\) we have
\begin{equation}\label{intro:main-correction-support}
 (u,p)=(u_B,p_B)\quad\text{whenever }R\notin[R_a,R_h].
\end{equation}
On $Q_\varepsilon$, the force $f$ is flat in $\lambda$ in the sense of \eqref{intro:flat-background-remainder}.

\item \textit{Finite energy and blowup.}
The velocity satisfies
\begin{equation}
 \sup_{0\le t<1}\|u(t)\|_{L^2(\mathbb R^3)} + \|\nabla u\|_{L^2(\mathbb R^3 \times [0,1])}<\infty
 \end{equation}
 and
 \begin{equation}\label{intro:main-blowup}
 \|u(t)\|_{L^\infty(\mathbb R^3)}
 \ge c(1-t)^{-(1+\delta)/2}
\end{equation}
for some \(c>0\) and all \(t<1\) sufficiently close to \(1\).
\end{enumerate}
\end{thm}

\subsection{Strategy and the pulse mechanism}
\label{intro:linear-mechanism}

The main building blocks are localized, oscillatory,
divergence-free velocities called \emph{pulses}. Their averaged
quadratic momentum flux supplies the stress missing from the
background equation. There are two types. A \emph{primary pulse}
solves an unforced principal amplitude equation and provides the
initial positive stress matching. A \emph{forced pulse} solves the
corresponding equation with a prescribed linear source and zero
initial data; it cancels nonzero angular harmonics left by earlier
corrections. The primary pulse construction determines the admissible
cone condition used in Part~I; see Sections~\ref{s7n:frequencies}
and~\ref{s6n:unforced-pulse}. Further pressure, stress, auxiliary and moment corrections are
needed to control all parts of the residual and to repeat the
improvement; see Sections~\ref{sec:initialization-ppt} and~\ref{sec:iteration-ppt}.

\paragraph{The principal linear dynamics.}
Adding a velocity \(w\) and pressure \(p_w\) to the background gives the exact identity
\begin{equation}\label{intro:linear-quadratic}
 \mathsf{NS}(u_B+w,p_B+p_w)
 =\mathsf{NS}(u_B,p_B)+\mathsf L_B(w,p_w)+(w\cdot\nabla)w,
\end{equation}
where
\[
 \mathsf L_B(w,p_w)
 =\partial_tw+(u_B\cdot\nabla)w+(w\cdot\nabla)u_B
                  -\Delta w+\nabla p_w
\]
is the linearization around the background \(u_B\).
We choose the primary pulse so that the leading terms in \(\mathsf L_B(w,p_w)\) cancel.
We use its averaged quadratic flux to match the background stress.

The precise framework of pulse construction is presented in Section~\ref{sec:pulse-framework}; here we only give a rough description. In a band $\lambda \asymp 2^{-k}$, we use pulses with spatial frequency $\sim 2^{(1+\frac{\delta}{2})k}$. The leading part of a pulse is given by the pre-pulse ansatz:
\begin{equation}\label{intro:pre-pulse}
 \begin{aligned}
 w^{\mathrm{pre}}
   &=2^{(1+\delta)k}\chi\mathbf h\cos(N\varphi),\\
 p_w&=2^{2(1+\delta)k}\chi
                     \operatorname{Re}(\pi e^{\mathrm iN\varphi}),
 \qquad N=\lceil2^{\delta k/2}\rceil.
 \end{aligned}
\end{equation}
Here \(\chi\) localizes the pulse, and \(\varphi\) is a phase function with $|\nabla \varphi| \sim 2^k$. For $w^{\mathrm{pre}}$ to be nearly divergence free, we require that $\mathbf h \cdot \nabla \varphi = 0$. Using a vector potential, it is easy to complete \(w^{\mathrm{pre}}\) to an
exactly divergence-free pulse \(w=w^{\mathrm{pre}}+w^{\mathrm{corr}}\).  The velocity amplitude $\mathbf h$ of a primary pre-pulse will be of order $2^{-\frac{\delta k}{2}}$ up to powers of $k$, so that $w^{\mathrm{pre}} \otimes w^{\mathrm{pre}}$ has the order $2^{(2+\delta)k}$ of the background stress. 

Take a local fast time $\tau$ with $\frac{d\tau}{dt} = 2^{(2+\delta)k}$. We choose the phase \(\varphi\) to be transported by the background \(u_B\) to leading order.
With the slow variables held fixed, we take the phase \(\varphi\) to be affine in the fast time \(\tau\).
With this phase choice, extracting the leading terms of the linearization $\mathsf L_B(w^{\mathrm{pre}}, p_w)$ gives the so-called \emph{principal amplitude equations}:
\begin{equation}\label{intro:amplitude-equation}
 \partial_\tau\mathbf h+\mathsf K_B\mathbf h+d\mathbf h
       +\mathrm iN\mathbf n_\varphi\pi=0,
 \qquad \mathbf n_\varphi\cdot\mathbf h=0,
\end{equation}
where
\[
 \begin{gathered}
 \mathbf n_\varphi=2^{-k}\nabla\varphi,\qquad
 d=2^{-\delta k}N^2|\mathbf n_\varphi|^2,\\
 \mathsf K_B=2^{-(2+\delta)k}
 \begin{pmatrix}
 0&-2u_B^\theta/r&0\\
 \partial_ru_B^\theta+u_B^\theta/r&0&0\\
 \partial_ru_B^z&0&0
 \end{pmatrix}.
 \end{gathered}
\]
Pressure can be solved from $\mathbf h$ using the orthogonality constraint $\mathbf n_\varphi\cdot\mathbf h=0$ as
\[
 \pi=\frac{\mathrm i}{N}
   \frac{\mathbf n_\varphi\cdot\mathsf K_B\mathbf h
          -(\partial_\tau\mathbf n_\varphi)\cdot\mathbf h}
        {|\mathbf n_\varphi|^2}.
\]

The competition between the shear term $\mathsf K_B \mathbf h$ and the viscosity term $d \mathbf h$ in the principal dynamics is crucial. The phase $\varphi$
and the initial polarization can be chosen such that the shear first
amplifies the pulse and then the viscous damping dominates. Consequently, the pulse amplitudes are controlled by 
a Gaussian envelope concentrated near the midpoint of a fast-time
interval of length \(\tau_{\rm max}\asymp k^2\). This delicate growth-and-decay mechanism allows us to construct pulses that are well localized in space-time.

\paragraph{The primary pulses determine a cone.}
The relevant quadratic quantity is the averaged radial transport
of tangential momentum:
\begin{equation}\label{intro:flux-map}
 \mathsf C(w)=
 \begin{pmatrix}
 \langle w^rw^\theta\rangle_{\theta,Y}\\
 \langle w^rw^z\rangle_{\theta,Y}
 \end{pmatrix}\in\RR^2.
\end{equation}
Here \(Y\in\TT^2\) is the independent auxiliary variable introduced
in Section~\ref{s6n:auxiliary}, and $\langle \rangle_{\theta, Y}$ stands for the average in both $\theta$ and $Y$ variables. The auxiliary variable also helps to separate the supports
of distinct local pulses, so that their cross products vanish. 

The directions relevant to stress matching are those of the
two-component flux \(\mathsf C(w)\). The unforced principal equation
limits the flux directions available to a pulse. Locally, we choose
two families whose fluxes point in linearly independent directions
on either side of the negative background shear. Multiplying a pulse
by a real amplitude \(a\) multiplies its flux by \(a^2\). Since the
supports of distinct pulses are separated, their fluxes add. The
stresses obtainable this way therefore form the cone of nonnegative
combinations of the two directions.

The cone condition \eqref{intro:admissible-cone} allows us to choose
these families so that the nonzero background stress lies strictly
inside this local cone. The uniform margin accommodates the slow
variation of the background and the small errors in the approximation
of the flux directions. Section~\ref{sec:matching-background} chooses the positive
amplitudes using the actual averaged fluxes. After summing over all
local pulses, the assembled pre-pulse \(w_0^{\mathrm{pre}}\) satisfies
\[
 \mathsf C(w_0^{\mathrm{pre}})=\mathsf T_B.
\]
Its averaged radial stress divergence cancels
\(-\mathcal D\mathsf T_B\) in \eqref{intro:background}.
Here the matching uses the complete background stress \(\mathsf T_B\).
The construction in \cite[Section~7.3, Proposition~7.5 and equation~(7.30)]{OpenAI}
uses the leading stress \(\mathsf T_{(0)}=\lambda^{-2-\delta}\mathcal T_{(0)}(R,Z)\) at this step.
\label{intro:matching-choice}

\paragraph{Initialization.}
This matching controls only the averaged radial flux. Oscillatory
terms, radial forces and errors from localization and divergence-free
completion still remain in the equation. The initialization in
Section~\ref{sec:initialization-ppt} first adds the primary pulses,
then an auxiliary velocity correction to reduce the part of the
tangential angular mean that still oscillates in \(Y\). Compact radial
support imposes integral conditions on pressure and stress primitives.
A moment correction reduces the defects in these conditions, which
control the remainders in pressure reconstruction and in later stress
corrections. Finally, we reconstruct the axisymmetric pressure to
control the radial angular mean. These steps establish the fixed
primary pulses and the bounds on the remaining force and integral
defects needed to begin the iteration.

\paragraph{Subsequent correction cycles.}
Later stress increments may point outside the initial cone. We keep
the primary pulses fixed and add small pulses with the same phases
and localizations as the corresponding primary pulses. Write \(v\)
for their assembled pre-pulse. The quadratic identity
\[
 \mathsf C(w_0^{\mathrm{pre}}+v)
 =\mathsf C(w_0^{\mathrm{pre}})
  +\mathsf B(w_0^{\mathrm{pre}},v)+\mathsf C(v)
\]
defines the cross flux \(\mathsf B(w_0^{\mathrm{pre}},v)\), which is
linear in \(v\). The primary scalar matching coefficients are positive, while the new
amplitudes may have either sign. Since the two flux directions are
linearly independent, their cross fluxes can match stress increments
of arbitrary direction, as shown in Section~\ref{subsec:signed-covariance}.
The new pulse's own quadratic flux \(\mathsf C(v)\) is left as a
smaller residual.

Each cycle first uses forced pulses to reduce the nonzero angular
modes of the force. It then takes a radial primitive of the tangential
double mean as the target for the cross-flux correction. The integral
defects control the compensating terms needed to keep this primitive
compactly supported. The cycle ends with the auxiliary and moment corrections.
At each stage, including initialization, we reconstruct the axisymmetric pressure
only once, after all velocity corrections are complete, without performing the
intermediate reconstructions used in \cite[Propositions~9.5 and~9.6]{OpenAI}. Each correction introduces new errors.
Over the full cycle, however, the bounds on both the force and the integral defects improve.
Figure~\ref{fig:intro-correction-cycle} shows how initialization leads into the iteration.

\begin{center}
\begin{minipage}{\linewidth}
\makeatletter\def\@captype{figure}\makeatother
\centering
\definecolor{ink}{HTML}{25394A}
\definecolor{teal}{HTML}{1B747C}
\definecolor{edge}{HTML}{9DAFB4}
\definecolor{mist}{HTML}{F3F8F7}
\tikzset{
  stage/.style={draw=edge,rounded corners=2pt,line width=.6pt,
    align=center,inner sep=0pt,text=ink},
  flow/.style={draw=ink,line width=.7pt,-{Latex[length=1.8mm]}},
  merge path/.style={draw=ink,line width=.7pt,rounded corners=2pt},
  cycle path/.style={flow,draw=teal,rounded corners=3pt},
  heading/.style={font=\bfseries,anchor=west,inner sep=0pt,text=ink},
  down/.style={text=teal,inner sep=0pt,font=\small}
}
\begin{tikzpicture}[font=\fontencoding{T1}\fontfamily{ntxtlf}\fontsize{10}{12}\selectfont]
  \node[heading] at (0,3.85) {Initialization};
  \node[stage,minimum width=8.6cm,minimum height=.84cm,anchor=west]
    (primary) at (0,3.12) {Primary pulses match the full stress \(\mathsf T_B\).};
  \node[heading] at (0,2.02) {Each iteration cycle};
  \node[stage,fill=mist,minimum width=3.62cm,minimum height=1.40cm,anchor=west]
    (forced) at (0,.98) {Forced pulses reduce\\nonzero angular modes.};
  \node[stage,fill=mist,minimum width=4.41cm,minimum height=1.40cm,anchor=west]
    (cross) at (4.19,.98) {Cross flux with the fixed\\primary pulses matches\\the stress increment.};
  \draw[flow] (forced.east) -- (cross.west);
  \begin{scope}[xshift=1.75cm]
  \node[heading] at (8.36,3.85) {Shared correction steps};
  \node[stage,minimum width=4.7cm,minimum height=2.52cm,anchor=west]
    (shared) at (8.36,2.28) {};
  \node[text=ink] at (10.71,3.11) {Auxiliary velocity correction};
  \node[down] at (10.71,2.70) {\(\downarrow\)};
  \node[text=ink] at (10.71,2.28) {Moment correction};
  \node[down] at (10.71,1.86) {\(\downarrow\)};
  \node[text=ink] at (10.71,1.45) {Pressure reconstruction};
  \end{scope}
  \draw[merge path] (cross.east) -- (9.31,.98) -- (9.31,3.12);
  \draw[flow] (primary.east) -- (10.11,3.12);
  \fill[ink] (9.31,3.12) circle (1pt);
  \draw[cycle path] (shared.south) -- (12.46,-.53) -- (1.81,-.53) -- (forced.south);
  \node[fill=white,text=teal,font=\small,inner xsep=4pt] at (7.135,-.53)
    {Improved force and defect bounds for the next cycle};
\end{tikzpicture}
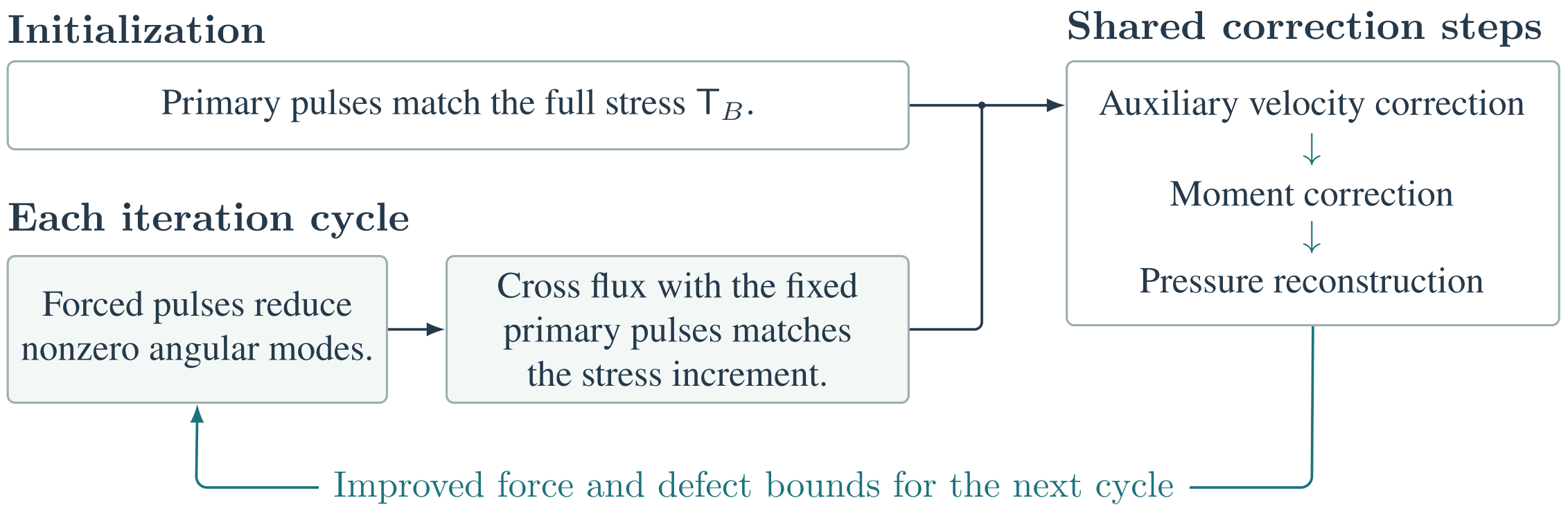
\caption{Initialization and one iteration cycle. The primary pulses chosen
during initialization remain fixed throughout the iteration.}
\label{fig:intro-correction-cycle}
\end{minipage}
\end{center}

Repeating the cycle gives arbitrarily high orders of decay as
\(\lambda\to0\). Section~\ref{sec:iteration-ppt} realizes these
finite-stage constructions by summing corrections with cutoffs
increasingly concentrated near \(\lambda=0\), obtaining a residual
flat to every order. All velocity corrections remain in the annulus,
preserving the background blowup in the heat exterior.

\subsection{Related constructions}
\label{intro:related-work}

Lifschitz and Hameiri~\cite{LH91} and Friedlander and
Vishik~\cite{FV91} study
short-wave perturbations through transported wavevectors and
velocity polarizations. For the columnar model obtained by fixing
\(\lambda\) and \(Z\), the condition \(\kappa_B>2\)
is equivalent to the sufficient instability
condition of Leibovich and Stewartson~\cite{LeibovichStewartson}.

Craik and Criminale~\cite{CraikCriminale} construct finite-amplitude
waves on affine background flows, with cancellation of the wave's
quadratic self-interaction. For a linear shear flow, Singh and
Sridhar~\cite{SinghSridhar2017} construct exact Navier--Stokes solutions
by superposing Kelvin modes whose wavevectors remain parallel as they evolve.
The pulses used here combine
this phase and amplitude viewpoint with localization and error
control for a spatially varying background.

The use of oscillations to realize a prescribed averaged stress is
central to convex integration for fluid equations. De~Lellis and
Sz\'ekelyhidi~\cite{DLS09,DLS13} develop the differential-inclusion
formulation and iterative stress-correction constructions for Euler.
Daneri and Sz\'ekelyhidi~\cite{DS2017} construct Mikado flows with
separated periodic supports whose averaged quadratic flux realizes
a prescribed positive-definite stress. For Navier--Stokes, Buckmaster
and Vicol~\cite{BV19} establish nonuniqueness of finite-energy weak
solutions using intermittent Beltrami waves.
Buckmaster, Colombo and Vicol~\cite{BCV2022} use intermittent jets
with disjoint supports, together with incompressibility and temporal correctors.
Cheskidov and Luo~\cite{CL2022} combine stationary Mikado flows with
prescribed temporal concentration and a temporal corrector.
In the present construction, the
oscillations are constrained by the linear dynamics of the fixed
background; the admissible stress cone records that constraint.
Temporal concentration here results from shear amplification
followed by viscous damping.

Tao~\cite{Tao2016} constructs finite-time blowup for an averaged
Navier--Stokes equation retaining the energy cancellation of the
nonlinearity.

Building on the forced three-dimensional Euler construction of
C\'ordoba and Mart\'inez-Zoroa~\cite{CM2023}, C\'ordoba,
Mart\'inez-Zoroa and Zheng~\cite{CMZ} obtain finite-time blowup
for hypodissipative Navier--Stokes with small dissipation orders
and forcing in \(L^1_tC_x^{1,\epsilon}\cap L^\infty_tL^2_x\).
For the incompressible porous media equation, C\'ordoba and
Mart\'inez-Zoroa~\cite{CMIPM} construct blowup with a source whose
spatial derivatives of every order remain uniformly bounded, and
Alp\"oge, Buckmaster and Coiculescu~\cite{ABCoIPM} adapt this
construction to \(\mathbb T^2\) with a source smooth in both space
and time up to the blowup time.
Alp\"oge and Buckmaster~\cite{ABEuler,ABBoussinesq} construct
finite-time blowup for Euler and Boussinesq with forcing smooth
in space and time through the blowup time.
These constructions amplify successively finer oscillatory layers
using the previously constructed flow. In OpenAI's
construction~\cite{OpenAI}, amplified pulses supply the collapsing
background's missing momentum flux, and further corrections make
the residual smooth through the singular time.

\subsection{Organization of the paper}

Section~\ref{sec:pulse-framework} sets up the coordinates,
localizations, principal equations and compatibility between auxiliary representations.
Section~\ref{sec:pulse-dynamics} constructs the unforced and forced
pre-pulses, matches the complete background stress, and gives a linear inverse for
lower-order stress increments.
Section~\ref{sec:initialization-ppt} states the initialization targets
and constructs the initial state by the primary pulse, the auxiliary
correction and the moment correction, followed by a single pressure
reconstruction. Section~\ref{sec:iteration-ppt} develops the
subsequent correction cycles and their estimates, constructs a smooth
realization of the iteration, and completes the proof by localization,
following \cite[Section~10]{OpenAI}.

\section{A framework for pulse construction}
\label{sec:pulse-framework}

\makeatletter
\newcommand{\pulseLabelAlias}[2]{\AtEndDocument{\@ifundefined{r@#2}{}{\begingroup
      \expandafter\let\expandafter\pulse@labeldata\csname r@#2\endcsname
      \immediate\write\@auxout{\string\newlabel{#1}{\unexpanded\expandafter{\pulse@labeldata}}}\endgroup
    }}}
\makeatother
\pulseLabelAlias{eq:covered-derivatives}{s6n:physical-derivatives}
\pulseLabelAlias{eq:covering-index}{s6n:covering-level}
\pulseLabelAlias{eq:covering-size}{eq:ck-size}
\pulseLabelAlias{eq:local-pulse-time}{s6n:pulse-time}
\pulseLabelAlias{eq:physical-derivatives}{s6n:physical-derivatives}
\pulseLabelAlias{eq:radial-exponent}{s6n:auxiliary-map}
\pulseLabelAlias{lem:pulse-separation}{s6n:separation}
\pulseLabelAlias{prop:local-pulse-time}{s6n:time-prop}
\pulseLabelAlias{prop:rectangle-preimages}{s6n:time-prop}
\pulseLabelAlias{sec:pulse-localization}{s6n:localizations}
\pulseLabelAlias{sec:pulse-scales}{s6n:auxiliary}
\pulseLabelAlias{eq:common-rectangle-radius}{s6n:common-rectangle-radius}
\pulseLabelAlias{eq:coordinates}{s6n:coordinates}
\pulseLabelAlias{eq:local-partition}{s6n:physical-cutoffs}
\pulseLabelAlias{eq:overlapping-covering-indices}{s6n:overlapping-covering-indices}
\pulseLabelAlias{eq:pulse-separation}{s6n:separated-supports}

A \emph{pulse} is a localized, exactly divergence-free oscillatory
velocity field, obtained by completing a \emph{pre-pulse} with a curl
correction. A primary pulse supplies an averaged quadratic momentum
flux; forced pulses absorb nonzero angular modes (also called harmonics
in \cite{OpenAI}) of the residual. This section fixes their local
coordinates and supports, and records the principal amplitude equation
for the primary pulse. We also define the forced pulse ansatz and
record its principal equation with a linear source. Finally, we describe
compatibility between auxiliary representations.

\subsection{The normalized coordinates and fields}
\label{s6n:slow}

Fix \(0<\delta<1/200\) and a sufficiently large integer \(k_0\).
A label \(\gamma=(k,\alpha,\iota)\) records three choices:
\(k\) specifies the scale band \(\mathcal B_k = \{(x,t):t<1, \ \lambda(t,z)\in(2^{-k-1},2^{-k+1})\} \subset \mathbb R^{3+1}\);
\(\alpha\in\ZZ^3\) labels a slow box, to be defined in
Section~\ref{s6n:localizations}; and \(\iota\in\{+,-\}\)
distinguishes the two pulse families used below.
We work on the domain \(0<\lambda(t,z)\le2^{-k_0}\), so that only bands
\(k\ge k_0\) are involved. In each band, define the \emph{normalized
coordinates} (also called \emph{slow coordinates} in \cite{OpenAI}) by
\begin{equation}\label{s6n:coordinates}
 \hat{R}=2^{2k-1}r^2,\qquad
 \hat{Z}=2^{(1-\delta)k}z,\qquad
 \hat{T}=2^{2k}(1-t).
\end{equation}
These coordinates use a fixed dyadic scale; recall that the similarity
coordinates of Part~I use the continuously varying scale \(\lambda\):
\[
 R=\frac{r^2}{2\lambda^2},\qquad Z=\frac{z}{\lambda^{1-\delta}},
 \qquad \lambda = \sqrt{\frac{1-t}{1-Z^2}}.
\]
Consequently,
\begin{equation}\label{s6n:coordinate-change}
 \hat{R}=(2^k\lambda)^2R,\qquad
 \hat{Z}=(2^k\lambda)^{1-\delta}Z,\qquad
 \hat{T}=(2^k\lambda)^2(1-Z^2).
\end{equation}
Note that on $\mathcal B_k$, \(2^k\lambda\) is a bounded positive function of
\((\hat{Z},\hat{T})\), independent of \(\hat{R}\).
On the active stress annulus
\((R,Z)\in[R_a,R_h]\times[-1,1]\), where
\(0<R_a<R_h<\infty\), this gives \(\hat{R}\asymp1\) and
\(|\hat{Z}|+|\hat{T}|\lesssim1\).
Holding \(k\) fixed, the chain rule for the normalized coordinates gives
\begin{equation}\label{s6n:slow-derivatives}
 \begin{aligned}
 2^{-k}\partial_r&=\sqrt{2\hat{R}}\,\partial_{\hat{R}},
 &2^{-(1-\delta)k}\partial_z&=\partial_{\hat{Z}},\\
 2^{-2k}\partial_t&=-\partial_{\hat{T}},
 &2^{-k}r^{-1}\partial_\theta&=(2\hat{R})^{-1/2}\partial_\theta .
 \end{aligned}
\end{equation}

\medskip
We normalize the background angular velocity, shear and stress in these
coordinates. Let \(u_B\) be the fixed smooth axisymmetric background
constructed in Part~I. Set
\begin{equation}\label{s6n:background-coefficients}
 \begin{gathered}
 \hat u_B=2^{-(1+\delta)k}u_B,\qquad
 \hat F_B=\frac{\hat u_B^\theta}{\sqrt{2\hat{R}}},\\
 \hat{\mathcal S}_B=
 \left(2\hat{R}\partial_{\hat{R}}\hat F_B,
       \sqrt{2\hat{R}}\,\partial_{\hat{R}}\hat u_B^z\right),\qquad
 \hat{\mathsf T}_B=2^{-(2+\delta)k}\mathsf T_B.
 \end{gathered}
\end{equation}
Here \(\hat F_B\) is the normalized angular velocity,
\(\hat{\mathcal S}_B\) records radial shear in the \((\theta,z)\)
components, and
\(\hat{\mathsf T}_B=(\hat{\mathsf T}_B^\theta,\hat{\mathsf T}_B^z)\)
is the normalized physical stress.

Since \(2^k\lambda\) does not depend on radius, the quantities
defined in \eqref{intro:full-background-shear} satisfy
\[
 (\hat F_B,\hat{\mathcal S}_B,\hat{\mathsf T}_B)
 =(2^k\lambda)^{-2-\delta}
   (F_B,\mathcal S_B,\lambda^{2+\delta}\mathsf T_B).
\]
The common positive multiplier preserves the stress direction and
the cone inequalities in \eqref{intro:admissible-cone}. Thus the
normalized background satisfies these inequalities with the uniform
positive margin \(c_{\rm cone}>0\) on
\((R_a,R_h)\times[-1,1]\).

For later estimates, \eqref{intro:background-velocity-bounds} and
the coordinate change \eqref{s6n:coordinate-change} give
\begin{equation}\label{s6n:background-velocity-bounds}
 \|\hat u_B^{\theta,z}\|_{C^m(\mathcal B_k \cap \{R_a \le R \le R_h\})}\le C_m,\qquad
 \|\hat u_B^r\|_{C^m(\mathcal B_k \cap \{R_a \le R \le R_h\})}\le C_m2^{-\delta k},\qquad m\ge0.
\end{equation}
Here \(C^m\) uses derivatives in the normalized coordinates
\((\hat R,\hat Z,\hat T)\). The constants depend on \(m\) and the fixed
background data, but not on \(k\).

\subsection{The auxiliary variables}
\label{s6n:auxiliary}

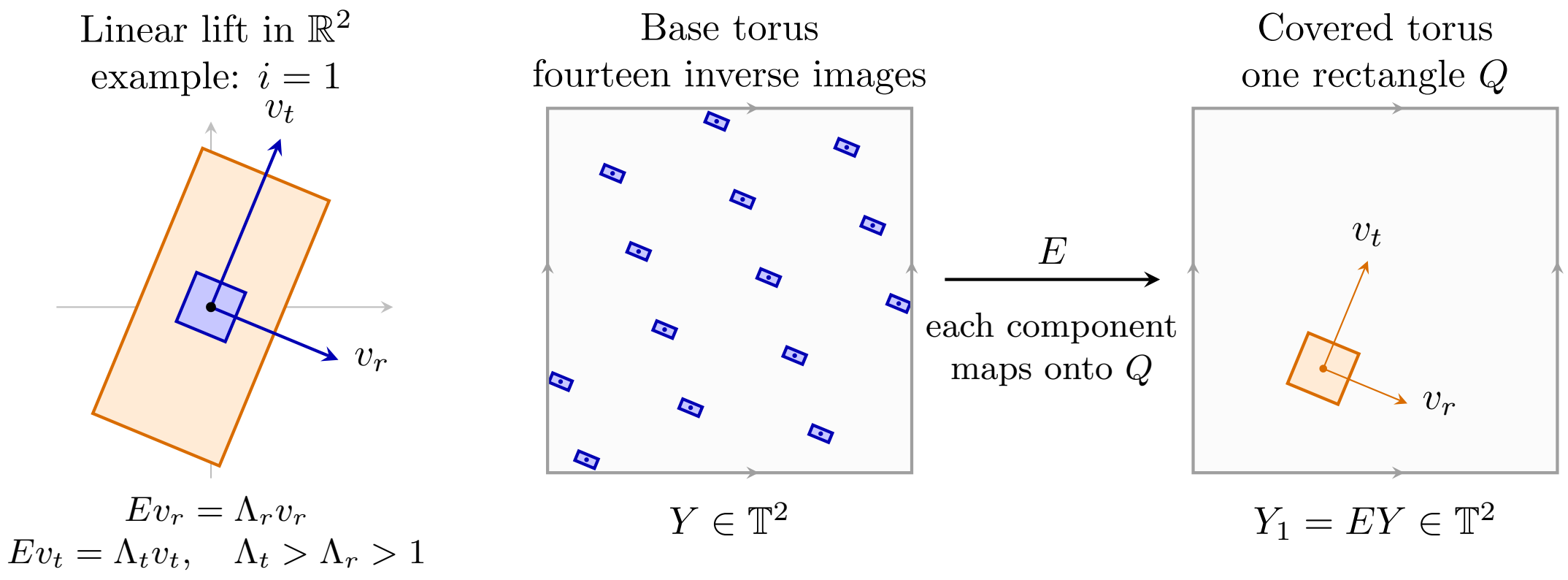
\begin{figure}[!ht]
\centering
\begin{tikzpicture}[>=stealth,font=\small,
 s6eigsource/.style={draw=blue!70!black,fill=blue!22,thick},
 s6eigtarget/.style={draw=orange!85!black,fill=orange!16,thick}]
 \node[align=center] at (1.55,3.8) {Linear lift in \(\RR^2\)\\example: \(i=1\)};
 \draw[->,gray!50] (.1,1.5)--(3.15,1.5);
 \draw[->,gray!50] (1.5,-.05)--(1.5,3.18);
 \begin{scope}[shift={(1.5,1.5)},rotate=-22.5]
  \pgfmathsetmacro{\er}{.24*(4-sqrt(2))}
  \pgfmathsetmacro{\et}{.24*(4+sqrt(2))}
  \path[s6eigtarget] (-\er,-\et) rectangle (\er,\et);
  \path[s6eigsource] (-.24,-.24) rectangle (.24,.24);
  \draw[->,blue!70!black,thick] (0,0)--(1.25,0)
       node[right,black] {\(v_r\)};
  \draw[->,blue!70!black,thick] (0,0)--(0,1.65)
       node[above,black] {\(v_t\)};
  \fill (0,0) circle (1.3pt);
 \end{scope}
 \node[align=center,font=\footnotesize] at (1.55,-.55)
   {\(Ev_r=\Lambda_rv_r\)\\\(Ev_t=\Lambda_tv_t,\quad\Lambda_t>\Lambda_r>1\)};
 \begin{scope}[xshift=4.55cm,x=3.3cm,y=3.3cm]
  \node[align=center] at (.5,1.15) {Base torus\\fourteen inverse images};
  \fill[gray!3] (0,0) rectangle (1,1);
  \foreach \j in {0,...,13} {
   \pgfmathsetmacro{\qx}{mod(.5+5*\j,14)/14}
   \pgfmathsetmacro{\qy}{mod(17.5-\j,14)/14}
   \begin{scope}[shift={(\qx,\qy)},rotate=-22.5]
    \pgfmathsetmacro{\ar}{.075/(4-sqrt(2))}
    \pgfmathsetmacro{\at}{.075/(4+sqrt(2))}
    \path[s6eigsource] (-\ar,-\at) rectangle (\ar,\at);
    \fill[blue!70!black] (0,0) circle (.6pt);
   \end{scope}
  }
  \draw[gray!75,thick] (0,0) rectangle (1,1);
  \draw[->,gray!75] (.40,0)--(.58,0);
  \draw[->,gray!75] (.40,1)--(.58,1);
  \draw[->,gray!75] (0,.40)--(0,.58);
  \draw[->,gray!75] (1,.40)--(1,.58);
  \node[below] at (.5,-.04) {\(Y\in\TT^2\)};
 \end{scope}
 \draw[->,thick] (8.15,1.75)--(10.1,1.75)
      node[midway,above] {\(E\)};
 \node[align=center,font=\footnotesize] at (9.12,1.12)
      {each component\\maps onto \(Q\)};
 \begin{scope}[xshift=10.4cm,x=3.3cm,y=3.3cm]
  \node[align=center] at (.5,1.15) {Covered torus\\one rectangle \(Q\)};
  \fill[gray!3] (0,0) rectangle (1,1);
  \begin{scope}[shift={(5/14,2/7)},rotate=-22.5]
   \path[s6eigtarget] (-.075,-.075) rectangle (.075,.075);
   \fill[orange!85!black] (0,0) circle (1pt);
   \draw[->,orange!85!black] (0,0)--(.25,0)
        node[right,black] {\(v_r\)};
   \draw[->,orange!85!black] (0,0)--(0,.32)
        node[above,black] {\(v_t\)};
  \end{scope}
  \draw[gray!75,thick] (0,0) rectangle (1,1);
  \draw[->,gray!75] (.40,0)--(.58,0);
  \draw[->,gray!75] (.40,1)--(.58,1);
  \draw[->,gray!75] (0,.40)--(0,.58);
  \draw[->,gray!75] (1,.40)--(1,.58);
  \node[below] at (.5,-.04) {\(Y_1=EY\in\TT^2\)};
 \end{scope}
\end{tikzpicture}
\caption{Expansion and covering, illustrated for \(i=1\).
Left: in the linear lift, \(E\) stretches a blue rectangle along the
orthogonal eigendirections to the orange rectangle, with the greater
expansion ratio along \(v_t\).
Middle and right: each of the fourteen blue rectangles maps onto the
orange rectangle \(Q\). Both squares represent unit tori, with opposite
sides identified. A sufficiently small target rectangle has
\(14^i\) inverse-image components.}
\label{fig:s6-torus-covering}
\end{figure}

The two-dimensional auxiliary torus is introduced to help separate the physical supports
of the pulses and to place the pulses suitably in physical space. Choose
\begin{equation}\label{s6n:eigendata}
 \begin{gathered}
 E=\begin{pmatrix}3&1\\1&5\end{pmatrix},\qquad
 \Lambda_r=4-\sqrt2,\quad\Lambda_t=4+\sqrt2,\\
 v_r=(1,1-\sqrt2)^{\mathsf T},\qquad
 v_t=(\sqrt2-1,1)^{\mathsf T},\qquad
 Ev_r=\Lambda_rv_r,\quad Ev_t=\Lambda_tv_t .
 \end{gathered}
\end{equation}
Both eigenvalues exceed one, \(\det E=14\), and \(v_r\perp v_t\).
Set
\begin{equation}\label{s6n:auxiliary-map}
 \begin{gathered}
 \rho_0=\frac{\log\Lambda_r}{\log\Lambda_t},\qquad
 d_0=(2+\delta)\rho_0-10^{-5}\delta>0,\\
 Y(t,r)=t\,v_t+r^{d_0}v_r\pmod{\ZZ^2},\qquad
 Y_i=E^iY\pmod{\ZZ^2}.
 \end{gathered}
\end{equation}
Since \(E\) is an integer matrix, \(Y\mapsto Y_i\) is well defined on
\(\TT^2=\RR^2/\ZZ^2\). The physical evaluation of \(Y_i\) reads
\[
 Y_i(t,r)=\Lambda_t^it\,v_t+\Lambda_r^ir^{d_0}v_r\pmod{\ZZ^2}.
\]
The linear lift \(E^i\) expands the two eigendirections by
\(\Lambda_r^i\) and \(\Lambda_t^i\). Modulo \(\ZZ^2\), it is a
\(14^i\)-to-one covering of the unit torus. A small target rectangle
therefore has many smaller inverse-image components.
For a set \(Q\subset\TT^2\), we write its full inverse image as
\[
 E^{-i}(Q):=\{Y\in\TT^2:E^iY\in Q\}.
\]
Figure~\ref{fig:s6-torus-covering} illustrates the expansion and the
covering of $E^i$. 

\medskip
\noindent\textbf{Physical-evaluation convention and physical derivatives.}
We write \(Y_i\) for either an independent variable in \(\TT^2\) or
the physical map \(Y_i(t,r)\), depending on the context. An
\emph{extended field} is a function \(f(t,r,\theta,z,Y_i)\) with
\(Y_i\) independent of the physical coordinates. Its \emph{physical
realization} is obtained by substituting \(Y_i=Y_i(t,r)\), and depends
only on \((t,r,\theta,z)\). Whenever the meaning is clear, we use the
same symbol \(f\) for both; this convention will be applied to the
phases, cutoffs and amplitudes that define the pulses.
The partial derivatives \(\nabla_{Y_i},\partial_t,\partial_r\)
act before physical evaluation, with the other independent variables
held fixed. The operators \(D_t,D_r\) represent physical differentiation
through the chain rule on extended fields; after evaluation, their
values agree with the ordinary derivatives of the physical realization.
When this distinction matters, we display the arguments or state
explicitly which variables are held fixed. Averages in \(Y\) or \(Y_i\)
apply to the extended field specified by the construction, before
physical evaluation.

On extended fields, the chain-rule operators are
\begin{equation}\label{s6n:physical-derivatives}
 \begin{aligned}
 D_t&= \partial_t + v_t \cdot \nabla_Y = \partial_t+\Lambda_t^iv_t\cdot\nabla_{Y_i},\\
 D_r&= \partial_r + d_0 r^{d_0 - 1} v_r \cdot \nabla_Y = \partial_r+\Lambda_r^id_0r^{d_0-1}v_r\cdot\nabla_{Y_i}.
 \end{aligned}
\end{equation}
Here \(\partial_t,\partial_r\) hold the independent variable
\(Y_i\) fixed. The chain rule is simply
\[
 \partial_q\bigl[f(t,r,\theta,z,Y_i(t,r))\bigr]
   =(D_q f)(t,r,\theta,z,Y_i(t,r)),\qquad q=t,r.
\]
Since \(Y_i(t,r)\) has no \(\theta,z\) dependence, differentiation
in \(\theta,z\) commutes with physical substitution. It is also
straightforward to check that
\[
 [D_t,D_r]=[D_t,\partial_z]=[D_r,\partial_z]
  =[D_t,\partial_\theta]=[D_r,\partial_\theta]=0.
\]
We only use these definitions away from the axis.

For a fixed scale band \(k\), choose the covering level
\begin{equation}\label{s6n:covering-level}
 i=i(k)=\left\lfloor
 \log_{\Lambda_t}\left(\frac{2^{(2+\delta)k}}{k^2}\right)
 \right\rfloor.
\end{equation}
With this choice understood, we suppress the superscript \((i)\) and
simply write \(D_t,D_r\). For a sufficiently large lower bound \(k_0\)
on \(k\), we have \(i(k)\ge1\).
Combining \eqref{s6n:physical-derivatives} with
\eqref{s6n:slow-derivatives}, for functions in
\((\hat R,\hat Z,\hat T,Y_i)\) we obtain
\begin{equation}\label{s6n:normalized-derivatives}
 \begin{aligned}
 D_t&=-2^{2k}\partial_{\hat T}
       +\Lambda_t^i v_t\cdot\nabla_{Y_i},\\
 D_r&=2^k\sqrt{2\hat R}\,\partial_{\hat R}
     +2^{(1-d_0)k}\Lambda_r^i d_0(2\hat R)^{(d_0-1)/2}
                           v_r\cdot\nabla_{Y_i}.
 \end{aligned}
\end{equation}
The choice of \(i(k)\) ensures
\begin{equation}\label{eq:ck-size}
 \begin{aligned}
 \frac1{\Lambda_t k^2}
 &<2^{-(2+\delta)k}\Lambda_t^{i(k)}\le\frac1{k^2},\\
 2^{-kd_0}\Lambda_r^{i(k)}
 &=2^{10^{-5}\delta k}
   \bigl(2^{-(2+\delta)k}\Lambda_t^{i(k)}\bigr)^{\rho_0}
   \asymp2^{10^{-5}\delta k}k^{-2\rho_0}.
 \end{aligned}
\end{equation}
The relation \(2^{-(2+\delta)k}\Lambda_t^{i(k)}\asymp k^{-2}\)
turns a fixed auxiliary interval into a fast pulse-time interval of
length \(O(k^2)\).
Here radial expansion means the growth of the covering in the
\(v_r\)-direction:
\[
 E^{i(k)}v_r=\Lambda_r^{i(k)}v_r,\qquad
 \Lambda_r^{i(k)}\asymp
 2^{(2+\delta)\rho_0 k}k^{-2\rho_0}.
\]
The factor \(2^{-kd_0}\), with
\(d_0=(2+\delta)\rho_0-10^{-5}\delta\), compensates for this growth.
In the normalized radial derivative, the auxiliary coefficient is
\[
 2^{-k}d_0r^{d_0-1}\Lambda_r^{i(k)}
 =d_0(2\hat R)^{(d_0-1)/2}
   2^{-kd_0}\Lambda_r^{i(k)}
 \asymp2^{10^{-5}\delta k}k^{-2\rho_0}
\]
on the active stress annulus. Thus the exponential loss is the extra
factor \(2^{10^{-5}\delta k}\) in a coefficient bound such as
\[
 |2^{-k}D_r f|\le
 C2^{10^{-5}\delta k}
 \bigl(|\partial_{\hat R}f|+|\nabla_{Y_i}f|\bigr).
\]

\medskip
\noindent\textbf{Diophantine estimates.}
Both eigenvectors $v_r$, $v_t$ are irrational with respect to $\mathbb T^2$, and they satisfy

\begin{lem}
For all \(n=(n_1,n_2)\in\ZZ^2\setminus\{0\}\), we have
\begin{equation}
|v_r\cdot n|\ge\frac{1}{(1+2\sqrt2)(1+|n|)},\qquad
|v_t\cdot n|\ge\frac{1}{(1+2\sqrt2)(1+|n|)},
\label{eq:diophantine-directions}
\end{equation}
where \(|n|=(n_1^2+n_2^2)^{1/2}\).
\end{lem}

\begin{proof}
  Let us prove the lower bound for 
\[
v_r\cdot n=(n_1+n_2)-\sqrt2\,n_2.
\]
If \(n_2=0\), then \(|v_r\cdot n|=|n_1|\ge1\).
If \(n_2\ne0\), then
\[
|v_r\cdot n|=|n_2|\left|\frac{n_1+n_2}{n_2}-\sqrt2\right|,
\]
so the question is how closely the rational number \((n_1+n_2)/n_2\) can approximate
\(\sqrt2\). Consider the quadratic polynomial
\[
P(x)=x^2-2,\qquad
P\!\left(\frac{n_1+n_2}{n_2}\right)
=\frac{(n_1+n_2)^2-2n_2^2}{n_2^2}.
\]
Since \(\sqrt2\) is irrational, \((n_1+n_2)^2-2n_2^2\) is a nonzero integer.
Its absolute value is therefore at least one. Factoring \(P\) gives
\[
\left|\frac{n_1+n_2}{n_2}-\sqrt2\right|
\left|\frac{n_1+n_2}{n_2}+\sqrt2\right|
=\left|P\!\left(\frac{n_1+n_2}{n_2}\right)\right|
\ge\frac1{n_2^2}.
\]
If \(\left|(n_1+n_2)/n_2-\sqrt2\right|\le1\), then
\(\left|(n_1+n_2)/n_2+\sqrt2\right|\le1+2\sqrt2\), and hence
\[
\left|\frac{n_1+n_2}{n_2}-\sqrt2\right|\ge\frac{1}{(1+2\sqrt2)n_2^2}.
\]
If \(\left|(n_1+n_2)/n_2-\sqrt2\right|>1\), the same bound holds because
\(|n_2|\ge1\) and \((1+2\sqrt2)^{-1}<1\).
Thus, in either case,
\[
\begin{aligned}
|v_r\cdot n|
&=|n_2|\left|\frac{n_1+n_2}{n_2}-\sqrt2\right|\\
&\ge\frac{1}{(1+2\sqrt2)|n_2|}
 \ge\frac{1}{(1+2\sqrt2)|n|}
 \ge\frac{1}{(1+2\sqrt2)(1+|n|)}.
\end{aligned}
\]
For the second estimate, the same rational-approximation argument applies to
\[
v_t\cdot n=\sqrt2\,n_1+(n_2-n_1).
\]
If \(n_1=0\), its absolute value is at least one. Otherwise, use the
rational number \((n_1-n_2)/n_1\) and the same polynomial \(P(x)=x^2-2\).
This proves the second bound as well.
\end{proof}

These bounds allow us to invert directional derivatives on the torus
with a finite loss of derivatives: an \(H^{s+1}\) bound on the source
controls the \(H^s\) norm of the solution. In particular, smooth
zero-mean sources still give smooth zero-mean solutions.

\begin{lem}[Inverse of a torus directional derivative]
\label{s6n:directional-inverse}
Let \(v\in\{v_r,v_t\}\) and let \(f\in C^\infty(\TT^2)\) have
zero torus mean. The equation \(v\cdot\nabla_{Y}g=f\) has a unique
smooth zero-mean solution, given by
\begin{equation}\label{s6n:torus-inverse}
 g(Y)=\sum_{n\ne0}
       \frac{\widehat f(n)}{2\pi\mathrm i\,v\cdot n}
       e^{2\pi\mathrm i n\cdot Y}.
\end{equation}
For every \(s\ge0\), it satisfies
\[
 \|g\|_{H^s(\TT^2)}\le C_s\|f\|_{H^{s+1}(\TT^2)}, \qquad \|g\|_{C^s(\TT^2)} \le C_s \|f\|_{C^{s+3}(\TT^2)}.
\]
\end{lem}

\begin{proof}

The Fourier coefficients of a smooth \(f\) decay faster than any
power, so the series and all its derivatives converge uniformly and
solve the equation. 
Since \(v\cdot n\ne0\) for every \(n\ne0\), the equation determines
all nonzero Fourier coefficients of \(g\); the zero-mean condition
determines the remaining coefficient $\widehat g(0)=0$. This proves uniqueness.

For \(n\ne0\),
\eqref{eq:diophantine-directions} gives
\[
|\widehat g(n)|
\le\frac{1+2\sqrt2}{2\pi}(1+|n|)|\widehat f(n)|.
\]
Multiplying the squared bound by
\((1+|n|^2)^s\) and summing over $n$ proves the Sobolev estimate. The $C^s$ estimate then follows from Sobolev embedding on the torus.
\end{proof}

\subsection{Localization in normalized and auxiliary variables}
\label{s6n:localizations}

In this subsection, we introduce various families of localizations, either directly in physical spacetime or in the auxiliary variables.

\smallskip
\noindent\textbf{Scale bands and slow boxes.}
We use nonnegative smooth cutoffs $\zeta_k$ in the scale parameter $\lambda$ with the following properties:
\begin{equation}\label{s6n:physical-cutoffs}
 \zeta_k(\lambda)=\zeta(2^k\lambda),\qquad
 \operatorname{supp}\zeta\Subset(1/2,2),\qquad
 \sum_{k\in\ZZ}\zeta_k(\lambda)^2=1.
\end{equation}
Further, for each $k$ we introduce a family of spacetime cut-offs $\zeta_{k, \alpha} = \zeta_{k, \alpha}(\hat{R},\hat{Z},\hat{T})$ satisfying 
\begin{equation}\label{s6n:slow-box-cutoffs}
\begin{gathered}
 Q_{k,\alpha}:=\left\{(\hat{R},\hat{Z},\hat{T}):
 (\sqrt{2\hat{R}},\hat{Z},\hat{T})
 \in k^{-6}\bigl(\alpha+(-1,1)^3\bigr)\right\},\\
 \operatorname{supp}(\zeta_k\zeta_{k,\alpha})
 \subset\operatorname{supp}\zeta_{k,\alpha}\subset Q_{k,\alpha},\\
 \sum_{\alpha\in\ZZ^3}\zeta_{k,\alpha}^2=1,\qquad
 \|\{\sqrt{2\hat{R}}\,\partial_{\hat{R}},
 \partial_{\hat{Z}},\partial_{\hat{T}}\}^\beta
 \zeta_{k,\alpha}\|_\infty\le C_\beta k^{6|\beta|}.
\end{gathered}
\end{equation}
Square partitions are a convenient choice here because momentum
flux will be quadratic in velocity. For later convenience, we duplicate the second family of cutoffs and denote 
$$\zeta_{k, \alpha, \pm} = \zeta_{k, \alpha}.$$

\begin{figure}[!ht]
\centering
\begin{tikzpicture}[>=stealth,font=\small]
 \def\sSixBoxCluster#1#2{\foreach \i in {0,1,2} {
   \foreach \j in {0,1,2} {
    \path[fill=gray!10] (1.05*\i,.75*\j) rectangle ++(1.4,1);
   }
  }
  \path[fill=blue!15] (1.05,.75) rectangle ++(1.4,1);
  \path[fill=orange!17] (2.10,1.50) rectangle ++(1.4,1);
  \path[fill=purple!30] (2.10,1.50) rectangle (2.45,1.75);
  \foreach \i in {0,1,2} {
   \foreach \j in {0,1,2} {
    \draw[gray!65,line width=#1]
      (1.05*\i,.75*\j) rectangle ++(1.4,1);
   }
  }
  \draw[blue!70!black,line width=#2]
    (1.05,.75) rectangle ++(1.4,1);
  \draw[orange!85!black,line width=#2]
    (2.10,1.50) rectangle ++(1.4,1);
 }
 \node at (2.25,4.12) {Band \(k\) on the slice \(z=0\)};
 \fill[gray!9] (0,.4) rectangle (4.5,3.4);
 \path[fill=blue!7]
   plot[domain=.25:4,samples=60] ({sqrt(\x)},{3.6-.8*\x})
   --plot[domain=4:.25,samples=60] ({2*sqrt(\x)},{3.6-.8*\x})
   --cycle;
 \draw[gray!55] (0,.4)--(4.5,.4);
 \draw[gray!55] (0,3.4)--(4.5,3.4);
 \draw[gray!60,dashed] (0,3.6)--(4.5,3.6);
 \draw[blue!45!black,semithick,domain=.25:4,samples=60]
   plot ({sqrt(\x)},{3.6-.8*\x});
 \draw[blue!45!black,semithick,domain=.25:4,samples=60]
   plot ({2*sqrt(\x)},{3.6-.8*\x});
 \draw[->] (0,0)--(4.75,0) node[right] {\(r\)};
 \draw[->] (0,0)--(0,3.92) node[above] {\(t\)};
 \draw (-.06,.4)--(.06,.4) node[left,xshift=-3pt] {\(t_-\)};
 \draw (-.06,3.4)--(.06,3.4) node[left,xshift=-3pt] {\(t_+\)};
 \node[above left] at (0,3.6) {\(1\)};
 \node[font=\footnotesize] at (3.74,3.08) {\(k\)-band};
 \node[font=\footnotesize,fill=blue!7,inner sep=1pt]
   at (1.73,.7) {\(R=R_a\)};
 \node[font=\footnotesize,fill=blue!7,inner sep=1pt]
   at (4.02,.7) {\(R=R_h\)};
 \begin{scope}[shift={(2.15,1.60)},scale=.12]
  \sSixBoxCluster{.20pt}{.35pt}
 \end{scope}
 \draw[red!60!black,dashed,semithick] (2.10,1.55) rectangle (2.62,1.95);
 \draw[->,blue!70!black,thin] (1.98,2.12)--(2.36,1.75);
 \node[above left,text=blue!70!black,inner sep=1pt] at (1.98,2.12) {\(\alpha\)};
 \node[font=\footnotesize] at (2.25,.16) {\(r\asymp2^{-k}\)};
 \node[align=center,font=\footnotesize] at (2.25,-.78)
   {\(t_-=1-4\,2^{-2k}\)\\[2pt]\(t_+=1-\tfrac14\,2^{-2k}\)};
 \draw[gray!55,dashed] (2.62,1.95)--(6.3,3.4);
 \draw[gray!55,dashed] (2.62,1.55)--(6.3,.4);
 \node[fill=white,inner sep=2pt,font=\footnotesize] at (5.23,1.92) {zoom};
 \begin{scope}[xshift=6.3cm]
  \node at (2.2,4.12) {Zoom: \(\alpha\) and its neighbors};
  \filldraw[fill=gray!3,draw=gray!55] (0,.4) rectangle (4.4,3.4);
  \draw[->] (0,0)--(4.65,0) node[right] {\(r\)};
  \draw[->] (-.13,.4)--(-.13,3.65) node[above] {\(t\)};
  \begin{scope}[shift={(.45,.65)}]
   \sSixBoxCluster{.45pt}{.85pt}
  \end{scope}
  \node[text=blue!70!black] at (2.20,1.90) {\(\alpha\)};
  \node[text=orange!85!black] at (3.25,2.65) {\(\alpha'\)};
  \draw[<->,blue!70!black] (1.50,.48)--(2.90,.48)
      node[midway,below,font=\footnotesize] {\(\Delta r\)};
  \draw[<->,blue!70!black] (4.15,1.40)--(4.15,2.40)
      node[midway,right,font=\footnotesize] {\(\Delta t\)};
  \node[align=center,font=\footnotesize] at (2.2,-.78)
    {\(\Delta r=O(2^{-k}k^{-6})\)\\[2pt]
     \(\Delta t=O(2^{-2k}k^{-6})\)};
 \end{scope}
\end{tikzpicture}
\caption{A band and a cluster of overlapping slow boxes in physical
coordinates. On \(z=0\), \(\lambda=\sqrt{1-t}\), so the \(k\)-band is
\(t_-<t<t_+\); the pale blue region is the stress annulus
\(R_a\le r^2/[2(1-t)]\le R_h\).
Both panels show the same cluster: the central box \(\alpha\) and its
eight neighbors on this slice. Boxes farther away are omitted.
The blue and orange outlines distinguish \(\alpha\) and one neighbor
\(\alpha'\), with their overlap shaded purple.
Each box is smaller than the radial and time scales \(2^{-k}\) and
\(2^{-2k}\) by the additional factor \(k^{-6}\).
The support rectangles and magnification are schematic. }
\label{fig:s6-band-boxes}
\end{figure}
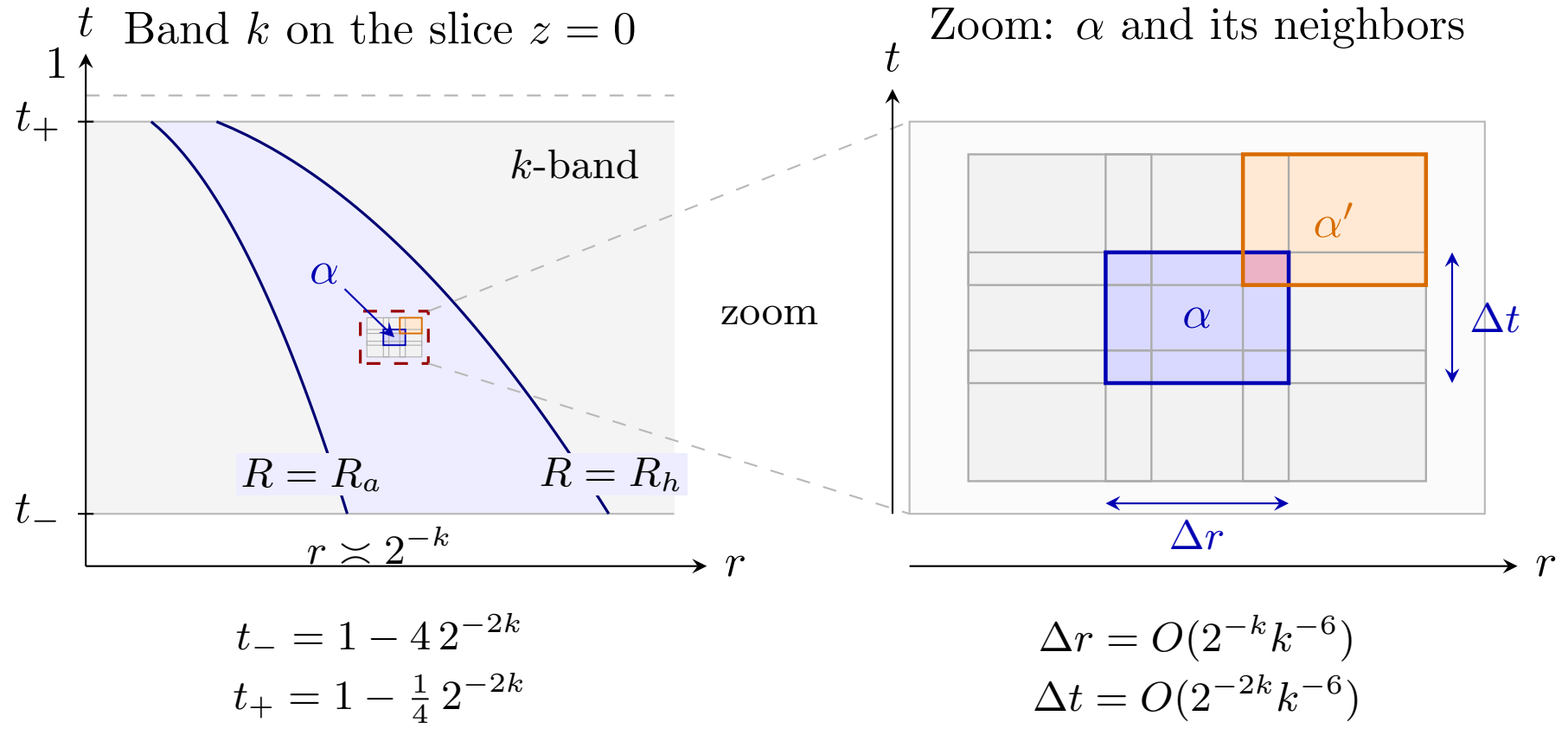

We refer to $\zeta_k$ as the band cutoffs and $\zeta_{k, \alpha}$ as the box cutoffs. The supports of $\zeta_{k, \alpha}$ will be called the slow boxes. In our region \(\lambda\le2^{-k_0}\),
only bands \(k\ge k_0\) contribute. On the stress annulus \(r\asymp2^{-k}\), a slow box has physical widths
\begin{equation}\label{s6n:box-widths}
 \Delta r=O(2^{-k}k^{-6}),\qquad
 \Delta z=O(2^{-(1-\delta)k}k^{-6}),\qquad
 \Delta t=O(2^{-2k}k^{-6}).
\end{equation}  
At a point, at most two band
cutoffs and at most eight box cutoffs in each band are nonzero.  
Figure~\ref{fig:s6-band-boxes} shows the band and box scales on the
physical \((r,t)\)-plane, together with a magnified view of the overlaps. Overlapping bands and boxes still permit unwanted products between distinct
pulses. The subsequent auxiliary localization will make those products identically zero.

\medskip
\noindent\textbf{Auxiliary rectangles and fast pulse time.}
For each label \(\gamma=(k,\alpha,\iota = \pm)\), we will choose a center \(O_\gamma\in\TT^2\).
A common width parameter \(r_0>0\) and all the centers will be chosen together
in Proposition~\ref{s6n:separation}. Using \(Y_i\) as the independent
coordinate on the target torus for band \(k\), set
\begin{equation}\label{s6n:rectangles}
 \mathcal R_\gamma
   =\{O_\gamma+l_rv_r+l_tv_t: |l_r|,|l_t|<r_0\}
 \qquad \pmod{\ZZ^2}.
\end{equation}
We require \(r_0<1/(4\sqrt2)\). Then the eigenvector coordinates
\((l_r,l_t)\), defined by \(Y_i=O_\gamma+l_rv_r+l_tv_t\), are
well defined and unique on the rectangle \(\mathcal R_\gamma\).
Here \(r_0\) is the half-width of \(\mathcal R_\gamma\) in each
eigenvector coordinate. These coordinates are local to the rectangles.

Define the fast pulse time by
\begin{equation}\label{s6n:pulse-time}
\begin{aligned}
 \tau_\gamma(Y_i)&=2^{(2+\delta)k}\Lambda_t^{-i}(l_t(Y_i)+r_0),\\
 \tau_{\rm max}=\tau_{\rm max}(k)&=2r_0\,2^{(2+\delta)k}\Lambda_t^{-i}.
\end{aligned}
\end{equation}
Thus \(\tau_\gamma\) runs from \(0\) to \(\tau_{\rm max}\) across
\(\mathcal R_\gamma\). By \eqref{eq:ck-size},
\[
 k^2\le2^{(2+\delta)k}\Lambda_t^{-i}<\Lambda_t k^2,
\]
and hence
\[
 2r_0k^2\le\tau_{\rm max}<2r_0\Lambda_t k^2.
\]
We omit \(\gamma\) from \(\tau\) when working with one label.

\begin{prop}[Pulse coordinates and time]\label{s6n:time-prop}
For the rectangle in \eqref{s6n:rectangles}, its full inverse image
under \(Y\mapsto E^{i(k)}Y\) consists of \(14^{i(k)}\) disjoint
rectangles, each mapped bijectively onto the given target rectangle.
On each component of the physical inverse image of
\(\mathcal R_\gamma\), the physical evaluation of \(\tau_\gamma(Y_i)\)
is a well-defined function of \(t\) of the form
\[
 \tau_\gamma(Y_i(t,r))=2^{(2+\delta)k}t+\text{constant},
\]
where the constant is determined by the lift of the target center.
The physical derivatives
of the local pulse coordinates satisfy
\begin{equation}\label{s6n:time-derivatives}
\begin{gathered}
 2^{-(2+\delta)k}D_t^{(i)}\tau_\gamma=1,\qquad
 D_r^{(i)}\tau_\gamma=\partial_z\tau_\gamma=0,\\
 D_t^{(i)}l_r=\partial_zl_r=0.
\end{gathered}
\end{equation}
\end{prop}

\begin{proof}
The inverse-image branches are
\[
 E^{-i}(O_\gamma+n)+\Lambda_r^{-i}l_rv_r
                   +\Lambda_t^{-i}l_tv_t\pmod{\ZZ^2},
 \qquad n\in\ZZ^2/E^i\ZZ^2.
\]
The quotient has \((\det E)^i=14^i\) elements. Each branch maps onto
the rectangle with the same coordinates \((l_r,l_t)\). If two branches
meet, uniqueness of those coordinates implies \(n-n'\in E^i\ZZ^2\),
so they are the same branch. The same observation proves injectivity.

On a fixed physical lift, write
\(Y_i=O_\gamma+n+l_rv_r+l_tv_t\), with a representative of
\(O_\gamma\) fixed and \(n\in\ZZ^2\). Orthogonality gives
\[
 l_t=\Lambda_t^it-\frac{v_t\cdot(O_\gamma+n)}{|v_t|^2}.
\]
Substituting this identity into \eqref{s6n:pulse-time} gives
the stated affine form. In local coordinates,
\[
 (v_t\cdot\nabla_{Y_i})l_t=(v_r\cdot\nabla_{Y_i})l_r=1,\qquad
 (v_r\cdot\nabla_{Y_i})l_t=(v_t\cdot\nabla_{Y_i})l_r=0.
\]
Hence \(D_tl_t=\Lambda_t^i\) and
\(D_rl_t=D_tl_r=\partial_zl_t=\partial_zl_r=0\).
Multiplication by \(2^{(2+\delta)k}\Lambda_t^{-i}\) proves \eqref{s6n:time-derivatives}.
\end{proof}

\begin{figure}[!ht]
\centering
\begin{tikzpicture}[
  >=stealth,
  font=\small,
  component/.style={
    draw=blue!70!black,fill=blue!22,line width=.45pt
  },
  target/.style={
    draw=orange!85!black,fill=orange!22,thick
  }
]
\def\torusframe{\draw[gray!65] (0,0) rectangle (3.2,3.2);
  \foreach \s in {0,3.2} {
    \draw[->,gray!65] (1.35,\s)--(1.85,\s);
    \draw[->,gray!65] (\s,1.35)--(\s,1.85);
  }
}

\node at (2,3.5) {Physical slow box, fixed \(z\)};
\filldraw[fill=gray!4,draw=gray!65]
  (0,0) rectangle (3.2,3.2);

\draw[->] (0,-.13)--(3.45,-.13) node[right] {\(r\)};
\draw[->] (-.13,0)--(-.13,3.4) node[above] {\(t\)};

\foreach \p in {-2,...,2} {
  \foreach \q in {-2,...,2} {
    \pgfmathsetmacro{\px}{1.20+.14*\p-.07*\q}
    \pgfmathsetmacro{\py}{1.80+.035*\p+.12*\q}
    \path[component]
      ({\px-.043},{\py-.022})
      rectangle
      ({\px+.043},{\py+.022});
  }
}
\node[text=gray!65] at (2.25,1.8) {\(\cdots\)};
\node[text=gray!65] at (1.2,2.5) {\(\vdots\)};

\node[font=\footnotesize] at (1.6,-.62)
  {\(Y_i(t,r)\in\mathcal R_\gamma\)};

\draw[->,thick] (3.5,1.6)--(4.5,1.6)
  node[midway,above] {\(Y(t,r)\)};

\begin{scope}[xshift=4.8cm]
  \node at (1.6,3.65) {Base torus};
  \fill[gray!3] (0,0) rectangle (3.2,3.2);

  \foreach \p in {-2,...,2} {
    \foreach \q in {-2,...,2} {
      \pgfmathsetmacro{\qx}{1.6+.26*\p-.08*\q}
      \pgfmathsetmacro{\qy}{1.6-.08*\p+.10*\q}
      \begin{scope}[shift={(\qx,\qy)},rotate=-22.5]
        \path[component]
          (-.058,-.021) rectangle (.058,.021);
      \end{scope}
    }
  }
  \node[text=gray!65] at (2.75,1.6) {\(\cdots\)};
  \node[text=gray!65] at (1.6,2.55) {\(\vdots\)};

  \torusframe
  \node[font=\footnotesize] at (1.6,-.62)
    {\(E^{-i}(\mathcal R_\gamma)\subset\mathbb T^2_Y\)};
\end{scope}

\draw[->,thick] (8.3,1.6)--(9.3,1.6)
  node[midway,above] {\(E^i\)};

\begin{scope}[xshift=9.6cm]
  \node at (1.6,3.65) {Covered target torus};
  \fill[gray!3] (0,0) rectangle (3.2,3.2);

  \begin{scope}[shift={(1.6,1.6)},rotate=-22.5]
    \path[target] (-.43,-.43) rectangle (.43,.43);

    \node at (0,0) {\(\mathcal R_\gamma\)};
  \end{scope}

  \torusframe
  \node[font=\footnotesize] at (1.6,-.62)
    {\(Y_i=E^iY\in\mathbb T^2\)};
\end{scope}
\end{tikzpicture}

\caption{Physical and auxiliary localization for one label
\(\gamma=(k,\alpha,\iota)\), with \(i=i(k)\).
Only a small cluster of physical pulse rectangles inside the slow box
and selected components of \(E^{-i}(\mathcal R_\gamma)\) are drawn;
the remaining components are omitted.
The physical pulse rectangles correspond to
\(Y(t,r)\in E^{-i}(\mathcal R_\gamma)\).
Their centers form a lattice in \((r^{d_0},t)\), so their placement
in the \((r,t)\)-panel is locally schematic.}
\label{fig:s6-physical-auxiliary}
\end{figure}
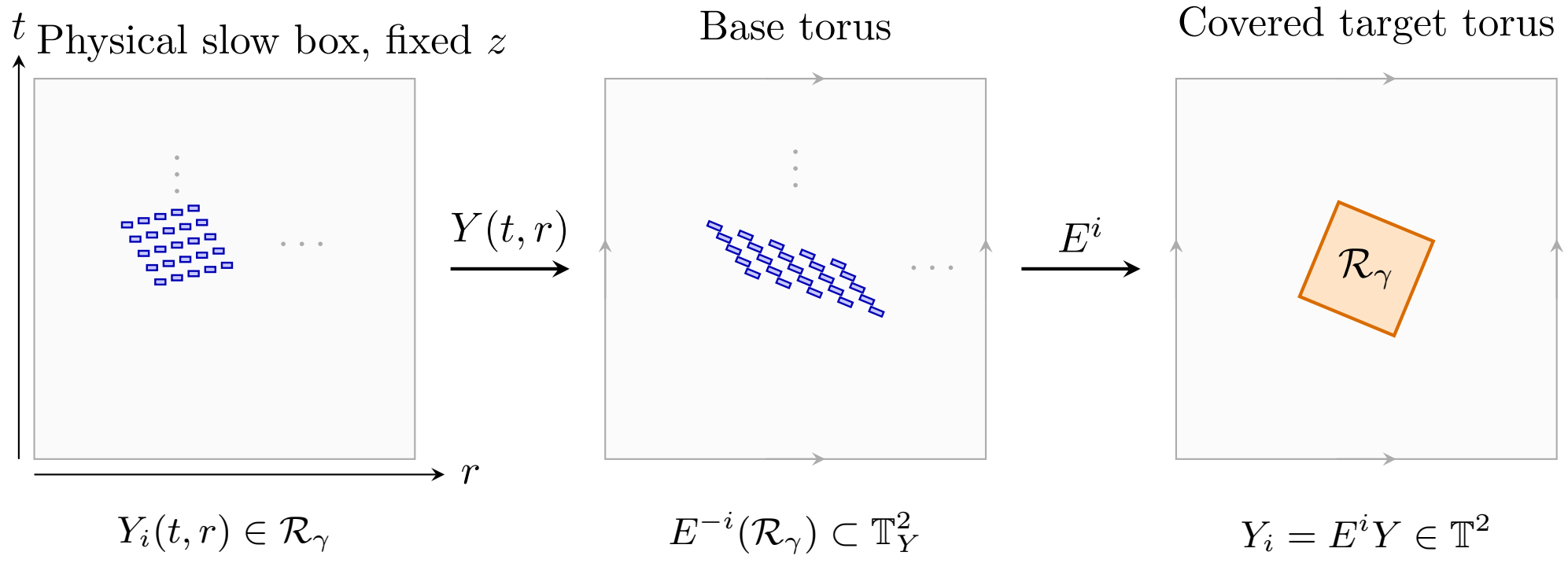

The rectangle \(\mathcal R_\gamma\) contains the support of
the auxiliary cutoff and hence of the
pulses to be constructed. The fast-time equations below,
including the forced principal equation, are solved for
\(0\le\tau\le\tau_{\rm max}\) before the temporal cutoff is applied.

We call the connected components of
\(\{(t,r):Y_i(t,r)\in\mathcal R_\gamma\}\) the
\emph{physical pulse rectangles}.
Thus, within one physical pulse rectangle and at fixed \(r,z\), an increment
\(\Delta\tau\) corresponds exactly to
\(\Delta t=2^{-(2+\delta)k}\Delta\tau\). See Figure \ref{fig:s6-physical-auxiliary} for an illustration.

Fix a representative of \(O_\gamma\) and a physical lift
\(Y_i=O_\gamma+n+l_rv_r+l_tv_t\), where the integer vector
\(n\in\ZZ^2\) selects the physical pulse rectangle. The coordinates are then
\[
\begin{aligned}
 l_t&=\Lambda_t^it-\frac{v_t\cdot(O_\gamma+n)}{|v_t|^2},\\
 l_r&=\Lambda_r^ir^{d_0}-\frac{v_r\cdot(O_\gamma+n)}{|v_r|^2}.
\end{aligned}
\]
The two constant offsets are the eigenvector-coordinate projections
of the lifted center \(O_\gamma+n\).
Thus each physical pulse rectangle is a product of intervals in \((t,r)\), with widths
\begin{equation}\label{s6n:aux-widths}
 \begin{aligned}
 \Delta t_{\rm aux}
   &\asymp(\Lambda_t^i)^{-1}
    \asymp 2^{-(2+\delta)k}k^2,\\
 \Delta r_{\rm aux}
   &\asymp(\Lambda_r^id_0r^{d_0-1})^{-1}
    \asymp 2^{-(1+10^{-5}\delta)k}k^{2\rho_0}
    \qquad(r\asymp2^{-k}).
 \end{aligned}
\end{equation}
Here we used \eqref{eq:ck-size}; the implicit constants may depend on
\(r_0,d_0\), which are fixed, but are independent of \(k\).

\medskip
\noindent\textbf{Physical separation via a choice of centers.}
Two distinct labels $\gamma$, $\gamma'$ are called neighbors when their physical slow supports
intersect, that is, 
$$\operatorname{supp}(\zeta_k\zeta_{k,\alpha}) \cap \operatorname{supp}(\zeta_{k'}\zeta_{k',\alpha'}) \neq \varnothing.$$
We now choose the centers $O_\gamma$ so that the physical realizations of the rectangles on neighboring slow boxes are separated.

\begin{prop}[Simultaneous separation]\label{s6n:separation}
For \(0<\delta<1/200\) and \(k_0\ge20\), the centers \(O_\gamma\)
and one absolute width parameter \(0<r_0<1/(4\sqrt2)\) can be chosen
so that any two neighboring labels \(\gamma,\gamma'\) satisfy
\begin{equation}\label{s6n:separated-supports}
 E^{-i(k)}(\mathcal R_\gamma)
       \cap E^{-i(k')}(\mathcal R_{\gamma'})=\varnothing.
\end{equation}
Consequently, the physical pulse rectangles for neighboring
labels are also separated:
\begin{equation}
 \{ (t,r): Y(t,r) \in E^{-i(k)}(\mathcal R_\gamma) \}
       \cap  \{ (t,r): Y(t,r) \in E^{-i(k')}(\mathcal R_{\gamma'})\}=\varnothing.
\end{equation}
\end{prop}

\begin{proof}
To choose only finitely many centers, we first bound the number of labels
whose physical cutoffs can intersect a fixed cutoff.
Suppose the two labels \(\gamma=(k,\alpha,\iota)\) and
\(\gamma'= (k',\alpha',\iota')\) are neighbors.
At a common point, \eqref{s6n:physical-cutoffs} gives
\begin{equation}
|k-k'|\le1.
\label{s6n:overlapping-band-indices}
\end{equation}
Thus, for fixed \(k\), there are at most three possible values of \(k'\).

We next bound the possible grid indices \(\alpha'\).
The two bands use different rescalings of the same physical point.
For any such \(k'\), a common physical point \((t,r,z)\) of the two
supports satisfies, by \eqref{s6n:coordinates},
\begin{equation}
\hat{R}'=2^{2(k'-k)}\hat{R},\qquad
\hat{Z}'=2^{(1-\delta)(k'-k)}\hat{Z},\qquad
\hat{T}'=2^{2(k'-k)}\hat{T}.
\label{s6n:overlapping-coordinate-scales}
\end{equation}
Using \eqref{s6n:overlapping-coordinate-scales} and the support bounds
for the box cutoffs, we first obtain
\begin{equation}
\begin{aligned}
\bigl|(k')^{-6}\alpha_1'-2^{k'-k}k^{-6}\alpha_1\bigr|
&\le\bigl|(k')^{-6}\alpha_1'-\sqrt{2\hat{R}'}\bigr|
+\bigl|\sqrt{2\hat{R}'}-2^{k'-k}k^{-6}\alpha_1\bigr|\\
&\le(k')^{-6}+2^{k'-k}k^{-6}.
\end{aligned}
\label{s6n:overlapping-grid-indices}
\end{equation}
Then \(\bigl|\alpha_1'-2^{k'-k}(k'/k)^6\alpha_1\bigr|
\le1+2^{k'-k}(k'/k)^6<9\).
Similarly, one obtains
\[
\begin{aligned}
\left|\alpha_2'-2^{(1-\delta)(k'-k)}\left(\frac{k'}k\right)^6\alpha_2\right|
&<9,\\
\left|\alpha_3'-2^{2(k'-k)}\left(\frac{k'}k\right)^6\alpha_3\right|
&<30.
\end{aligned}
\]
Consequently, there are only finitely many choices for \(\alpha'\).
Their number is bounded uniformly in \(k,\alpha\) and \(k'\).
Note that we cannot directly bound \(|\alpha_1-\alpha_1'|\) uniformly
when \(k\ne k'\).
For \(k'=k\), one simply has \(\|\alpha'-\alpha\|_\infty\le2\).
There are also two possible values of \(\iota'\).
Thus there is an integer \(n_0\ge2\), independent
of the label, such that each label has at most \(n_0-1\) neighbors.

The center choice will need to distinguish only finitely many relative
coverings. We verify that intersecting bands satisfy \(|i(k')-i(k)|\le2\).
Moreover, by \eqref{s6n:covering-level} and the mean value theorem, one also has
\[
\begin{aligned}
|i(k')-i(k)|
&\le1+\frac{(2+\delta)|k'-k|\log2+2|\log k'-\log k|}
{\log\Lambda_t}\\
&\le1+\frac{2(2+\delta)\log2+4/k_0}{\log\Lambda_t}
<3.
\end{aligned}
\]
The additional \(1\) accounts for the difference of the fractional parts. Since the
left-hand side is an integer, we conclude that
\begin{equation}
|i(k')-i(k)|\le2.
\label{s6n:overlapping-covering-indices}
\end{equation}

To obtain separated auxiliary supports, we assign different colors to neighboring labels.
We exclude centers that would coincide after applying \(I,E\) or \(E^2\).
Enumerate the labels and successively assign each a color in \(\{1,\ldots,n_0\}\) and a center.
For each label, exclude the colors of its previously assigned neighbors
and choose the smallest remaining color \(\nu\).
Whenever a color \(\nu\) is used for the first time, choose its center
\(O_\nu\in\TT^2\) so that
\begin{equation}
O_\nu\notin
\bigcup_{\mu<\nu}\ \bigcup_{i-i'\in\{0,1,2\}}
\left(\{E^{i-i'}O_\mu\}\cup E^{-(i-i')}(\{O_\mu\})\right).
\label{s6n:color-center-choice}
\end{equation}
If the color has already been used, reuse its center.
In either case, set \(O_\gamma=O_\nu\) for the current label
of color \(\nu\). Continue in the same way.
At each step, at most \(n_0-1\) colors and finitely many centers are excluded.
Thus the choices are possible. Since only finitely many colors are used,
\eqref{s6n:color-center-choice} gives a common \(\varepsilon_0>0\) such that,
for the colors used,
\begin{equation}
\bigl\|(O_\mu-E^{i-i'}O_\nu)\bmod\ZZ^2\bigr\|_\infty
\ge\varepsilon_0,\qquad 1\le\mu\ne\nu\le n_0,\quad i-i'\in\{0,1,2\}.
\label{s6n:color-center-separation}
\end{equation}
Here the remainder is taken in \([-\tfrac12,\tfrac12)^2\), so
\[
\bigl\|X\bmod\ZZ^2\bigr\|_\infty
=\min_{n\in\ZZ^2}\|X-n\|_\infty\le\|X\|_\infty.
\]
Choose
\begin{equation}
0<r_0<\min\left\{\frac1{4\sqrt2},\,
\frac{\varepsilon_0}{2\sqrt2(1+\Lambda_t^2)}\right\}.
\label{s6n:common-rectangle-radius}
\end{equation}

\needspace{10\baselineskip}
We now verify \eqref{s6n:separated-supports} by contradiction.
A point in both auxiliary supports would force the center distance to be at least
\(\varepsilon_0\) and strictly less than \(\varepsilon_0\).
Take two distinct labels with intersecting physical supports.
Their colors \(\nu,\mu\) are different. After interchanging the labels
if necessary, take \(i=i(k)\), \(i'=i(k')\) with
\(i-i'\in\{0,1,2\}\).
If a point \(Y\) belonged to both inverse images in
\eqref{s6n:separated-supports}, there would be coordinates satisfying
\begin{equation}
\begin{aligned}
E^iY&=O_\nu+l_tv_t+l_rv_r\pmod{\ZZ^2},\\
E^{i'}Y&=O_\mu+l_t'v_t+l_r'v_r\pmod{\ZZ^2},\\
&\hspace{5mm}|l_r|,|l_t|,|l_r'|,|l_t'|<r_0.
\end{aligned}
\label{s6n:hypothetical-common-phase}
\end{equation}
Applying \(E^{i-i'}\) to the second line and subtracting the first gives
\[
O_\nu-E^{i-i'}O_\mu
\equiv(\Lambda_r^{i-i'}l_r'-l_r)v_r
+(\Lambda_t^{i-i'}l_t'-l_t)v_t\pmod{\ZZ^2}.
\]
Using \(\|av_r+bv_t\|_\infty\le\sqrt2\max\{|a|,|b|\}\), we obtain
\begin{equation}
\begin{aligned}
\varepsilon_0
&\le\bigl\|(O_\nu-E^{i-i'}O_\mu)\bmod\ZZ^2\bigr\|_\infty\\
&\le\bigl\|(\Lambda_r^{i-i'}l_r'-l_r)v_r
+(\Lambda_t^{i-i'}l_t'-l_t)v_t\bigr\|_\infty\\
&\le\sqrt2\max\bigl\{|\Lambda_r^{i-i'}l_r'-l_r|,
|\Lambda_t^{i-i'}l_t'-l_t|\bigr\}\\
&<2\sqrt2(1+\Lambda_t^2)r_0<\varepsilon_0,
\end{aligned}
\label{s6n:rectangle-separation-contradiction}
\end{equation}
a contradiction.

Finally, the physical realizations are inverse images of the disjoint
sets in \eqref{s6n:separated-supports} under the same map
\((t,r)\mapsto Y(t,r)\). They are therefore disjoint as well.
\end{proof}

\medskip
\noindent\textbf{The auxiliary cutoff.}
Choose a nonnegative, nonzero function
\(\Phi_r\in C_c^\infty((-r_0,r_0))\) and a function
\(\Phi_t\in C_c^\infty((1/6,5/6))\) such that
\(0\le\Phi_t\le1\) and \(\Phi_t=1\) on \([3/10,7/10]\).
For each label \(\gamma\), define the \emph{auxiliary cutoff} by
\begin{equation}\label{s6n:aux-cutoff}
 \Phi_\gamma(Y_i)
 =\Phi_r(l_r(Y_i))
  \Phi_t\!\left(\frac{l_t(Y_i)+r_0}{2r_0}\right),
\end{equation}
using the eigenvector coordinates centered at \(O_\gamma\), and extend
it by zero to \(\TT^2\). The temporal factor in \(\Phi_\gamma\)
is precisely \(\Phi_t(\tau_\gamma/\tau_{\rm max})\).

Each connected component of the physical support of the cutoff lies in one of
these rectangles, so the widths in \eqref{s6n:aux-widths} bound its extent
in \(r\) and \(t\).
Table~\ref{tab:s6-prepulse-directional-scales} collects the directional localization
scales. 

\begin{center}
\small
\setlength{\tabcolsep}{4pt}
\renewcommand{\arraystretch}{1.35}
\begin{longtable}{@{}p{2.5cm}p{2.5cm}p{2.85cm}p{3.25cm}p{3cm}@{}}
\caption{Anisotropic physical cutoff widths. The auxiliary cutoff row gives upper bounds. ``None'' means no additional cutoff.}
\label{tab:s6-prepulse-directional-scales}\\
\toprule
Layer & \(\theta\) & \(z\) & \(r\) & \(t\)\\
\midrule
Slow box
 & None
 & \(2^{-(1-\delta)k}k^{-6}\)
 & \(2^{-k}k^{-6}\)
 & \(2^{-2k}k^{-6}\)\\[8pt]
Auxiliary cutoff
 & None
 & None
 & \(2^{-(1+10^{-5}\delta)k}k^{2\rho_0}\)
 & \(2^{-(2+\delta)k}k^2\)\\[8pt]
\bottomrule
\end{longtable}
\end{center}

\subsection{The pulse ansatz and principal amplitude equation}
\label{s6n:ansatz}

\medskip
\noindent\textbf{Pre-pulses.}
Fix a geometric label \(\gamma=(k,\alpha,\iota)\),
\(\iota\in\{+,-\}\). In this subsection we consider the fundamental
real cosine carrier and the unforced principal amplitude equation.
We suppress the label \(\gamma\), writing
\(\mathbf h,\pi,\varphi,\chi,w^{\mathrm{pre}},w\), etc., instead of
\(\mathbf h_\gamma,\pi_\gamma,\varphi_\gamma\), etc.

Put \(i=i(k)\) and \(N=\lceil2^{\delta k/2}\rceil\), and define
the full cut-off function by
\[
 \chi(t,r,z,Y_i)
 =(\zeta_k\zeta_{k,\alpha})(t,r,z)\Phi_\gamma(Y_i).
\]
We consider an axisymmetric real vector amplitude function of the form
\[
 \mathbf h=\mathbf h(\hat R,\hat Z,\hat T,\tau(Y_i)).
\]
Its components in the orthonormal cylindrical frame
\((e_r,e_\theta,e_z)\) are independent of \(\theta\).
Using the physical-evaluation convention of Section~\ref{s6n:auxiliary},
define the pre-pulse by
\begin{equation}\label{s6n:pre-pulse}
 w^{\mathrm{pre}}=2^{(1+\delta)k}\chi\mathbf h\cos(N\varphi).
\end{equation}
Here \(w^{\mathrm{pre}}\) may be viewed as an extended field or, after
substituting \(Y_i=Y_i(t,r)\), as its physical realization.
The phase function \(\varphi\) will be discussed
below. Both \(\mathbf h\) and \(\varphi\) are allowed to depend on
the label \((k,\alpha,\pm)\).

We can choose three real frequency parameters \(p_\theta,p_r,p_z\)
for the given label \(\gamma\), with
\(Np_\theta\in\ZZ\setminus\{0\}\). They are constant for this label.
The phase is needed only on the auxiliary rectangle
\(\mathcal R_\gamma\), where we define it by
\begin{equation}\label{s6n:phase}
 \begin{aligned}
 \varphi
 &=p_\theta\theta+2^{\delta k}p_z\hat Z
   +p_r\sqrt{2\hat R}
   -\tau p_\theta\hat F_B-\tau p_z\hat u_B^z.
 \end{aligned}
\end{equation}
The first three terms prescribe the wavevector at \(\tau=0\).
The last two terms are designed to compensate the leading angular
and axial transport by the fixed background \(u_B\). Using
\[
 2^{\delta k}p_z\hat Z=2^kp_z z,\qquad
 p_r\sqrt{2\hat R}=2^kp_r r,
\]
the normalized spatial phase gradient is
\begin{equation}\label{s6n:wavevector}
\begin{aligned}
 \mathbf n_\varphi
 &:=2^{-k}\nabla\varphi
 =2^{-k}(D_r\varphi,\ r^{-1}\partial_\theta\varphi,
                               \partial_z\varphi)\\
 &=\begin{pmatrix}
 p_r-\tau\sqrt{2\hat R}
       (p_\theta\partial_{\hat R}\hat F_B
        +p_z\partial_{\hat R}\hat u_B^z)\\
 p_\theta/\sqrt{2\hat R}\\
 p_z-2^{-\delta k}\tau
       (p_\theta\partial_{\hat Z}\hat F_B
        +p_z\partial_{\hat Z}\hat u_B^z)
 \end{pmatrix}.
\end{aligned}
\end{equation}
The components \(n_\varphi^r,n_\varphi^\theta,n_\varphi^z\) are in
the orthonormal cylindrical frame. The cosine carrier is
\[
 \cos(N\varphi)
 =\frac{e^{\mathrm iN\varphi}+e^{-\mathrm iN\varphi}}2,
\]
with spatial frequencies \(\pm N2^k\mathbf n_\varphi\).

The time derivative and tangential transport give
\[
\begin{aligned}
 D_t\varphi
 &=\partial_t\varphi+\Lambda_t^i v_t\cdot\nabla_{Y_i}\varphi\\
 &=-2^{(2+\delta)k}(p_\theta\hat F_B+p_z\hat u_B^z)\\
 &\quad+2^{2k}\tau
       (p_\theta\partial_{\hat T}\hat F_B
        +p_z\partial_{\hat T}\hat u_B^z),\\[3pt]
 u_B^{\theta,z}\cdot\nabla_{\theta,z}\varphi
 &=2^{(2+\delta)k}
   (\hat F_B\partial_\theta+\hat u_B^z2^{-k}\partial_z)\varphi\\
 &=2^{(2+\delta)k}(p_\theta\hat F_B+p_z\hat u_B^z)\\
 &\quad-2^{2k}\tau\hat u_B^z
       (p_\theta\partial_{\hat Z}\hat F_B
        +p_z\partial_{\hat Z}\hat u_B^z).
\end{aligned}
\]
Here, on functions before physical evaluation,
\[
 u_B\cdot\nabla
 =u_B^rD_r+\frac{u_B^\theta}{r}\partial_\theta+u_B^z\partial_z.
\]
Consequently, the \emph{physical transport defect of the phase} is
given by
\begin{equation}\label{s6n:transport-defect}
\begin{aligned}
 &2^{-(2+\delta)k}(D_t+u_B\cdot\nabla)\varphi=\mathfrak d,\\
 \mathfrak d
 &:=-2^{-\delta k}\tau\Bigl[
    p_\theta(-\partial_{\hat T}
              +\hat u_B^z\partial_{\hat Z})\hat F_B
   +p_z(-\partial_{\hat T}
              +\hat u_B^z\partial_{\hat Z})\hat u_B^z
       \Bigr]+\hat u_B^r n_\varphi^r.
\end{aligned}
\end{equation}
For bounded frequency parameters \(p_\theta,p_r,p_z\), the estimates
on the background solution imply
\[
 |\mathfrak d|\le C2^{-\delta k}(1+\tau)
       \le C2^{-\delta k}k^2,\qquad
 0\le\tau\le\tau_{\rm max}\asymp k^2.
\]
The same bound holds after any fixed composition of the normalized
physical derivatives \(2^{-(2+\delta)k}D_t\), \(2^{-k}D_r\), and
\(2^{-k}\partial_z\), on each physical pulse rectangle.

\medskip
\noindent\textbf{The principal amplitude equation.}
Adding a velocity \(w\) and pressure \(p_w\) to the background produces
the linear terms
\[
 \mathsf L_B(w,p_w)
 :=\partial_tw+(u_B\cdot\nabla)w+(w\cdot\nabla)u_B
                       -\Delta w+\nabla p_w.
\]
For scalar cylindrical components, write
\[
 \Delta
 =D_r^2+\frac1rD_r+\frac1{r^2}\partial_\theta^2+\partial_z^2.
\]
Since the background is axisymmetric and independent of the auxiliary
variable, the three components are
\begin{equation}\label{s6n:linearized-components}
\begin{aligned}
\mathsf L_B(w,p_w)^r
={}&\left(D_t+u_B^rD_r+\frac{u_B^\theta}{r}\partial_\theta
                       +u_B^z\partial_z\right)w^r\\
&+w^r\partial_ru_B^r+w^z\partial_zu_B^r
  -\frac{2u_B^\theta}{r}w^\theta\\
&-\left(\Delta-\frac1{r^2}\right)w^r
  +\frac2{r^2}\partial_\theta w^\theta+D_rp_w,\\[3pt]
\mathsf L_B(w,p_w)^\theta
={}&\left(D_t+u_B^rD_r+\frac{u_B^\theta}{r}\partial_\theta
                       +u_B^z\partial_z\right)w^\theta\\
&+\left(\partial_ru_B^\theta+\frac{u_B^\theta}{r}\right)w^r
  +w^z\partial_zu_B^\theta+\frac{u_B^r}{r}w^\theta\\
&-\left(\Delta-\frac1{r^2}\right)w^\theta
  -\frac2{r^2}\partial_\theta w^r+\frac1r\partial_\theta p_w,\\[3pt]
\mathsf L_B(w,p_w)^z
={}&\left(D_t+u_B^rD_r+\frac{u_B^\theta}{r}\partial_\theta
                       +u_B^z\partial_z\right)w^z\\
&+w^r\partial_ru_B^z+w^z\partial_zu_B^z
  -\Delta w^z+\partial_zp_w.
\end{aligned}
\end{equation}
Here \(D_t,D_r\) represent physical differentiation through the chain
rule of Section~\ref{s6n:auxiliary}.

Corresponding to \eqref{s6n:pre-pulse}, take the pressure ansatz
\begin{equation}\label{s6n:pressure-ansatz}
 p_w=2^{2(1+\delta)k}\chi
                     \operatorname{Re}(\pi e^{\mathrm iN\varphi}).
\end{equation}
The velocity amplitude \(\mathbf h\) considered earlier is real,
while \(\pi\) may be complex here. We suppress the cutoffs.
Keeping the following principal terms and dividing by \(2^{(3+2\delta)k}\), we obtain
\begin{equation}\label{s6n:principal-ode}
\begin{gathered}
 \underbrace{\partial_\tau\mathbf h}_{\text{fast time derivative}}
 +\underbrace{\mathsf K_B\mathbf h}_{\text{background shear}}
 +\underbrace{d\mathbf h}_{\text{viscosity}}
 +\underbrace{\mathrm iN\mathbf n_\varphi\pi}_{\text{pressure}}=0,\\
 \underbrace{\mathbf n_\varphi\cdot\mathbf h=0}_{\text{incompressibility}},\qquad
 d=2^{-\delta k}N^2|\mathbf n_\varphi|^2,\qquad
 \mathsf K_B=
 \begin{pmatrix}
 0&-2\hat F_B&0\\
 2\hat F_B+\hat{\mathcal S}_B^\theta&0&0\\
 \hat{\mathcal S}_B^z&0&0
 \end{pmatrix}.
\end{gathered}
\end{equation}
More explicitly, the background shear coefficients satisfy
\[
\begin{aligned}
 2^{-(2+\delta)k}\left(-\frac{2u_B^\theta}{r}\right)
   &=-2\hat F_B,\\
 2^{-(2+\delta)k}\left(\partial_ru_B^\theta
                           +\frac{u_B^\theta}{r}\right)
   &=\sqrt{2\hat{R}}\,\partial_{\hat{R}}\hat u_B^\theta
           +\frac{\hat u_B^\theta}{\sqrt{2\hat{R}}}
     =2\hat F_B+\hat{\mathcal S}_B^\theta,\\
 2^{-(2+\delta)k}\partial_ru_B^z
   &=\sqrt{2\hat{R}}\,\partial_{\hat{R}}\hat u_B^z
     =\hat{\mathcal S}_B^z.
\end{aligned}
\]
The fast time derivative gives
\[
 2^{(1+\delta)k}(D_t\tau)\partial_\tau\mathbf h
 =2^{(3+2\delta)k}\partial_\tau\mathbf h.
\]
The principal viscous damping gives
\[
 2^{(1+\delta)k}N^2|\nabla\varphi|^2\mathbf h
 =2^{(3+2\delta)k}d\mathbf h.
\]
The principal pressure term gives
\[
 2^{2(1+\delta)k}\mathrm iN\nabla\varphi\pi
 =2^{(3+2\delta)k}\mathrm iN\mathbf n_\varphi\pi.
\]
Thus the sum of these principal terms is
\[
 2^{(3+2\delta)k}\chi\operatorname{Re}\!\left[
  (\partial_\tau\mathbf h+\mathsf K_B\mathbf h+d\mathbf h+\mathrm iN\mathbf n_\varphi\pi)
                      e^{\mathrm iN\varphi}\right]=0.
\]

The phase transport defect, for example, contributes
\[
 2^{(3+2\delta)k}\chi\operatorname{Re}\!\left(
      \mathrm iN\mathfrak d\,\mathbf h e^{\mathrm iN\varphi}\right)
\]
to the residual. Cutoff derivatives, transport derivatives in
\((\hat{R},\hat{Z},\hat{T})\), the smaller background terms
and the other viscous terms remain to be estimated.

\medskip
\noindent\textbf{Eliminating pressure.}
The principal amplitude equation~\eqref{s6n:principal-ode} determines
the linear dynamics of \(\mathbf h\). Differentiating the
orthogonality constraint gives
\[
 \mathbf n_\varphi\cdot\partial_\tau\mathbf h=-(\partial_\tau\mathbf n_\varphi)\cdot\mathbf h.
\]
Taking the scalar product of \eqref{s6n:principal-ode} with
\(\mathbf n_\varphi\) then yields
\[
 -(\partial_\tau\mathbf n_\varphi)\cdot\mathbf h+\mathbf n_\varphi\cdot\mathsf K_B\mathbf h
                       +\mathrm iN|\mathbf n_\varphi|^2\pi=0.
\]
Solving for the pressure gives
\begin{equation}\label{s6n:projected-ode}
\begin{aligned}
 \partial_\tau\mathbf h&=(\mathsf A-d)\mathbf h, \qquad \mathsf A=-\mathsf K_B+
 \frac{\mathbf n_\varphi(\mathbf n_\varphi^{\mathsf T}\mathsf K_B
                         -(\partial_\tau\mathbf n_\varphi)^{\mathsf T})}{|\mathbf n_\varphi|^2},\\
 \pi&=\frac{\mathrm i}{N}
       \frac{\mathbf n_\varphi\cdot\mathsf K_B\mathbf h
                          -(\partial_\tau\mathbf n_\varphi)\cdot\mathbf h}{|\mathbf n_\varphi|^2}.
\end{aligned}
\end{equation}
In particular,
\[
 \mathbf n_\varphi^{\mathsf T}\mathsf A=-(\partial_\tau\mathbf n_\varphi)^{\mathsf T},\qquad
 \partial_\tau(\mathbf n_\varphi\cdot\mathbf h)=-d(\mathbf n_\varphi\cdot\mathbf h).
\]
Thus transverse initial data remain transverse. The matrix
\(\mathsf A\) is real and preserves real amplitudes; \(\pi\) is then
purely imaginary. Taking the scalar product with \(\mathbf h\) gives
\begin{equation}\label{s6n:energy}
 \frac12\partial_\tau|\mathbf h|^2
 =-h^r\hat{\mathcal S}_B\cdot(h^\theta,h^z)-d|\mathbf h|^2.
\end{equation}
Shear can supply energy, while viscosity removes it.
Section~\ref{sec:pulse-dynamics} chooses the frequency parameters so that shear
dominates first and viscosity dominates later, and constructs the
unforced amplitude from suitable nonzero initial data. We only
formulate the unforced principal equation here; its solution and
estimates are given in Section~\ref{s6n:unforced-pulse}.

\medskip
\noindent\textbf{Pulses and exact divergence-free completion.}
The orthogonality constraint removes only the leading divergence
of the pre-pulse. For \(\mathbf n_\varphi\ne0\) on the support of
\(\chi\mathbf h\), define a physical vector potential \(A\)
and a pulse \(w\) by
\begin{equation}\label{s6n:curl-completion}
 \begin{aligned}
 A&=-\frac{2^{\delta k}}{N}\chi
       \frac{\mathbf n_\varphi\times\mathbf h}{|\mathbf n_\varphi|^2}\sin(N\varphi),\\
 w&=\nabla\times A.
 \end{aligned}
\end{equation}
The curl here is the ordinary physical curl, taken after substitution
\(Y_i=Y_i(t,r)\); the potential is extended by zero outside its support.
Since
\[
 \mathbf n_\varphi\times(\mathbf n_\varphi\times\mathbf h)=-|\mathbf n_\varphi|^2\mathbf h,
\]
the derivative falling on the sine gives exactly \(w^{\mathrm{pre}}\).
The remaining terms form \(w^{\mathrm{corr}}\), and
\[
 \begin{gathered}
 w=w^{\mathrm{pre}}+w^{\mathrm{corr}},\qquad \nabla\cdot w=0,\\
 \operatorname{supp}w\subset\operatorname{supp}(\chi\mathbf h)
                           \subset\operatorname{supp}\chi.
 \end{gathered}
\]
These are physical supports, so the curl preserves the support
separation of distinct geometric labels \(\gamma\). The pre-pulse determines the principal amplitude
dynamics and background stress matching. For the primary pulses used later, we apply \eqref{s6n:curl-completion} with \(\mathbf h\) replaced by the matched amplitude \(\mathbf h_{\mathrm{match},\gamma}=2^{-\delta k/2}g_{\mathrm{match},\gamma}\mathbf h_\gamma\) defined in \eqref{s6n:amplitude-data}. The matching factor $g_{\mathrm{match}, \gamma}$ is inserted in the vector potential before taking the curl, so its spatial derivatives are included in the curl correction. The resulting complete pulses are the velocity building blocks added to the background and used in the full residual.

\paragraph{General harmonics with linear sources.}
\label{s6n:general-harmonics}

Fix a geometric label \(\gamma\), continue to suppress it from the
notation, and use
the phase \(\varphi\), wavevector \(\mathbf n_\varphi\), fast time \(\tau\)
and coefficients \(d,\mathsf K_B\) as above.
For a nonzero harmonic index \(m\in\ZZ\setminus\{0\}\), consider a prescribed linear source
\[
 2^{(3+2\delta)k}\mathbf f_m e^{\mathrm imN\varphi}.
\]
Here \(\mathbf f_m\) is in fast-time units: after multiplying the force
by \(2^{-(3+\delta)k}\), its coefficient is \(2^{\delta k}\mathbf f_m\).
The vector coefficient $\mathbf f_m$ may be complex and need not be orthogonal
to \(\mathbf n_\varphi\). We assume that the source coefficient
\(\mathbf f_m\) is already supported in
\(\operatorname{supp}(\zeta_k\zeta_{k,\alpha})\) in the normalized
variables and in \(\operatorname{supp}\Phi_r\) in \(l_r\). Thus no
further spatial cutoff functions are needed here. The coefficient is
independent of \(\theta\) in the cylindrical frame and may depend on
\((\hat R,\hat Z,\hat T,l_r,\tau)\).
Only the temporal cutoff is to be applied. Set
\[
 \psi(\tau)=\Phi_t(\tau/\tau_{\rm max}),\qquad
 \begin{aligned}
 w_m^{\mathrm{pre}}&=2^{(1+\delta)k}\psi\mathbf h_m
                                           e^{\mathrm imN\varphi},\\
 p_{w,m}&=2^{2(1+\delta)k}\psi\pi_m e^{\mathrm imN\varphi}.
 \end{aligned}
\]
Cancellation of the prescribed source requires the forced principal
equation
\begin{equation}\label{s6n:general-principal-ode}
 \begin{aligned}
 \partial_\tau\mathbf h_m+\mathsf K_B\mathbf h_m+m^2d\mathbf h_m
                   +\mathrm imN\mathbf n_\varphi\pi_m&=-\mathbf f_m,\\
 \mathbf n_\varphi\cdot\mathbf h_m&=0.
 \end{aligned}
\end{equation}
The derivative \(\partial_\tau\) keeps \((\hat R,\hat Z,\hat T,l_r)\) fixed.
The factors \(m^2\) and \(m\) come from two carrier derivatives in
the Laplacian and one in the pressure gradient, respectively.
After the temporal cutoff, the principal equation gives
\[
 (\partial_\tau+\mathsf K_B+m^2d)(\psi\mathbf h_m)
 +\mathrm i mN\mathbf n_\varphi\,\psi\pi_m+\mathbf f_m
 =\psi'\mathbf h_m+(1-\psi)\mathbf f_m.
\]
For the sources used in Section~\ref{iter:nonzero-modes}, the Gaussian
bounds for \(\mathbf h_m\) and \(\mathbf f_m\) make both terms on the
right flat on the cutoff region near the time endpoints, as shown
in the estimate following \eqref{iter:nonzero-flat-error}.
We only formulate the forced principal equation here; its solution and
estimates, with zero initial data, are
given in Section~\ref{s7n:forced-pulse}.

As in the unforced case, differentiate the orthogonality constraint
and take the scalar product of \eqref{s6n:general-principal-ode}
with \(\mathbf n_\varphi\). With \(\mathsf A\) as in
\eqref{s6n:projected-ode}, eliminating the pressure gives
\begin{equation}\label{s6n:general-projected-ode}
 \begin{aligned}
 \pi_m&=\frac{\mathrm i}{mN}
 \frac{\mathbf n_\varphi\cdot\mathsf K_B\mathbf h_m
       -(\partial_\tau\mathbf n_\varphi)\cdot\mathbf h_m
       +\mathbf n_\varphi\cdot\mathbf f_m}{|\mathbf n_\varphi|^2},\\
 \partial_\tau\mathbf h_m
   &=(\mathsf A-m^2dI_3)\mathbf h_m-\Pi\mathbf f_m,\\
 \Pi&=I_3-\frac{\mathbf n_\varphi\mathbf n_\varphi^{\mathsf T}}
                      {|\mathbf n_\varphi|^2}.
 \end{aligned}
\end{equation}

The corresponding potential and exactly divergence-free pulse are
\[
 A_m=\frac{\mathrm i2^{\delta k}}{mN}\psi
        \frac{\mathbf n_\varphi\times\mathbf h_m}{|\mathbf n_\varphi|^2}
                             e^{\mathrm imN\varphi},\qquad
 w_m=\nabla\times A_m.
\]
As in \eqref{s6n:curl-completion}, the carrier derivative gives
\(w_m^{\mathrm{pre}}\); the remaining terms are the curl correction.
Finite sums with
\(\mathbf h_{-m}=\overline{\mathbf h_m}\),
\(\pi_{-m}=\overline{\pi_m}\), and
\(\mathbf f_{-m}=\overline{\mathbf f_m}\) give real velocity, pressure,
and source fields.

\subsection{Common auxiliary variables and compatibility}
\label{sec:pulse-common}

For each reference band
\(\mathcal B_k=\{(t,r,z):1/2<2^k\lambda(t,z)<2\}\), \(k\ge k_0\), set
\[
\begin{aligned}
 \mathcal K_k&=\{\ell\ge k_0:
 \mathcal B_k\cap\operatorname{supp}[\zeta_\ell(\lambda)]\ne\varnothing\},\\
 i_0&=i_0(k)=\min_{\ell\in\mathcal K_k}i(\ell).
\end{aligned}
\]
The support condition in \eqref{s6n:physical-cutoffs} gives
\(\mathcal K_k\subseteq\{k-1,k,k+1\}\). The argument for
\eqref{s6n:overlapping-covering-indices} therefore gives
\(0\le i(\ell)-i_0\le2\) for every \(\ell\in\mathcal K_k\).
Use \(Y_{i_0}\in\TT^2\) as the independent common variable throughout
\(\mathcal B_k\), with physical evaluation
\[
 Y_{i_0}(t,r)=\Lambda_t^{i_0}t\,v_t
             +\Lambda_r^{i_0}r^{d_0}v_r\pmod{\ZZ^2}
\]
under the convention of Section~\ref{s6n:auxiliary}.
For a label \(\gamma=(\ell,\alpha,\iota)\) with \(\ell\in\mathcal K_k\),
 note that
\begin{equation}
 Y_{i(\ell)}=E^{i(\ell)-i_0}Y_{i_0},\qquad 0\le i(\ell)-i_0\le2.
\label{eq:common-auxiliary-variable}
\end{equation}
This relation holds on the level of the independent auxiliary variables
and also, after substitution, for their physical realizations.
For neighboring labels \(\gamma=(\ell,\alpha,\iota)\) and
\(\gamma'=(\ell',\alpha',\iota')\) with \(\ell,\ell'\in\mathcal K_k\),
the separation in \eqref{s6n:separated-supports} reads
\[
 E^{-(i(\ell)-i_0)}(\mathcal R_\gamma)
 \cap E^{-(i(\ell')-i_0)}(\mathcal R_{\gamma'})=\varnothing.
\]

\medskip
\noindent\textbf{Compatibility on overlaps.}
Consider $|k-k'| \le 1$ so that $\mathcal B_k \cap \mathcal B_{k'} \neq \emptyset$. Let $(\hat R,\hat Z,\hat T)$ be the normalized coordinates with respect to scale $k$. On \(\mathcal B_k \cap \mathcal B_{k'}\), consider two functions \(f(\hat R,\hat Z,\hat T,Y_{i_0(k)})\) and \(g(\hat R,\hat Z,\hat T,Y_{i_0(k')})\), which are smooth in all variables
and periodic in $Y_{i_0(k)}$ and $Y_{i_0(k')}$ respectively. 
\emph{Compatibility} means that
\[
 f(\hat R,\hat Z,\hat T,E^{i_0(k)}Y)
 =g(\hat R,\hat Z,\hat T,E^{i_0(k')}Y),
 \qquad Y\in\TT^2,
\]
on \(\mathcal B_k\cap\mathcal B_{k'}\). Their physical evaluations
then agree on \(\mathcal B_k\cap\mathcal B_{k'}\).
All torus averages and directional inverses below keep
\((\hat R,\hat Z,\hat T)\) fixed.

\begin{lem}[Averages and directional inverses on overlaps]
\label{lem:common-compatibility}
Let \(f,g\) be compatible smooth periodic functions as above.
At every fixed
\((\hat R,\hat Z,\hat T)\), their torus averages with normalized
Haar measure agree:
\[
 \langle f\rangle_{Y_{i_0(k)}}
 =\langle g\rangle_{Y_{i_0(k')}}.
\]
If \((f_1,g_1)\) and \((f_2,g_2)\) are compatible pairs in these
same auxiliary variables, then \((f_1f_2,g_1g_2)\) is compatible.

Suppose \(f,g\) have zero torus mean. For \(\nu\in\{r,t\}\), set
\[
\begin{aligned}
 \widetilde f_\nu
 &=\Lambda_\nu^{-i_0(k)}
       (v_\nu\cdot\nabla_{Y_{i_0(k)}})^{-1}f,\\
 \widetilde g_\nu
 &=\Lambda_\nu^{-i_0(k')}
       (v_\nu\cdot\nabla_{Y_{i_0(k')}})^{-1}g.
\end{aligned}
\]
For each \(\nu\), the functions \(\widetilde f_\nu\) and
\(\widetilde g_\nu\) are compatible. The inverses act on the
respective auxiliary variables and are the zero-mean inverses of
Lemma~\ref{s6n:directional-inverse}.
\end{lem}

\begin{proof}
Since \(E\) is a nonsingular integer matrix, every nonconstant
Fourier mode of \(h(E^jY)\) has zero average. Thus, for every smooth
periodic \(h\),
\begin{equation}\label{eq:common-torus-average}
 \int_{\TT^2}h(E^jY)\,dY=\int_{\TT^2}h(Y)\,dY,
 \qquad j\in\ZZ_{\ge0}.
\end{equation}
Applying this identity to the compatibility relation for \(f,g\)
proves equality of their averages. Multiplying the compatibility
relations for the two pairs proves compatibility of their products.

Suppress the fixed normalized coordinates. The chain rule gives
\[
\begin{aligned}
 (v_\nu\cdot\nabla_Y)
       \bigl[\widetilde f_\nu(E^{i_0(k)}Y)\bigr]
 &=f(E^{i_0(k)}Y),\\
 (v_\nu\cdot\nabla_Y)
       \bigl[\widetilde g_\nu(E^{i_0(k')}Y)\bigr]
 &=g(E^{i_0(k')}Y).
\end{aligned}
\]
The two right-hand sides agree, and both pulled-back solutions have
zero mean by \eqref{eq:common-torus-average}. Uniqueness in
Lemma~\ref{s6n:directional-inverse} proves compatibility.
\end{proof}

In Section~\ref{init:auxiliary-correction}, the time-direction inverse
is applied (under local common auxiliary variables) to the fields
\(\langle f^a\rangle_\theta-\langle f^a\rangle_{\theta,Y_{i_0}}\),
\(a=\theta,z\). The lemma ensures that the resulting fields
\(A_{0,\mathrm{aux}}\) and \(B_{0,\mathrm{aux}}\) agree on overlaps.

\section{Principal dynamics and stress matching}
\label{sec:pulse-dynamics}
\label{sec:leading-pulse-construction}

We first choose the frequency parameters of the phase function in
\eqref{s6n:phase}, construct the unforced pre-pulse, and solve the
forced principal equation with zero initial data. We then use the
unforced pre-pulse amplitudes to match the complete background stress and to
realize  stress increments in arbitrary directions.
The positive matching is used to construct the primary pulse in
Section~\ref{sec:initialization-ppt}; the forced inverse and linearized
matching are used for the pulse corrections in Section~\ref{sec:iteration-ppt}.

In the frequency and amplitude calculations of
Subsections~\ref{s7n:frequencies}--\ref{s7n:forced-pulse},
we fix the band number \(k\) and a slow box \(\alpha\) such that
\(\operatorname{supp}(\zeta_k\zeta_{k,\alpha})\) meets the active annulus
\(R_a<R<R_h\), while retaining the sign \(\iota=\pm\) in the label
\(\gamma=(k,\alpha,\iota)\). Put \(i=i(k)\) and use the local
auxiliary variable \(Y_i=(Y_{i,1},Y_{i,2})\), independently of
\((\hat R,\hat Z,\hat T)\).
We suppress the dependence on \(k,\alpha\), writing, for example,
\(\tau_\pm=\tau_{(k,\alpha,\pm)}\). The derivatives \(\partial_{\tau_\pm}\)
hold \((\hat R,\hat Z,\hat T,l_r)\) fixed.
For \(k\ge k_0\), the estimates in these three subsections hold
throughout the corresponding slow box, with \(|l_r|<r_0\) and
\(0\le\tau_\pm\le\tau_{\rm max}\).
Here \((l_r,\tau_\pm)\) are the local coordinates associated with
the auxiliary rectangles \(\mathcal R_{(k,\alpha,\pm)}\)
in \eqref{s6n:rectangles}--\eqref{s6n:pulse-time}.

\subsection{Frequency choice and wavevector bounds}
\label{s7n:frequencies}

Recall the phase ansatz \eqref{s6n:phase} from Section~\ref{s6n:ansatz},
which makes the material derivative of the phase small.
For each slow box and each sign, it depends on three real parameters
\(p_{\theta,\pm},p_{r,\pm},p_{z,\pm}\):
\begin{equation}\label{s7n:phase-ansatz}
\begin{gathered}
 \varphi_\pm=p_{\theta,\pm}\theta
       +2^{\delta k}p_{z,\pm}\hat Z
       +p_{r,\pm}\sqrt{2\hat R}
       -\tau_\pm(p_{\theta,\pm}\hat F_B+p_{z,\pm}\hat u_B^z),\\
 Np_{\theta,\pm}\in\ZZ\setminus\{0\}.
\end{gathered}
\end{equation}
For bounded parameters, the transport cancellation in
\eqref{s6n:transport-defect} gives
\[
 2^{-(2+\delta)k}(D_t+u_B\cdot\nabla)\varphi_\pm
       =O\bigl(2^{-\delta k}(1+|\tau_\pm|)\bigr).
\]
The two signs label two pulse families whose averaged flux directions
span a cone. This cone determines the quantitative cone condition imposed
in Part~I: a stress outside it cannot be canceled by a nonnegative
combination of these two pre-pulse fluxes.

On each physical pulse rectangle, as defined in
Section~\ref{s6n:localizations}, recall from \eqref{s6n:pulse-time} that
\[
 \tau_\pm=\tau_{(k,\alpha,\pm)}(Y_i)=2^{(2+\delta)k}\Lambda_t^{-i}(l_t+r_0).
\]
By Proposition~\ref{s6n:separation}, the physical pulse
rectangles for distinct neighboring labels are disjoint. In
particular, \(\tau_+\) and \(\tau_-\) are used on disjoint physical
regions. Recall from \eqref{s6n:wavevector} that
\[
\begin{aligned}
 \mathbf n_{\varphi,\pm}
 &:=2^{-k}\nabla\varphi_\pm\\
 &=2^{-k}(D_r\varphi_\pm,\ r^{-1}\partial_\theta\varphi_\pm,
                                      \partial_z\varphi_\pm)\\
 &=\begin{pmatrix}
 p_{r,\pm}-\tau_\pm\sqrt{2\hat R}
       (p_{\theta,\pm}\partial_{\hat R}\hat F_B
        +p_{z,\pm}\partial_{\hat R}\hat u_B^z)\\
 p_{\theta,\pm}/\sqrt{2\hat R}\\
 p_{z,\pm}-2^{-\delta k}\tau_\pm
       (p_{\theta,\pm}\partial_{\hat Z}\hat F_B
        +p_{z,\pm}\partial_{\hat Z}\hat u_B^z)
 \end{pmatrix}.
\end{aligned}
\]

The fields \(\hat F_B\) and \(\hat u_B^z\) in
\(\mathbf n_{\varphi,\pm}\) are not frozen: they vary pointwisely within
the slow box. The three parameters
\(p_{\theta,\pm},p_{r,\pm},p_{z,\pm}\) must, however, be constant
for each \((k,\alpha)\) and each sign \(\pm\), with
\(Np_{\theta,\pm}\) being a nonzero integer.
We therefore choose them using frozen background values. Fix any
reference point \((\hat R_0,\hat Z_0,\hat T_0)\) in this slow box,
within the domain of the background coefficients.

\paragraph{Frozen quantities and the frozen cone.}
We use subscript \(0\) to denote evaluation at the reference point
\((\hat R_0,\hat Z_0,\hat T_0)\). Set
\[
\begin{gathered}
 \hat F_{B,0}=\hat F_B(\hat R_0,\hat Z_0,\hat T_0),\qquad
 \hat{\mathcal S}_{B,0}=\hat{\mathcal S}_B(\hat R_0,\hat Z_0,\hat T_0),\\
 \mathbf s_{B,0}=\frac{\hat{\mathcal S}_{B,0}}{|\hat{\mathcal S}_{B,0}|},\qquad
 \mathbf s_{B,0}^\perp=(-s_{B,0}^z,s_{B,0}^\theta),\qquad
 \kappa_{B,0}=-\frac{|\hat{\mathcal S}_{B,0}|^2}
                    {\hat F_{B,0}\hat{\mathcal S}_{B,0}^\theta}.
\end{gathered}
\]
These are constants associated with the given box \((k,\alpha)\).
The uniform derivative bounds in \eqref{intro:background-velocity-bounds}
and the positive lower bounds for \(F_B\) and
\(-\mathcal S_B^{\theta}\) in Section~\ref{intro:flux-cone}
allow these lower bounds to persist on a slightly enlarged radial annulus
\[
 R_a'\le R\le R_h',\qquad
 0<R_a'<R_a<R_h<R_h'<\infty,
\]
where \(R_a',R_h'\) are fixed independently of \(k,\alpha\).
Together with the \(O(k^{-6})\) diameter of the slow box defined by
\eqref{s6n:slow-box-cutoffs}, these bounds give
\[
 |\hat F_B-\hat F_{B,0}|
 +|\hat{\mathcal S}_B-\hat{\mathcal S}_{B,0}|\le Ck^{-6},
 \qquad
 \hat F_{B,0}\asymp-\hat{\mathcal S}_{B,0}^\theta\asymp1.
\]
Write
\[
 \kappa_B=-\frac{|\hat{\mathcal S}_B|^2}
                   {\hat F_B\hat{\mathcal S}_B^\theta}.
\]
On the active stress annulus, the uniform cone condition  gives, wherever \(\hat{\mathsf T}_B\ne0\),
\begin{equation}\label{s7n:uniform-cone-margins}
\begin{gathered}
 \kappa_B\ge2+c_{\rm cone},\qquad
 -\hat{\mathsf T}_B\cdot
       \frac{\hat{\mathcal S}_B}{|\hat{\mathcal S}_B|}
       \ge c_{\rm cone}|\hat{\mathsf T}_B|,\\
 \sqrt{\frac{\kappa_B-2}{2}}
 \left|\frac{\hat{\mathsf T}_B\cdot\hat{\mathcal S}_B^\perp}
                  {\hat{\mathsf T}_B\cdot\hat{\mathcal S}_B}\right|
       \le1-c_{\rm cone},
\end{gathered}
\end{equation}
as already mentioned in Section~\ref{s6n:slow}. Here \(c_{\rm cone}>0\) is independent of \(k,\alpha\) and the position within
the box, and is chosen small enough to work in all three inequalities
of \eqref{s7n:uniform-cone-margins}.
To quantify the effect of freezing, the preceding bounds imply
\[
 |\kappa_B-\kappa_{B,0}|
 +\left|\frac{\hat{\mathcal S}_B}{|\hat{\mathcal S}_B|}
                     -\mathbf s_{B,0}\right|\le Ck^{-6}.
\]
In particular, for sufficiently large \(k_0\),
\[
 -\hat{\mathsf T}_B\cdot\mathbf s_{B,0}
       \ge\frac{c_{\rm cone}}2|\hat{\mathsf T}_B|,\qquad
 \kappa_{B,0}\ge2+\frac{c_{\rm cone}}2.
\]
Thus, wherever \(\hat{\mathsf T}_B\ne0\), both denominators below
are bounded away from zero after division by the stress magnitude,
and
\[
\begin{aligned}
 &\left|
 \sqrt{\frac{\kappa_{B,0}-2}{2}}
   \frac{\hat{\mathsf T}_B\cdot\mathbf s_{B,0}^\perp}
        {\hat{\mathsf T}_B\cdot\mathbf s_{B,0}}
 -\sqrt{\frac{\kappa_B-2}{2}}
   \frac{\hat{\mathsf T}_B\cdot\hat{\mathcal S}_B^\perp}
        {\hat{\mathsf T}_B\cdot\hat{\mathcal S}_B}
 \right|\le Ck^{-6}.
\end{aligned}
\]
The second term inside the absolute value has magnitude at most
\(1-c_{\rm cone}\) by the last inequality in \eqref{s7n:uniform-cone-margins}.
Choosing \(k_0\) large  so that \(Ck^{-6}\le c_{\rm cone}/2\),
we obtain the explicit bound
\[
 \sqrt{\frac{\kappa_{B,0}-2}{2}}
 \left|\frac{\hat{\mathsf T}_B\cdot\mathbf s_{B,0}^\perp}
                 {\hat{\mathsf T}_B\cdot\mathbf s_{B,0}}\right|
       \le1-c_{\rm cone}+Ck^{-6}\le1-\frac{c_{\rm cone}}2.
\]
These estimates are uniform up to the radial edges
\(R=R_a\) and \(R=R_h\), where the stress tends to zero. Choose one constant \(M_{\rm cone}>0\),
independent of the box, such that
\[
 1-\frac{c_{\rm cone}}4<\frac{M_{\rm cone}}{\sqrt{1+M_{\rm cone}^2}}<1.
\]
Then, wherever \(\hat{\mathsf T}_B\ne0\) in the box,
\begin{equation}\label{s6n:cone-aperture}
 \hat{\mathsf T}_B\cdot\hat{\mathcal S}_{B,0}<0,\qquad
 \sqrt{\frac{\kappa_{B,0}-2}{2}}
 \left|\frac{\hat{\mathsf T}_B\cdot\mathbf s_{B,0}^\perp}
                 {\hat{\mathsf T}_B\cdot\mathbf s_{B,0}}\right|
 <\frac{M_{\rm cone}}{\sqrt{1+M_{\rm cone}^2}}<1.
\end{equation}
Moreover, these two inequalities hold with uniform margins. This is the frozen cone condition in each box, satisfied by the
varying stress direction.

\paragraph{Slow radial drift.}
By the definition of \(\hat{\mathcal S}_B\), the wavevector is exactly
\begin{equation}\label{s7n:wavevector-drift}
\begin{aligned}
 \mathbf n_{\varphi,\pm}(\tau_\pm)
 ={}&\left(p_{r,\pm},\frac{p_{\theta,\pm}}{\sqrt{2\hat R}},p_{z,\pm}\right)
 -\tau_\pm\left[
 \left(\frac{p_{\theta,\pm}}{\sqrt{2\hat R}},p_{z,\pm}\right)
       \cdot\hat{\mathcal S}_B\right]e_r\\
 &-2^{-\delta k}\tau_\pm
   (p_{\theta,\pm}\partial_{\hat Z}\hat F_B
        +p_{z,\pm}\partial_{\hat Z}\hat u_B^z)e_z.
\end{aligned}
\end{equation}
Recall from Proposition~\ref{s6n:time-prop} that the fast pulse-time
interval \(0\le\tau_\pm\le\tau_{\rm max}\) has length
\(\tau_{\rm max}\asymp k^2\).
For bounded frequency parameters, the last term in
\eqref{s7n:wavevector-drift} is \(O(2^{-\delta k}\tau_{\rm max})\).
The principal part of \(\mathbf n_{\varphi,\pm}(\tau_\pm)\) is a
constant vector plus a radial component linear in \(\tau_\pm\);
its slope is minus the inner product of the \(\theta,z\) frequency
pair with the shear.

Recall the principal equation \eqref{s6n:principal-ode} and its
projected form \eqref{s6n:projected-ode}:
\[
\begin{gathered}
 \partial_{\tau_\pm}\mathbf h_\pm+\mathsf K_B\mathbf h_\pm+d_\pm\mathbf h_\pm
                  +\mathrm iN\mathbf n_{\varphi,\pm}\pi_\pm=0,
 \qquad \mathbf n_{\varphi,\pm}\cdot\mathbf h_\pm=0,\\
 \partial_{\tau_\pm}\mathbf h_\pm=(\mathsf A_\pm-d_\pm I_3)\mathbf h_\pm,\qquad
 d_\pm=2^{-\delta k}N^2|\mathbf n_{\varphi,\pm}|^2,\qquad
 N=\lceil2^{\delta k/2}\rceil,
\end{gathered}
\]
where
\[
\begin{aligned}
 \mathsf A_\pm&=-\mathsf K_B+
 \frac{\mathbf n_{\varphi,\pm}
       (\mathbf n_{\varphi,\pm}^{\mathsf T}\mathsf K_B
                         -(\partial_{\tau_\pm}\mathbf n_{\varphi,\pm})^{\mathsf T})}
      {|\mathbf n_{\varphi,\pm}|^2},\\
 \mathsf K_B&=
 \begin{pmatrix}
 0&-2\hat F_B&0\\
 2\hat F_B+\hat{\mathcal S}_B^\theta&0&0\\
 \hat{\mathcal S}_B^z&0&0
 \end{pmatrix}.
\end{aligned}
\]
We seek frequency parameters for which
\[
 |\mathbf n_{\varphi,\pm}|\asymp1,\qquad
 |(\partial_{\tau_\pm}\mathbf n_{\varphi,\pm})|\le C\tau_{\rm max}^{-1}=O(k^{-2}).
\]
This would allow an adiabatic separation of the instantaneous spectral
states in the analysis of the next subsection (see \cite{Wasow1965}
for mathematical background on asymptotic decoupling in linear systems
with slowly varying coefficients and separated eigenvalues).

Formula \eqref{s7n:wavevector-drift} suggests imposing
\[
 \left|\left(\frac{p_{\theta,\pm}}{\sqrt{2\hat R_0}},p_{z,\pm}\right)
              \cdot\hat{\mathcal S}_{B,0}\right|
       =O(\tau_{\rm max}^{-1}).
\]
Thus the \(\theta,z\) frequency pair should be almost perpendicular
to the frozen shear. Exact perpendicularity would remove the desired
radial drift. Instead, we add a shear component of size
\(\tau_{\rm max}^{-1}\), with opposite signs for the two sign families
\(\iota=\pm\). This produces an order-one change in the radial
component of the wavevector over the full pulse-time interval.

\paragraph{The frequency choices.}
Introduce the positive reference growth parameter
\begin{equation}\label{s6n:growth-parameters}
\begin{gathered}
 \Lambda_{B,0}^2=-2\hat F_{B,0}s_{B,0}^\theta
       (2\hat F_{B,0}s_{B,0}^\theta+|\hat{\mathcal S}_{B,0}|),\\
 \sqrt{\frac{\kappa_{B,0}-2}{2}}
       =-\frac{\Lambda_{B,0}}{2\hat F_{B,0}s_{B,0}^\theta},\qquad
 \Lambda_{B,0}\asymp\sqrt{\frac{\kappa_{B,0}-2}{2}}\asymp1,
\end{gathered}
\end{equation}
and the frequency scale \(B_p>0\), which will be used to
define the three frequency parameters:
\begin{equation}\label{s6n:frequency-choice}
 B_p^2=\frac{\Lambda_{B,0}}
       {2^{-\delta k}N^2(1+M_{\rm cone}^2)^{3/2}},\qquad B_p\asymp1.
\end{equation}
This choice balances the reference growth and viscous damping at the
pulse-time midpoint, as verified in Section~\ref{s6n:unforced-pulse}.
First set
\begin{equation}\label{s7n:unrounded-frequencies}
 \left(\frac{\widetilde p_{\theta,\pm}}{\sqrt{2\hat R_0}},p_{z,\pm}\right)
 =B_p\left(\mathbf s_{B,0}^\perp
       \mp\frac{M_{\rm cone}}{\tau_{\rm max}}
              \frac{\hat{\mathcal S}_{B,0}}{|\hat{\mathcal S}_{B,0}|^2}\right).
\end{equation}
The first summand \(\mathbf s_{B,0}^\perp\) is perpendicular to the
frozen shear $\hat{\mathcal{S}}_{B, 0}$ in the \((\theta,z)\)-plane. The second gives the desired
radial drift with the corresponding sign \(\pm\). To ensure angular
periodicity, we round the angular frequency parameter and
choose \(p_{\theta,\pm}\) such that
\begin{equation}\label{s7n:frequency-rounding}
 Np_{\theta,\pm}\in\ZZ\setminus\{0\},\qquad
 |p_{\theta,\pm}-\widetilde p_{\theta,\pm}|\le N^{-1}.
\end{equation}
Finally, choose the radial frequency parameter
\begin{equation}\label{s6n:pr-pm}
 p_{r,\pm}=\pm B_pM_{\rm cone}/2.
\end{equation}
The three frequency parameters are now fixed for each box
\((k,\alpha)\) and sign \(\iota=\pm\). Introduce the linear functions
\begin{equation}\label{s7n:radial-paths}
\begin{gathered}
 s_\pm(\tau_\pm)=\pm M_{\rm cone}
       \left(\frac12+\frac{\tau_\pm}{\tau_{\rm max}}\right),\qquad
 \partial_{\tau_\pm}s_\pm=\pm\frac{M_{\rm cone}}{\tau_{\rm max}},\\
 s_\pm(0)=\pm\frac{M_{\rm cone}}2,\qquad
 s_\pm(\tau_{\rm max})=\pm\frac{3M_{\rm cone}}2,\qquad
 \frac{M_{\rm cone}}2\le|s_\pm(\tau_\pm)|\le\frac{3M_{\rm cone}}2.
\end{gathered}
\end{equation}
At the respective midpoints \(\tau_\pm=\tau_{\rm max}/2\),
\(s_\pm=\pm M_{\rm cone}\).
The frequency parameters above are chosen to achieve the wavevector
approximation and slow-variation estimates in the following lemma.

For partial derivatives in the five independent variables
\((\hat R,\hat Z,\hat T,Y_{i,1},Y_{i,2})\), we use multi-indices
\(\beta\in\ZZ_{\ge0}^5\). The label, reference point and chosen
frequency parameters are held fixed when taking such derivatives. These partial derivatives, as well as \(\partial_{\tau_\pm}\), keep \(\theta\) fixed, hence do not act on the cylindrical frame.

\begin{lem}[Wavevector bounds]
\label{s7n:phase-bounds}
With the frequency parameters chosen in
\eqref{s6n:frequency-choice}--\eqref{s6n:pr-pm}, the following bounds
hold throughout the slow box corresponding to \((k,\alpha)\), for
\(0\le\tau_\pm\le\tau_{\rm max}\):
\begin{equation}\label{s6n:normal-approximation}
\begin{gathered}
 |\mathbf n_{\varphi,\pm}-B_p(s_\pm(\tau_\pm),\mathbf s_{B,0}^\perp)|
       +|(\partial_{\tau_\pm}\mathbf n_{\varphi,\pm})|\le Ck^{-2},\\
 0<c\le|(n_{\varphi,\pm}^\theta,n_{\varphi,\pm}^z)|\le|\mathbf n_{\varphi,\pm}|\le C.
\end{gathered}
\end{equation}
In addition,
\[
 |(n_{\varphi,\pm}^\theta,n_{\varphi,\pm}^z)\cdot\hat{\mathcal S}_B|
       +|\partial_{\tau_\pm}\mathsf A_\pm|+|\partial_{\tau_\pm}d_\pm|\le Ck^{-2}.
\]
For every \(\beta\), there exist constants
\(C_{\beta}>0\) and \(b_{\beta}\ge0\) such that
\begin{equation}\label{s7n:phase-derivatives}
 \begin{aligned}
 &|\partial_{\hat R,\hat Z,\hat T,Y_{i,1},Y_{i,2}}^{\beta}\mathbf n_{\varphi,\pm}|
 +\left|\partial_{\hat R,\hat Z,\hat T,Y_{i,1},Y_{i,2}}^{\beta}
       \frac{\mathbf n_{\varphi,\pm}}{|\mathbf n_{\varphi,\pm}|^2}\right|\\
 &\quad+|\partial_{\hat R,\hat Z,\hat T,Y_{i,1},Y_{i,2}}^{\beta}\mathsf A_\pm|
 +|\partial_{\hat R,\hat Z,\hat T,Y_{i,1},Y_{i,2}}^{\beta}d_\pm|
 \le C_{\beta} k^{b_{\beta}}.
 \end{aligned}
\end{equation}
The constants are uniform in the band index \(k\), the slow box
label \(\alpha\), and the signs \(\pm\) in \(\gamma\).
\end{lem}

\begin{proof}
\emph{Step 1: the drift coefficient throughout the slow box.}
Since \(1\le2^{-\delta k}N^2\le4\),
\eqref{s6n:growth-parameters} and the fixed choice of \(M_{\rm cone}\)
give \(c\le B_p\le C\). Also
\(\hat R_0\asymp1\) and \(|\hat{\mathcal S}_{B,0}|\asymp1\), so
\[
 |\widetilde p_{\theta,\pm}|+|p_{\theta,\pm}|
             +|p_{z,\pm}|+|p_{r,\pm}|\le C.
\]
Taking the inner product of \eqref{s7n:unrounded-frequencies} with
\(\hat{\mathcal S}_{B,0}\) gives, exactly,
\[
\begin{aligned}
 \left(\frac{\widetilde p_{\theta,\pm}}{\sqrt{2\hat R_0}},p_{z,\pm}\right)
       \cdot\hat{\mathcal S}_{B,0}
 &=\mp\frac{B_pM_{\rm cone}}{\tau_{\rm max}},\\
 p_{r,\pm}-\tau_\pm
 \left(\frac{\widetilde p_{\theta,\pm}}{\sqrt{2\hat R_0}},p_{z,\pm}\right)
       \cdot\hat{\mathcal S}_{B,0}
 &=B_ps_\pm(\tau_\pm).
\end{aligned}
\]
Moreover,
\[
 \left|\left(\frac{\widetilde p_{\theta,\pm}}{\sqrt{2\hat R_0}},p_{z,\pm}\right)
                 -B_p\mathbf s_{B,0}^\perp\right|
 =\frac{B_pM_{\rm cone}}{\tau_{\rm max}|\hat{\mathcal S}_{B,0}|}
 \le C\tau_{\rm max}^{-1}.
\]

On the active annulus, \(\hat R\) is bounded away from zero. The
slow-box diameter and the background derivative bounds imply
\[
 \left|\frac1{\sqrt{2\hat R}}-\frac1{\sqrt{2\hat R_0}}\right|
       +|\hat{\mathcal S}_B-\hat{\mathcal S}_{B,0}|
       \le Ck^{-6}.
\]
To separate the rounding error from the spatial variation, write
\[
\begin{aligned}
 &\left(\frac{p_{\theta,\pm}}{\sqrt{2\hat R}},p_{z,\pm}\right)
                \cdot\hat{\mathcal S}_B
 -\left(\frac{\widetilde p_{\theta,\pm}}{\sqrt{2\hat R_0}},p_{z,\pm}\right)
                \cdot\hat{\mathcal S}_{B,0}\\
 &\quad=\frac{p_{\theta,\pm}-\widetilde p_{\theta,\pm}}
                    {\sqrt{2\hat R_0}}\hat{\mathcal S}_{B,0}^\theta
 +p_{\theta,\pm}\left(\frac1{\sqrt{2\hat R}}
                         -\frac1{\sqrt{2\hat R_0}}\right)
                           \hat{\mathcal S}_{B,0}^\theta\\
 &\qquad\quad+
 \left(\frac{p_{\theta,\pm}}{\sqrt{2\hat R}},p_{z,\pm}\right)
               \cdot(\hat{\mathcal S}_B-\hat{\mathcal S}_{B,0}).
\end{aligned}
\]
The three terms are bounded by \(CN^{-1}\), \(Ck^{-6}\), and
\(Ck^{-6}\), respectively. Hence
\[
 \left|\left(\frac{p_{\theta,\pm}}{\sqrt{2\hat R}},p_{z,\pm}\right)
              \cdot\hat{\mathcal S}_B
          \pm\frac{B_pM_{\rm cone}}{\tau_{\rm max}}\right|
       \le C(k^{-6}+N^{-1}).
\]
Multiplying by \(\tau_\pm\le\tau_{\rm max}\) and using the identity
above for \(B_ps_\pm(\tau_\pm)\) yields
\[
 |n_{\varphi,\pm}^r-B_ps_\pm(\tau_\pm)|
       \le C\tau_{\rm max}(k^{-6}+N^{-1}).
\]

\emph{Step 2: the \(\theta,z\) components and the time derivative.}
The \(\theta,z\) components differ from the prescribed vector
\((\widetilde p_{\theta,\pm}/\sqrt{2\hat R_0},p_{z,\pm})\) by
\[
\begin{aligned}
 &(n_{\varphi,\pm}^\theta,n_{\varphi,\pm}^z)
 -\left(\frac{\widetilde p_{\theta,\pm}}{\sqrt{2\hat R_0}},p_{z,\pm}\right)\\
 &\quad=\left(
 \frac{p_{\theta,\pm}-\widetilde p_{\theta,\pm}}{\sqrt{2\hat R_0}}
 +p_{\theta,\pm}\left(\frac1{\sqrt{2\hat R}}-\frac1{\sqrt{2\hat R_0}}\right),\right.\\
 &\hspace{42mm}\left.
 -2^{-\delta k}\tau_\pm
  (p_{\theta,\pm}\partial_{\hat Z}\hat F_B
                       +p_{z,\pm}\partial_{\hat Z}\hat u_B^z)\right).
\end{aligned}
\]
The background derivatives \(\partial_{\hat Z}\hat F_B\) and
\(\partial_{\hat Z}\hat u_B^z\) are uniformly bounded. Combining this
identity with Step 1 gives
\[
 |(n_{\varphi,\pm}^\theta,n_{\varphi,\pm}^z)-B_p\mathbf s_{B,0}^\perp|
 \le C(\tau_{\rm max}^{-1}+k^{-6}+N^{-1}
                                  +2^{-\delta k}\tau_{\rm max}).
\]
Since the background coefficients and frequency parameters are independent
of \(\tau_\pm\), differentiating \eqref{s7n:wavevector-drift} gives
\[
 (\partial_{\tau_\pm}\mathbf n_{\varphi,\pm})
 =-\left[\left(\frac{p_{\theta,\pm}}{\sqrt{2\hat R}},p_{z,\pm}\right)
             \cdot\hat{\mathcal S}_B\right]e_r
 -2^{-\delta k}(p_{\theta,\pm}\partial_{\hat Z}\hat F_B
                   +p_{z,\pm}\partial_{\hat Z}\hat u_B^z)e_z.
\]
Step 1 therefore implies
\[
 |(\partial_{\tau_\pm}\mathbf n_{\varphi,\pm})|
 \le C(\tau_{\rm max}^{-1}+k^{-6}+N^{-1}+2^{-\delta k}).
\]
Now \(\tau_{\rm max}\asymp k^2\) and \(N^{-1}\le2^{-\delta k/2}\).
For sufficiently large \(k_0\),
\[
 \tau_{\rm max}k^{-6}\le Ck^{-4},\qquad
 \tau_{\rm max}N^{-1}+2^{-\delta k}\tau_{\rm max}\le Ck^{-2}.
\]
The preceding three bounds prove the first line of
\eqref{s6n:normal-approximation}. Since
\(|B_p\mathbf s_{B,0}^\perp|=B_p\ge c\), its \(\theta,z\) estimate
also gives
\[
 |(n_{\varphi,\pm}^\theta,n_{\varphi,\pm}^z)|\ge B_p-Ck^{-2}\ge c/2.
\]
The upper bound follows from \(B_p\le C\) and
\(|s_\pm|\le3M_{\rm cone}/2\). This proves the second line.
Finally, write
\[
\begin{aligned}
 (n_{\varphi,\pm}^\theta,n_{\varphi,\pm}^z)\cdot\hat{\mathcal S}_B
 ={}&[(n_{\varphi,\pm}^\theta,n_{\varphi,\pm}^z)-B_p\mathbf s_{B,0}^\perp]
                      \cdot\hat{\mathcal S}_B\\
 &+B_p\mathbf s_{B,0}^\perp
                 \cdot(\hat{\mathcal S}_B-\hat{\mathcal S}_{B,0}).
\end{aligned}
\]
The two terms are \(O(k^{-2})\) and \(O(k^{-6})\), proving the
near-perpendicularity estimate.

\emph{Step 3: derivatives of the principal equation coefficients.}
For the time derivatives, put \(\widetilde d_\pm=|\mathbf n_{\varphi,\pm}|^2\) in this
calculation. Then
\[
 \partial_{\tau_\pm}\widetilde d_\pm=2\mathbf n_{\varphi,\pm}\cdot(\partial_{\tau_\pm}\mathbf n_{\varphi,\pm})=O(k^{-2}),\qquad
 \partial_{\tau_\pm}d_\pm=2^{-\delta k}N^2\partial_{\tau_\pm}\widetilde d_\pm=O(k^{-2}).
\]
Using \(\partial_{\tau_\pm}\mathsf K_B=0\) and \(\partial_{\tau_\pm}^2\mathbf n_{\varphi,\pm}=0\),
\[
\begin{aligned}
 \partial_{\tau_\pm}\mathsf A_\pm
 ={}&\frac{(\partial_{\tau_\pm}\mathbf n_{\varphi,\pm})
        (\mathbf n_{\varphi,\pm}^{\mathsf T}\mathsf K_B
                    -(\partial_{\tau_\pm}\mathbf n_{\varphi,\pm})^{\mathsf T})
       +\mathbf n_{\varphi,\pm}(\partial_{\tau_\pm}\mathbf n_{\varphi,\pm})^{\mathsf T}\mathsf K_B}{\widetilde d_\pm}\\
 &-\frac{\partial_{\tau_\pm}\widetilde d_\pm}{\widetilde d_\pm^2}\,
       \mathbf n_{\varphi,\pm}(\mathbf n_{\varphi,\pm}^{\mathsf T}\mathsf K_B
                              -(\partial_{\tau_\pm}\mathbf n_{\varphi,\pm})^{\mathsf T})
       =O(k^{-2}).
\end{aligned}
\]
All denominators are uniformly bounded away from zero by Step 2.

For derivatives in these five variables, the label, reference point
and chosen frequency parameters are held fixed. The derivatives in the normalized variables of the
background and of \((2\hat R)^{-1/2}\) are uniformly bounded.
For a derivative \(\partial_{\hat R,\hat Z,\hat T}^{\beta}\), the
independent auxiliary variable is fixed, so \(\tau_\pm\) is held constant. Thus
\eqref{s7n:wavevector-drift} gives
\[
\begin{aligned}
 \partial_{\hat R,\hat Z,\hat T}^{\beta}\mathbf n_{\varphi,\pm}
 ={}&\partial_{\hat R,\hat Z,\hat T}^{\beta}\left(p_{r,\pm},
                \frac{p_{\theta,\pm}}{\sqrt{2\hat R}},p_{z,\pm}\right)\\
 &-\tau_\pm\partial_{\hat R,\hat Z,\hat T}^{\beta}\left[
 \left(\frac{p_{\theta,\pm}}{\sqrt{2\hat R}},p_{z,\pm}\right)
                  \cdot\hat{\mathcal S}_B\right]e_r\\
 &-2^{-\delta k}\tau_\pm\partial_{\hat R,\hat Z,\hat T}^{\beta}
      (p_{\theta,\pm}\partial_{\hat Z}\hat F_B
                  +p_{z,\pm}\partial_{\hat Z}\hat u_B^z)e_z.
\end{aligned}
\]
Each differentiated coefficient on the right is uniformly bounded,
so
\[
 |\partial_{\hat R,\hat Z,\hat T}^{\beta}\mathbf n_{\varphi,\pm}|
 \le C_{\beta}(1+|\tau_\pm|)
 \le C_{\beta}(1+\tau_{\rm max})\le C_{\beta}'k^2.
\]
The factor \(1+\tau_{\rm max}\) comes from the affine dependence on
\(\tau_\pm\); the same factor suffices at every fixed
order of differentiation in the normalized variables.
The fast time is affine in \(Y_i\), with
\(|\nabla_{Y_i}\tau_\pm|=O(k^2)\). Consequently, all fixed partial
derivatives in \((\hat R,\hat Z,\hat T,Y_{i,1},Y_{i,2})\) of
\(\mathbf n_{\varphi,\pm}\),
\((\partial_{\tau_\pm}\mathbf n_{\varphi,\pm})\), and \(\mathsf K_B\)
have polynomial bounds.
Differentiating \(\widetilde d_\pm \widetilde d_\pm^{-1}=1\), for \(\beta\ne0\), gives
\[
 \partial_{\hat R,\hat Z,\hat T,Y_{i,1},Y_{i,2}}^{\beta}(\widetilde d_\pm^{-1})
 =-\widetilde d_\pm^{-1}\sum_{0<\beta'\le \beta}\binom{\beta}{\beta'}
            (\partial_{\hat R,\hat Z,\hat T,Y_{i,1},Y_{i,2}}^{\beta'}\widetilde d_\pm)\,\partial_{\hat R,\hat Z,\hat T,Y_{i,1},Y_{i,2}}^{\beta-\beta'}(\widetilde d_\pm^{-1}).
\]
Starting with \(\widetilde d_\pm^{-1}\le c^{-2}\), induction on \(|\beta|\)
gives polynomial bounds for all these derivatives. The product rule
then gives the same bounds for
\(\mathbf n_{\varphi,\pm}/\widetilde d_\pm\), \(\mathsf A_\pm\), and \(d_\pm\), which are
precisely \eqref{s7n:phase-derivatives}.
\end{proof}

\subsection{The unforced pre-pulse and its Gaussian envelope}
\label{s6n:unforced-pulse}

We first express the principal equation in an orthonormal basis of the
constraint plane \(\mathbf n_{\varphi,\pm}^{\perp}\). The wavevector
varies at rate \(O(k^{-2})\), while the two principal eigenvalues have
a uniformly positive gap. We choose initial data near the growing
eigenvector and show that the ratio of the two amplitude coordinates
remains \(O(k^{-2})\)-close to the corresponding eigenvector ratio
throughout the pulse-time interval; see \eqref{s6n:growing-ratio}.

Write
\[
 \mathbf n_{\varphi,\pm}^{\theta,z}
       =(n_{\varphi,\pm}^{\theta},n_{\varphi,\pm}^{z}),
 \qquad
 \mathbf n_{\varphi,\pm}
       =B_p(s_\pm,\mathbf s_{B,0}^{\perp})+O(k^{-2}).
\]
We introduce two orthonormal basis vectors
\(\mathbf n_{\varphi,\pm,x}^{\perp}\) and
\(\mathbf n_{\varphi,\pm,y}^{\perp}\), which are exactly orthogonal
to \(\mathbf n_{\varphi,\pm}\):
\[
\begin{aligned}
 \mathbf n_{\varphi,\pm,x}^{\perp}
 &=\frac{(0,(\mathbf n_{\varphi,\pm}^{\theta,z})^{\perp})}
              {|\mathbf n_{\varphi,\pm}^{\theta,z}|}
   =-(0,\mathbf s_{B,0})+O(k^{-2}),\\
 \mathbf n_{\varphi,\pm,y}^{\perp}
 &=\frac{(|\mathbf n_{\varphi,\pm}^{\theta,z}|^2,
               -n_{\varphi,\pm}^r\mathbf n_{\varphi,\pm}^{\theta,z})}
          {|\mathbf n_{\varphi,\pm}|\,
                         |\mathbf n_{\varphi,\pm}^{\theta,z}|}
   =\frac{(1,-s_\pm\mathbf s_{B,0}^{\perp})}
              {\sqrt{1+s_\pm^2}}+O(k^{-2}).
\end{aligned}
\]
Then, due to the constraint
\(\mathbf n_{\varphi,\pm}\cdot\mathbf h_\pm=0\), the amplitude vector
\(\mathbf h_\pm\) can be uniquely decomposed as
\begin{equation}\label{s6n:moving-frame}
\begin{gathered}
 \mathbf h_\pm=\bar x_\pm\mathbf n_{\varphi,\pm,x}^{\perp}
                  +\bar y_\pm\mathbf n_{\varphi,\pm,y}^{\perp},\\
 \bar x_\pm=\mathbf h_\pm\cdot\mathbf n_{\varphi,\pm,x}^{\perp},
 \qquad
 \bar y_\pm=\mathbf h_\pm\cdot\mathbf n_{\varphi,\pm,y}^{\perp}
       =\frac{|\mathbf n_{\varphi,\pm}|}
                       {|\mathbf n_{\varphi,\pm}^{\theta,z}|}h_\pm^r,
 \\
 |\mathbf h_\pm|^2=\bar x_\pm^2+\bar y_\pm^2.
\end{gathered}
\end{equation}

By a direct computation and Lemma~\ref{s7n:phase-bounds}, projecting
the principal equation \eqref{s6n:principal-ode} onto these two basis
vectors gives
\begin{equation}\label{s6n:frame-equation}
\begin{gathered}
 \partial_{\tau_\pm}\binom{\bar x_\pm}{\bar y_\pm}
 =\left[
 \frac1{\sqrt{1+s_\pm^2}}
 \begin{pmatrix}
 0&2\hat F_{B,0}s_{B,0}^{\theta}+|\hat{\mathcal S}_{B,0}|\\
 -2\hat F_{B,0}s_{B,0}^{\theta}&0
 \end{pmatrix}
 +\mathcal E_\pm-d_\pm I_2
 \right]\binom{\bar x_\pm}{\bar y_\pm},\\
 |\mathcal E_\pm|
 +\left|d_\pm-2^{-\delta k}N^2B_p^2(1+s_\pm^2)\right|
 \le Ck^{-2}.
\end{gathered}
\end{equation}
Here the error includes the changes in the background and basis
vectors. For every fixed \(\beta\in\ZZ_{\ge0}^5\), we have
\[
 |\partial_{\hat R,\hat Z,\hat T,Y_{i,1},Y_{i,2}}^\beta\mathcal E_\pm|
 +\left|\partial_{\hat R,\hat Z,\hat T,Y_{i,1},Y_{i,2}}^\beta\bigl(d_\pm
       -2^{-\delta k}N^2B_p^2(1+s_\pm^2)\bigr)\right|
 \le C_\beta k^{b_\beta},
\]
with constants \(C_\beta>0\) and \(b_\beta\ge0\) uniform in
\(k,\alpha\), both signs, and \(0\le\tau_\pm\le\tau_{\rm max}\).
Denote the leading matrix by
\begin{equation}\label{s6n:leading-matrix}
 \mathsf B_\pm
 :=\frac1{\sqrt{1+s_\pm^2}}
 \begin{pmatrix}
 0&2\hat F_{B,0}s_{B,0}^{\theta}+|\hat{\mathcal S}_{B,0}|\\
 -2\hat F_{B,0}s_{B,0}^{\theta}&0
 \end{pmatrix}.
\end{equation}
Its action on the two eigenvectors is
\[
\begin{aligned}
 \mathsf B_\pm\binom{1}{\sqrt{2/(\kappa_{B,0}-2)}}
 &=\frac{\Lambda_{B,0}}{\sqrt{1+s_\pm^2}}
                       \binom{1}{\sqrt{2/(\kappa_{B,0}-2)}},\\[6pt]
 \mathsf B_\pm\binom{1}{-\sqrt{2/(\kappa_{B,0}-2)}}
 &=-\frac{\Lambda_{B,0}}{\sqrt{1+s_\pm^2}}
                       \binom{1}{-\sqrt{2/(\kappa_{B,0}-2)}}.
\end{aligned}
\]
Note that we used the definition
\[
 \kappa_{B,0}
 =-\frac{|\hat{\mathcal S}_{B,0}|^2}
         {\hat F_{B,0}\hat{\mathcal S}_{B,0}^\theta} = -\frac{|\hat{\mathcal S}_{B,0}|}
         {\hat F_{B,0}s_{B,0}^\theta}
\]
in the eigenvector calculation above.
The scalar damping term \(-d_\pm I_2\) in
\eqref{s6n:frame-equation} shifts both eigenvalues by \(-d_\pm\) and
leaves these eigenvectors unchanged. We will choose positive initial
coordinates close to the growing eigenvector
\((1,\sqrt{2/(\kappa_{B,0}-2)})^{\mathsf T}\), so that the ratio satisfies
\begin{equation}\label{s6n:growing-ratio}
 \frac{\bar y_\pm}{\bar x_\pm}
       =\sqrt{\frac2{\kappa_{B,0}-2}}+O(k^{-2}),
 \qquad 0\le\tau_\pm\le\tau_{\rm max},
\end{equation}
uniformly throughout the pulse-time interval, as proved below.
By \eqref{s6n:moving-frame}, the components of \(\mathbf h_\pm\) are
given by
\[
\begin{aligned}
 h_\pm^r
 &=\frac{|\mathbf n_{\varphi,\pm}^{\theta,z}|}
              {|\mathbf n_{\varphi,\pm}|}\bar y_\pm,\\
 (h_\pm^\theta,h_\pm^z)
 &=\frac{\bar x_\pm(\mathbf n_{\varphi,\pm}^{\theta,z})^{\perp}}
              {|\mathbf n_{\varphi,\pm}^{\theta,z}|}
   -\frac{n_{\varphi,\pm}^r\bar y_\pm\mathbf n_{\varphi,\pm}^{\theta,z}}
              {|\mathbf n_{\varphi,\pm}|\,
                         |\mathbf n_{\varphi,\pm}^{\theta,z}|}.
\end{aligned}
\]
Together with \eqref{s6n:growing-ratio} and
Lemma~\ref{s7n:phase-bounds}, these identities show that the radial
quadratic flux \(h_\pm^r(h_\pm^\theta,h_\pm^z)\) has direction close to
the vectors \(\mathbf q_\pm(\tau_\pm)\) in the \((\theta,z)\)-plane;
more precisely,
\[
\begin{aligned}
 h_\pm^r(h_\pm^\theta,h_\pm^z)
 &=(h_\pm^r)^2\bigl[\mathbf q_\pm(\tau_\pm)+O(k^{-2})\bigr],\\
 \mathbf q_\pm(\tau_\pm)
 &:=-s_\pm\mathbf s_{B,0}^{\perp}
       -\sqrt{\frac{(\kappa_{B,0}-2)(1+s_\pm^2)}2}\,
                                             \mathbf s_{B,0}.
\end{aligned}
\]
At the fast time \emph{midpoints} \(\tau_\pm=\tau_{\rm max}/2\), we have
\(s_\pm=\pm M_{\rm cone}\). The cone spanned by these two midpoint
reference directions contains a nonzero stress \(\hat{\mathsf T}_B\)
in its interior precisely when
\[
 \hat{\mathsf T}_B\cdot\mathbf s_{B,0}<0,\qquad
 \sqrt{\frac{\kappa_{B,0}-2}2}
 \left|\frac{\hat{\mathsf T}_B\cdot\mathbf s_{B,0}^{\perp}}
                  {\hat{\mathsf T}_B\cdot\mathbf s_{B,0}}\right|
 <\frac{M_{\rm cone}}{\sqrt{1+M_{\rm cone}^2}}.
\]
This is the frozen cone condition \eqref{s6n:cone-aperture}:
to leading order, the growing eigenvector of \(\mathsf A_\pm\)
determines the accessible stress directions.
Figure~\ref{s7n:midpoint-cone-figure} illustrates the stress cone at the pulse-time midpoint.
\begin{center}
\begin{minipage}{\linewidth}
\makeatletter\def\@captype{figure}\makeatother
\centering
\begin{tikzpicture}[>=stealth,font=\small,x=1.25cm,y=1cm]
 \fill[blue!7] (0,0)--(-2.08,2.6)--(2.08,2.6)--cycle;
 \draw[->,draw=black!45,text=black!80] (-2.35,0)--(2.55,0)
   node[right] {\(\mathbf s_{B,0}^\perp\)};
 \draw[->,draw=black!45,text=black!80] (0,-.2)--(0,3)
   node[above] {\(-\mathbf s_{B,0}\)};
 \draw[blue!70!black,thick,->] (0,0)--(-1.6,2)
   node[above left] {\(\mathbf q_+\)};
 \draw[blue!70!black,thick,->] (0,0)--(1.6,2)
   node[above right] {\(\mathbf q_-\)};
 \draw[orange!85!black,very thick,->] (0,0)--(.42,1.7)
   node[right] {\(\hat{\mathsf T}_B\)};
 \fill (0,0) circle (1pt);
 \node[below left] at (0,0) {\(0\)};
\end{tikzpicture}
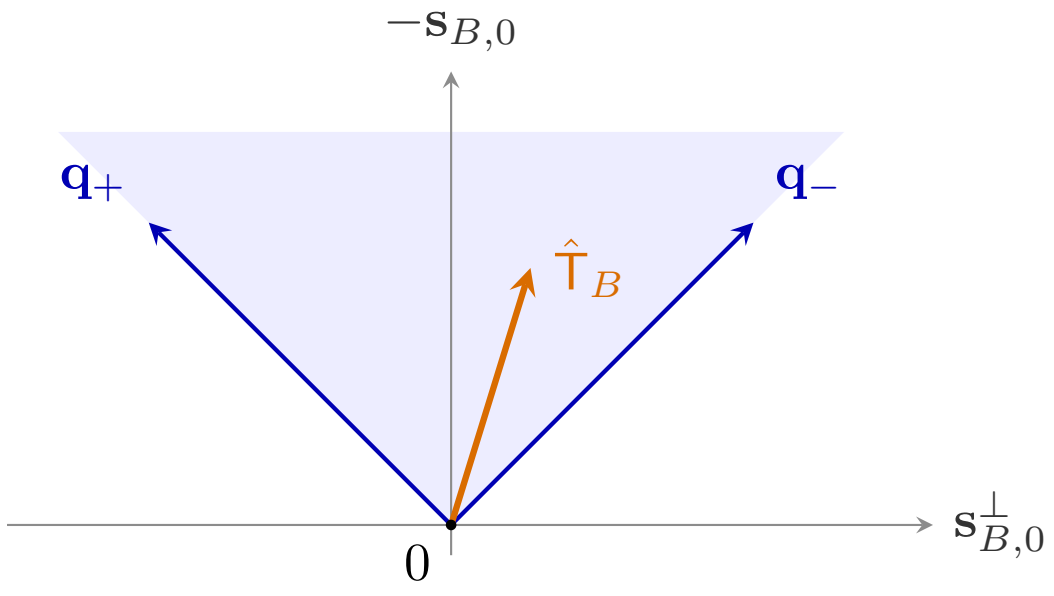
\caption{The midpoint reference directions
\(\mathbf q_\pm=\mathbf q_\pm(\tau_{\rm max}/2)\) enclose the
stress direction. The actual flux directions
differ by \(O(k^{-2})\). Directions are drawn schematically.}
\label{s7n:midpoint-cone-figure}
\end{minipage}
\end{center}

The choice of \(B_p\) in \eqref{s6n:frequency-choice} makes the
reference growth and damping balance at the midpoint
\(\tau_{\rm max}/2\). Indeed, the positive eigenvalue
\(\Lambda_{B,0}/\sqrt{1+s_\pm^2}\) of \(\mathsf B_\pm\) and the
damping coefficient \(d_\pm\) almost cancel at \(\tau_{\rm max}/2\):
\[
\begin{aligned}
 \left|\frac{\Lambda_{B,0}}{\sqrt{1+M_{\rm cone}^2}}
                 -d_\pm(\tau_{\rm max}/2)\right|
 &=\left|2^{-\delta k}N^2B_p^2(1+M_{\rm cone}^2)
                 -d_\pm(\tau_{\rm max}/2)\right|\\
 &\le Ck^{-2}.
\end{aligned}
\]
Define
\begin{equation}\label{s6n:gaussian-envelope}
\begin{aligned}
 a(\tau)&=\frac{\Lambda_{B,0}}{\sqrt{1+s_+(\tau)^2}}
       -\frac{\Lambda_{B,0}(1+s_+(\tau)^2)}{(1+M_{\rm cone}^2)^{3/2}},\\
 \mathcal G(\tau)&=\exp\left(\int_{\tau_{\rm max}/2}^{\tau}a(\tau')\,d\tau'\right).
\end{aligned}
\end{equation}
Since \(s_+(\tau)^2=s_-(\tau)^2\), the same functions \(a\) and
\(\mathcal G\) apply to both phase families, evaluated at their
respective fast times \(\tau_\pm\). In particular,
\[
 a(\tau)>0\quad(\tau<\tau_{\rm max}/2),\qquad
 a(\tau_{\rm max}/2)=0,\qquad
 a(\tau)<0\quad(\tau>\tau_{\rm max}/2).
\]
The reference envelope \(\mathcal G(\tau)\) first grows and then decays.
We choose the exact initial data
\begin{equation}\label{s6n:homogeneous-data}
\begin{aligned}
 \bar x_\pm(\tau=0)
 &=\sqrt{\frac{(\kappa_{B,0}-2)(1+s_\pm(0)^2)}2}\,\mathcal G(0),\\
 \bar y_\pm(\tau=0)
 &=\frac{|\mathbf n_{\varphi,\pm}(0)|}
               {|\mathbf n_{\varphi,\pm}^{\theta,z}(0)|}\,\mathcal G(0),
\end{aligned}
\end{equation}
and then solve the principal equation \eqref{s6n:principal-ode}
on \(\tau\in[0,\tau_{\rm max}]\).
These data give \(h_\pm^r(0)=\mathcal G(0)\), and the ratio
\eqref{s6n:growing-ratio} is satisfied at \(\tau_\pm=0\).
Figure~\ref{fig:s7n-growth-decay} shows the Gaussian bounds and scales
proved below.

\begin{center}
\begin{minipage}{\linewidth}
\makeatletter\def\@captype{figure}\makeatother
\centering
\begin{tikzpicture}[x=1.02cm,y=1.02cm,>=stealth,font=\small]
 \draw[->] (0,0)--(10.5,0) node[right] {\(\tau_\pm\)};
 \draw[->] (0,0)--(0,3.9) node[above] {\(|\mathbf h_\pm|\)};
 \draw[gray!60,densely dotted] (5,0)--(5,3.45);
 \draw[blue!65!black,dashed,thick,domain=0:10,samples=100]
       plot (\x,{3.2*exp(-.30*(\x-5)^2)});
 \draw[orange!85!black,very thick,domain=0:10,samples=100]
       plot (\x,{2.45*exp(-.56*(\x-5)^2)});
 \draw[blue!65!black,dashed,thick,domain=0:10,samples=100]
       plot (\x,{1.65*exp(-.90*(\x-5)^2)});
 \node[anchor=west,blue!65!black] at (6.45,3.35)
       {\(C e^{-c(\tau_\pm-\tau_{\rm max}/2)^2/\tau_{\rm max}}\)};
 \draw[blue!65!black,thin] (6.5,3.1)--(6.15,2.17);
 \node[anchor=west,blue!65!black] at (.35,1.65)
       {\(c e^{-C(\tau_\pm-\tau_{\rm max}/2)^2/\tau_{\rm max}}\)};
 \draw[blue!65!black,thin] (3.55,1.42)--(4.12,.82);
 \node[orange!85!black,above right] at (5,2.45) {\(\asymp1\)};
 \node at (2.35,.50) {growth};
 \node at (7.7,.50) {decay};
 \node[below] at (0,0) {\(0\)};
 \node[below] at (5,0) {\(\tau_{\rm max}/2\)};
 \node[below] at (10,0) {\(\tau_{\rm max}\)};
 \draw[<->] (4,-.65)--(6,-.65)
       node[midway,below] {\(\asymp\sqrt{\tau_{\rm max}}\asymp k\)};
 \draw[<->] (0,-1.35)--(10,-1.35)
       node[midway,below] {\(\tau_{\rm max}\asymp k^2\)};
 \node[align=center,font=\footnotesize] at (1.15,2.65)
       {endpoint size\\\(e^{-\Theta(k^2)}\)};
\end{tikzpicture}
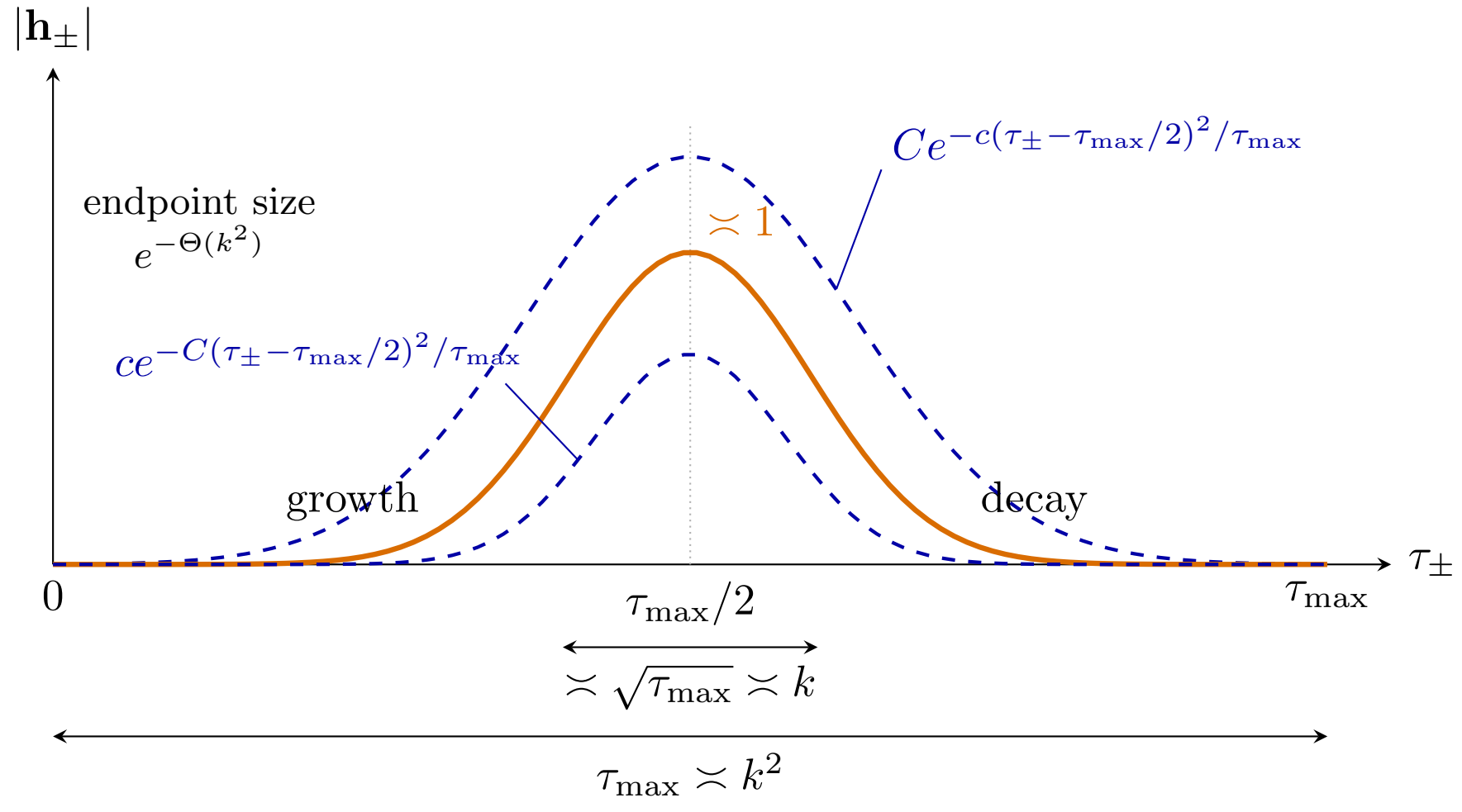
\caption{Schematic pre-pulse amplitude between Gaussian upper and lower
bounds. The central pulse-time width is \(O(k)\), corresponding to
physical time \(O(2^{-(2+\delta)k}k)\). The dashed curves include
constant prefactors; the schematic shape does not assert exact
monotonicity of the perturbed amplitude.}
\label{fig:s7n-growth-decay}
\end{minipage}
\end{center}

\begin{prop}[Unforced pre-pulses and Gaussian bounds]
\label{prop:homogeneous-pulse}
The data \eqref{s6n:homogeneous-data} determine a unique solution of
\eqref{s6n:principal-ode} on \(0\le\tau_\pm\le\tau_{\rm max}\), with
\[
 \mathbf h_\pm\in\RR^3,\qquad \pi_\pm\in\mathrm i\RR,\qquad
 \mathbf n_{\varphi,\pm}\cdot\mathbf h_\pm=0.
\]
Both \(\mathbf h_\pm\) and \(\pi_\pm\) are independent of the local
pulse coordinate \(l_r\). The envelope satisfies
\begin{equation}\label{s6n:gaussian-bounds}
 e^{-C(\tau_\pm-\tau_{\rm max}/2)^2/\tau_{\rm max}}
 \le\mathcal G(\tau_\pm)
 \le e^{-c(\tau_\pm-\tau_{\rm max}/2)^2/\tau_{\rm max}}\le1,
\end{equation}
and the size and direction of \(\mathbf h_\pm\) satisfy
\begin{equation}\label{s6n:homogeneous-direction}
\begin{gathered}
 c\mathcal G\le h_\pm^r\le|\mathbf h_\pm|\le C\mathcal G,\\
 \frac{(h_\pm^\theta,h_\pm^z)}{h_\pm^r}
 =-s_\pm\mathbf s_{B,0}^{\perp}
  -\sqrt{\frac{(\kappa_{B,0}-2)(1+s_\pm^2)}2}\,
                                  \mathbf s_{B,0}+O(k^{-2}).
\end{gathered}
\end{equation}
At the fast time midpoints \(\tau_\pm=\tau_{\rm max}/2\), the cone
generated by the vectors \(h_\pm^r(h_\pm^\theta,h_\pm^z)\)
contains \(\hat{\mathsf T}_B\) (with its direction uniformly separated from the boundary rays) wherever \(\hat{\mathsf T}_B\ne0\).
For every multi-index \(\beta\) in the normalized variables and \(j\ge0\),
\begin{equation}\label{s6n:homogeneous-derivatives}
\begin{aligned}
 |\partial_{\hat R,\hat Z,\hat T}^{\beta}
                        \partial_{\tau_\pm}^j\mathbf h_\pm|
 &\le C_{\beta,j}k^{b_{\beta,j}}\mathcal G(\tau_\pm),\\
 |\partial_{\hat R,\hat Z,\hat T}^{\beta}
                        \partial_{\tau_\pm}^j\pi_\pm|
 &\le C_{\beta,j}2^{-\delta k/2}k^{b_{\beta,j}}\mathcal G(\tau_\pm).
\end{aligned}
\end{equation}
\end{prop}

\begin{proof}
The exact basis in \eqref{s6n:moving-frame} and its derivatives are
uniformly bounded by powers of \(k\); its fast-time derivative is
\(O(k^{-2})\) by Lemma~\ref{s7n:phase-bounds}. In the frozen model,
the two off-diagonal entries of \(\mathsf B_\pm\) are
\[
 \frac{2\hat F_{B,0}s_{B,0}^{\theta}+|\hat{\mathcal S}_{B,0}|}
      {\sqrt{1+s_\pm^2}},\qquad
 \frac{-2\hat F_{B,0}s_{B,0}^{\theta}}{\sqrt{1+s_\pm^2}},
\]
and the diagonal entries vanish. The estimates
\[
 \mathsf K_B-\mathsf K_{B,0}=O(k^{-6}),\qquad
 \mathbf n_{\varphi,\pm}
          -B_p(s_\pm,\mathbf s_{B,0}^{\perp})=O(k^{-2}),\qquad
 |(\partial_{\tau_\pm}\mathbf n_{\varphi,\pm})|=O(k^{-2})
\]
therefore give \eqref{s6n:frame-equation}, including the polynomial
bounds for all fixed partial derivatives in
\((\hat R,\hat Z,\hat T,Y_{i,1},Y_{i,2})\). Here
\(\mathsf K_{B,0}=\mathsf K_B(\hat R_0,\hat Z_0,\hat T_0)\).
The projected equation \eqref{s6n:projected-ode} preserves the
constraint in \eqref{s6n:principal-ode}.

Direct differentiation of \eqref{s6n:gaussian-envelope} gives
\[
\begin{gathered}
 (\partial_\tau a)(\tau_\pm)
 =-\frac{\Lambda_{B,0}M_{\rm cone}|s_\pm|}{\tau_{\rm max}}
 \left[\frac1{(1+s_\pm^2)^{3/2}}
       +\frac2{(1+M_{\rm cone}^2)^{3/2}}\right],\\
 a(\tau_{\rm max}/2)=0,\qquad
 -\frac C{\tau_{\rm max}}\le\partial_\tau a\le-\frac c{\tau_{\rm max}}.
\end{gathered}
\]
Integration then proves \eqref{s6n:gaussian-bounds}.

As long as \(\bar x_\pm>0\), put
\[
 \rho_\pm=\frac{\bar y_\pm}{\bar x_\pm}
                     -\sqrt{\frac2{\kappa_{B,0}-2}}.
\]
The scalar damping cancels in the evolution of \(\rho_\pm\), and
\eqref{s6n:frame-equation} gives
\[
\begin{aligned}
 \partial_{\tau_\pm}\rho_\pm
 ={}&-\frac{2\Lambda_{B,0}}{\sqrt{1+s_\pm^2}}\rho_\pm
 -\frac{2\hat F_{B,0}s_{B,0}^{\theta}+|\hat{\mathcal S}_{B,0}|}
              {\sqrt{1+s_\pm^2}}\rho_\pm^2\\
 &+O(k^{-2})(1+|\rho_\pm|+\rho_\pm^2),
 \qquad \rho_\pm(0)=O(k^{-2}).
\end{aligned}
\]
Since \(\Lambda_{B,0}/\sqrt{1+s_\pm^2}\ge c>0\), for fixed
sufficiently large \(K\) and \(k_0\), the interval
\([-Kk^{-2},Kk^{-2}]\) is forward invariant under the evolution of
\(\rho_\pm\) with the prescribed initial data: \(\rho_\pm(0)\) lies
in this interval, and at its two endpoints,
\[
 \left.\partial_{\tau_\pm}\rho_\pm\right|_{\rho_\pm=Kk^{-2}}
 \le k^{-2}(-2cK+C+O(k^{-2}))<0,\qquad
 \left.\partial_{\tau_\pm}\rho_\pm\right|_{\rho_\pm=-Kk^{-2}}>0.
\]
Using the first coordinate equation, we obtain
\[
\begin{aligned}
 \partial_{\tau_\pm}\left(\log\frac{\bar x_\pm}{\mathcal G}\right)
 ={}&\frac{2\hat F_{B,0}s_{B,0}^{\theta}+|\hat{\mathcal S}_{B,0}|}
            {\sqrt{1+s_\pm^2}}\rho_\pm
 +(\mathcal E_\pm)_{11}
 +(\mathcal E_\pm)_{12}\frac{\bar y_\pm}{\bar x_\pm}\\
 &-\left[d_\pm-2^{-\delta k}N^2B_p^2(1+s_\pm^2)\right]
 =O(k^{-2}).
\end{aligned}
\]
The initial ratio \(\bar x_\pm(0)/\mathcal G(0)\asymp1\) and
\(\tau_{\rm max}k^{-2}\asymp1\) imply
\[
 \left|\log\frac{\bar x_\pm}{\mathcal G}\right|\le C,
 \qquad |\rho_\pm|\le Ck^{-2},\qquad
 \bar x_\pm\asymp\bar y_\pm\asymp\mathcal G.
\]
This prevents the appearance of any zero for \(\bar x_\pm\). Substitution into
\eqref{s6n:moving-frame} proves \eqref{s6n:homogeneous-direction}.
At the midpoint, its two reference directions are
\[
 \mathbf q_\pm(\tau_{\rm max}/2)
 =\mp M_{\rm cone}\mathbf s_{B,0}^{\perp}
  -\sqrt{\frac{(\kappa_{B,0}-2)(1+M_{\rm cone}^2)}2}\,
                                                   \mathbf s_{B,0}.
\]
The strict margin in \eqref{s6n:cone-aperture} persists under the
\(O(k^{-2})\) change in these directions.

For any solution, the weighted energy satisfies
\[
 \frac12\partial_{\tau_\pm}
 \left(\frac2{\kappa_{B,0}-2}|\bar x_\pm|^2+|\bar y_\pm|^2\right)
 \le\bigl(a(\tau_\pm)+Ck^{-2}\bigr)
       \left(\frac2{\kappa_{B,0}-2}|\bar x_\pm|^2
                         +|\bar y_\pm|^2\right).
\]
This energy is uniformly equivalent to \(|\mathbf h_\pm|^2\).
Integration gives the intermediate estimate
\begin{equation}\label{s7n:forward-envelope}
 |\mathbf h_\pm(\tau_2)|
 \le C\frac{\mathcal G(\tau_2)}{\mathcal G(\tau_1)}
                     |\mathbf h_\pm(\tau_1)|,
 \qquad 0\le\tau_1\le\tau_2\le\tau_{\rm max}.
\end{equation}
For a source \((g_x,g_y)\), let
\(\bar x_\pm^g,\bar y_\pm^g\) solve the inhomogeneous equation
\[
\begin{gathered}
 \partial_{\tau_\pm}\binom{\bar x_\pm^g}{\bar y_\pm^g}
 =\bigl(\mathsf B_\pm+\mathcal E_\pm-d_\pm I_2\bigr)
       \binom{\bar x_\pm^g}{\bar y_\pm^g}
       +\binom{g_x}{g_y},\\
 \binom{\bar x_\pm^g}{\bar y_\pm^g}(0)=\binom00.
\end{gathered}
\]
By Duhamel's principle, the corresponding estimate is
\[
 |(\bar x_\pm^g,\bar y_\pm^g)(\tau_\pm)|
 \le C\int_0^{\tau_\pm}
       \frac{\mathcal G(\tau_\pm)}{\mathcal G(\tau)}
                       |(g_x,g_y)(\tau)|\,d\tau.
\]
For derivatives in the normalized coordinates
\((\hat R,\hat Z,\hat T)\), differentiate \eqref{s6n:frame-equation}
at fixed \(\tau_\pm\). The principal linear operator is unchanged;
the source is a sum of coefficient derivatives times lower
derivatives of \((\bar x_\pm,\bar y_\pm)\). With the label,
reference point, and frequency parameters held fixed,
\(\bar x_\pm(0)\) depends only on frozen quantities. Hence
\[
 \partial_{\hat R,\hat Z,\hat T}^{\beta}\bar x_\pm(0)=0,
 \qquad |\beta|\ge1.
\]
For \(\bar y_\pm(0)\), the formula in \eqref{s6n:homogeneous-data}
and the lower bound for \(|\mathbf n_{\varphi,\pm}^{\theta,z}(0)|\)
give
\[
 |\partial_{\hat R,\hat Z,\hat T}^{\beta}\bar y_\pm(0)|
 \le C_\beta\mathcal G(0),
 \qquad \beta\in\ZZ_{\ge0}^3.
\]
Each induction step costs at most further powers of \(k\), since
\(\tau_{\rm max}=O(k^2)\).

Fast-time derivatives of $(\bar{x}_\pm, \bar{y}_\pm)$ follow from the equation. The basis
\eqref{s6n:moving-frame} and pressure formula \eqref{s6n:projected-ode},
together with \(N^{-1}\le2^{-\delta k/2}\), then give
\eqref{s6n:homogeneous-derivatives}.
\end{proof}

\subsection{The forced principal equation}
\label{s7n:forced-pulse}

We now fix a nonzero harmonic index \(m\in\ZZ\setminus\{0\}\).
As before, we suppress
\(k,\alpha\) in the coefficient subscripts.

On \(0\le\tau_\pm\le\tau_{\rm max}\), we consider the forced
principal equation with zero initial data
\begin{equation}\label{s7n:forced-system}
\begin{aligned}
 \partial_{\tau_\pm}\mathbf h_{\pm,m}+\mathsf K_B\mathbf h_{\pm,m}
 +m^2d_\pm\mathbf h_{\pm,m}
 +\mathrm iNm\mathbf n_{\varphi,\pm}\pi_{\pm,m}
 &=-\mathbf f_{\pm,m},\\
 \mathbf n_{\varphi,\pm}\cdot\mathbf h_{\pm,m}&=0,
 \qquad \mathbf h_{\pm,m}(0)=0.
\end{aligned}
\end{equation}
Let \(s_0\in\RR\). On each inverse-image rectangle, we use the
coordinates \(Y_i\) of its image and suppress the copy index;
sources on distinct copies may differ.
We assume that the source \(\mathbf f_{\pm,m}\) is axisymmetric and smooth in the independent
normalized slow and auxiliary variables. Its support in the normalized variables lies in
\(\operatorname{supp}(\zeta_k\zeta_{k,\alpha})\), its transverse auxiliary support
lies in \(\operatorname{supp}\Phi_r\), and its radial support lies in
the closed active annulus, with smooth zero extension across these
boundaries. For every \(\beta\in\ZZ_{\ge0}^5\), we assume
the following bound uniformly over the source family, with
constants and exponents independent of \(k,\alpha,\pm\) and the
rectangle copy:
\begin{equation}\label{s7n:forced-source-bounds}
 |\partial_{\hat R,\hat Z,\hat T,Y_{i,1},Y_{i,2}}^{\beta}\mathbf f_{\pm,m}|
 \le C_{\beta}2^{-s_0\delta k}k^{b_{\beta}}
       \sqrt\omega\,d_{\mathrm{edge}}^{-d_{\beta}}\mathcal G(\tau_\pm).
\end{equation}
We establish the basic estimates for the equation \eqref{s7n:forced-system} that will be used
to construct the forced pulse corrections. We follow the method of
the previous subsection, using the same \(\bar x,\bar y\) decomposition
and weighted energy estimates.

\begin{prop}[The zero-data forced inverse]
\label{s7n:forced-inverse}
Under the above assumptions, \eqref{s7n:forced-system} has a unique
solution, and
\begin{equation}\label{s7n:forced-solution-bounds}
\begin{aligned}
 |\partial_{\hat R,\hat Z,\hat T,Y_{i,1},Y_{i,2}}^{\beta}\mathbf h_{\pm,m}|
 &\le C'_{\beta}2^{-s_0\delta k}k^{b'_{\beta}}
       \sqrt\omega\,d_{\mathrm{edge}}^{-d'_{\beta}}\mathcal G(\tau_\pm),\\
 |\partial_{\hat R,\hat Z,\hat T,Y_{i,1},Y_{i,2}}^{\beta}\pi_{\pm,m}|
 &\le C'_{\beta}2^{-(s_0+1/2)\delta k}k^{b'_{\beta}}
       \sqrt\omega\,d_{\mathrm{edge}}^{-d'_{\beta}}\mathcal G(\tau_\pm).
\end{aligned}
\end{equation}
The support of the solution, projected onto
\((\hat{R},\hat{Z},\hat{T},l_r)\), is contained in the
corresponding projection of the source support. 
The solution map is linear. It also satisfies
\[
 \mathbf f_{\pm,-m}=\overline{\mathbf f_{\pm,m}}
 \quad\Longrightarrow\quad
 \mathbf h_{\pm,-m}=\overline{\mathbf h_{\pm,m}},\qquad
 \pi_{\pm,-m}=\overline{\pi_{\pm,m}}.
\]
\end{prop}

\begin{proof}
Recall the pressure formula and projected equation from
\eqref{s6n:general-projected-ode}:
\begin{equation}\label{s7n:forced-projection}
\begin{aligned}
 \pi_{\pm,m}
 &=\frac{\mathrm i}{Nm}
   \frac{\mathbf n_{\varphi,\pm}\cdot\mathsf K_B\mathbf h_{\pm,m}
       -(\partial_{\tau_\pm}\mathbf n_{\varphi,\pm})\cdot\mathbf h_{\pm,m}
       +\mathbf n_{\varphi,\pm}\cdot\mathbf f_{\pm,m}}
        {|\mathbf n_{\varphi,\pm}|^2},\\
 \partial_{\tau_\pm}\mathbf h_{\pm,m}
 &=(\mathsf A_\pm-m^2d_\pm I_3)\mathbf h_{\pm,m}
      -\Pi_\pm\mathbf f_{\pm,m},\\
 \Pi_\pm&=I_3-\frac{\mathbf n_{\varphi,\pm}
                         \mathbf n_{\varphi,\pm}^{\mathsf T}}
                         {|\mathbf n_{\varphi,\pm}|^2}.
\end{aligned}
\end{equation}
Decompose \(\mathbf h_{\pm,m}\) as in \eqref{s6n:moving-frame}.
In the non-homogeneous setting, equation~\eqref{s6n:frame-equation} becomes
\begin{equation} \label{eq:two-coord}
\begin{aligned}
 \partial_{\tau_\pm}\binom{\bar x_{\pm,m}}{\bar y_{\pm,m}}
 &=\left[
 \mathsf B_\pm
 +\mathcal E_\pm-m^2d_\pm I_2\right]
 \binom{\bar x_{\pm,m}}{\bar y_{\pm,m}}\\
 &\quad-\binom{\mathbf n_{\varphi,\pm,x}^{\perp}\cdot\mathbf f_{\pm,m}}
                  {\mathbf n_{\varphi,\pm,y}^{\perp}\cdot\mathbf f_{\pm,m}},
 \qquad
 \binom{\bar x_{\pm,m}(0)}{\bar y_{\pm,m}(0)}=0.
\end{aligned}
\end{equation}
In the energy
\(2|\bar x_{\pm,m}|^2/(\kappa_{B,0}-2)+|\bar y_{\pm,m}|^2\),
the additional scalar term \(-(m^2-1)d_\pm I_2\) is nonpositive.
The argument proving \eqref{s7n:forward-envelope} therefore gives,
uniformly for \(m\ne0\),
\begin{equation}\label{s7n:forced-forward-envelope}
 |\mathbf h_{\pm,m}(\tau_\pm)|
 \le C\frac{\mathcal G(\tau_\pm)}{\mathcal G(\tau)}
                    |\mathbf h_{\pm,m}(\tau)|
   +C\int_\tau^{\tau_\pm}\frac{\mathcal G(\tau_\pm)}{\mathcal G(\tau')}
                    |\mathbf f_{\pm,m}(\tau')|\,d\tau',
 \qquad 0\le \tau\le\tau_\pm\le\tau_{\rm max}.
\end{equation}
In particular, zero initial data imply
\begin{equation}\label{s7n:forced-duhamel}
 |\mathbf h_{\pm,m}(\tau_\pm)|
 \le C\int_0^{\tau_\pm}\frac{\mathcal G(\tau_\pm)}{\mathcal G(\tau')}
                    |\mathbf f_{\pm,m}(\tau')|\,d\tau'.
\end{equation}
In this integral, \((\hat R,\hat Z,\hat T,l_r)\) are fixed. The radial
weights \(\omega(R)\) and \(d_{\mathrm{edge}}(R)\) are therefore
independent of the integration variable. Using
\eqref{s7n:forced-source-bounds} and
\(\tau_{\rm max}\asymp k^2\), we obtain
\[
\begin{aligned}
 |\mathbf h_{\pm,m}(\tau_\pm)|
 &\le C2^{-s_0\delta k}k^b\sqrt\omega\,
             d_{\mathrm{edge}}^{-d_1}
       \int_0^{\tau_\pm}\frac{\mathcal G(\tau_\pm)}{\mathcal G(\tau')}
                                     \mathcal G(\tau')\,d\tau'\\
 &\le C2^{-s_0\delta k}k^{b+2}\sqrt\omega\,
             d_{\mathrm{edge}}^{-d_1}\mathcal G(\tau_\pm).
\end{aligned}
\]

For derivatives in \((\hat R,\hat Z,\hat T,l_r)\), keep
\(\tau_\pm\) fixed and differentiate the displayed two-coordinate
equation \eqref{eq:two-coord}. Its principal linear operator is unchanged.
The remaining terms are derivatives of the projected source and
products of coefficient derivatives with lower-order derivatives of
\(\bar x_{\pm,m},\bar y_{\pm,m}\). All fixed derivatives of those
coefficients and of the orthonormal basis are polynomially bounded
in \(k\), for fixed \(m\).
The same energy estimate gives the first bound in
\eqref{s7n:forced-solution-bounds} by induction on \(|\beta|\).
Fast-time derivatives follow from the equation \eqref{eq:two-coord}. The change from
\((l_r,\tau_\pm)\) to \(Y_i\) costs only powers of \(k\). The pressure bound
follows from \eqref{s7n:forced-projection}, the nondegeneracy in
Lemma~\ref{s7n:phase-bounds}, and
\(|Nm|^{-1}\le2^{-\delta k/2}\).

If the source vanishes for all fast times at fixed
\((\hat R,\hat Z,\hat T,l_r)\), estimate~\eqref{s7n:forced-duhamel}
gives a zero solution there. Applying the differentiated
equations on the support boundaries in the normalized and transverse variables gives smooth
zero extension there; the weighted bounds give it at the radial
edges. Conjugate symmetry follows from uniqueness.
\end{proof}

\subsection{Matching the complete background stress}
\label{sec:matching-background}
\label{s8n:matching-data}

\pulseLabelAlias{subsec:pulse-phase}{s7n:frequencies}
\pulseLabelAlias{subsec:homogeneous-pulse}{s6n:unforced-pulse}
\pulseLabelAlias{eq:pulse-growth-parameters}{s6n:growth-parameters}
\pulseLabelAlias{eq:pulse-frequency-choice}{s6n:frequency-choice}
\pulseLabelAlias{eq:pulse-phase-choice}{s6n:phase}
\pulseLabelAlias{eq:pulse-phase-defect}{s6n:transport-defect}
\pulseLabelAlias{eq:pulse-principal-ode}{s6n:principal-ode}
\pulseLabelAlias{eq:pulse-frame-equation}{s6n:frame-equation}
\pulseLabelAlias{eq:pulse-homogeneous-data}{s6n:homogeneous-data}
\pulseLabelAlias{prop:pulse-forced-inverse}{s7n:forced-inverse}

We use the complete background from Section~\ref{intro:flux-cone}.
The scalar weights of the unforced pre-pulses from
Section~\ref{s6n:unforced-pulse} will be chosen so that their averaged
radial tangential flux equals \(\mathsf T_B\). In band \(k\), the
normalized target is
\[
 \hat{\mathsf T}_B=2^{-(2+\delta)k}\mathsf T_B.
\]
We use the uniform cone condition \eqref{intro:admissible-cone} and,
to prove smoothness of the matching amplitudes at the annular edges,
the bounds \eqref{intro:weighted-background-stress}.
\pulseLabelAlias{eq:completed-background-residual}{intro:background}

The following update and cylindrical divergence identities are written
for physical fields. On extended fields, their derivatives are represented
by the chain-rule operators of Section~\ref{s6n:auxiliary}.
For later corrections we record the exact update identity:
\begin{equation}
\begin{aligned}
\mathsf{NS}(u+\widetilde u,p+\widetilde p)
={}&\mathsf{NS}(u,p)+\partial_t\widetilde u-\Delta\widetilde u
+\nabla\widetilde p\\
&+(u\cdot\nabla)\widetilde u+(\widetilde u\cdot\nabla)u
+(\widetilde u\cdot\nabla)\widetilde u.
\end{aligned}
\label{eq:physical-residual-update}
\end{equation}

For a real vector field \(w\), its radial flux pair
\(w^r(w^\theta,w^z)\) enters the tangential components of
\(\nabla\cdot(w\otimes w)\) through the operators
\(\partial_r+2/r\) and \(\partial_r+1/r\), respectively. These are
the quadratic fluxes used for stress matching. Indeed, in the cylindrical orthonormal frame
\((e_r,e_\theta,e_z)\),
\[
\begin{aligned}
\bigl[\nabla\cdot(w\otimes w)\bigr]^r
&=\left(\partial_r+\frac1r\right)(w^r)^2
+\frac1r\partial_\theta(w^rw^\theta)
+\partial_z(w^rw^z)-\frac{(w^\theta)^2}{r},\\[4pt]
\bigl[\nabla\cdot(w\otimes w)\bigr]^\theta
&=\left(\partial_r+\frac2r\right)(w^rw^\theta)
+\frac1r\partial_\theta(w^\theta)^2
+\partial_z(w^zw^\theta),\\[4pt]
\bigl[\nabla\cdot(w\otimes w)\bigr]^z
&=\left(\partial_r+\frac1r\right)(w^rw^z)
+\frac1r\partial_\theta(w^\theta w^z)
+\partial_z(w^z)^2.
\end{aligned}
\]
The \(\theta\)-derivatives vanish under angular averaging.
The radial fluxes \(w^rw^\theta\) and \(w^rw^z\) transport azimuthal
and axial momentum, respectively. They can potentially match and cancel the background stress
because their divergences have precisely the operators and components
of the background stress divergence in
\eqref{eq:completed-background-residual}.
The radial component and the axial flux terms remain in the residual
for later corrections.
All auxiliary averages below precede physical substitution. In the
local construction, they use \(Y_i\), with \(i=i(k)\).

\paragraph{The matched amplitudes.}
For each \((k,\alpha)\) whose cutoff support
\(\operatorname{supp}(\zeta_k\zeta_{k,\alpha})\) meets the active
stress annulus, use the fixed box \(Q_{k,\alpha}\) in
\eqref{s6n:slow-box-cutoffs}. In this and the next subsection,
statements on \(Q_{k,\alpha}\) refer to its intersection with
\(\mathcal B_k\) and the domain of the background coefficients,
expressed in normalized coordinates.
Let \(g_{\mathrm{match},\gamma}
=g_{\mathrm{match},\gamma}(\hat R,\hat Z,\hat T)\ge0\)
be smooth on \(Q_{k,\alpha}\), with polynomially bounded
fixed derivatives in these variables.
For every label \(\gamma=(k,\alpha,\iota)\), with \(\iota=\pm1\),
let \((\mathbf h_\gamma,\pi_\gamma)\) be the unforced solution of
Proposition~\ref{prop:homogeneous-pulse}, and set
\begin{equation}\label{s6n:amplitude-data}
\begin{aligned}
 \mathbf h_{\mathrm{match},\gamma}
 &=2^{-\delta k/2}g_{\mathrm{match},\gamma}\mathbf h_\gamma,\\
 \pi_{\mathrm{match},\gamma}
 &=2^{-\delta k/2}g_{\mathrm{match},\gamma}\pi_\gamma.
\end{aligned}
\end{equation}
The scale \(\lambda\) depends only on \((\hat Z,\hat T)\), which are
held fixed when differentiating in the auxiliary fast time. Thus
\(\partial_{\tau_\gamma}g_{\mathrm{match},\gamma}=0\), including when
\(g_{\mathrm{match},\gamma}\) depends on \(\lambda\). This pair \eqref{s6n:amplitude-data} therefore
solves the unforced principal equation and satisfies
\begin{equation}\label{s6n:envelope}
 |\mathbf h_{\mathrm{match},\gamma}|
 \le C2^{-\delta k/2}g_{\mathrm{match},\gamma}\mathcal G(\tau_\gamma)
 \le C2^{-\delta k/2}g_{\mathrm{match},\gamma}
       \exp\!\left(-c\frac{(\tau_\gamma-\tau_{\rm max}/2)^2}{\tau_{\rm max}}\right).
\end{equation}
The corresponding pre-pulse is
\[
 w_\gamma^{\rm pre}
 =2^{(1+\delta)k}\chi_\gamma\mathbf h_{\mathrm{match},\gamma}
   \cos(N\varphi_\gamma).
\]
The factor \(g_{\mathrm{match},\gamma}\) remains free for quadratic
stress matching.
We choose \(g_{\mathrm{match},\gamma}>0\) in the active annulus
\[
 \hat R\in\bigl((2^k\lambda)^2R_a,(2^k\lambda)^2R_h\bigr)
\]
such that the sum of their averaged radial fluxes over both signs
and all boxes and bands equals
\[
 \mathsf T_B
 =\begin{pmatrix}\mathsf T_B^\theta\\\mathsf T_B^z\end{pmatrix}.
\]
All averages in this subsection keep \((\hat{R},\hat{Z},\hat{T})\) fixed.
They are taken before physical evaluation.
For a real vector field \(w\), recall \eqref{intro:flux-map}:
\begin{equation}
\begin{aligned}
\mathsf C(w)
&=\left\langle
w^r(w^\theta,w^z)^{\mathsf T}
\right\rangle_{\theta,{Y_i}}\\
&=\int_{\TT^2}\frac1{2\pi}\int_0^{2\pi}
\begin{pmatrix}w^rw^\theta\\w^rw^z\end{pmatrix}
\,d\theta\,dY_i,\qquad \int_{\TT^2}dY_i=1.
\end{aligned}
\label{eq:pulse-covariance-pair}
\end{equation}
The vector \(\mathsf C(w)\in\RR^2\) records radial
transport of azimuthal and axial momentum.

On the auxiliary rectangle \(\mathcal R_\gamma\) associated with
each label \(\gamma=(k,\alpha,\iota)\),
use the unforced solution of Section~\ref{s6n:unforced-pulse} and define
\begin{equation}
\begin{aligned}
b_\gamma
&=\Phi_r(l_r)\Phi_t(\tau_\gamma/\tau_{\rm max})
\mathbf h_\gamma\cos(N\varphi_\gamma),\\
H_\gamma&=\mathsf C(b_\gamma),\qquad
H_{k,\alpha}=\bigl[H_{(k,\alpha,+)}\mid H_{(k,\alpha,-)}\bigr].
\end{aligned}
\label{eq:basic-pulses-covariance-matrix}
\end{equation}
Define each \(b_\gamma\) on the \(Y_i\)-torus, where \(i=i(k)\),
and extend it by zero off its assigned rectangle.
Here \(\tau_\gamma=\tau_\gamma(Y_i)\) on the corresponding rectangle.
These normalized pre-pulses include the auxiliary cutoffs.
They do not yet include \(\zeta_k\zeta_{k,\alpha}\), the cutoff in \((\hat{R},\hat{Z},\hat{T})\).
The scalar matching amplitude is also absent.
The two signs in each box have disjoint auxiliary supports, so
\[
b_\gamma^r
(b_{(k,\alpha,-\iota)}^\theta,b_{(k,\alpha,-\iota)}^z)^{\mathsf T}=0,
\qquad \gamma=(k,\alpha,\iota).
\]
Since \(g_{\mathrm{match},\gamma}\) depend only on
\((\hat{R},\hat{Z},\hat{T})\), 
\begin{equation}
\mathsf C\!\left(\sum_{\iota=\pm1}
g_{\mathrm{match},(k,\alpha,\iota)}b_{(k,\alpha,\iota)}\right)
=H_{k,\alpha}\begin{pmatrix}g_{\mathrm{match},(k,\alpha,+)}^{\,2}\\
                  g_{\mathrm{match},(k,\alpha,-)}^{\,2}\end{pmatrix}.
\label{eq:two-pulse-positive-combination}
\end{equation}

\begin{prop}[Positive matching of the complete background stress]
\label{prop:leading-positive-matching}
After choosing \(k_0\) large, \(H_{k,\alpha}\) is invertible
throughout \(Q_{k,\alpha}\), with the domain convention above.
On its intersection with the active annulus
\((2^k\lambda)^2R_a<\hat R<(2^k\lambda)^2R_h\), both components of
\(H_{k,\alpha}^{-1}(2^{-(2+\delta)k}\mathsf T_B)\) are positive. Define
\begin{equation}
\begin{gathered}
\begin{pmatrix}g_{\mathrm{match},(k,\alpha,+)}^{\,2}\\
                g_{\mathrm{match},(k,\alpha,-)}^{\,2}\end{pmatrix}
=H_{k,\alpha}^{-1}\bigl(2^{-(2+\delta)k}\mathsf T_B\bigr),
\qquad g_{\mathrm{match},\gamma}>0,\\
W_{k,\alpha}
=2^{-\delta k/2}
\sum_{\iota=\pm1}g_{\mathrm{match},(k,\alpha,\iota)}
b_{(k,\alpha,\iota)}.
\end{gathered}
\label{eq:positive-matching-amplitudes}
\end{equation}
Set \(g_{\mathrm{match},\gamma}=0\) outside the active annulus; these amplitudes extend smoothly across both radial edges. For every
fixed multi-index \(\beta\) in \((\hat R,\hat Z,\hat T)\),
\begin{equation}
 \bigl|\partial_{\hat R,\hat Z,\hat T}^{\beta}
                  g_{\mathrm{match},\gamma}\bigr|
 +\bigl|\partial_{\hat R,\hat Z,\hat T}^{\beta}
                  (g_{\mathrm{match},\gamma}^{\,2})\bigr|
 \le C_\beta k^{b_\beta}.
\label{eq:matching-square-root-bounds}
\end{equation}
The normalized field has averaged radial flux
\begin{equation}
\mathsf C(W_{k,\alpha})
=2^{-2(1+\delta)k}\mathsf T_B.
\label{eq:local-leading-covariance}
\end{equation}
\end{prop}

\begin{proof}
Each relevant \(Q_{k,\alpha}\) has diameter \(O(k^{-6})\) in
\((\sqrt{2\hat R},\hat Z,\hat T)\) and meets the active annulus.
For sufficiently large \(k_0\), its part under consideration lies
in the fixed enlarged radial annulus of Section~\ref{s7n:frequencies}.
Using the same frozen parameters, the wavevector and unforced amplitude
constructions and estimates of Lemma~\ref{s7n:phase-bounds} and
Proposition~\ref{prop:homogeneous-pulse} apply there by the same argument
using the \(O(k^{-6})\) variation and background bounds; the cone
estimates are used only on the active annulus.

Since
\(\partial_\theta\varphi_\gamma=p_{\theta,\gamma}\) and
\(Np_{\theta,\gamma}\in\ZZ\setminus\{0\}\),
\[
\frac1{2\pi}\int_0^{2\pi}\cos^2(N\varphi_\gamma)\,d\theta
=\frac12+\frac1{4\pi}\int_0^{2\pi}\cos(2N\varphi_\gamma)\,d\theta
=\frac12.
\]
The unforced amplitude is independent of \(l_r\), and
\[
dY_i=|\det(v_r,v_t)|\,dl_r\,dl_t,
\qquad
dl_t=2^{-(2+\delta)k}\Lambda_t^i\,d\tau.
\]
It follows that each column \(H_\gamma\) of \(H_{k,\alpha}\) satisfies
\begin{equation}
H_\gamma
=\frac{|\det(v_r,v_t)|}{2}
\left(\int_{\mathbb R}\Phi_r(l_r)^2\,dl_r\right)
2^{-(2+\delta)k}\Lambda_t^i
\int_0^{\tau_{\rm max}}\Phi_t(\tau/\tau_{\rm max})^2
h_\gamma^r(h_\gamma^\theta,h_\gamma^z)^{\mathsf T}
\,d\tau.
\label{eq:covariance-column-integral}
\end{equation}
Here we used that the covering maps preserve normalized Haar averages.
Define the positive numbers
\begin{equation}
J_\gamma
:=\frac{|\det(v_r,v_t)|}{2}
\left(\int_{\mathbb R}\Phi_r(l_r)^2\,dl_r\right)
2^{-(2+\delta)k}\Lambda_t^i
\int_0^{\tau_{\rm max}}\Phi_t(\tau/\tau_{\rm max})^2(h_\gamma^r)^2
\,d\tau.
\label{eq:covariance-column-mass}
\end{equation}
The bounds from \eqref{s6n:homogeneous-direction} and
\eqref{s6n:gaussian-bounds}, namely \(h_\gamma^r\asymp\mathcal G\) and
\[
e^{-C(\tau-\tau_{\rm max}/2)^2/\tau_{\rm max}}
\le\mathcal G\le
e^{-c(\tau-\tau_{\rm max}/2)^2/\tau_{\rm max}}
\]
give
\begin{equation}
\begin{gathered}
\int_0^{\tau_{\rm max}}\Phi_t(\tau/\tau_{\rm max})^2(h_\gamma^r)^2
\,d\tau\asymp\sqrt{\tau_{\rm max}},\\
\int_0^{\tau_{\rm max}}\frac{|\tau-\tau_{\rm max}/2|}{\tau_{\rm max}}
\Phi_t(\tau/\tau_{\rm max})^2(h_\gamma^r)^2\,d\tau\le C.
\end{gathered}
\label{eq:covariance-gaussian-concentration}
\end{equation}
Indeed, set \(s=(\tau-\tau_{\rm max}/2)/\sqrt{\tau_{\rm max}}\).
For large \(k\), \(|s|\le1\) implies
\(\tau/\tau_{\rm max}\in[3/10,7/10]\), where \(\Phi_t=1\). Hence
\[
\begin{aligned}
c\sqrt{\tau_{\rm max}}\int_{-1}^1 e^{-Cs^2}\,ds
&\le\int_0^{\tau_{\rm max}}\Phi_t^2(h_\gamma^r)^2\,d\tau
\le C\sqrt{\tau_{\rm max}}\int_{\RR}e^{-cs^2}\,ds,\\[4pt]
\int_0^{\tau_{\rm max}}\frac{|\tau-\tau_{\rm max}/2|}{\tau_{\rm max}}
\Phi_t^2(h_\gamma^r)^2\,d\tau
&\le C\int_{\RR}|s|e^{-cs^2}\,ds\le C.
\end{aligned}
\]
Since \(2^{-(2+\delta)k}\Lambda_t^i=2r_0/\tau_{\rm max}\), we have
\[
J_\gamma\asymp \tau_{\rm max}^{-1}\sqrt{\tau_{\rm max}}
=\tau_{\rm max}^{-1/2}\asymp k^{-1},\qquad
\frac{\displaystyle\int_0^{\tau_{\rm max}}\frac{|\tau-\tau_{\rm max}/2|}{\tau_{\rm max}}
\Phi_t^2(h_\gamma^r)^2\,d\tau}
{\displaystyle\int_0^{\tau_{\rm max}}\Phi_t^2(h_\gamma^r)^2\,d\tau}
\le\frac{C}{\sqrt{\tau_{\rm max}}}\le\frac Ck.
\]

Let \(\mathbf s_{B,0,k,\alpha}\) and \(\kappa_{B,0,k,\alpha}\)
denote the frozen values at the \((k,\alpha)\)-reference point used in
Section~\ref{s7n:frequencies}, and write
\(s_\gamma(\tau)=\iota M_{\rm cone}(1/2+\tau/\tau_{\rm max})\)
as in \eqref{s7n:radial-paths}.
Equation~\eqref{s6n:homogeneous-direction} gives
\[
\begin{gathered}
\frac{(h_\gamma^\theta,h_\gamma^z)^{\mathsf T}}{h_\gamma^r}
=-s_\gamma(\tau)\mathbf s_{B,0,k,\alpha}^\perp
-\sqrt{\frac{(\kappa_{B,0,k,\alpha}-2)(1+s_\gamma(\tau)^2)}2}\,\mathbf s_{B,0,k,\alpha}+O(k^{-2}),\\
s_\gamma(\tau_{\rm max}/2)=\iota M_{\rm cone}.
\end{gathered}
\]
Moreover,
\[
\begin{aligned}
|s_\gamma(\tau)-\iota M_{\rm cone}|
&=\frac{M_{\rm cone}}{\tau_{\rm max}}|\tau-\tau_{\rm max}/2|,\\
\left|\sqrt{1+s_\gamma(\tau)^2}-\sqrt{1+M_{\rm cone}^2}\right|
&\le |s_\gamma(\tau)-\iota M_{\rm cone}|.
\end{aligned}
\]
Set
\(A_{\mathrm{cone},k,\alpha}=\sqrt{(\kappa_{B,0,k,\alpha}-2)(1+M_{\rm cone}^2)/2}>0\).
By \eqref{eq:covariance-column-integral} and
\eqref{eq:covariance-column-mass}, the normalized column is the weighted average
\[
\frac{H_\gamma}{J_\gamma}
=\frac{\displaystyle\int_0^{\tau_{\rm max}}
\Phi_t(\tau/\tau_{\rm max})^2(h_\gamma^r)^2
\frac{(h_\gamma^\theta,h_\gamma^z)^{\mathsf T}}{h_\gamma^r}\,d\tau}
{\displaystyle\int_0^{\tau_{\rm max}}
\Phi_t(\tau/\tau_{\rm max})^2(h_\gamma^r)^2\,d\tau}.
\]
The preceding direction bounds and Gaussian concentration show that
this average differs from the frozen midpoint direction
\(-A_{\mathrm{cone},k,\alpha}\mathbf s_{B,0,k,\alpha}
-\iota M_{\rm cone}\mathbf s_{B,0,k,\alpha}^\perp\)
by at most \(C(\tau_{\rm max}^{-1/2}+k^{-2})\le Ck^{-1}\).
Thus
\begin{equation}
\begin{gathered}
H_\gamma=J_\gamma\bigl(-A_{\mathrm{cone},k,\alpha}\mathbf s_{B,0,k,\alpha}-\iota M_{\rm cone}\mathbf s_{B,0,k,\alpha}^\perp+\mathbf r_\gamma\bigr),\\
\mathbf r_\gamma\in\RR^2,\qquad |\mathbf r_\gamma|\le Ck^{-1}.
\end{gathered}
\label{eq:two-covariance-directions}
\end{equation}
To see why the matching coefficients are positive, first consider
the reference matching problem obtained by setting \(\mathbf r_\gamma=0\).
In the following reference calculation, \(g_{\mathrm{match},\gamma}\)
denotes the solution of that problem. Write
\[
2^{-(2+\delta)k}\mathsf T_B
=T_{\parallel,k,\alpha}\mathbf s_{B,0,k,\alpha}+T_{\perp,k,\alpha}\mathbf s_{B,0,k,\alpha}^\perp.
\]
The two components of the matching equation are
\[
\begin{aligned}
-A_{\mathrm{cone},k,\alpha}\sum_{\iota=\pm1}
J_{(k,\alpha,\iota)}g_{\mathrm{match},(k,\alpha,\iota)}^2
&=T_{\parallel,k,\alpha},\\
-M_{\rm cone}\sum_{\iota=\pm1}\iota
J_{(k,\alpha,\iota)}g_{\mathrm{match},(k,\alpha,\iota)}^2
&=T_{\perp,k,\alpha}.
\end{aligned}
\]
Adding and subtracting gives, for \(\gamma=(k,\alpha,\iota)\),
\begin{equation}
J_\gamma g_{\mathrm{match},\gamma}^{\,2}
=\frac12\left(-\frac{T_{\parallel,k,\alpha}}{A_{\mathrm{cone},k,\alpha}}
-\iota\frac{T_{\perp,k,\alpha}}{M_{\rm cone}}\right).
\label{eq:reference-positive-squares}
\end{equation}
The frozen cone condition \eqref{s6n:cone-aperture} holds pointwise
for \(\hat{\mathsf T}_B\), with a fixed margin throughout the active
annulus. By the definition of \(A_{\mathrm{cone},k,\alpha}\), it gives a fixed
\(c_{\rm match}>0\) such that
\[
\begin{gathered}
T_{\parallel,k,\alpha}<0,\qquad
\frac{A_{\mathrm{cone},k,\alpha}|T_{\perp,k,\alpha}|}{M_{\rm cone}(-T_{\parallel,k,\alpha})}\le1-c_{\rm match},\\
\frac12\left(-\frac{T_{\parallel,k,\alpha}}{A_{\mathrm{cone},k,\alpha}}-\iota\frac{T_{\perp,k,\alpha}}{M_{\rm cone}}\right)
\ge\frac{c_{\rm match}(-T_{\parallel,k,\alpha})}{2A_{\mathrm{cone},k,\alpha}}.
\end{gathered}
\]
We now return to the actual matching amplitudes. Restoring \(\mathbf r_\gamma\),
the actual matrix satisfies
\[
[\mathbf s_{B,0,k,\alpha}\mid\mathbf s_{B,0,k,\alpha}^\perp]^{\mathsf T}
H_{k,\alpha}\begin{pmatrix}J_{(k,\alpha,+)}^{-1}&0\\0&J_{(k,\alpha,-)}^{-1}\end{pmatrix}
=\begin{pmatrix}-A_{\mathrm{cone},k,\alpha}&-A_{\mathrm{cone},k,\alpha}\\-M_{\rm cone}&M_{\rm cone}\end{pmatrix}+O(k^{-1}),
\]
and therefore
\[
\begin{pmatrix}J_{(k,\alpha,+)}g_{\mathrm{match},(k,\alpha,+)}^2\\J_{(k,\alpha,-)}g_{\mathrm{match},(k,\alpha,-)}^2\end{pmatrix}
=\left[
\begin{pmatrix}-(2A_{\mathrm{cone},k,\alpha})^{-1}&-(2M_{\rm cone})^{-1}\\
-(2A_{\mathrm{cone},k,\alpha})^{-1}&(2M_{\rm cone})^{-1}\end{pmatrix}
+O(k^{-1})\right]\begin{pmatrix}T_{\parallel,k,\alpha}\\T_{\perp,k,\alpha}\end{pmatrix}.
\]
The error is at most \(Ck^{-1}(|T_{\parallel,k,\alpha}|+|T_{\perp,k,\alpha}|)\le C'k^{-1}(-T_{\parallel,k,\alpha})\).
This is smaller than \(c_{\rm match}(-T_{\parallel,k,\alpha})/(4A_{\mathrm{cone},k,\alpha})\) for large \(k_0\).
Also \(|\det H_{k,\alpha}|=J_{(k,\alpha,+)}J_{(k,\alpha,-)}|{-2A_{\mathrm{cone},k,\alpha}M_{\rm cone}}+O(k^{-1})|\).
Using \(J_\gamma\asymp k^{-1}\) gives
\begin{equation}
\begin{gathered}
|\det H_{k,\alpha}|\ge ck^{-2},\qquad \|H_{k,\alpha}^{-1}\|\le Ck,
\\
ck\,\bigl|2^{-(2+\delta)k}\mathsf T_B\bigr|
\le g_{\mathrm{match},\gamma}^{\,2}
\le Ck\,\bigl|2^{-(2+\delta)k}\mathsf T_B\bigr|.
\end{gathered}
\label{eq:positive-matching-inverse-bounds}
\end{equation}

We next prove the derivative bounds and smoothness at the radial edges.
The target satisfies
\[
 \hat{\mathsf T}_B=(2^k\lambda)^{-(2+\delta)}
                (\lambda^{2+\delta}\mathsf T_B)(\lambda,R,Z).
\]
The similarity-normalized stress \(\lambda^{2+\delta}\mathsf T_B\)
may depend on \(\lambda\).
The coordinate change \eqref{s6n:coordinate-change} gives
\[
 (2^k\lambda)^2-(2^k\lambda)^{2\delta}\hat Z^2=\hat T,\qquad
 R=\hat R/(2^k\lambda)^2,\qquad Z=\hat Z/(2^k\lambda)^{1-\delta}.
\]
Implicit differentiation yields
\[
 \partial_{\hat T}(2^k\lambda)=\frac1{2(2^k\lambda)(1-\delta Z^2)},\qquad
 \partial_{\hat Z}(2^k\lambda)=\frac{(2^k\lambda)^\delta Z}{1-\delta Z^2},\qquad
 \partial_{\hat R}(2^k\lambda)=0.
\]
For a function of the independent variables \((\lambda,R,Z)\),
the chain rule after pullback to band coordinates takes the form
\begin{equation}
\begin{aligned}
 \partial_{\hat R}
   &=(2^k\lambda)^{-2}\partial_R,\\
 \partial_{\hat T}
   &=\frac{\lambda\partial_\lambda-2R\partial_R-(1-\delta)Z\partial_Z}
           {2(2^k\lambda)^2(1-\delta Z^2)},\\
 \partial_{\hat Z}
   &=\frac{(2^k\lambda)^{\delta-1}}{1-\delta Z^2}
       \left[Z\lambda\partial_\lambda-2RZ\partial_R+(1-Z^2)\partial_Z\right].
\end{aligned}
\label{eq:matching-target-chain-rule}
\end{equation}
Since \(1/2<2^k\lambda<2\), \(|Z|\le1\), and \(R\) lies in a fixed
enlarged annulus, all fixed derivatives of these coefficients are
bounded uniformly in \(k\). Applying
\eqref{intro:weighted-background-stress}, including its
\(\lambda\partial_\lambda\) derivatives, gives
\begin{equation}
 |\hat{\mathsf T}_B|\ge c\omega,\qquad
 \bigl|\partial_{\hat R,\hat Z,\hat T}^{\beta}
                    \hat{\mathsf T}_B\bigr|
 \le C_\beta\omega d_{\mathrm{edge}}^{-d_\beta}.
\label{eq:matching-target-edge-bounds}
\end{equation}
Here \(\omega\) and \(d_{\mathrm{edge}}\) are evaluated at
\(R=\hat R/(2^k\lambda)^2\).
Differentiation of \eqref{eq:covariance-column-integral}, using
\eqref{s6n:homogeneous-derivatives}, gives
\(\lVert\partial_{\hat R,\hat Z,\hat T}^{\beta}H_{k,\alpha}\rVert\le C_{\beta} k^{b_{\beta}}\). For \(|\beta|>0\),
\[
\partial_{\hat R,\hat Z,\hat T}^{\beta}H_{k,\alpha}^{-1}
=-H_{k,\alpha}^{-1}\sum_{0<\beta'\le \beta}\binom{\beta}{\beta'}
(\partial_{\hat R,\hat Z,\hat T}^{\beta'}H_{k,\alpha})(\partial_{\hat R,\hat Z,\hat T}^{\beta-\beta'}H_{k,\alpha}^{-1}),
\]
so induction gives \(\lVert\partial_{\hat R,\hat Z,\hat T}^{\beta}H_{k,\alpha}^{-1}\rVert\le C_{\beta} k^{b_{\beta}}\).
Finally,
\[
\partial_{\hat R,\hat Z,\hat T}^{\beta}
\begin{pmatrix}g_{\mathrm{match},(k,\alpha,+)}^2\\g_{\mathrm{match},(k,\alpha,-)}^2\end{pmatrix}
=\sum_{\beta'\le \beta}\binom{\beta}{\beta'}(\partial_{\hat R,\hat Z,\hat T}^{\beta'}H_{k,\alpha}^{-1})
\partial_{\hat R,\hat Z,\hat T}^{\beta-\beta'}\bigl(2^{-(2+\delta)k}\mathsf T_B\bigr).
\]
Together with \eqref{eq:positive-matching-inverse-bounds} and
\eqref{eq:matching-target-edge-bounds}, this yields
\begin{equation}
 g_{\mathrm{match},\gamma}^{\,2}\ge ck\omega,\qquad
 \bigl|\partial_{\hat R,\hat Z,\hat T}^{\beta}
             (g_{\mathrm{match},\gamma}^{\,2})\bigr|
 \le C_\beta k^{b_\beta}\omega d_{\mathrm{edge}}^{-d_\beta}.
\label{eq:matching-squares-edge-bounds}
\end{equation}

For a nonzero multi-index \(\beta\), each term in a derivative of the
positive square root $g_{\mathrm{match},\gamma} = \sqrt{g_{\mathrm{match},\gamma}^2}$ has the form
\[
C\,g_{\mathrm{match},\gamma}^{\,1-2j}
\prod_{\nu=1}^j
\partial_{\hat R,\hat Z,\hat T}^{\beta^{(\nu)}}(g_{\mathrm{match},\gamma}^{\,2}),
\qquad
\sum_{\nu=1}^j\beta^{(\nu)}=\beta,\quad |\beta^{(\nu)}|\ge1.
\]
For \(j\ge1\), its absolute value is at most
\[
C_{\beta}(k\omega)^{1/2-j}
\prod_{\nu=1}^j\bigl(k^{b_{\beta^{(\nu)}}}\omega d_{\mathrm{edge}}^{-d_{\beta^{(\nu)}}}\bigr)
\le C_{\beta} k^{b'_{\beta}}\sqrt\omega d_{\mathrm{edge}}^{-d'_{\beta}}.
\]
The case \(\beta=0\) follows by taking the square root of its upper
bound. For every fixed exponent \(b\ge0\), both
\(\omega d_{\mathrm{edge}}^{-b}\) and
\(\sqrt\omega d_{\mathrm{edge}}^{-b}\) are bounded and tend to zero
at either edge. Thus every slow derivative of \(g_{\mathrm{match},\gamma}\)
extends by zero, and \eqref{eq:matching-square-root-bounds} follows.
This also proves smoothness across the moving boundaries
\(\hat R=(2^k\lambda)^2R_a,(2^k\lambda)^2R_h\).

Finally,
\eqref{eq:two-pulse-positive-combination} gives
\[
\mathsf C(W_{k,\alpha})
=2^{-\delta k}H_{k,\alpha}
\begin{pmatrix}g_{\mathrm{match},(k,\alpha,+)}^{\,2}\\
                g_{\mathrm{match},(k,\alpha,-)}^{\,2}\end{pmatrix}
=2^{-\delta k}2^{-(2+\delta)k}\mathsf T_B,
\]
which is \eqref{eq:local-leading-covariance}.
\end{proof}

For the primary pre-pulses, the amplitude used in the ansatz of
Section~\ref{s6n:ansatz} is the matched amplitude from
\eqref{s6n:amplitude-data}:
\begin{equation}
\begin{aligned}
\mathbf h_{\mathrm{match},\gamma}
 &=2^{-\delta k/2}g_{\mathrm{match},\gamma}\mathbf h_\gamma,\\
(\mathbf h_{\mathrm{match},\gamma})^r
 &=2^{-\delta k/2}g_{\mathrm{match},\gamma}h_\gamma^r.
\end{aligned}
\label{eq:model-radial-amplitude-identified}
\end{equation}
The scalar \(g_{\mathrm{match},\gamma}\) chooses how much of each flux
direction is needed.  Multiplying each summand by the cutoffs in the normalized variables and restoring the physical
velocity scale gives the pre-pulse:
\[
 w_\gamma^{\mathrm{pre}}
 =2^{(1+\delta)k}\zeta_k\zeta_{k,\alpha}
   2^{-\delta k/2}g_{\mathrm{match},\gamma}b_\gamma.
\]
Thus the cutoff pair \(\zeta_k\zeta_{k,\alpha}W_{k,\alpha}\)
is the normalized sum of the two pre-pulses with labels
\((k,\alpha,+)\) and \((k,\alpha,-)\).
The two \(m=\pm1\) harmonic coefficients of
\(2^{-(1+\delta)k}w_\gamma^{\mathrm{pre}}\) are
\[
 \frac12\chi_\gamma\mathbf h_{\mathrm{match},\gamma}
 =\frac12\,2^{-\delta k/2}\zeta_k\zeta_{k,\alpha}
   \Phi_r\Phi_t\,g_{\mathrm{match},\gamma}\mathbf h_\gamma.
\]
The profiles \(\mathbf h_\gamma\) are independent of \(l_r\), and
\(\tau_\gamma\) is affine in \(Y_i\) with gradient \(O(k^2)\).
Using \eqref{s6n:homogeneous-derivatives}, the weighted square-root
estimate above, and the cutoff bounds, the product rule gives, for
every multi-index \(\beta\in\ZZ_{\ge0}^5\),
\begin{equation}
\left|\partial_{\hat R,\hat Z,\hat T,Y_{i,1},Y_{i,2}}^{\beta}\!\left[
2^{-\delta k/2}\zeta_k\zeta_{k,\alpha}\Phi_r\Phi_t
g_{\mathrm{match},\gamma}\mathbf h_\gamma\right]\right|
\le C_{\beta}2^{-\delta k/2}k^{b_{\beta}}
\sqrt\omega\,d_{\mathrm{edge}}^{-d_\beta}\mathcal G.
\label{eq:matched-prepulse-coefficient-bounds}
\end{equation}
These coefficients have the prescribed cutoff support, vanish for
\(R\notin[R_a,R_h]\), and extend smoothly by zero across all support
boundaries.
The factor \(2^{-\delta k/2}\) gives the required physical stress scale:
\[
2^{2(1+\delta)k}\bigl(2^{-\delta k/2}\bigr)^2
 =2^{(2+\delta)k}\asymp\lambda^{-2-\delta}.
\]
The factor \(k\) in \(g_{\mathrm{match},\gamma}^2\) compensates
for the auxiliary mass \(J_\gamma\asymp k^{-1}\) from
\eqref{eq:covariance-column-mass}.
To assemble different bands, we can use the common auxiliary variable
\(Y_{i_0}\) of Section~\ref{sec:pulse-common}. A native profile
\(f(Y_i)\) is pulled back as \(f(E^{i-i_0}Y_{i_0})\); in particular,
\(b_\gamma\) is pulled back on the full preimage of its rectangle.
Equation~\eqref{eq:common-torus-average} gives the same normalized
average in \(Y_i\) and \(Y_{i_0}\). Lemma~\ref{lem:common-compatibility}
shows that the averaged quadratic fluxes also agree on overlaps.
For assembled fields, \(\mathsf C\) denotes the average in
\(Y_{i_0}\); this agrees with the native average for each summand.
We assemble the pre-pulses in physical velocity units over all labels,
keeping this auxiliary variable independent:

\begin{equation}
\begin{aligned}
w_0^{\mathrm{pre}}
&=\sum_\gamma w_\gamma^{\mathrm{pre}}\\
&=\sum_{k,\alpha}2^{(1+\delta)k}
\zeta_k(\lambda)\zeta_{k,\alpha}(\hat{R},\hat{Z},\hat{T})
W_{k,\alpha}.
\end{aligned}
\label{eq:assembled-leading-pulse}
\end{equation}
Each summand uses its own band coordinates. The sum is locally finite
for \(\lambda>0\). The support separation in
Proposition~\ref{s6n:separation} eliminates every product between distinct labels, including
the two signs. For \(0<\lambda\le2^{-k_0}\), the squared partition gives
\begin{equation}
\begin{aligned}
\mathsf C(w_0^{\mathrm{pre}})
&=\sum_{k,\alpha}2^{2(1+\delta)k}
\zeta_k^2\zeta_{k,\alpha}^2\,
\mathsf C(W_{k,\alpha})\\
&=\sum_{k,\alpha}\zeta_k^2\zeta_{k,\alpha}^2\,
\mathsf T_B
=\mathsf T_B.
\end{aligned}
\label{eq:assembled-leading-stress}
\end{equation}
The squared partition counts each \((\hat{R},\hat{Z},\hat{T})\)-box once.
The two signs have already been combined inside \(W_{k,\alpha}\).

By \eqref{eq:assembled-leading-stress} and the cylindrical divergence
formulas above, the averaged radial tangential fluxes of
\(w_0^{\mathrm{pre}}\) cancel the complete background stress divergence:
\begin{equation}
\begin{aligned}
-\left(\partial_r+\frac2r\right)\mathsf T_B^\theta
+\left(\partial_r+\frac2r\right)
 \mathsf C(w_0^{\mathrm{pre}})^\theta&=0,\\
-\left(\partial_r+\frac1r\right)\mathsf T_B^z
+\left(\partial_r+\frac1r\right)
 \mathsf C(w_0^{\mathrm{pre}})^z&=0.
\end{aligned}
\label{eq:leading-stress-cancellation}
\end{equation}
Equation~\eqref{eq:leading-stress-cancellation} concerns the averaged
radial pair. The remaining terms in
\(\nabla\cdot(w_0^{\mathrm{pre}}\otimes w_0^{\mathrm{pre}})\) include, \emph{e.g.},
the axial flux divergences
\[
\partial_z\langle
 (w_0^{\mathrm{pre}})^z(w_0^{\mathrm{pre}})^\theta
 \rangle_{\theta, Y_{i_0}},\qquad
\partial_z\langle((w_0^{\mathrm{pre}})^z)^2
 \rangle_{\theta, Y_{i_0}},
\]
the radial component, and the parts with zero angular or auxiliary mean.

\subsection{Linearized matching of arbitrary stress increments}
\label{subsec:signed-covariance}

After \(\mathsf C(w_0^{\mathrm{pre}})=\mathsf T_B\),
we allow an arbitrary two-component stress increment. The prescribed
increment need not satisfy the cone condition used for positive
matching. We keep the positive amplitudes
\(g_{\mathrm{match},\gamma}\) fixed and linearize the quadratic
map around the fixed primary pre-pulse.
On each reference band, use the common auxiliary variable \(Y_{i_0}\)
from Section~\ref{sec:pulse-common}. For two real extended vector
fields \(U,V\) expressed in this variable, define
\begin{equation}
\begin{aligned}
\mathsf B(U,V)
&=\left\langle
U^r(V^\theta,V^z)+V^r(U^\theta,U^z)
\right\rangle_{\theta, Y_{i_0}},\\
\mathsf C(U+V)&=\mathsf C(U)+\mathsf B(U,V)+\mathsf C(V).
\end{aligned}
\label{eq:covariance-polarization}
\end{equation}
(The order of the two entries is always \((\theta,z)\).) The averages
fix \((\hat{R},\hat{Z},\hat{T})\). For inputs defined on the same native
torus, these averages may equivalently be taken in \(Y_i\), with
\(i=i(k)\), by \eqref{eq:common-torus-average}.

\begin{prop}[A linear inverse for an arbitrary stress increment]
\label{prop:pulse-signed-covariance}
Fix \(k\ge k_0\) and \(s_0\in\RR\). On the band \(\mathcal B_k\),
let \(\Sigma=\Sigma(\hat R,\hat Z,\hat T)\in\RR^2\) be smooth,
supported in the closed active radial annulus with smooth zero extension. Assume that, for every \(\beta\),
there are \(C_\beta>0\) and \(b_\beta,d_\beta\ge0\), independent of
\(k\), such that
\begin{equation}\label{eq:signed-stress-input-bounds}
 |\partial_{\hat R,\hat Z,\hat T}^{\beta}\Sigma|
 \le C_\beta2^{-s_0\delta k}k^{b_\beta}
       \omega d_{\mathrm{edge}}^{-d_\beta}.
\end{equation}
On each \(Q_{k,\alpha}\) from \eqref{s6n:slow-box-cutoffs}, for
\(\gamma=(k,\alpha,\iota)\), \(\iota=\pm1\), define
\begin{equation}
\begin{aligned}
\begin{pmatrix}
d_{\Sigma,(k,\alpha,+)}\\d_{\Sigma,(k,\alpha,-)}
\end{pmatrix}
&=H_{k,\alpha}^{-1}\Sigma,\\
\delta g_{\mathrm{match},\gamma}
&=\frac{d_{\Sigma,\gamma}}{2g_{\mathrm{match},\gamma}},\\
\mathcal L_{k,\alpha}\Sigma
&=2^{-\delta k/2}\sum_{\iota=\pm1}
\delta g_{\mathrm{match},(k,\alpha,\iota)}b_{(k,\alpha,\iota)}.
\end{aligned}
\label{eq:signed-covariance-inverse}
\end{equation}
The quotient is defined in the open active annulus, where
\(g_{\mathrm{match},\gamma}>0\), and extends smoothly by zero across
both radial edges.
Then
\[
2^{\delta k}\mathsf B(W_{k,\alpha},\mathcal L_{k,\alpha}\Sigma)=\Sigma.
\]
\end{prop}

\begin{proof}
The supports of \(b_{(k,\alpha,+)}\) and
\(b_{(k,\alpha,-)}\) are disjoint.
All scalar matching coefficients are independent of the variables being averaged. Thus
\[
\begin{aligned}
&2^{\delta k}\mathsf B(W_{k,\alpha},\mathcal L_{k,\alpha}\Sigma)\\
&\quad=H_{k,\alpha}
\begin{pmatrix}
2g_{\mathrm{match},(k,\alpha,+)}\delta g_{\mathrm{match},(k,\alpha,+)}\\
2g_{\mathrm{match},(k,\alpha,-)}\delta g_{\mathrm{match},(k,\alpha,-)}
\end{pmatrix}\\
&\quad=H_{k,\alpha}H_{k,\alpha}^{-1}\Sigma=\Sigma.
\end{aligned}
\]
The factor \(1/2\) from the average of \(\cos^2(N\varphi_\gamma)\)
is already in \(H_{k,\alpha}\). The new factor \(2\) comes from the derivative of the
square \(g_{\mathrm{match},\gamma}^2\).

The inverse-matrix derivative bounds established in the proof of
Proposition~\ref{prop:leading-positive-matching}, together with
\eqref{eq:signed-stress-input-bounds}, imply
\[
|\partial_{\hat R,\hat Z,\hat T}^{\beta}d_{\Sigma,\gamma}|
\le C_{\beta}2^{-s_0\delta k}k^{b_{\beta}}
\omega d_{\mathrm{edge}}^{-d_{\beta}}.
\]
The zeroth-order bound for the quotient follows directly from
\(g_{\mathrm{match},\gamma}^2\ge ck\omega\).
Write \(y_\gamma=g_{\mathrm{match},\gamma}^2\). Every nonzero derivative
of \(y_\gamma^{-1/2}\) is a sum of terms
\[
C\,y_\gamma^{-j-1/2}
\prod_{\nu=1}^j\partial_{\hat R,\hat Z,\hat T}^{\beta^{(\nu)}}y_\gamma,
\qquad \sum_{\nu=1}^j\beta^{(\nu)}=\beta,\quad |\beta^{(\nu)}|\ge1.
\]
By \eqref{eq:matching-squares-edge-bounds},
\(y_\gamma\ge ck\omega\), and each numerator derivative carries
one factor \(\omega\). Hence these derivatives have weight \(\omega^{-1/2}\)
and only polynomial and edge-distance losses. Multiplying by
\(d_{\Sigma,\gamma}\), whose weight is \(\omega\), gives
\[
|\partial_{\hat R,\hat Z,\hat T}^{\beta}\delta g_{\mathrm{match},\gamma}|
\le C_{\beta}2^{-s_0\delta k}k^{b'_{\beta}}
\sqrt\omega\,d_{\mathrm{edge}}^{-d'_{\beta}}.
\]
Every slow derivative of the quotient tends to zero at both radial edges.
Thus the quotient extends smoothly by zero across these edges.
\end{proof}

The quadratic identity for the normalized pre-pulses is
\[
\begin{aligned}
2^{\delta k}\mathsf C(W_{k,\alpha}+\mathcal L_{k,\alpha}\Sigma)
&=2^{-(2+\delta)k}\mathsf T_B+\Sigma\\
&\quad+2^{\delta k}\mathsf C(\mathcal L_{k,\alpha}\Sigma).
\end{aligned}
\]
Thus the inverse supplies the prescribed stress increment to first order.

We now assemble the local pre-pulses in physical units.
Let \(\mathsf S(r,z,t)\) be the stress pair to be cancelled.
Assume it is defined on the working domain \(0<\lambda\le2^{-k_0}\)
and that, on each band \(k\ge k_0\), the normalized target
\(2^{-(2+\delta)k}\mathsf S\) satisfies the assumptions of
Proposition~\ref{prop:pulse-signed-covariance}, with constants and exponents
independent of \(k\). This region is covered
by the squared partition for \(k\ge k_0\). Set
\begin{equation}
\mathcal L\mathsf S
=\sum_{k,\alpha}2^{(1+\delta)k}\zeta_k\zeta_{k,\alpha}
\mathcal L_{k,\alpha}\bigl(2^{-(2+\delta)k}\mathsf S\bigr).
\label{eq:assembled-signed-inverse}
\end{equation}
The normalized input is recomputed with the scale of each summand.
Support separation and the squared partition give
\[
\begin{aligned}
\mathsf B(w_0^{\mathrm{pre}},
\mathcal L\mathsf S)
&=\sum_{k,\alpha}2^{(2+\delta)k}
\zeta_k^2\zeta_{k,\alpha}^2
\bigl(2^{-(2+\delta)k}\mathsf S\bigr)\\
&=\mathsf S.
\end{aligned}
\]
On every reference band, the same physical input \(\mathsf S\) gives
\begin{equation}
2^{-(2+\delta)k}
 \mathsf B(w_0^{\mathrm{pre}},\mathcal L\mathsf S)=2^{-(2+\delta)k}\mathsf S.
\label{eq:assembled-signed-covariance}
\end{equation}
These identities hold before physical evaluation and agree on overlapping
bands under the coordinate changes of Section~\ref{sec:pulse-common}.

For \(\Sigma\) as in Proposition~\ref{prop:pulse-signed-covariance},
with \(i=i(k)\) and \(\gamma=(k,\alpha,\iota)\), the localized cosine
coefficients and the averaged radial flux \(\mathsf C\) defined in
\eqref{eq:pulse-covariance-pair} satisfy, uniformly over all bands and boxes,
\begin{equation}
\begin{aligned}
&\left|\partial_{\hat R,\hat Z,\hat T,Y_{i,1},Y_{i,2}}^{\beta}\left[
2^{-\delta k/2}\chi_\gamma
\delta g_{\mathrm{match},\gamma}\mathbf h_\gamma\right]\right|\\
&\qquad\le C_\beta2^{-(s_0+1/2)\delta k}k^{b_\beta}
\sqrt\omega\,d_{\mathrm{edge}}^{-d_\beta}\mathcal G_\gamma(\tau_\gamma),\\
&2^{\delta k}\left|\partial_{\hat R,\hat Z,\hat T}^{\beta}
\mathsf C(\zeta_k\zeta_{k,\alpha}\mathcal L_{k,\alpha}\Sigma)\right|\\
&\qquad\le C_\beta2^{-2s_0\delta k}k^{b_\beta}
\omega d_{\mathrm{edge}}^{-d_\beta}.
\end{aligned}
\label{eq:signed-covariance-bounds}
\end{equation}
The coefficients inherit the supports of the cutoffs in
\(\chi_\gamma\) and extend smoothly by zero across their boundaries.

The derivative bound for \(\delta g_{\mathrm{match},\gamma}\) proved in
Proposition~\ref{prop:pulse-signed-covariance}, together with the cutoff
and homogeneous amplitude bounds, gives the first estimate in
\eqref{eq:signed-covariance-bounds} by the product rule.
Applying the product rule to the quadratic flux and averaging gives
the second estimate, since \(2(s_0+1/2)-1=2s_0\).
For \(s_0>0\), the quadratic stress is of lower order than the input
scale \(2^{-s_0\delta k}\) as \(k\to\infty\).
Bounded label overlap and the uniform bounds on the coordinate changes
preserve these derivative estimates in the common auxiliary variable \(Y_{i_0}\).
For the labelled harmonic coefficients of \(2^{-(1+\delta)k}w_0^{\mathrm{pre}}\)
and \(2^{-(1+\delta)k}\mathcal L\mathsf S\), the
\(\partial_{\hat R,\hat Z,\hat T,Y_{i_0,1},Y_{i_0,2}}^\beta\) derivatives
carry the decay factors \(2^{-\delta k/2}\) and
\(2^{-(s_0+1/2)\delta k}\), respectively.
The edge weights and the Gaussian envelope of each label are retained.
Here each native summand uses the target
\(2^{-(2+\delta)\ell}\mathsf S\) at its own scale \(\ell\); the normalization uses the reference scale \(k\),
as in Section~\ref{sec:pulse-common}.

\section{Initialization}
\label{sec:initialization-ppt}

We will correct the background \((u_B,p_B)\) in successive \emph{stages}.
Stage zero is the initialization: it creates the primary pulse \(w_0\)
and then adds \(u_{0,\mathrm{aux}}\), \(u_{0,\mathrm{moment}}\) and
pressure corrections. The later stages use four types of
velocity correction and will be constructed in
Section~\ref{sec:iteration-ppt}.
Many corrections are defined by summing over labels
\(\gamma=(\ell,\alpha,\iota)\) from overlapping scale bands.
We estimate these sums on each reference band \(\mathcal B_k\),
where only terms with \(|\ell-k|\le1\) can contribute.
In these estimates, \((\hat R,\hat Z,\hat T)\) are the normalized
slow coordinates at scale \(k\), unless explicitly stated otherwise.
We first specify the fields and
estimates that these constructions must produce. 

Write
\[
\mathsf{NS}(u,p)
 =\partial_tu+(u\cdot\nabla)u-\Delta u+\nabla p,
\qquad \mathrm{tan}=(\theta,z).
\]
In this section, \(u^{[j]},p^{[j]}\) denote the full updated velocity and
pressure after stage \(j\), including \(u_B,p_B\). Throughout the finite-stage
construction, fields retain the independent auxiliary variable unless
physical evaluation is explicitly stated. On \(\lambda\asymp2^{-k}\), the
reference sizes are
\[
\begin{array}{c|ccc}
\text{quantity}&u_B^{\mathrm{tan}}&\mathsf T_B&f_B^{\mathrm{tan}}\\ \hline
\text{reference size}&2^{(1+\delta)k}&2^{(2+\delta)k}&2^{(3+\delta)k},
\end{array}
\]
where
\[
\begin{aligned}
f_B&:=\mathsf{NS}(u_B,p_B)-\mathsf E_B,\\
f_B^r&=0,\\
f_B^\theta&=-r^{-2}\partial_r(r^2\mathsf T_B^\theta),\\
f_B^z&=-r^{-1}\partial_r(r\mathsf T_B^z).
\end{aligned}
\]
For every subsequent correction step, \(E\) denotes the flat error that may update from step to step,
written as \(E^{[j]}\) after finishing stage \(j\), and \(f\) includes the entire
force left after subtracting \(E\):
\begin{equation}
f:=\mathsf{NS}(u,p)-E,\qquad
f^{[j]}:=\mathsf{NS}(u^{[j]},p^{[j]})-E^{[j]}.
\label{init:force-definition}
\end{equation}
Note that after adding the primary pulse \(w_0\) to the background
\(u_B\), we will estimate the full \(f\) without extracting a stress divergence.

All angular and auxiliary averages below are taken in the common
variable \(Y_{i_0}\) of Section~\ref{sec:pulse-common}, with
\(i_0=i_0(k)\) on \(\mathcal B_k\) ($k \ge k_0$), before physical evaluation.
Explicitly, with normalized Haar measure on \(\mathbb T^2\),
\[
\langle g\rangle_\theta=\frac1{2\pi}\int_0^{2\pi}g\,d\theta,
\qquad
\langle g\rangle_{\theta,Y_{i_0}}
 =\int_{\mathbb T^2}\langle g\rangle_\theta\,dY_{i_0}.
\]
Vector averages are taken componentwise in the cylindrical frame.
For the flat-error estimates, \(D_{t,x}^{\beta}\) denotes compositions
of the physical space--time derivatives on extended fields:
\[
D_qg=\partial_qg+(\partial_qY_{i_0})\cdot\nabla_{Y_{i_0}}g,
\qquad q\in\{t,x_1,x_2,x_3\}.
\]
Here \(\partial_q\) holds \(Y_{i_0}\) fixed, and vector fields are
differentiated in Cartesian components. Thus
\[
(D_{t,x}^{\beta}g)\big|_{Y_{i_0}=Y_{i_0}(t,r)}
 =\partial_{t,x}^{\beta}
       \bigl(g\big|_{Y_{i_0}=Y_{i_0}(t,r)}\bigr).
\]
At each stage we will retain
\[
\begin{gathered}
|D_{t,x}^{\beta}E^{[j]}|\le C_{j,\beta,N}\lambda^N
 \quad\text{for every fixed }\beta\text{ and }N>0,\\
\langle E^{[j]}\rangle_{\theta,Y_{i_0}}=\mathsf E_B,\qquad
E^{[j]}=\mathsf E_B\quad\text{if }R\notin[R_a,R_h].
\end{gathered}
\]
The same properties are required of \(E\) at intermediate steps:
the flat terms added to \(\mathsf E_B\) always have zero double average and are
supported in the annulus.
The derivatives in \(\mathsf{NS}\) are physical derivatives:
when acting on extended fields they include differentiation of
\(Y_{i_0}(t,r)\), as in \eqref{s6n:physical-derivatives} and
Section~\ref{sec:pulse-common}. For a smooth extended field \(g\),
periodicity in \(Y_{i_0}\) gives
\(\int_{\mathbb T^2}\nabla_{Y_{i_0}}g\,dY_{i_0}=0\), so
\[
\begin{aligned}
\langle D_tg\rangle_{Y_{i_0}}
 &=\partial_t\langle g\rangle_{Y_{i_0}},&
\langle D_rg\rangle_{Y_{i_0}}
 &=\partial_r\langle g\rangle_{Y_{i_0}},\\
\langle\partial_\theta g\rangle_{Y_{i_0}}
 &=\partial_\theta\langle g\rangle_{Y_{i_0}},&
\langle\partial_zg\rangle_{Y_{i_0}}
 &=\partial_z\langle g\rangle_{Y_{i_0}}.
\end{aligned}
\]
Pointwise estimates that we perform below are understood as uniform
in the angular and auxiliary variables.

\paragraph{The corrections and their expected sizes.}
The special stage \(j=0\) is called the \emph{initialization}. We seek
\begin{equation}
\begin{aligned}
u^{[0]}&=u_B+w_0+u_{0,\mathrm{aux}}+u_{0,\mathrm{moment}},\\
p^{[0]}&=p_B+p_{w_0}+\widetilde p^{[0]}.
\end{aligned}
\label{init:stage-zero-fields}
\end{equation}
The primary pulse has zero angular mean. The fields
\(u_{0,\mathrm{aux}},u_{0,\mathrm{moment}}\) and
\(\widetilde p^{[0]}\) are axisymmetric. Compared with later stages
\(j\ge1\), no additional \(w_{0,m\ne0}\) or \(w_{0,\mathrm{stress}}\)
is needed in this initial stage.
With \(\omega,d_{\mathrm{edge}}\) from \eqref{intro:background-edge-weight}, the velocity corrector bounds will be
\begin{equation}
\begin{aligned}
&2^{-(1+\delta)k}|w_0|
 \lesssim k^C2^{-0.5\delta k}\sqrt\omega\,d_{\mathrm{edge}}^{-C},\\
&2^{-(1+\delta)k}
 \bigl(|u_{0,\mathrm{aux}}^{\mathrm{tan}}| + |u_{0,\mathrm{moment}}^{\mathrm{tan}}| + 2^{\delta k}|u_{0,\mathrm{aux}}^r| + 2^{\delta k}|u_{0,\mathrm{moment}}^r|\bigr)\\
&\qquad\lesssim k^C2^{-(1-10^{-5})\delta k}
                  \omega d_{\mathrm{edge}}^{-C}.
\end{aligned}
\label{init:initial-increment-sizes}
\end{equation}
All such normalized estimates in this section are understood to hold
on the band \(\mathcal B_k\), uniformly in \(k\).
For \(j\ge1\), the four new velocity corrections and the pressure
increments have the form
\begin{equation}
\begin{aligned}
u^{[j]}&=u^{[j-1]}+w_{j,m\ne0}+w_{j,\mathrm{stress}}
                         +u_{j,\mathrm{aux}}+u_{j,\mathrm{moment}},\\
p^{[j]}&=p^{[j-1]}+p_{j,m\ne0}+p_{j,\mathrm{stress}}
                         +\widetilde p^{[j]}-\widetilde p^{[j-1]}.
\end{aligned}
\label{init:later-stage-fields}
\end{equation}
Here \(w_{j,m\ne0}\) corrects the nonzero angular modes
\(f-\langle f\rangle_\theta\), based on Section~\ref{s7n:forced-pulse}.
The correction \(w_{j,\mathrm{stress}}\) corrects
\(\langle f^{\mathrm{tan}}\rangle_{\theta,Y_{i_0}}\) by its averaged
quadratic products with the primary pulse, based on
Section~\ref{subsec:signed-covariance}. The next correction
\(u_{j,\mathrm{aux}}\) is aimed at reducing
\[
\langle f^{\mathrm{tan}}\rangle_\theta
       -\langle f^{\mathrm{tan}}\rangle_{\theta,Y_{i_0}},
\]
and \(u_{j,\mathrm{moment}}\) reduces the three moment defects defined
below. In each of these correction steps \(f\) is recomputed from the
updated velocity, pressure and flat error by \eqref{init:force-definition}.
The pulse pressures have zero angular mean, whereas
\(\widetilde p^{[j]}\) is the accumulated axisymmetric correction:
\[
\langle p_{w_0}\rangle_\theta
 =\langle p_{j,m\ne0}\rangle_\theta
 =\langle p_{j,\mathrm{stress}}\rangle_\theta=0,
\qquad
\widetilde p^{[j]}=\langle p^{[j]}-p_B\rangle_\theta.
\]
The expected estimates for \(j\ge1\) are
\begin{equation}
\begin{aligned}
&2^{-(1+\delta)k}|w_{j,m\ne0}|
 \lesssim k^{C_j}2^{-(0.6+0.1j)\delta k}
                  \sqrt\omega\,d_{\mathrm{edge}}^{-C_j},\\
&2^{-(1+\delta)k}|w_{j,\mathrm{stress}}|
 \lesssim k^{C_j}2^{-(0.6+0.1j)\delta k}2^{10^{-5}\delta k}
                  \sqrt\omega\,d_{\mathrm{edge}}^{-C_j},\\
&2^{-(1+\delta)k}\sum_{a=\mathrm{aux},\mathrm{moment}}
 \bigl(|u_{j,a}^{\mathrm{tan}}|+2^{\delta k}|u_{j,a}^r|\bigr)\\
&\qquad\lesssim k^{C_j}2^{-(0.6+0.1j)\delta k}
                    2^{-(0.5-2\cdot10^{-5})\delta k}
                    \omega d_{\mathrm{edge}}^{-C_j}.
\end{aligned}
\label{init:later-increment-sizes}
\end{equation}
Implicit constants may depend on $j$, but not on $k$. Each velocity increment is required to be divergence-free. The velocity and
pressure increments will be constructed with compact radial support in the
fixed active annulus;
in particular, at each completed stage,
\begin{equation}
\begin{gathered}
\nabla\cdot u^{[j]}=0,\qquad
(u^{[j]},p^{[j]})=(u_B,p_B)\quad\text{if }R\notin[R_a,R_h],\\
\operatorname{supp}_R f^{[j]}\subset[R_a,R_h].
\end{gathered}
\label{init:stage-support}
\end{equation}
On the working domain
\(0<\lambda\le2^{-k_0}\), \(|Z|\le1\), these support conditions imply
\[
r\le\sqrt{2R_h}\,2^{-k_0},\qquad
|z|=\lambda^{1-\delta}|Z|\le2^{-(1-\delta)k_0}.
\]
Hence all velocity and pressure increments and the forces \(f^{[j]}\)
will have compact spatial support in the same fixed set
\(K\Subset\mathbb R^3\).

\paragraph{The force bounds after each stage.}
For every \(j\ge0\), our targets are
\begin{equation}
\begin{aligned}
&2^{-(3+\delta)k}
 \bigl(|f^{[j]}|+|f^{[j]}-\langle f^{[j]}\rangle_\theta|\bigr)
 \lesssim k^{C_j}2^{(0.3-0.1j)\delta k}
                    \sqrt\omega\,d_{\mathrm{edge}}^{-C_j},\\
&2^{-(3+\delta)k}\Bigl(
 \bigl|\langle(f^{[j]})^{\mathrm{tan}}\rangle_\theta
       -\langle(f^{[j]})^{\mathrm{tan}}\rangle_{\theta,Y_{i_0}}\bigr|
       +|\langle f^{[j]}\rangle_{\theta,Y_{i_0}}|\Bigr)
 \lesssim k^{C_j}2^{-(0.2+0.1j)\delta k}
                    \omega d_{\mathrm{edge}}^{-C_j}.
\end{aligned}
\label{init:iteration-force-targets}
\end{equation}
The implicit constants here may depend on \(j\), but not on \(k\).
The double average in the second line includes all three components:
\[
|\langle f^{[j]}\rangle_{\theta,Y_{i_0}}|^2
 =\sum_{\mu\in\{r,\theta,z\}}
       |\langle(f^{[j]})^\mu\rangle_{\theta,Y_{i_0}}|^2.
\]
Since \(\sqrt\omega\,d_{\mathrm{edge}}^{-C_j}\) is bounded, the first line gives
\[
|f^{[j]}|\lesssim k^{C_j}2^{[3+(1.3-0.1j)\delta]k}
 \asymp k^{C_j}\lambda^{(0.1j-1.3)\delta-3}.
\]
The exponent of \(\lambda\) tends to infinity with \(j\).
Thus, for any prescribed power, a sufficiently large finite stage
will give that decay as \(\lambda\to0\), with a constant depending
on the chosen stage. Passing from these finite stages to the final solution
requires the construction in Section~\ref{sec:iteration-ppt}.

\paragraph{After the primary pulse.}
The first construction will provide \((w_0,p_{w_0})\) with
\begin{equation}
\begin{gathered}
\mathsf{NS}(u_B+w_0,p_B+p_{w_0})=E+f,
 \qquad\langle p_{w_0}\rangle_\theta=0,\\
(u_B+w_0,p_B+p_{w_0})=(u_B,p_B)\quad\text{if }R\notin[R_a,R_h],
 \qquad\operatorname{supp}_R f\subset[R_a,R_h],\\
2^{-(3+\delta)k}|f|\lesssim k^C2^{0.01\delta k}
                    \sqrt\omega\,d_{\mathrm{edge}}^{-C},\\
2^{-(3+\delta)k}|\langle f^{\theta,z}\rangle_{\theta,Y_{i_0}}|
                     \lesssim k^C2^{-0.99\delta k}
                    \omega d_{\mathrm{edge}}^{-C}.
\end{gathered}
\label{init:primary-starting-point}
\end{equation}
This is the input to the remaining initialization steps, before
\(u_{0,\mathrm{aux}},u_{0,\mathrm{moment}},\widetilde p^{[0]}\)
have been added. For example, a term contributing to the tangential
double mean of \(f\) is
\[
2^{-(3+\delta)k}
\left|\partial_z\langle w_0^zw_0^\theta\rangle_{\theta,Y_{i_0}}\right|
 \lesssim k^C2^{-\delta k}.
\]
The target in \eqref{init:primary-starting-point} allows the slightly
weaker factor \(2^{-0.99\delta k}\). The radial force can still have
reference size \(2^{(3+\delta)k}\): it contains, for example,
\[
-\frac{\langle(w_0^\theta)^2\rangle_\theta}{r},\qquad
\frac{\langle|w_0|^2\rangle_\theta}{r}
 \lesssim k^C2^{2(1+\delta)k-\delta k+k}
 =k^C2^{(3+\delta)k}.
\]

\paragraph{The moments and their defects.}
For an updated divergence-free velocity \(u\), with
\((u-u_B,p-p_B)\) compactly supported in the radial annulus, use
\(f=\mathsf{NS}(u,p)-E\) and define
\begin{equation}
\mathcal P=\int_0^\infty
                      \langle f^r\rangle_{\theta,Y_{i_0}}\,dr.
\label{init:pressure-moment}
\end{equation}
We will call \(\mathcal P\) the pressure defect, which quantitatively measures the obstruction in cancelling the radial $\theta$-mean force with pressure reconstruction (see Section \ref{init:pressure-reconstruction}).

We also need two physical flux moments with the following
formal expressions. Here ``physical'' and the superscript \(\mathrm{phys}\)
refer to the unnormalized variables and units; auxiliary averages are
still taken before physical evaluation. (In fact, we use only their differences from the background, without requiring convergence
of the individual integrals.) 
\begin{equation}
\begin{aligned}
\langle M^{\theta z,\mathrm{phys}}\rangle_{\theta,Y_{i_0}}
 &=\int_0^\infty r^2\langle u^\theta u^z\rangle_{\theta,Y_{i_0}}\,dr,\\
\langle M^{z\theta,\mathrm{phys}}\rangle_{\theta,Y_{i_0}}
 &=\int_0^\infty r\left\langle
          (u^z)^2+p+\frac r2 f^r\right\rangle_{\theta,Y_{i_0}}\,dr.
\end{aligned}
\label{init:flux-moments}
\end{equation}
The formal background moments
\(M_B^{\theta z,\mathrm{phys}},M_B^{z\theta,\mathrm{phys}}\)
are given by \eqref{init:flux-moments} with
\((u,p,f)=(u_B,p_B,f_B)\). We point out that after removing the respective similarity factors \(\lambda^{1-2\delta}\) and \(\lambda^{-2\delta}\), the formal expansions of these background flux moments recover the leading total moments \(M^{\theta z}(\infty,Z),M^{z\theta}(\infty,Z)\) and the lower-order moments used in Part~I \cite{PartI}.

We define the physical moment defects directly by
\[
\begin{aligned}
\mathcal J^{\theta z}
 &:=
       \langle M^{\theta z,\mathrm{phys}}-M_B^{\theta z,\mathrm{phys}}
                                      \rangle_{\theta,Y_{i_0}}\\
 &=\int_0^\infty r^2
       \langle u^\theta u^z-u_B^\theta u_B^z\rangle_{\theta,Y_{i_0}}\,dr,\\
\mathcal J^{z\theta}
 &:=
       \langle M^{z\theta,\mathrm{phys}}-M_B^{z\theta,\mathrm{phys}}
                                      \rangle_{\theta,Y_{i_0}}\\
 &=\int_0^\infty r\Bigl\langle
       (u^z)^2-(u_B^z)^2+p-p_B
       +\frac r2(f^r-f_B^r)\Bigr\rangle_{\theta,Y_{i_0}}\,dr.
\end{aligned}
\]
The differences \(u-u_B,p-p_B,f-f_B\) are supported in
\(R\in[R_a,R_h]\), so both integrands are compactly supported and
these integrals are finite. These physical defects do not depend on
the reference band $k$. On \(\mathcal B_k\), the estimates below use the
normalized quantities
\[
2^{-2(1+\delta)k}\mathcal P,\qquad
2^{(1-2\delta)k}\mathcal J^{\theta z},\qquad
2^{-2\delta k}\mathcal J^{z\theta},
\]
with the factors written explicitly. Derivatives of these defects are
written as \(\partial_{\hat T,\hat Z}^{\beta}\), with \(k\) fixed.
The bounds for these defects below are also required after every fixed
such derivative, with constants and powers of \(k\) allowed to depend
on \(\beta\).

During the construction, we require the two linear velocity moments
to remain equal to their background values at each \(t,z\):
\begin{equation}
\int_0^\infty r^2\langle u^\theta-u_B^\theta\rangle_{\theta,Y_{i_0}}\,dr=0,
\qquad
\int_0^\infty r\langle u^z-u_B^z\rangle_{\theta,Y_{i_0}}\,dr=0.
\label{init:linear-moments}
\end{equation}
To derive the corresponding weighted force identities, first write the
two tangential components of \(\mathsf{NS}(u,p)=f+E\) after
double averaging. Incompressibility puts the transport terms in
divergence form, and periodicity gives
\(\langle\partial_\theta g\rangle_{\theta,Y_{i_0}}=0\).
Together with the averaging identities above, this yields
\[
\begin{aligned}
\langle f^\theta+E^\theta\rangle_{\theta,Y_{i_0}}
={}&\partial_t\langle u^\theta\rangle_{\theta,Y_{i_0}}
 +\frac1{r^2}\partial_r
       \bigl(r^2\langle u^ru^\theta\rangle_{\theta,Y_{i_0}}\bigr)
 +\partial_z\langle u^\theta u^z\rangle_{\theta,Y_{i_0}}\\
&-\left(\partial_r^2+\frac1r\partial_r+\partial_z^2-\frac1{r^2}\right)
                 \langle u^\theta\rangle_{\theta,Y_{i_0}},\\
\langle f^z+E^z\rangle_{\theta,Y_{i_0}}
={}&\partial_t\langle u^z\rangle_{\theta,Y_{i_0}}
 +\frac1r\partial_r
       \bigl(r\langle u^ru^z\rangle_{\theta,Y_{i_0}}\bigr)
 +\partial_z\langle(u^z)^2+p\rangle_{\theta,Y_{i_0}}\\
&-\left(\partial_r^2+\frac1r\partial_r+\partial_z^2\right)
                 \langle u^z\rangle_{\theta,Y_{i_0}}.
\end{aligned}
\]
Subtract the corresponding background equations. The flat terms
cancel because
\[
\langle E-\mathsf E_B\rangle_{\theta,Y_{i_0}}=0.
\]
Multiply the \(\theta\) equation by \(r^2\), the \(z\) equation by
\(r\), and integrate over \(0<r<\infty\).
The radial transport terms give
\[
\begin{aligned}
&\int_0^\infty\partial_r
 \bigl(r^2\langle u^ru^\theta-u_B^ru_B^\theta
                         \rangle_{\theta,Y_{i_0}}\bigr)\,dr\\
&\qquad=
 \bigl[r^2\langle u^ru^\theta-u_B^ru_B^\theta
                         \rangle_{\theta,Y_{i_0}}\bigr]_0^\infty=0,\\
&\int_0^\infty\partial_r
 \bigl(r\langle u^ru^z-u_B^ru_B^z
                         \rangle_{\theta,Y_{i_0}}\bigr)\,dr\\
&\qquad=
 \bigl[r\langle u^ru^z-u_B^ru_B^z
                         \rangle_{\theta,Y_{i_0}}\bigr]_0^\infty=0.
\end{aligned}
\]
The boundary terms vanish because \(u-u_B\) is supported in the
radial annulus. For the radial viscosity, integration by parts gives
\[
\begin{aligned}
&\int_0^\infty r^2
 \left(\partial_r^2+\frac1r\partial_r-\frac1{r^2}\right)
          \langle u^\theta-u_B^\theta\rangle_{\theta,Y_{i_0}}\,dr\\
&\qquad=
 \bigl[r^2\partial_r\langle u^\theta-u_B^\theta\rangle_{\theta,Y_{i_0}}
       -r\langle u^\theta-u_B^\theta\rangle_{\theta,Y_{i_0}}
 \bigr]_0^\infty=0,\\
&\int_0^\infty r
 \left(\partial_r^2+\frac1r\partial_r\right)
          \langle u^z-u_B^z\rangle_{\theta,Y_{i_0}}\,dr\\
&\qquad=
 \bigl[r\partial_r\langle u^z-u_B^z\rangle_{\theta,Y_{i_0}}
 \bigr]_0^\infty=0.
\end{aligned}
\]
The time derivatives and the axial viscosity terms vanish by
\eqref{init:linear-moments}, which holds at every \(t,z\):
\[
\begin{aligned}
(\partial_t-\partial_z^2)
 \int_0^\infty r^2
       \langle u^\theta-u_B^\theta\rangle_{\theta,Y_{i_0}}\,dr&=0,\\
(\partial_t-\partial_z^2)
 \int_0^\infty r
       \langle u^z-u_B^z\rangle_{\theta,Y_{i_0}}\,dr&=0.
\end{aligned}
\]
The background force contributes no weighted integral, since
\(\mathsf T_B\) is compactly supported and
\[
\begin{aligned}
\int_0^\infty r^2f_B^\theta\,dr
 &=-\int_0^\infty\partial_r(r^2\mathsf T_B^\theta)\,dr
   =-[r^2\mathsf T_B^\theta]_0^\infty=0,\\
\int_0^\infty rf_B^z\,dr
 &=-\int_0^\infty\partial_r(r\mathsf T_B^z)\,dr
   =-[r\mathsf T_B^z]_0^\infty=0.
\end{aligned}
\]
Only the axial derivatives of the flux differences remain.
Using the definitions of the moment defects and \(f_B^r=0\), we obtain
\begin{equation}
\begin{aligned}
\int_0^\infty r^2\langle f^\theta\rangle_{\theta,Y_{i_0}}\,dr
 &=\partial_z\int_0^\infty r^2
      \langle u^\theta u^z-u_B^\theta u_B^z\rangle_{\theta,Y_{i_0}}\,dr\\
 &=\partial_z\mathcal J^{\theta z},\\
\int_0^\infty r\langle f^z\rangle_{\theta,Y_{i_0}}\,dr
 &=\partial_z\int_0^\infty r
      \langle(u^z)^2-(u_B^z)^2+p-p_B\rangle_{\theta,Y_{i_0}}\,dr\\
 &=\partial_z\Bigl[\mathcal J^{z\theta}
       -\frac12\int_0^\infty r^2
                  \langle f^r\rangle_{\theta,Y_{i_0}}\,dr\Bigr].
\end{aligned}
\label{init:weighted-force-moments}
\end{equation}
The last term in \eqref{init:weighted-force-moments} is useful:
it cancels with the corresponding contribution from pressure
reconstruction.
This cancellation will be used in Section~\ref{sec:iteration-ppt}. The defects $\mathcal J^{\theta z}$ and $\mathcal J^{z \theta}$ quantitatively measures the obstruction in writing the tangential mean force as stress divergence, after pressure reconstruction.

At the end of stage \(j\), write \(\mathcal P^{[j]}\) and
\(\mathcal J^{a,[j]}\),
\(a\in\{\theta z,z\theta\}\), for the pressure and moment defects
defined above with \((u,p,f)=(u^{[j]},p^{[j]},f^{[j]})\).
In addition to the force bounds in \eqref{init:iteration-force-targets},
we also require that after each stage \(j\),
\begin{equation}
\begin{aligned}
&2^{-2(1+\delta)k}|\mathcal P^{[j]}|
 +2^{(1-2\delta)k}|\mathcal J^{\theta z,[j]}|
 +2^{-2\delta k}|\mathcal J^{z\theta,[j]}|\\
&\qquad\lesssim k^{C_j}2^{-(1.2+0.1j)\delta k},\qquad j\ge0.
\end{aligned}
\label{init:iteration-moment-target}
\end{equation}
At the end of each stage, pressure reconstruction also gives
\(\langle(f^{[j]})^r\rangle_\theta=\rho_r\mathcal P^{[j]}\),
where \(\rho_r\) is the fixed normalized profile from
\eqref{init:pressure-bump}.

\paragraph{The target of initialization.}
The first task is to complete \eqref{init:stage-zero-fields} so that
\[
\mathsf{NS}(u^{[0]},p^{[0]})=E^{[0]}+f^{[0]},\qquad
\nabla\cdot u^{[0]}=0,
\]
with equality to the background outside \([R_a,R_h]\) and
\(\operatorname{supp}_R f^{[0]}\subset[R_a,R_h]\).
The three required estimates, on \(\lambda\asymp2^{-k}\), are
\begin{equation}
\begin{aligned}
&2^{-(3+\delta)k}|f^{[0]}-\langle f^{[0]}\rangle_\theta|
                         \lesssim k^C2^{0.3\delta k}
                         \sqrt\omega\,d_{\mathrm{edge}}^{-C},\\
&2^{-(3+\delta)k}\Bigl(
 |\langle f^{[0]}\rangle_{\theta,Y_{i_0}}|
 +\bigl|\langle(f^{[0]})^{\theta,z}\rangle_\theta
       -\langle(f^{[0]})^{\theta,z}\rangle_{\theta,Y_{i_0}}\bigr|\Bigr)
                         \lesssim k^C2^{-0.2\delta k}
                         \omega d_{\mathrm{edge}}^{-C},\\
&2^{-2(1+\delta)k}|\mathcal P^{[0]}|
 +2^{(1-2\delta)k}|\mathcal J^{\theta z,[0]}|
 +2^{-2\delta k}|\mathcal J^{z\theta,[0]}|\\
&\qquad\lesssim k^C2^{-1.2\delta k}.
\end{aligned}
\label{init:stage-zero-target}
\end{equation}
We will first add the primary pulses and establish
\eqref{init:primary-starting-point}. These estimates imply stronger
bounds than those required in \eqref{init:stage-zero-target} for the
nonzero angular part \(f-\langle f\rangle_\theta\) and the tangential
double mean \(\langle f^{\mathrm{tan}}\rangle_{\theta,Y_{i_0}}\), since
\[
0.01<0.3,\qquad 0.99>0.2.
\]
This is why we do not need the corrections \(w_{0,m\ne0}\) and
\(w_{0,\mathrm{stress}}\) for the initialization.
It remains to construct \(u_{0,\mathrm{aux}}\),
\(u_{0,\mathrm{moment}}\) and \(\widetilde p^{[0]}\) so that the full
mean and moment bounds in \eqref{init:stage-zero-target} hold, while
preserving the bounds for \(f-\langle f\rangle_\theta\) and
\(\langle f^{\mathrm{tan}}\rangle_{\theta,Y_{i_0}}\) required in
\eqref{init:stage-zero-target}, together with the support conditions.

\subsection{Estimates after adding the primary pulse}
\label{init:primary-pulse}

We now prove the bounds in \eqref{init:primary-starting-point}.
Use the unforced amplitudes \((\mathbf h_\gamma,\pi_\gamma)\) of
Proposition~\ref{prop:homogeneous-pulse} and the weights
\(g_{\mathrm{match},\gamma}\) of
Proposition~\ref{prop:leading-positive-matching}. For a label
\(\gamma=(\ell,\alpha,\iota)\), all coefficients for the pulse associated with this label are defined using its
own scale \(\ell\). With the cutoff \(\chi_\gamma\) from
Section~\ref{s6n:ansatz}, set
\begin{equation}
\begin{gathered}
N_\ell=\lceil2^{\delta\ell/2}\rceil,\qquad
\mathbf n_\gamma=2^{-\ell}\nabla\varphi_\gamma,\\
\mathbf a_\gamma=2^{-\delta\ell/2}\chi_\gamma
                   g_{\mathrm{match},\gamma}\mathbf h_\gamma,
\qquad
\pi_{0,\gamma}=2^{-\delta\ell/2}\chi_\gamma
                   g_{\mathrm{match},\gamma}\pi_\gamma,\\
\mathbf b_\gamma=-\frac{2^{-\ell}}{N_\ell}\nabla\times
       \left(\frac{\mathbf n_\gamma\times\mathbf a_\gamma}
                        {|\mathbf n_\gamma|^2}\right).
\end{gathered}
\label{init:primary-coefficients}
\end{equation}
Here and below the curl and gradient are physical derivatives,
including the auxiliary chain rule and the action on cylindrical frame.
The construction \eqref{s6n:curl-completion} with matching factors included gives
\begin{equation}
\begin{aligned}
w_0
 &=\sum_\gamma\nabla\times\left[
 -\frac{2^{\delta\ell}}{N_\ell}
 \frac{\mathbf n_\gamma\times\mathbf a_\gamma}{|\mathbf n_\gamma|^2}
                 \sin(N_\ell\varphi_\gamma)\right]\\
 &=\sum_\gamma2^{(1+\delta)\ell}
       \bigl(\mathbf a_\gamma\cos(N_\ell\varphi_\gamma)
             +\mathbf b_\gamma\sin(N_\ell\varphi_\gamma)\bigr),\\
p_{w_0}
 &=\sum_\gamma2^{2(1+\delta)\ell}
             \operatorname{Re}\bigl(\pi_{0,\gamma}
                                    e^{\mathrm iN_\ell\varphi_\gamma}\bigr).
\end{aligned}
\label{init:primary-fields}
\end{equation}
The cosine sum is the assembled pre-pulse
\(w_0^{\mathrm{pre}}\) defined in \eqref{eq:assembled-leading-pulse},
and the sine sum is the curl correction.
The sum over \(\gamma\) means
\[
\sum_\gamma
 =\sum_{\ell\ge k_0}\sum_{\alpha\in\mathbb Z^3}
                         \sum_{\iota\in\{+,-\}},
\]
where \(\ell\) labels the native scale band, \(\alpha\) the slow box, and
\(\iota\) the two pulse families. Each summand is extended smoothly
by zero off its cutoff support. On \(\mathcal B_k\), only neighboring scales \(|\ell-k|\le1\) occur.

\begin{prop}[The force after the primary pulse]
\label{init:primary-pulse-estimates}
For the fields in \eqref{init:primary-fields}, there is a flat error
\(E\) such that
\[
\mathsf{NS}(u_B+w_0,p_B+p_{w_0})=E+f,\qquad
\nabla\cdot w_0=0,\qquad
\langle w_0\rangle_\theta=\langle p_{w_0}\rangle_\theta=0.
\]
The fields \(w_0,p_{w_0},f,E-\mathsf E_B\) are real, smooth and supported
in \(R\in[R_a,R_h]\). Moreover, on each fixed compact profile range,
\[
\langle E\rangle_\theta=\mathsf E_B,\qquad
|D_{t,x}^{\beta}E|\le C_{\beta,N}\lambda^N
\quad\text{for every fixed }\beta\text{ and }N>0.
\]
Uniformly on \(\mathcal B_k\),
\begin{equation}
\begin{aligned}
2^{-(1+\delta)k}|w_0|&\lesssim k^C2^{-0.5\delta k}
                         \sqrt\omega\,d_{\mathrm{edge}}^{-C},\\
2^{-(3+\delta)k}|f|&\lesssim k^C2^{0.01\delta k}
                         \sqrt\omega\,d_{\mathrm{edge}}^{-C},\\
2^{-(3+\delta)k}
 |\langle f^{\mathrm{tan}}\rangle_{\theta,Y_{i_0}}|
 &\lesssim k^C2^{-0.99\delta k}\omega d_{\mathrm{edge}}^{-C}.
\end{aligned}
\label{init:primary-pointwise-bounds}
\end{equation}
The angular mean, needed for the subsequent auxiliary and pressure
corrections, is given exactly by
\begin{equation}
\begin{aligned}
\langle f^r\rangle_\theta
 &=\left(D_r+\frac1r\right)\langle(w_0^r)^2\rangle_\theta
   +\partial_z\langle w_0^rw_0^z\rangle_\theta
   -\frac{\langle(w_0^\theta)^2\rangle_\theta}{r},\\
\langle f^\theta\rangle_\theta
 &=\left(D_r+\frac2r\right)
        \bigl(\langle w_0^rw_0^\theta\rangle_\theta-\mathsf T_B^\theta\bigr)
   +\partial_z\langle w_0^zw_0^\theta\rangle_\theta,\\
\langle f^z\rangle_\theta
 &=\left(D_r+\frac1r\right)
        \bigl(\langle w_0^rw_0^z\rangle_\theta-\mathsf T_B^z\bigr)
   +\partial_z\langle(w_0^z)^2\rangle_\theta.
\end{aligned}
\label{init:primary-angular-force}
\end{equation}
\end{prop}

\begin{proof}
\emph{Coefficient bounds.}
Fix a reference band \(\mathcal B_k\). For every active label
\(\gamma=(\ell,\alpha,\iota)\) contributing to this band, \(\ell\) is its native scale and
\(|\ell-k|\le1\). All bounds below are stated on \(\mathcal B_k\)
in terms of \(k\); the exact construction formulas retain \(\ell\).
The coefficient bounds for a pulse associated with $\gamma$ are
\begin{equation}
\begin{aligned}
|\partial_{\hat R,\hat Z,\hat T,Y_{i_0,1},Y_{i_0,2}}^\beta\mathbf a_\gamma|
 &\lesssim k^{C_\beta}2^{-\delta k/2}
                  \sqrt\omega\,d_{\mathrm{edge}}^{-C_\beta}
                  \mathcal G_\gamma(\tau_\gamma),\\
|\partial_{\hat R,\hat Z,\hat T,Y_{i_0,1},Y_{i_0,2}}^\beta\pi_{0,\gamma}|
 &\lesssim k^{C_\beta}2^{-\delta k}
                  \sqrt\omega\,d_{\mathrm{edge}}^{-C_\beta}
                  \mathcal G_\gamma(\tau_\gamma),\\
|\partial_{\hat R,\hat Z,\hat T,Y_{i_0,1},Y_{i_0,2}}^\beta\mathbf b_\gamma|
 &\lesssim k^{C_\beta}2^{-(1-10^{-5})\delta k}
                  \sqrt\omega\,d_{\mathrm{edge}}^{-C_\beta}
                  \mathcal G_\gamma(\tau_\gamma).
\end{aligned}
\label{init:primary-coefficient-bounds}
\end{equation}
Here \(\beta\in\mathbb Z_{\ge0}^5\), and differentiation acts before
physical evaluation, with the slow and auxiliary variables independent
and \(\ell,k,i_0(k)\) fixed. The auxiliary variable is \(Y_{i_0(k)}\),
while the normalized slow coordinates use scale \(\ell\).
Since \(|\ell-k|\le1\), using scale \(k\) instead changes only the
implicit constants.
The envelope \(\mathcal G_\gamma(\tau_\gamma)\) from
\eqref{s6n:gaussian-envelope} belongs to the individual label;
the carrier is left undifferentiated. We sometimes suppress this factor
in the intermediate scale calculations below, but retain it in the
coefficient bounds at every fixed derivative order.
The bound for \(\mathbf a_\gamma\) follows from
\eqref{eq:matched-prepulse-coefficient-bounds}.
The pressure estimate \eqref{s6n:homogeneous-derivatives} has an
additional factor \(2^{-\delta\ell/2}\); Leibniz' rule, the same weighted
matching factor bounds, and the formulas in
\eqref{init:primary-coefficients} give the bound for \(\pi_{0,\gamma}\).

For \(\mathbf b_\gamma\), Lemma~\ref{s7n:phase-bounds} gives polynomial
bounds for the derivatives of \(\mathbf n_\gamma/|\mathbf n_\gamma|^2\)
in these variables.
By \eqref{s6n:normalized-derivatives} and \eqref{eq:ck-size}, the
largest loss in \(2^{-\ell}\nabla\times\) is \(2^{10^{-5}\delta k}\)
from the radial derivative. The bounded difference
\(0\le i(\ell)-i_0\le2\) and Leibniz' rule give the same loss after
any fixed order of differentiation in these variables. Thus \eqref{init:primary-coefficients} and
\(N_\ell^{-1}\lesssim2^{-\delta k/2}\) give
\[
\begin{aligned}
|\partial_{\hat R,\hat Z,\hat T,Y_{i_0,1},Y_{i_0,2}}^\beta\mathbf b_\gamma|
 &\lesssim k^{C_\beta}N_\ell^{-1}
       2^{10^{-5}\delta k}2^{-\delta k/2}
       \sqrt\omega\,d_{\mathrm{edge}}^{-C_\beta}\\
 &\lesssim k^{C_\beta}2^{-(1-10^{-5})\delta k}
       \sqrt\omega\,d_{\mathrm{edge}}^{-C_\beta},
\end{aligned}
\]
for every fixed multi-index \(\beta\).
Summing over \(|\ell-k|\le1\) gives the velocity bound in
\eqref{init:primary-pointwise-bounds}.
Clearly, \(w_0\) is divergence-free since it is a curl.
Also, \(N_\ell p_{\theta,\gamma}\in\mathbb Z\setminus\{0\}\) gives
the zero angular means $\langle w_0\rangle_\theta=\langle p_{w_0}\rangle_\theta=0$ stated in
Proposition~\ref{init:primary-pulse-estimates}.

\emph{The flat error.}
The matching factor $g_{\mathrm{match}, \gamma}$ depends only on the normalized coordinates $\hat{R}, \hat{Z}, \hat{T}$ and satisfies
\[
\partial_{\tau_\gamma}g_{\mathrm{match},\gamma}=0.
\]
Consequently, the principal equation \eqref{s6n:principal-ode} implies
\begin{equation}
\begin{aligned}
&(\partial_{\tau_\gamma}+\mathsf K_{B,\ell}+d_\ell)\mathbf a_\gamma
       +\mathrm iN_\ell\mathbf n_\gamma\pi_{0,\gamma}\\
&\qquad=2^{-\delta\ell/2}(\partial_{\tau_\gamma}\chi_\gamma)
                    g_{\mathrm{match},\gamma}\mathbf h_\gamma,
\qquad d_\ell=2^{-\delta\ell}N_\ell^2|\mathbf n_\gamma|^2.
\end{aligned}
\label{init:primary-principal-cancellation}
\end{equation}
The matrix \(\mathsf K_{B,\ell}\) is the matrix in
\eqref{s6n:principal-ode} at the native scale \(\ell\) of \(\gamma\).
Assign the right-hand
side to the flat error by setting
\begin{equation}
E=\mathsf E_B+\sum_\gamma2^{(3+2\delta)\ell-\delta\ell/2}
 (\partial_{\tau_\gamma}\chi_\gamma)
 g_{\mathrm{match},\gamma}\mathbf h_\gamma
                         \cos(N_\ell\varphi_\gamma).
\label{init:primary-flat-error}
\end{equation}
On the support of \(\partial_{\tau_\gamma}\chi_\gamma\),
\[
\left|\tau_\gamma-\frac{\tau_{\max}(\ell)}2\right|
 \ge\frac{\tau_{\max}(\ell)}5,
\qquad \mathcal G_\gamma(\tau_\gamma)\le e^{-ck^2},
\]
by \eqref{s6n:gaussian-bounds} and \(\tau_{\max}(\ell)\asymp k^2\).
Combining \eqref{s6n:homogeneous-derivatives} and
\eqref{eq:matching-square-root-bounds} with the cutoff bounds and
the physical chain rule in
\eqref{s6n:physical-derivatives}--\eqref{s6n:normalized-derivatives}
gives, on \(\mathcal B_k\),
\[
|D_{t,x}^\beta(E-\mathsf E_B)|
 \le C_\beta2^{C_\beta k}k^{C_\beta}e^{-ck^2}
 \le C_{\beta,N}\lambda^N\qquad(N>0).
\]
It has zero angular mean and the same radial support as the pulse.
This proves the assertions about \(E\).

\emph{The full force.}
Using \(\mathsf L_B\) defined immediately after
\eqref{intro:linear-quadratic}, the exact update identity gives
\begin{equation}
\begin{aligned}
f={}&-\left(\partial_r+\frac2r\right)\mathsf T_B^\theta e_\theta
     -\left(\partial_r+\frac1r\right)\mathsf T_B^z e_z\\
 &+\mathsf L_B(w_0,p_{w_0})-(E-\mathsf E_B)+(w_0\cdot\nabla)w_0.
\end{aligned}
\label{init:primary-force-expansion}
\end{equation}
We next estimate
\(\mathsf L_B(w_0,p_{w_0})-(E-\mathsf E_B)\) and \((w_0\cdot\nabla)w_0\).
The first expression is the remaining linear force after the cutoff
contribution \(E-\mathsf E_B\) has been assigned to the flat error $E$; the second is the quadratic self-interaction of the pulse.

First consider the linear term contributed by one label
\(\gamma=(\ell,\alpha,\iota)\) with \(|\ell-k|\le1\). In the following
calculation, hatted coordinates and quantities use the native scale
\(\ell\). Put
\[
\mathbf c_{\gamma,m}
 =\frac{\mathbf a_\gamma-\mathrm im\mathbf b_\gamma}{2},
\qquad m=\pm1.
\]
The pressure formula \eqref{s6n:projected-ode} gives
\(\pi_{0,\gamma}\in\mathrm i\mathbb R\). We have
\[
\begin{aligned}
\mathbf a_\gamma\cos(N_\ell\varphi_\gamma)
 +\mathbf b_\gamma\sin(N_\ell\varphi_\gamma)
 &=\sum_{m=\pm1}\mathbf c_{\gamma,m}
                         e^{\mathrm imN_\ell\varphi_\gamma},\\
\operatorname{Re}(\pi_{0,\gamma}e^{\mathrm iN_\ell\varphi_\gamma})
 &=\sum_{m=\pm1}\frac m2\pi_{0,\gamma}
                         e^{\mathrm imN_\ell\varphi_\gamma}.
\end{aligned}
\]
Let \(\mathbf v\) denote either \(\mathbf a_\gamma\) or \(\mathbf b_\gamma\).
The chain rule for the coefficients gives
\[
2^{-(2+\delta)\ell}D_t\mathbf v
 =\partial_{\tau_\gamma}\mathbf v
            -2^{-\delta\ell}\partial_{\hat T}\mathbf v,
\qquad
2^{-\ell}\partial_z\mathbf v
 =2^{-\delta\ell}\partial_{\hat Z}\mathbf v,
\]
\[
\widehat D_r\mathbf v:=2^{-\ell}D_r\mathbf v
 =\sqrt{2\hat R}\,\partial_{\hat R}\mathbf v
  +2^{-d_0\ell}\Lambda_r^{i(\ell)}d_0(2\hat R)^{(d_0-1)/2}
          v_r\cdot\nabla_{Y_{i(\ell)}}\mathbf v.
\]
Note that here \(\partial_{\tau_\gamma}\) costs only a power of \(k\)
for pointwise estimates, whereas each \(2^{-\ell}D_r\) can cost
\(2^{10^{-5}\delta k}\). Write the phase transport defect as
\[
\mathfrak d_\gamma
 =2^{-(2+\delta)\ell}(D_t+u_B\cdot\nabla)\varphi_\gamma,
\qquad
|\mathfrak d_\gamma|\lesssim k^C2^{-\delta k},
\]
as in \eqref{s6n:transport-defect}. Using the definition of $\mathsf L_B(w_0, p_{w_0})$
and the equation
\eqref{init:primary-principal-cancellation} gives
\begin{equation}
\mathsf L_B(w_0,p_{w_0})-(E-\mathsf E_B)
 =\sum_{\gamma=(\ell,\alpha,\iota)}2^{(3+\delta)\ell}
       \sum_{m=\pm1}\mathbf f_{\gamma,m}
                            e^{\mathrm imN_\ell\varphi_\gamma},
\label{init:primary-linear-harmonics}
\end{equation}
where, writing \(\hat r=2^\ell r=\sqrt{2\hat R}\) for the normalized radius, the coefficients are given exactly by
\begin{equation}
\begin{aligned}
2^{-\delta\ell}\mathbf f_{\gamma,m}
={}&-\frac{\mathrm im}{2}
  (\partial_{\tau_\gamma}+\mathsf K_{B,\ell}+d_\ell)\mathbf b_\gamma\\
&+\bigl(-2^{-\delta\ell}\partial_{\hat T}
       +\hat u_B^r\widehat D_r
       +\hat u_B^z2^{-\delta\ell}\partial_{\hat Z}\bigr)
                                             \mathbf c_{\gamma,m}\\
&+\mathrm imN_\ell\mathfrak d_\gamma\mathbf c_{\gamma,m}\\
&+c_{\gamma,m}^z2^{-\delta\ell}\partial_{\hat Z}\hat u_B
  +(c_{\gamma,m}^r\widehat D_r\hat u_B^r,
                  \hat r^{-1}\hat u_B^r c_{\gamma,m}^\theta,0)\\
&+\frac m2\bigl(\widehat D_r\pi_{0,\gamma},0,
                    2^{-\delta\ell}\partial_{\hat Z}\pi_{0,\gamma}\bigr)\\
&-2^{-\delta\ell}\Bigl[
  \bigl(\widehat D_r^2+\hat r^{-1}\widehat D_r
              +2^{-2\delta\ell}\partial_{\hat Z}^2\bigr)\mathbf c_{\gamma,m}
  -\hat r^{-2}(c_{\gamma,m}^r,c_{\gamma,m}^\theta,0)\Bigr]\\
&-\mathrm im2^{-\delta\ell}N_\ell\Bigl[
  2\bigl(n_\gamma^r\widehat D_r
              +n_\gamma^z2^{-\delta\ell}\partial_{\hat Z}\bigr)
                                             \mathbf c_{\gamma,m}\\
&\hspace{33mm}
 +\bigl(\widehat D_r n_\gamma^r+\hat r^{-1}n_\gamma^r
              +2^{-\delta\ell}\partial_{\hat Z}n_\gamma^z\bigr)
                                             \mathbf c_{\gamma,m}\\
&\hspace{33mm}
 +2\hat r^{-1}n_\gamma^\theta
                   (-c_{\gamma,m}^\theta,c_{\gamma,m}^r,0)\Bigr].
\end{aligned}
\label{init:primary-linear-coefficients}
\end{equation}
The seven groups in \eqref{init:primary-linear-coefficients} are all supported in the active annulus and have
the following pointwise upper bounds. Each right-hand entry is to be
multiplied by
\(k^C\sqrt\omega\,d_{\mathrm{edge}}^{-C}\mathcal G_\gamma(\tau_\gamma)\).
In the verification of these seven bounds below, we suppress the common
radial and Gaussian factors
\(\sqrt\omega\,d_{\mathrm{edge}}^{-C}\mathcal G_\gamma(\tau_\gamma)\).
\(\mathbf v\) denotes either \(\mathbf a_\gamma\) or
\(\mathbf b_\gamma\).
\[
\renewcommand{\arraystretch}{1.35}
\begin{array}{l|c}
\text{expression to be bounded}&\text{remaining factor}\\ \hline
|(\partial_{\tau_\gamma}+\mathsf K_{B,\ell}+d_\ell)\mathbf b_\gamma|
 &2^{-(1-10^{-5})\delta k}\\
2^{-\delta\ell}(|\partial_{\hat T}\mathbf v|
                  +|\widehat D_r\mathbf v|+|\partial_{\hat Z}\mathbf v|)
 &2^{-(3/2-10^{-5})\delta k}\\
N_\ell|\mathfrak d_\gamma|\,|\mathbf v|
 &2^{-\delta k}\\
2^{-\delta\ell}|\mathbf v|
 &2^{-3\delta k/2}\\
|\widehat D_r\pi_{0,\gamma}|
                +2^{-\delta\ell}|\partial_{\hat Z}\pi_{0,\gamma}|
 &2^{-(1-10^{-5})\delta k}\\
2^{-\delta\ell}\bigl(|\widehat D_r^2\mathbf v|
 +2^{-2\delta\ell}|\partial_{\hat Z}^2\mathbf v|
 +|\widehat D_r\mathbf v|+|\mathbf v|\bigr)
 &2^{-(3/2-2\cdot10^{-5})\delta k}\\
2^{-\delta\ell}N_\ell\bigl(|\widehat D_r\mathbf v|
             +2^{-\delta\ell}|\partial_{\hat Z}\mathbf v|+|\mathbf v|\bigr)
 &2^{-(1-10^{-5})\delta k}
\end{array}
\]
The rows account, respectively, for the curl coefficient's principal
terms, slow coefficient transport, phase transport, the lower-order
background derivatives, the pressure coefficient's slow derivatives,
viscosity acting on coefficients, and viscosity cross terms.

The largest bound is \(k^C2^{-(1-10^{-5})\delta k}\),
attained in the first, fifth and seventh rows. We detail the seventh
row. By Lemma~\ref{s7n:phase-bounds}, \(\mathbf n_\gamma\) and its
derivatives in these variables are bounded by powers of \(k\).
Its auxiliary dependence is only through \(\tau_\gamma\), with
\(D_r\tau_\gamma=0\); hence \(\widehat D_r\mathbf n_\gamma\) also
has a polynomial bound. The radial derivative of
\(\mathbf c_{\gamma,m}\) is the worst term in this group:
\[
\begin{aligned}
2^{-\delta\ell}N_\ell
       |n_\gamma^r\widehat D_r\mathbf c_{\gamma,m}|
 &\lesssim k^C
  {2^{-\delta k}}\quad
  \underbrace{2^{\delta k/2}}_{\text{carrier derivative}}\quad
  \underbrace{2^{10^{-5}\delta k}}_{\text{radial derivative}}\quad
  \underbrace{2^{-\delta k/2}}_{\text{coefficient}}\\
 &=k^C2^{-(1-10^{-5})\delta k}.
\end{aligned}
\]
The terms with \(\partial_{\hat Z}\mathbf c_{\gamma,m}\),
derivatives of \(\mathbf n_\gamma\), or the frame factor
\(\hat r^{-1}\) have the better bound \(k^C2^{-\delta k}\).
For the first row, the polynomial bounds on
\(\mathsf K_{B,\ell},d_\ell\) and the fast-time derivative give
\[
\begin{aligned}
|(\partial_{\tau_\gamma}+\mathsf K_{B,\ell}+d_\ell)\mathbf b_\gamma|
 &\le|\partial_{\tau_\gamma}\mathbf b_\gamma|
       +(|\mathsf K_{B,\ell}|+|d_\ell|)|\mathbf b_\gamma|\\
 &\lesssim k^C2^{-(1-10^{-5})\delta k}.
\end{aligned}
\]
This uses only \eqref{init:primary-coefficient-bounds}, and does
not require \(\mathbf b_\gamma\) to solve the principal equation.
The pressure gradient in the fifth row has the same bound:
\[
|\widehat D_r\pi_{0,\gamma}|
   +2^{-\delta\ell}|\partial_{\hat Z}\pi_{0,\gamma}|
 \lesssim k^C\bigl(2^{10^{-5}\delta k}2^{-\delta k}
                         +2^{-2\delta k}\bigr)
 \lesssim k^C2^{-(1-10^{-5})\delta k}.
\]
The other rows are estimated in the same way and are smaller.
For the second and fourth rows,
\eqref{s6n:background-velocity-bounds} supplies the factor
\(2^{-\delta k}\) in \(\hat u_B^r\) and the background
derivatives outside \(\mathsf K_{B,\ell}\). The third row uses
\(N_\ell|\mathfrak d_\gamma|\,|\mathbf v|
 \lesssim k^C2^{\delta k/2}2^{-\delta k}2^{-\delta k/2}
 =k^C2^{-\delta k}\).
The sixth row has the viscosity factor \(2^{-\delta\ell}\)
and at most two radial derivative losses, giving
\(k^C2^{-(3/2-2\cdot10^{-5})\delta k}\).
Summing the rows and restoring these factors gives
\[
\begin{aligned}
2^{-\delta\ell}|\mathbf f_{\gamma,m}|
 &\lesssim k^C2^{-(1-10^{-5})\delta k}
                  \sqrt\omega\,d_{\mathrm{edge}}^{-C}
                  \mathcal G_\gamma(\tau_\gamma),\\
|\mathbf f_{\gamma,m}|
 &\lesssim k^C2^{10^{-5}\delta k}
                  \sqrt\omega\,d_{\mathrm{edge}}^{-C}
                  \mathcal G_\gamma(\tau_\gamma),
\end{aligned}
\]
where the second line uses \(|\ell-k|\le1\). Using
\(\mathcal G_\gamma\le1\), restoring physical units in
\eqref{init:primary-linear-harmonics}, and summing over the active
labels gives, in particular,
\begin{equation}
|\mathsf L_B(w_0,p_{w_0})-(E-\mathsf E_B)|
 \lesssim k^C2^{(3+\delta)k}2^{3\cdot10^{-5}\delta k}
                  \sqrt\omega\,d_{\mathrm{edge}}^{-C}.
\label{init:primary-linear-bound}
\end{equation}

For the quadratic term, write \(w_{0,\gamma}\) for one summand in
\eqref{init:primary-fields}. Orthogonality and the curl formula give
\(\mathbf n_\gamma\cdot\mathbf a_\gamma=0\) and
\(\nabla\cdot\mathbf b_\gamma=0\). Expanding
\(\nabla\cdot w_{0,\gamma}=0\) therefore gives
\[
N_\ell\mathbf n_\gamma\cdot\mathbf b_\gamma
 =-2^{-\ell}\nabla\cdot\mathbf a_\gamma,
\qquad
\nabla\cdot\mathbf a_\gamma
 =(D_r+r^{-1})a_\gamma^r+\partial_za_\gamma^z.
\]
In particular, the carrier derivative in self-advection satisfies
\[
w_{0,\gamma}\cdot\nabla(N_\ell\varphi_\gamma)
 =-2^{(1+\delta)\ell}(\nabla\cdot\mathbf a_\gamma)
                           \sin(N_\ell\varphi_\gamma).
\]
Thus the large factor \(N_\ell\) is reduced thanks to the inherent structure of the pulse. Expanding \((w_{0,\gamma}\cdot\nabla)w_{0,\gamma}\), including
the cylindrical frame derivatives, now gives
\begin{align*}
|(w_{0,\gamma}\cdot\nabla)w_{0,\gamma}|
&\lesssim 2^{2(1+\delta)k}(|\mathbf a_\gamma|+|\mathbf b_\gamma|)
 \!\Bigl(|D_r\mathbf a_\gamma|+|D_r\mathbf b_\gamma|
       +|\partial_z\mathbf a_\gamma|+|\partial_z\mathbf b_\gamma|
       +r^{-1}(|\mathbf a_\gamma|+|\mathbf b_\gamma|)\Bigr)\\
&\lesssim k^C2^{2(1+\delta)k}
       2^{-\delta k/2}2^k2^{-(1/2-10^{-5})\delta k}
       \omega d_{\mathrm{edge}}^{-C}\\
&=k^C2^{(3+\delta)k}2^{10^{-5}\delta k}
       \omega d_{\mathrm{edge}}^{-C}.
\end{align*}
Products between different labels vanish by
the support separation in Proposition~\ref{s6n:separation}. The background force in stress divergence form in
\eqref{init:primary-force-expansion} is bounded by
\(k^C2^{(3+\delta)k}\omega d_{\mathrm{edge}}^{-C}\), by
\eqref{intro:weighted-background-stress}. Combining these estimates with
\eqref{init:primary-linear-bound}, \(\omega\le\sqrt\omega\), and
\(3\cdot10^{-5}<0.01\) gives
the full-force bound in \eqref{init:primary-pointwise-bounds}.
The bounds for \(\mathbf f_{\gamma,m}\) in
\eqref{init:primary-linear-harmonics}, where the flat error has been
subtracted, and for the nonzero angular harmonic coefficients of
\((w_{0,\gamma}\cdot\nabla)w_{0,\gamma}\) retain
\(\mathcal G_\gamma(\tau_\gamma)\), also after differentiation to any fixed
order in these variables: each term contains a pulse coefficient or a derivative
of one.

\emph{The angular mean and the tangential double average.}
Both \(\mathsf L_B(w_0,p_{w_0})\) and \(E-\mathsf E_B\) have zero angular
mean. Taking the angular mean of \eqref{init:primary-force-expansion}
and using \(\nabla\cdot w_0=0\) gives the three cylindrical
divergence formulas \eqref{init:primary-angular-force}.
For every \(\mu,\nu\in\{r,\theta,z\}\), support separation and
\(N_\ell p_{\theta,\gamma}\ne0\) give
\begin{equation}
\langle w_0^\mu w_0^\nu\rangle_\theta
 =\frac12\sum_\gamma2^{2(1+\delta)\ell}
          (a_\gamma^\mu a_\gamma^\nu+b_\gamma^\mu b_\gamma^\nu).
\label{init:primary-angular-products}
\end{equation}
Indeed, the cylindrical components of \(\mathbf a_\gamma,\mathbf b_\gamma\)
are independent of \(\theta\), and
\[
\langle\cos^2(N_\ell\varphi_\gamma)\rangle_\theta
 =\langle\sin^2(N_\ell\varphi_\gamma)\rangle_\theta=\tfrac12,
\qquad
\langle\cos(N_\ell\varphi_\gamma)\sin(N_\ell\varphi_\gamma)\rangle_\theta=0.
\]
In particular, there is no mixed \(\mathbf a_\gamma\mathbf b_\gamma\)
term in \eqref{init:primary-angular-products}.
The matching identity
\eqref{eq:assembled-leading-stress} yields the exact radial flux
\begin{equation}
\langle w_0^rw_0^{\mathrm{tan}}\rangle_{\theta,Y_{i_0}}-\mathsf T_B
 =\frac12\sum_\gamma2^{2(1+\delta)\ell}
              \langle b_\gamma^rb_\gamma^{\mathrm{tan}}\rangle_{Y_{i_0}}.
\label{init:primary-matched-flux}
\end{equation}
By \eqref{init:primary-coefficient-bounds}, for \(\nu=\theta,z\),
\[
\begin{aligned}
\left|\frac12\sum_\gamma2^{2(1+\delta)\ell}
 \langle b_\gamma^rb_\gamma^\nu\rangle_{Y_{i_0}}\right|
 &\lesssim k^C2^{2(1+\delta)k}2^{-2(1-10^{-5})\delta k}
                  \omega d_{\mathrm{edge}}^{-C}\\
 &=k^C2^{(2+2\cdot10^{-5}\delta)k}\omega d_{\mathrm{edge}}^{-C},\\
\left|\partial_r\left[\frac12\sum_\gamma2^{2(1+\delta)\ell}
 \langle b_\gamma^rb_\gamma^\nu\rangle_{Y_{i_0}}\right]\right|
 &\lesssim k^C2^{(3+2\cdot10^{-5}\delta)k}\omega d_{\mathrm{edge}}^{-C}.
\end{aligned}
\]
For the derivative estimate, averaging removes the auxiliary part of
\(D_r\), so that only the slow radial derivative remains; it costs
\(2^k\) and a power of \(k\), with the power of
\(d_{\mathrm{edge}}^{-1}\) allowed to increase. The same bound holds for the
 terms with \(r^{-1}\) in \eqref{init:primary-angular-force}.
For the axial derivative,
\eqref{init:primary-angular-products} gives, for \(\nu=\theta,z\),
\[
\begin{aligned}
\partial_z\langle w_0^zw_0^\nu\rangle_\theta
 =\frac12\sum_\gamma2^{2(1+\delta)\ell}\Bigl[
 &(\partial_za_\gamma^z)a_\gamma^\nu
                 +a_\gamma^z\partial_za_\gamma^\nu\\
 &+(\partial_zb_\gamma^z)b_\gamma^\nu
                 +b_\gamma^z\partial_zb_\gamma^\nu\Bigr].
\end{aligned}
\]
On the coefficients at native scale \(\ell\),
\(\partial_z=2^{(1-\delta)\ell}\partial_{\hat Z}\).
Thus \eqref{init:primary-coefficient-bounds} and \(|\ell-k|\le1\)
imply
\[
\begin{aligned}
|\partial_z\mathbf a_\gamma|
 &\lesssim k^C2^{(1-\delta)k}2^{-\delta k/2}
                  \sqrt\omega\,d_{\mathrm{edge}}^{-C},\\
|\partial_z\mathbf b_\gamma|
 &\lesssim k^C2^{(1-\delta)k}2^{-(1-10^{-5})\delta k}
                  \sqrt\omega\,d_{\mathrm{edge}}^{-C}.
\end{aligned}
\]
Combining these bounds with finite overlap and then taking the
auxiliary average yields
\begin{equation}
\begin{aligned}
|\partial_z\langle w_0^zw_0^\nu\rangle_{\theta,Y_{i_0}}|
 &\le\left\langle
       |\partial_z\langle w_0^zw_0^\nu\rangle_\theta|
       \right\rangle_{Y_{i_0}}\\
 &\lesssim k^C2^{2(1+\delta)k}2^{(1-\delta)k}
       \Bigl(2^{-\delta k}+2^{-2(1-10^{-5})\delta k}\Bigr)
       \omega d_{\mathrm{edge}}^{-C}\\
 &\lesssim k^C2^{3k}\omega d_{\mathrm{edge}}^{-C}.
\end{aligned}
\label{init:primary-axial-flux}
\end{equation}
Taking the auxiliary average in \eqref{init:primary-angular-force}
and combining the radial covariance bounds above with
\eqref{init:primary-axial-flux} proves
\[
2^{-(3+\delta)k}|\langle f^{\mathrm{tan}}\rangle_{\theta,Y_{i_0}}|
 \lesssim k^C2^{-(1-2\cdot10^{-5})\delta k}\omega d_{\mathrm{edge}}^{-C}
 \le k^C2^{-0.99\delta k}\omega d_{\mathrm{edge}}^{-C}.
\]
This completes the proof of \eqref{init:primary-pointwise-bounds}.
\end{proof}

\subsection{The auxiliary velocity correction}
\label{init:auxiliary-correction}

We now reduce the part of the tangential angular mean that still
oscillates in the auxiliary variable. The input is the state of
Proposition~\ref{init:primary-pulse-estimates}:
\[
u=u_B+w_0,\qquad p=p_B+p_{w_0},\qquad
\mathsf{NS}(u,p)=E+f.
\]
The pressure remains unchanged in this subsection. We will construct
an axisymmetric, divergence-free increment \(u_{0,\mathrm{aux}}\)
with zero auxiliary mean, whose fast time derivative cancels
\[
\langle f^a\rangle_\theta-\langle f^a\rangle_{\theta,Y_{i_0}},
\qquad a=\theta,z,
\]
up to a flat term.

\emph{Bounds for the source.}
Fix a reference band \(\mathcal B_k\), and use
\[
i_0=i_0(k)=\min_{\ell\in\mathcal K_k}i(\ell)
\]
as in Section~\ref{sec:pulse-common}. In the following estimates, the
normalized coordinates use the reference scale \(k\), and
\(k,i_0(k)\) are held fixed when differentiating.
The exact angular products \eqref{init:primary-angular-products}
and coefficient bounds \eqref{init:primary-coefficient-bounds} give
\[
2^{-2(1+\delta)k}
 \left|\partial_{\hat R,\hat Z,\hat T,Y_{i_0,1},Y_{i_0,2}}^\beta\langle w_0^a w_0^b\rangle_\theta\right|
 \lesssim k^{C_\beta}\omega d_{\mathrm{edge}}^{-C_\beta}2^{-\delta k},
 \qquad a,b\in\{r,\theta,z\}.
\]
Since the background stress is independent of the auxiliary
variable, \eqref{init:primary-angular-force} gives
\[
\begin{aligned}
&\langle f^\theta\rangle_\theta
               -\langle f^\theta\rangle_{\theta,Y_{i_0}}\\
&\quad=\left(D_r+\frac2r\right)
 \Bigl(\langle w_0^rw_0^\theta\rangle_\theta
           -\langle w_0^rw_0^\theta\rangle_{\theta,Y_{i_0}}\Bigr)\\
&\qquad+\partial_z
 \Bigl(\langle w_0^zw_0^\theta\rangle_\theta
           -\langle w_0^zw_0^\theta\rangle_{\theta,Y_{i_0}}\Bigr),\\[3pt]
&\langle f^z\rangle_\theta-\langle f^z\rangle_{\theta,Y_{i_0}}\\
&\quad=\left(D_r+\frac1r\right)
 \Bigl(\langle w_0^rw_0^z\rangle_\theta
           -\langle w_0^rw_0^z\rangle_{\theta,Y_{i_0}}\Bigr)\\
&\qquad+\partial_z
 \Bigl(\langle(w_0^z)^2\rangle_\theta
           -\langle(w_0^z)^2\rangle_{\theta,Y_{i_0}}\Bigr).
\end{aligned}
\]
The radial derivative bound from \eqref{s6n:normalized-derivatives}
and \eqref{eq:ck-size}, together with
\(2^{-k}\partial_z=2^{-\delta k}\partial_{\hat Z}\), implies,
for \(a=\theta,z\),
\[
\begin{aligned}
&2^{-(3+\delta)k}\Bigl(
 |\partial_{\hat R,\hat Z,\hat T,Y_{i_0,1},Y_{i_0,2}}^\beta D_r\langle w_0^rw_0^a\rangle_\theta|
 +\bigl|\partial_{\hat R,\hat Z,\hat T,Y_{i_0,1},Y_{i_0,2}}^\beta\bigl[r^{-1}
                   \langle w_0^rw_0^a\rangle_\theta\bigr]\bigr|\Bigr)\\
&\qquad\lesssim k^{C_\beta}\omega d_{\mathrm{edge}}^{-C_\beta}
             2^{\delta k}2^{10^{-5}\delta k}2^{-\delta k},\\
&2^{-(3+\delta)k}
 |\partial_{\hat R,\hat Z,\hat T,Y_{i_0,1},Y_{i_0,2}}^\beta\partial_z\langle w_0^zw_0^a\rangle_\theta|
 \lesssim k^{C_\beta}\omega d_{\mathrm{edge}}^{-C_\beta}2^{-\delta k}.
\end{aligned}
\]
The corresponding double-averaged terms satisfy the same bounds.
Consequently, at every fixed derivative order,
\begin{equation}
\begin{aligned}
2^{-(3+\delta)k}\Bigl(
 &\left|\partial_{\hat R,\hat Z,\hat T,Y_{i_0,1},Y_{i_0,2}}^\beta\!\left(
 \langle f^\theta\rangle_\theta
       -\langle f^\theta\rangle_{\theta,Y_{i_0}}\right)\right|\\
+&\left|\partial_{\hat R,\hat Z,\hat T,Y_{i_0,1},Y_{i_0,2}}^\beta\!\left(
 \langle f^z\rangle_\theta
       -\langle f^z\rangle_{\theta,Y_{i_0}}\right)\right|\Bigr)
 \lesssim k^{C_\beta}\omega d_{\mathrm{edge}}^{-C_\beta}2^{10^{-5}\delta k}.
\end{aligned}
\label{init:auxiliary-source-bound}
\end{equation}
The bound at higher derivative orders is needed because the
auxiliary inverse loses a finite number of torus derivatives.

\emph{Inverting the fast time derivative.}
In the current common coordinate, the physical time derivative is
\[
D_t=\partial_t+\Lambda_t^{i_0}v_t\cdot\nabla_{Y_{i_0}}.
\]
For a zero-mean function
\(h(Y_{i_0})=\sum_{n\ne0}h_n e^{2\pi\mathrm i n\cdot Y_{i_0}}\),
Lemma~\ref{s6n:directional-inverse} gives the unique zero-mean inverse
\[
(v_t\cdot\nabla_{Y_{i_0}})^{-1}h
 =\sum_{n\ne0}\frac{h_n}{2\pi\mathrm i(v_t\cdot n)}
                         e^{2\pi\mathrm i n\cdot Y_{i_0}}.
\]
The Diophantine bound \eqref{eq:diophantine-directions} controls the
divisor by \(C(1+|n|)\). Fourier decay and summation over
\(\mathbb Z^2\) therefore give
\begin{equation}
\|(v_t\cdot\nabla_{Y_{i_0}})^{-1}h\|_{C^m_{Y_{i_0}}}
 \le C_m\|h\|_{C^{m+4}_{Y_{i_0}}}.
\label{init:auxiliary-inverse-bound}
\end{equation}
Define two real scalar fields by
\begin{equation}
\begin{aligned}
A_{0,\mathrm{aux}}
 &=-\Lambda_t^{-i_0}(v_t\cdot\nabla_{Y_{i_0}})^{-1}
   \bigl(\langle f^\theta\rangle_\theta
                  -\langle f^\theta\rangle_{\theta,Y_{i_0}}\bigr),\\
B_{0,\mathrm{aux}}
 &=-\Lambda_t^{-i_0}(v_t\cdot\nabla_{Y_{i_0}})^{-1}
   \bigl(\langle f^z\rangle_\theta
                  -\langle f^z\rangle_{\theta,Y_{i_0}}\bigr).
\end{aligned}
\label{init:auxiliary-scalars}
\end{equation}
Both have zero auxiliary mean.
Lemma~\ref{lem:common-compatibility} shows that these definitions
agree on overlapping reference bands, so they are independent of $k$. Equivalently,
after pullback to the original variable \(Y\), they are the unique
zero-mean solutions of
\[
\begin{aligned}
(v_t\cdot\nabla_Y)A_{0,\mathrm{aux}}
 &=-\langle f^\theta\rangle_\theta
                  +\langle f^\theta\rangle_{\theta,Y_{i_0}},\\
(v_t\cdot\nabla_Y)B_{0,\mathrm{aux}}
 &=-\langle f^z\rangle_\theta
                  +\langle f^z\rangle_{\theta,Y_{i_0}}.
\end{aligned}
\]
Here both right-hand sides are understood after the same pullback.

By Section~\ref{sec:pulse-common},
\[
0\le i(\ell)-i_0\le2,\qquad \ell\in\mathcal K_k.
\]
Combining this with \eqref{eq:ck-size} gives
\(\Lambda_t^{-i_0}\asymp k^2 2^{-(2+\delta)k}\). In particular,
\eqref{init:auxiliary-inverse-bound} yields
\[
\begin{aligned}
\|A_{0,\mathrm{aux}}\|_{C^m_{Y_{i_0}}}
 +\|B_{0,\mathrm{aux}}\|_{C^m_{Y_{i_0}}}
 &\lesssim \Lambda_t^{-i_0}
       \sum_{a=\theta,z}
       \left\|\langle f^a\rangle_\theta
          -\langle f^a\rangle_{\theta,Y_{i_0}}\right\|_{C^{m+4}_{Y_{i_0}}}\\
 &\lesssim k^{C_m}\omega d_{\mathrm{edge}}^{-C_m}
       \frac{2^{(3+\delta)k}2^{0.01\delta k}}{2^{(2+\delta)k}}\\
 &=k^{C_m}\omega d_{\mathrm{edge}}^{-C_m}2^{(1+\delta)k}2^{-0.99\delta k}.
\end{aligned}
\]
Keeping the factors in \eqref{init:auxiliary-source-bound} gives,
at every fixed order of differentiation in these variables,
\begin{equation}
\begin{aligned}
|\partial_{\hat R,\hat Z,\hat T,Y_{i_0,1},Y_{i_0,2}}^\beta A_{0,\mathrm{aux}}|
                 +|\partial_{\hat R,\hat Z,\hat T,Y_{i_0,1},Y_{i_0,2}}^\beta B_{0,\mathrm{aux}}|
 &\lesssim k^{C_\beta}\omega d_{\mathrm{edge}}^{-C_\beta}\Lambda_t^{-i_0}
                    2^{(3+\delta)k}2^{10^{-5}\delta k}\\
 &\lesssim k^{C_\beta}\omega d_{\mathrm{edge}}^{-C_\beta}2^k2^{10^{-5}\delta k}\\
 &\le k^{C_\beta}\omega d_{\mathrm{edge}}^{-C_\beta}2^{(1+\delta)k}2^{-0.99\delta k}.
\end{aligned}
\label{init:auxiliary-scalar-bounds}
\end{equation}
Indeed the inverse commutes with derivatives in the normalized slow
variables, and its finite loss of auxiliary derivatives is supplied
by \eqref{init:auxiliary-source-bound}. The radial weight is independent
of \(Y_{i_0}\), so the inverse preserves it.

\emph{Completing the divergence-free velocity.}
The two scalar solutions give the tangential field
\[
A_{0,\mathrm{aux}}e_\theta+B_{0,\mathrm{aux}}e_z.
\]
Its divergence is \(\partial_zB_{0,\mathrm{aux}}\), which need not vanish.
To obtain a divergence-free field \(u_{0,\mathrm{aux}}\) supported in
\(R\in[R_a,R_h]\), we need to make a flat modification to the axial
component and add a radial component.
Choose a smooth function \(\chi_{m,*}(R)\) which is zero near
\(R_a\), one near \(R_h\), and constant outside \([R_a,R_h]\), and set
\[
\chi_m(t,r,z)=\chi_{m,*}\!\left(\frac{r^2}{2\lambda(t,z)^2}\right).
\]
For a scalar field \(h\) supported in the annulus, define a modified radial primitive operator
\begin{equation}
\begin{aligned}
\mathcal I_0h(t,r,z,Y_{i_0})
={}&\int_0^r h\bigl(t,s,z,
  Y_{i_0}+\Lambda_r^{i_0}(s^{d_0}-r^{d_0})v_r\bigr)\,ds\\
 &-\chi_m(t,r,z)\int_0^\infty h\bigl(t,s,z,
  Y_{i_0}+\Lambda_r^{i_0}(s^{d_0}-r^{d_0})v_r\bigr)\,ds.
\end{aligned}
\label{init:auxiliary-radial-primitive}
\end{equation}
The translation follows the physical radial map: after setting
\(Y_{i_0}=Y_{i_0}(t,r)\), the argument of \(h\) becomes
\(Y_{i_0}(t,s)\). Thus \(D_r\) differentiates this primitive as an
ordinary radial integral, except for the derivative of \(\chi_m\).
Below the annulus both terms vanish; above it the two integrals
agree and \(\chi_m=1\). This is why the subtraction is included.
Since \(E^jv_r=\Lambda_r^jv_r\), the operator \(\mathcal I_0\) commutes with the common-coordinate pullbacks.
For compatible inputs on overlapping reference bands, the resulting primitives therefore agree after pullback to the original torus variable \(Y\).

Set
\begin{equation}
\Psi_{0,\mathrm{aux}}=\frac1r\mathcal I_0(rB_{0,\mathrm{aux}}),
\qquad
u_{0,\mathrm{aux}}
 =A_{0,\mathrm{aux}}e_\theta
                   +\nabla\times(\Psi_{0,\mathrm{aux}}e_\theta).
\label{init:auxiliary-velocity}
\end{equation}
The curl uses physical derivatives. Its components are
\begin{equation}
\begin{aligned}
u_{0,\mathrm{aux}}^r&=-\partial_z\Psi_{0,\mathrm{aux}},&
u_{0,\mathrm{aux}}^\theta&=A_{0,\mathrm{aux}},\\
u_{0,\mathrm{aux}}^z
 &=\left(D_r+\frac1r\right)\Psi_{0,\mathrm{aux}}
   =B_{0,\mathrm{aux}}+\mathcal R_{0,\mathrm{aux}},\\
\mathcal R_{0,\mathrm{aux}}
 &=-\frac{\partial_r\chi_m}{r}\int_0^\infty sB_{0,\mathrm{aux}}
  \bigl(t,s,z,Y_{i_0}+\Lambda_r^{i_0}(s^{d_0}-r^{d_0})v_r\bigr)\,ds.
\end{aligned}
\label{init:auxiliary-components}
\end{equation}

\begin{lem}[Size and cancellation of the auxiliary increment]
\label{init:auxiliary-increment}
The field \(u_{0,\mathrm{aux}}\) is real and smooth, is supported in
\(R\in[R_a,R_h]\), and satisfies
\begin{equation}
\nabla\cdot u_{0,\mathrm{aux}}=0,\qquad
\partial_\theta u_{0,\mathrm{aux}}^{r,\theta,z}=0,\qquad
\langle u_{0,\mathrm{aux}}\rangle_{Y_{i_0}}=0.
\label{init:auxiliary-invariants}
\end{equation}
For every fixed \(\beta\in\mathbb Z_{\ge0}^5\), uniformly on
\(\mathcal B_k\),
\begin{equation}
\begin{aligned}
2^{-\delta k}|\partial_{\hat R,\hat Z,\hat T,Y_{i_0,1},Y_{i_0,2}}^\beta\Psi_{0,\mathrm{aux}}|
 +2^{-(1+\delta)k}
       |\partial_{\hat R,\hat Z,\hat T,Y_{i_0,1},Y_{i_0,2}}^\beta u_{0,\mathrm{aux}}^{\mathrm{tan}}|
 &\lesssim k^{C_\beta}\omega d_{\mathrm{edge}}^{-C_\beta}2^{-0.99\delta k},\\
2^{-(1+\delta)k}|\partial_{\hat R,\hat Z,\hat T,Y_{i_0,1},Y_{i_0,2}}^\beta u_{0,\mathrm{aux}}^r|
 &\lesssim k^{C_\beta}\omega d_{\mathrm{edge}}^{-C_\beta}2^{-1.99\delta k}.
\end{aligned}
\label{init:auxiliary-velocity-bounds}
\end{equation}
Both \(\mathcal R_{0,\mathrm{aux}}\) and
\(\Lambda_t^{i_0}(v_t\cdot\nabla_{Y_{i_0}})\mathcal R_{0,\mathrm{aux}}\)
have zero auxiliary mean and are flat with every fixed physical
derivative. The fast time cancellation is exactly
\begin{equation}
\begin{aligned}
\Lambda_t^{i_0}(v_t\cdot\nabla_{Y_{i_0}})u_{0,\mathrm{aux}}^\theta
 &=-\langle f^\theta\rangle_\theta
                     +\langle f^\theta\rangle_{\theta,Y_{i_0}},\\
\Lambda_t^{i_0}(v_t\cdot\nabla_{Y_{i_0}})u_{0,\mathrm{aux}}^z
 &=-\langle f^z\rangle_\theta
                     +\langle f^z\rangle_{\theta,Y_{i_0}}\\
 &\quad+\Lambda_t^{i_0}(v_t\cdot\nabla_{Y_{i_0}})\mathcal R_{0,\mathrm{aux}}.
\end{aligned}
\label{init:auxiliary-fast-cancellation}
\end{equation}
\end{lem}

\begin{proof}
We first estimate the primitive without differentiating a large
radial translation. In this proof abbreviate
\(A=A_{0,\mathrm{aux}}, B=B_{0,\mathrm{aux}},\Psi=\Psi_{0,\mathrm{aux}}\),
and put
\[
U=(2^kr)^{d_0},\qquad M=2^{-kd_0}\Lambda_r^{i_0},\qquad
b(U,Y_{i_0})=2^{-(1+\delta)k}
                  B(t,2^{-k}U^{1/d_0},z,Y_{i_0}).
\]
The variables \(\hat T,\hat Z\) are suppressed. On the active
annulus, \(U\) stays in a fixed compact subset of \((0,\infty)\).
Let
\[
F(U,Y_{i_0})=\frac1{d_0}U^{2/d_0-1}b(U,Y_{i_0}),
\]
extended smoothly by zero outside its radial support. Changing the
radial integration variable to \(U+\sigma=(2^ks)^{d_0}\) gives
\begin{equation}
\begin{aligned}
2^{-\delta k}\Psi
=U^{-1/d_0}\Bigl[&(1-\chi_m)\int_{-\infty}^0
 F(U+\sigma,Y_{i_0}+M\sigma v_r)\,d\sigma\\
 &-\chi_m\int_0^\infty
 F(U+\sigma,Y_{i_0}+M\sigma v_r)\,d\sigma\Bigr].
\end{aligned}
\label{init:auxiliary-fixed-shift}
\end{equation}
At fixed \(\sigma\), the auxiliary shift is independent of the
normalized slow variables. The nonzero integration length and the
fixed derivatives of \(U^{-1/d_0}\) and \(\chi_m\) are bounded.
Near the inner edge only the first integral in
\eqref{init:auxiliary-fixed-shift} occurs, and near the outer edge
only the second occurs. In the logarithmic distance \(x\) from the
corresponding edge, the estimate
\[
\int_0^x e^{-c/s^2}s^{-q}\,ds
 \lesssim e^{-c/x^2}x^{-q'},\qquad c=c_{\mathrm{collar}, a}\text{ or }c_{\mathrm{collar}, h},
\]
for a suitable \(q'\), preserves the same exponential factor.
The factor from the opposite edge is bounded above and below there,
and \(\omega\) is bounded below on fixed interior subannuli.
Thus \eqref{init:auxiliary-scalar-bounds} gives
\[
|\partial_{\hat R,\hat Z,\hat T,Y_{i_0,1},Y_{i_0,2}}^\beta\Psi|
 \lesssim k^{C_\beta}\omega d_{\mathrm{edge}}^{-C_\beta}
                     2^{10^{-5}\delta k}.
\]
This also proves the smooth zero extension of \(\Psi\).
Since \(\partial_z=2^{(1-\delta)k}\partial_{\hat Z}\),
\[
\begin{aligned}
|\partial_{\hat R,\hat Z,\hat T,Y_{i_0,1},Y_{i_0,2}}^\beta u_{0,\mathrm{aux}}^r|
 &=2^{(1-\delta)k}|\partial_{\hat R,\hat Z,\hat T,Y_{i_0,1},Y_{i_0,2}}^\beta\partial_{\hat Z}\Psi|\\
 &\lesssim k^{C_\beta}\omega d_{\mathrm{edge}}^{-C_\beta}
                     2^{(1-\delta)k}2^{10^{-5}\delta k}.
\end{aligned}
\]
This shows the extra factor \(2^{-\delta k}\) in the radial component.

We next estimate the full integral in
\(\mathcal R_{0,\mathrm{aux}}\). Its zero Fourier mode vanishes
because \(\langle B\rangle_{Y_{i_0}}=0\). Write
\(F=\sum_{n\ne0}F_n(U)e^{2\pi\mathrm i n\cdot Y_{i_0}}\).
For every integer \(J\ge1\), integration by parts gives
\[
\begin{aligned}
&\int_{\mathbb R}F_n(U+\sigma)
          e^{2\pi\mathrm i M(n\cdot v_r)\sigma}\,d\sigma\\
&\qquad=
 \left(\frac{-1}{2\pi\mathrm i M(n\cdot v_r)}\right)^J
 \int_{\mathbb R}F_n^{(J)}(U+\sigma)
          e^{2\pi\mathrm i M(n\cdot v_r)\sigma}\,d\sigma.
\end{aligned}
\]
There are no endpoint terms because \(F_n\) and all its radial
derivatives vanish at both ends of the support. By
\eqref{eq:diophantine-directions}, each inverse divisor costs at
most \(C(1+|n|)\).
Moreover, \eqref{eq:ck-size} gives
\[
M \asymp k^{-2\rho_0}2^{10^{-5}\delta k},
\qquad M^{-J}\lesssim k^{2J\rho_0}2^{-J\,10^{-5}\delta k}.
\]
The factor \((\partial_r\chi_m)/r\) in
\eqref{init:auxiliary-components} has size \(O(2^{2k})\), with the
same scaling after fixed normalized derivatives. Applying
\eqref{init:auxiliary-scalar-bounds} to the derivatives of \(F\)
and restoring the units of the radial integral yields
\begin{equation}
2^{-(1+\delta)k}|\partial_{\hat R,\hat Z,\hat T,Y_{i_0,1},Y_{i_0,2}}^\beta\mathcal R_{0,\mathrm{aux}}|
 \lesssim C_{\beta, J} k^{C_{\beta,J}}
                  2^{-(1-10^{-5}+J\,10^{-5})\delta k}.
\label{init:auxiliary-tail-bound}
\end{equation}
Every fixed physical derivative, and the extra fast time derivative
in \eqref{init:auxiliary-fast-cancellation}, costs at most
\(k^C2^{Ck}\). Choosing \(J\) after the derivative order and the
desired decay power proves, for every fixed \(\beta\) and \(N>0\),
\[
|D_{t,x}^\beta\mathcal R_{0,\mathrm{aux}}|
 +\left|D_{t,x}^\beta\!\left[
   \Lambda_t^{i_0}(v_t\cdot\nabla_{Y_{i_0}})
                         \mathcal R_{0,\mathrm{aux}}\right]\right|
 \le C_{\beta,N}\lambda^N.
\]
\(\mathcal R_{0,\mathrm{aux}}\) and its derivatives in these variables are
supported in \(\operatorname{supp}\chi_{m,*}'\Subset(R_a,R_h)\),
where \(\omega\) is bounded below. Hence the tail bound also holds
with a factor \(\omega\), after changing the constant.
Since \(u_{0,\mathrm{aux}}^\theta=A\) and
\(u_{0,\mathrm{aux}}^z=B+\mathcal R_{0,\mathrm{aux}}\), we have proved
\begin{equation}
\begin{aligned}
|\partial_{\hat R,\hat Z,\hat T,Y_{i_0,1},Y_{i_0,2}}^\beta u_{0,\mathrm{aux}}^{\mathrm{tan}}|
 &\lesssim k^{C_\beta}\omega d_{\mathrm{edge}}^{-C_\beta}2^k2^{10^{-5}\delta k},\\
|\partial_{\hat R,\hat Z,\hat T,Y_{i_0,1},Y_{i_0,2}}^\beta u_{0,\mathrm{aux}}^r|
 &\lesssim k^{C_\beta}\omega d_{\mathrm{edge}}^{-C_\beta}2^{(1-\delta)k}2^{10^{-5}\delta k}.
\end{aligned}
\label{init:auxiliary-physical-bounds}
\end{equation}
Dividing these estimates by \(2^{(1+\delta)k}\), and the potential
estimate by \(2^{\delta k}\), proves the bounds with exponents
\(0.99\) and \(1.99\) in \eqref{init:auxiliary-velocity-bounds}.

Translations preserve the auxiliary average. Hence
\(\langle\Psi\rangle_{Y_{i_0}}
 =\langle\mathcal R_{0,\mathrm{aux}}\rangle_{Y_{i_0}}=0\), and the
commutation identities at the start of Section~\ref{sec:initialization-ppt}
give \(\langle u_{0,\mathrm{aux}}\rangle_{Y_{i_0}}=0\).
The swirl component is divergence-free, and so is the curl.
More explicitly, using \([D_r,\partial_z]=0\),
\[
\left(D_r+\frac1r\right)u_{0,\mathrm{aux}}^r
                 +\partial_zu_{0,\mathrm{aux}}^z
 =-\left(D_r+\frac1r\right)\partial_z\Psi
                 +\partial_z\left(D_r+\frac1r\right)\Psi=0.
\]
This uses the actual potential, including derivatives of its
\((t,z)\)-dependent cutoff. Finally,
\eqref{init:auxiliary-scalars} and
\eqref{init:auxiliary-components} give
\eqref{init:auxiliary-fast-cancellation}.
The torus inverse preserves radial support, and the primitive
\(\mathcal I_0\) preserves support in the active annulus. The auxiliary
support may spread over the torus.
\end{proof}

\emph{The updated force.}
Define the actual new fields and flat error by
\begin{equation}
\begin{aligned}
u_{\mathrm{new}}&=u_B+w_0+u_{0,\mathrm{aux}},&
p_{\mathrm{new}}&=p_B+p_{w_0},\\
E_{\mathrm{new}}&=E+
 \Lambda_t^{i_0}(v_t\cdot\nabla_{Y_{i_0}})
                              \mathcal R_{0,\mathrm{aux}}\,e_z,&
f_{\mathrm{new}}&=\mathsf{NS}(u_{\mathrm{new}},p_{\mathrm{new}})
                                                      -E_{\mathrm{new}}.
\end{aligned}
\label{init:auxiliary-new-fields}
\end{equation}
The added flat term is axisymmetric and has zero auxiliary mean.
It therefore preserves
\(\langle E_{\mathrm{new}}\rangle_{\theta,Y_{i_0}}=\mathsf E_B\).

\begin{prop}[The force after the auxiliary correction]
\label{init:auxiliary-force}
The fields in \eqref{init:auxiliary-new-fields} satisfy
\mbox{\(\nabla\cdot u_{\mathrm{new}}=0\)}. Their velocity and pressure
increments over the background, and \(f_{\mathrm{new}}\) and
\(E_{\mathrm{new}}-\mathsf E_B\), are smooth and supported in
\(R\in[R_a,R_h]\). The flat
error satisfies
\[
|D_{t,x}^\beta E_{\mathrm{new}}|\le C_{\beta,N}\lambda^N,
\qquad \langle E_{\mathrm{new}}\rangle_{\theta,Y_{i_0}}=\mathsf E_B
\]
for every fixed \(\beta\) and \(N>0\). Uniformly on \(\mathcal B_k\),
\begin{equation}
\begin{aligned}
2^{-(1+\delta)k}|u_{0,\mathrm{aux}}|
 &\lesssim k^C\omega d_{\mathrm{edge}}^{-C}2^{-0.99\delta k},\\
2^{-(3+\delta)k}\bigl(|f_{\mathrm{new}}|
               +|f_{\mathrm{new}}-\langle f_{\mathrm{new}}\rangle_\theta|\bigr)
 &\lesssim k^C\sqrt\omega\,d_{\mathrm{edge}}^{-C}2^{0.01\delta k},\\
2^{-(3+\delta)k}
 \left|\langle f_{\mathrm{new}}^{\mathrm{tan}}\rangle_\theta
          -\langle f_{\mathrm{new}}^{\mathrm{tan}}\rangle_{\theta,Y_{i_0}}\right|
 &\lesssim k^C\omega d_{\mathrm{edge}}^{-C}2^{-0.99\delta k},\\
2^{-(3+\delta)k}
 \left|\langle f_{\mathrm{new}}^{\mathrm{tan}}\rangle_{\theta,Y_{i_0}}\right|
 &\lesssim k^C\omega d_{\mathrm{edge}}^{-C}2^{-0.99\delta k}.
\end{aligned}
\label{init:auxiliary-output-bounds}
\end{equation}
The two linear velocity moments in \eqref{init:linear-moments}
remain equal to their background values.
\end{prop}

\begin{proof}
For this proof write \(v=u_{0,\mathrm{aux}}\), and denote the
axisymmetric convective flux increment by
\[
\mathsf A^{ab}=u_B^av^b+u_B^bv^a+v^av^b,
\qquad a,b\in\{r,\theta,z\}.
\]
The exact update identity \eqref{eq:physical-residual-update} gives
\begin{equation}
\begin{aligned}
f_{\mathrm{new}}=f
 &+\mathsf L_B(v,0)+(v\cdot\nabla)v
  -(E_{\mathrm{new}}-E)\\
 &+(w_0\cdot\nabla)v+(v\cdot\nabla)w_0.
\end{aligned}
\label{init:auxiliary-force-update}
\end{equation}
The two mixed terms have zero angular mean: \(v\) is axisymmetric,
and \(w_0\) and each of its cylindrical component derivatives have
zero angular mean. All other added terms are axisymmetric.

\emph{Cancellation in the tangential angular mean.}
Put \(\Delta_0=D_r^2+r^{-1}D_r+\partial_z^2\).
The cylindrical momentum equations and
\eqref{init:auxiliary-fast-cancellation} give the exact identities
\begin{equation}
\begin{aligned}
\langle f_{\mathrm{new}}^\theta\rangle_\theta
={}&\langle f^\theta\rangle_{\theta,Y_{i_0}}
 +\partial_tv^\theta
 -\left(\Delta_0-\frac1{r^2}\right)v^\theta
 +\left(D_r+\frac2r\right)\mathsf A^{r\theta}
 +\partial_z\mathsf A^{z\theta},\\
\langle f_{\mathrm{new}}^z\rangle_\theta
={}&\langle f^z\rangle_{\theta,Y_{i_0}}
 +\partial_tv^z-\Delta_0v^z
 +\left(D_r+\frac1r\right)\mathsf A^{rz}
 +\partial_z\mathsf A^{zz}.
\end{aligned}
\label{init:auxiliary-tangential-force}
\end{equation}
Here \(\partial_t=D_t-\Lambda_t^{i_0}v_t\cdot\nabla_{Y_{i_0}}\)
holds the auxiliary variable fixed. In the second line the extra
fast derivative of \(\mathcal R_{0,\mathrm{aux}}\) is exactly the
term assigned to \(E_{\mathrm{new}}-E\).

We estimate the displayed terms in physical force units.
The background bounds \eqref{s6n:background-velocity-bounds} and
\eqref{init:auxiliary-physical-bounds} imply, for every fixed
\(\beta\in\mathbb Z_{\ge0}^5\),
\[
\begin{aligned}
2^{-2(1+\delta)k}|\partial_{\hat R,\hat Z,\hat T,Y_{i_0,1},Y_{i_0,2}}^\beta\mathsf A^{ra}|
 &\lesssim k^{C_\beta}\omega d_{\mathrm{edge}}^{-C_\beta}2^{-(2-10^{-5})\delta k},\\
2^{-2(1+\delta)k}|\partial_{\hat R,\hat Z,\hat T,Y_{i_0,1},Y_{i_0,2}}^\beta\mathsf A^{za}|
 &\lesssim k^{C_\beta}\omega d_{\mathrm{edge}}^{-C_\beta}2^{-(1-10^{-5})\delta k},
\qquad a=\theta,z.
\end{aligned}
\]
For the following estimates, put \(c_\theta=2\) and \(c_z=1\).
For example, \(u_B^rv^a\) has the factors
\(2^{-\delta k}2^{-(1-10^{-5})\delta k}\), and
\(u_B^av^r\) has the same bound. Their radial derivatives therefore
satisfy
\[
\begin{aligned}
2^{-(3+\delta)k}
 \left|\left(D_r+\frac{c_a}{r}\right)\mathsf A^{ra}\right|
 &\lesssim k^C\omega d_{\mathrm{edge}}^{-C}2^{\delta k}2^{10^{-5}\delta k}
                       2^{-(2-10^{-5})\delta k}\\
 &=k^C\omega d_{\mathrm{edge}}^{-C}2^{-(1-2\cdot10^{-5})\delta k}.
\end{aligned}
\]
The axial derivative has the better factor
\(2^{-k}\partial_z=2^{-\delta k}\partial_{\hat Z}\).
Also, \(\partial_t=-2^{2k}\partial_{\hat T}\), whereas two
physical radial derivatives cost at most
\(k^C2^{2k}2^{2\cdot10^{-5}\delta k}\), with a possible increase in
the power of \(d_{\mathrm{edge}}^{-1}\). Thus, for \(a=\theta,z\),
\[
\begin{array}{c|c}
\text{term divided by }2^{(3+\delta)k}
 &\text{upper bound, apart from }k^C\omega d_{\mathrm{edge}}^{-C}\\ \hline
|\partial_tv^a|&2^{-(1-10^{-5})\delta k}\\
|(D_r+c_a/r)\mathsf A^{ra}|&2^{-(1-2\cdot10^{-5})\delta k}\\
|\partial_z\mathsf A^{za}|&2^{-(1-10^{-5})\delta k}\\
|\Delta_0v^a|+r^{-2}|v^a|&2^{-(1-3\cdot10^{-5})\delta k}
\end{array}
\]
Subtracting the auxiliary average of
\eqref{init:auxiliary-tangential-force} removes its first term.
The remaining terms and their averages have the displayed bounds.
Since \(1-3\cdot10^{-5}>0.99\), this proves the third line of
\eqref{init:auxiliary-output-bounds}.

\emph{Preserving the tangential double mean.}
The zero auxiliary mean of \(v\) gives
\[
\langle\mathsf A^{ab}\rangle_{Y_{i_0}}
 =\langle v^av^b\rangle_{Y_{i_0}}.
\]
Averaging commutes with physical derivatives as described at the
start of Section~\ref{sec:initialization-ppt}. In particular the
time and viscous terms in
\eqref{init:auxiliary-tangential-force} have zero auxiliary average.
Consequently,
\begin{equation}
\langle f_{\mathrm{new}}^a-f^a\rangle_{\theta,Y_{i_0}}
 =\left(\partial_r+\frac{c_a}{r}\right)
                        \langle v^rv^a\rangle_{Y_{i_0}}
      +\partial_z\langle v^zv^a\rangle_{Y_{i_0}},
\quad a=\theta,z.
\label{init:auxiliary-double-mean-change}
\end{equation}
After averaging, the radial derivative has no auxiliary loss.
The radial and axial terms respectively satisfy
\[
\begin{aligned}
2^{-(3+\delta)k}
 \left|\left(\partial_r+\frac{c_a}{r}\right)
                      \langle v^rv^a\rangle_{Y_{i_0}}\right|
 &\lesssim k^C\omega d_{\mathrm{edge}}^{-C}2^{\delta k}
           2^{-(2-10^{-5})\delta k}2^{-(1-10^{-5})\delta k},\\
2^{-(3+\delta)k}
 \left|\partial_z\langle v^zv^a\rangle_{Y_{i_0}}\right|
 &\lesssim k^C\omega d_{\mathrm{edge}}^{-C}2^{-2(1-10^{-5})\delta k}.
\end{aligned}
\]
Both are bounded by \(k^C\omega d_{\mathrm{edge}}^{-C}2^{-(2-2\cdot10^{-5})\delta k}\), so the
incoming tangential double-mean bound in
\eqref{init:primary-pointwise-bounds} is preserved.

\emph{Preserving the nonzero angular modes.}
From \eqref{init:auxiliary-force-update}, exactly,
\begin{equation}
 f_{\mathrm{new}}-\langle f_{\mathrm{new}}\rangle_\theta
  =f-\langle f\rangle_\theta
          +(w_0\cdot\nabla)v+(v\cdot\nabla)w_0.
\label{init:auxiliary-nonzero-modes}
\end{equation}
The primary pulse and its physical spatial gradient satisfy
\[
\begin{aligned}
|w_0|&\lesssim k^C\sqrt\omega\,d_{\mathrm{edge}}^{-C}
                         2^{(1+\delta)k}2^{-\delta k/2},\\
|\nabla w_0|&\lesssim k^C\sqrt\omega\,d_{\mathrm{edge}}^{-C}2^{(2+\delta)k}.
\end{aligned}
\]
For the second bound, a carrier derivative contributes
\(2^\ell N_\ell\mathbf n_\gamma\), while the primary coefficient
has size \(k^C\sqrt\omega\,d_{\mathrm{edge}}^{-C}2^{-\delta k/2}\); the factors
\(N_\ell\lesssim2^{\delta k/2}\) and \(2^{-\delta k/2}\) cancel.
Coefficient derivatives and cylindrical frame derivatives are
smaller. Similarly,
\[
|\nabla v|\lesssim k^C\omega d_{\mathrm{edge}}^{-C}2^{(2+\delta)k}
                                    2^{-(1-2\cdot10^{-5})\delta k}.
\]
These gradients include the frame terms \(r^{-1}w_0\) and
\(r^{-1}v\). Since \(\omega^{3/2}\le\sqrt\omega\),
\[
\begin{aligned}
2^{-(3+\delta)k}|(v\cdot\nabla)w_0|
 &\lesssim k^C\sqrt\omega\,d_{\mathrm{edge}}^{-C}2^{\delta k}2^{-(1-10^{-5})\delta k}\\
 &=k^C\sqrt\omega\,d_{\mathrm{edge}}^{-C}2^{10^{-5}\delta k},\\
2^{-(3+\delta)k}|(w_0\cdot\nabla)v|
 &\lesssim k^C\sqrt\omega\,d_{\mathrm{edge}}^{-C}2^{\delta k}2^{-\delta k/2}
                         2^{-(1-2\cdot10^{-5})\delta k}\\
 &=k^C\sqrt\omega\,d_{\mathrm{edge}}^{-C}2^{-(1/2-2\cdot10^{-5})\delta k}.
\end{aligned}
\]
Together with \eqref{init:primary-pointwise-bounds} and
\(10^{-5}<0.01\), this proves the nonzero-mode bound.
The torus inverse need not preserve the auxiliary support of an individual pulse.
Each mixed term contains the corresponding primary pulse or a
derivative of it, so it is supported in that label's primary-pulse
support.

\emph{The radial angular mean and the remaining assertions.}
The radial component of \eqref{init:auxiliary-force-update} has
angular mean
\begin{equation}
\begin{aligned}
\langle f_{\mathrm{new}}^r-f^r\rangle_\theta
={}&D_tv^r+\left(D_r+\frac1r\right)\mathsf A^{rr}
 +\partial_z\mathsf A^{zr}-\frac{\mathsf A^{\theta\theta}}r
 -\left(\Delta_0-\frac1{r^2}\right)v^r.
\end{aligned}
\label{init:auxiliary-radial-change}
\end{equation}
Here the largest possible new term is the centrifugal contribution:
\[
\begin{aligned}
2^{-(3+\delta)k}\frac{|\mathsf A^{\theta\theta}|}{r}
 &\lesssim k^C\omega d_{\mathrm{edge}}^{-C}2^{\delta k}
       \left(2^{-(1-10^{-5})\delta k}
                         +2^{-2(1-10^{-5})\delta k}\right)\\
 &\lesssim k^C\omega d_{\mathrm{edge}}^{-C}2^{10^{-5}\delta k},\\
2^{-(3+\delta)k}|D_tv^r|
 &\lesssim k^C\omega d_{\mathrm{edge}}^{-C}2^{\delta k}2^{-(2-10^{-5})\delta k}\\
 &=k^C\omega d_{\mathrm{edge}}^{-C}2^{-(1-10^{-5})\delta k}.
\end{aligned}
\]
The radial and axial flux derivatives in
\eqref{init:auxiliary-radial-change} are bounded by
\(k^C\omega d_{\mathrm{edge}}^{-C}2^{-(2-2\cdot10^{-5})\delta k}\) in force units, and its
viscous term by \(k^C\omega d_{\mathrm{edge}}^{-C}2^{-(2-3\cdot10^{-5})\delta k}\).
The incoming radial angular mean has the same radial weight
by \eqref{init:primary-angular-force}. Combining these bounds with the
already proved tangential and nonzero-mode bounds proves the full
force estimate in \eqref{init:auxiliary-output-bounds}.

The velocity estimate, flatness and support follow from
Lemma~\ref{init:auxiliary-increment} and
\eqref{init:auxiliary-new-fields}. Finally, pointwise in \((t,r,z)\),
\[
\langle u_{\mathrm{new}}-u_B\rangle_{\theta,Y_{i_0}}
 =\langle w_0\rangle_{\theta,Y_{i_0}}
                 +\langle v\rangle_{Y_{i_0}}=0.
\]
Multiplication by \(r^2\) in the angular component and by \(r\)
in the axial component, followed by radial integration, proves
preservation of both linear velocity moments.
\end{proof}

Thus the auxiliary correction has reduced the remaining tangential
auxiliary oscillation to \(k^C\omega d_{\mathrm{edge}}^{-C}2^{-0.99\delta k}\) in force units,
which is stronger than the \(k^C\omega d_{\mathrm{edge}}^{-C}2^{-0.2\delta k}\) target in
\eqref{init:stage-zero-target}. The nonzero angular modes and the
tangential double mean retain their previous bounds. The radial
angular mean in \eqref{init:auxiliary-radial-change} still requires
an axisymmetric pressure correction. In the next two subsections, we first reduce the three
moment defects and then reconstruct the pressure for the updated
velocity.

\subsection{The moment correction}
\label{init:moment-correction}

We start from the state at the end of
Section~\ref{init:auxiliary-correction}, and write
\[
u=u_B+w_0+u_{0,\mathrm{aux}},\qquad
p=p_B+p_{w_0},\qquad
\mathsf{NS}(u,p)=E+f.
\]
For any axisymmetric pressure increment \(q\) supported in the annulus,
\[
\int_0^\infty\langle D_rq\rangle_{Y_{i_0}}\,dr
 =\int_0^\infty\partial_r\langle q\rangle_{Y_{i_0}}\,dr=0.
\]
Thus pressure reconstruction cannot reduce \(\mathcal P\).
The reconstruction in Section~\ref{init:pressure-reconstruction}
will bound the normalized radial angular mean by
\(C2^{-(2+\delta)k}|\mathcal P|\).
Reducing \(2^{-2(1+\delta)k}|\mathcal P|\) to \(k^C2^{-1.2\delta k}\) will therefore
give the required \(k^C2^{-0.2\delta k}\) radial bound.

Let \(I_3\Subset(R_a,R_h)\) be the interval in
\eqref{intro:background-moment-patch}, expressed in
\(R=r^2/(2\lambda^2)\).
The three defects \(\mathcal P,\mathcal J^{\theta z},
\mathcal J^{z\theta}\) need to be reduced to the size required
in \eqref{init:stage-zero-target}. We construct
\(u_{0,\mathrm{moment}}\), independent of \(\theta,Y_{i_0}\) and
supported in the reserved interval \(I_3\), by prescribing five
radial integrals. Two preserve the linear velocity moments; the
other three cancel the current defect values through the terms
linear in the new increment and the background. We then reconstruct
the pressure for the updated velocity.

\emph{The radial source and its integral.}
Fix a reference band \(\mathcal B_k\), and use the common
coordinates from Section~\ref{init:auxiliary-correction}.
For every fixed \(\beta\in\mathbb Z_{\ge0}^5\), we have
\begin{equation}
2^{-(3+\delta)k}
 \bigl|\partial_{\hat R,\hat Z,\hat T,Y_{i_0,1},Y_{i_0,2}}^\beta\langle f^r\rangle_\theta\bigr|
 \lesssim k^{C_\beta}2^{10^{-5}\delta k}
                  \omega d_{\mathrm{edge}}^{-C_\beta}.
\label{init:pressure-radial-source}
\end{equation}
To verify this estimate, first consider the three primary-pulse
terms in the radial identity \eqref{init:primary-angular-force}.
The product formula \eqref{init:primary-angular-products} and
the coefficient bounds \eqref{init:primary-coefficient-bounds} give
\[
\begin{aligned}
&2^{-(3+\delta)k}
 \left|\partial_{\hat R,\hat Z,\hat T,Y_{i_0,1},Y_{i_0,2}}^\beta\!\left[
 \left(D_r+\frac1r\right)\langle(w_0^r)^2\rangle_\theta\right]\right|
 \lesssim k^{C_\beta}2^{10^{-5}\delta k}
                  \omega d_{\mathrm{edge}}^{-C_\beta},\\
&2^{-(3+\delta)k}\left(
 \left|\partial_{\hat R,\hat Z,\hat T,Y_{i_0,1},Y_{i_0,2}}^\beta\!\left[\frac{\langle(w_0^\theta)^2\rangle_\theta}{r}
                    \right]\right|
 +\left|\partial_{\hat R,\hat Z,\hat T,Y_{i_0,1},Y_{i_0,2}}^\beta\partial_z
                    \langle w_0^rw_0^z\rangle_\theta\right|\right)
 \lesssim k^{C_\beta}
                  \omega d_{\mathrm{edge}}^{-C_\beta}.
\end{aligned}
\]
For the auxiliary contribution, put \(v=u_{0,\mathrm{aux}}\) in
\eqref{init:auxiliary-radial-change}. The background bounds and
\eqref{init:auxiliary-physical-bounds} yield
\[
\begin{aligned}
&2^{-(3+\delta)k}
 \left|\partial_{\hat R,\hat Z,\hat T,Y_{i_0,1},Y_{i_0,2}}^\beta\!\left[
             \frac{2u_B^\theta v^\theta+(v^\theta)^2}{r}\right]\right|\\
&\qquad\lesssim k^{C_\beta}2^{\delta k}
 \left(2^{-(1-10^{-5})\delta k}
            +2^{-2(1-10^{-5})\delta k}\right)
             \omega d_{\mathrm{edge}}^{-C_\beta}\\
&\qquad\lesssim k^{C_\beta}2^{10^{-5}\delta k}
             \omega d_{\mathrm{edge}}^{-C_\beta}.
\end{aligned}
\]
For the other terms in \eqref{init:auxiliary-radial-change},
the same derivative and product bounds give, with
\(\Delta_0=D_r^2+r^{-1}D_r+\partial_z^2\),
\[
\begin{aligned}
&2^{-(3+\delta)k}|\partial_{\hat R,\hat Z,\hat T,Y_{i_0,1},Y_{i_0,2}}^\beta D_tv^r|
 \lesssim k^{C_\beta}2^{-(1-10^{-5})\delta k}
                  \omega d_{\mathrm{edge}}^{-C_\beta},\\
&2^{-(3+\delta)k}\Bigl(
 \left|\partial_{\hat R,\hat Z,\hat T,Y_{i_0,1},Y_{i_0,2}}^\beta\!\left[
 \left(D_r+\frac1r\right)(2u_B^rv^r+(v^r)^2)\right]\right|\\
&\hspace{23mm}
 +|\partial_{\hat R,\hat Z,\hat T,Y_{i_0,1},Y_{i_0,2}}^\beta\partial_z(u_B^zv^r+u_B^rv^z+v^zv^r)|\Bigr)
 \lesssim k^{C_\beta}2^{-(2-2\cdot10^{-5})\delta k}
                  \omega d_{\mathrm{edge}}^{-C_\beta},\\
&2^{-(3+\delta)k}
 \left|\partial_{\hat R,\hat Z,\hat T,Y_{i_0,1},Y_{i_0,2}}^\beta\!\left[
             \left(\Delta_0-\frac1{r^2}\right)v^r\right]\right|
 \lesssim k^{C_\beta}2^{-(2-3\cdot10^{-5})\delta k}
                  \omega d_{\mathrm{edge}}^{-C_\beta}.
\end{aligned}
\]
Together these estimates prove \eqref{init:pressure-radial-source}.
Since \(\omega d_{\mathrm{edge}}^{-C_\beta}\) is bounded and the
radial support has length \(O(2^{-k})\), it follows that
\begin{equation}
\begin{aligned}
\mathcal P
 &=\int_0^\infty
                   \langle f^r\rangle_{\theta,Y_{i_0}}\,dr,\\
2^{-2(1+\delta)k}|\partial_{\hat T,\hat Z}^\beta\mathcal P|
 &\lesssim k^{C_\beta}2^{-2(1+\delta)k}2^{-k}
                      2^{(3+\delta)k}2^{10^{-5}\delta k}
  =k^{C_\beta}2^{-(1-10^{-5})\delta k}.
\end{aligned}
\label{init:pressure-defect-input}
\end{equation}
Here \(k\) is fixed when taking the derivatives in \((\hat T,\hat Z)\).

\emph{The incoming defects.}
The identity proved at the end of Section~\ref{init:auxiliary-correction}
is
\[
\langle u-u_B\rangle_{\theta,Y_{i_0}}=0.
\]
Using this identity, \(\langle w_0\rangle_\theta=0\), and
\(\langle p_{w_0}\rangle_\theta=0\), the definitions of the two
physical moment defects give
\begin{equation}
\begin{aligned}
\mathcal J^{\theta z}
 &=\int_0^\infty r^2
 \left\langle w_0^\theta w_0^z
       +u_{0,\mathrm{aux}}^\theta u_{0,\mathrm{aux}}^z
 \right\rangle_{\theta,Y_{i_0}}\,dr,\\
\mathcal J^{z\theta}
 &=\int_0^\infty r
 \left\langle (w_0^z)^2+(u_{0,\mathrm{aux}}^z)^2
 \right\rangle_{\theta,Y_{i_0}}\,dr
 +\frac12\int_0^\infty r^2
                  \langle f^r\rangle_{\theta,Y_{i_0}}\,dr.
\end{aligned}
\label{init:moment-input-formulas}
\end{equation}
The mixed products with the background vanish under double
averaging, and the mixed products with \(w_0\) vanish under angular
averaging. On \(\mathcal B_k\), the radial integration factors are
\(\int_{\mathrm{supp}}r^2\,dr\lesssim2^{-3k}\) and
\(\int_{\mathrm{supp}}r\,dr\lesssim2^{-2k}\).
Thus \eqref{init:primary-angular-products},
\eqref{init:auxiliary-physical-bounds} and
\eqref{init:pressure-radial-source} imply, after every fixed
derivative in \((\hat T,\hat Z)\),
\[
\begin{aligned}
&2^{(1-2\delta)k}\left|\partial_{\hat T,\hat Z}^\beta
 \mathcal J^{\theta z}\right|
 \lesssim k^{C_\beta}
       \left(2^{-\delta k}+2^{-2(1-10^{-5})\delta k}\right),\\
&2^{-2\delta k}\left|\partial_{\hat T,\hat Z}^\beta
 \mathcal J^{z\theta}\right|
 \lesssim k^{C_\beta}\left(2^{-\delta k}
       +2^{-2(1-10^{-5})\delta k}+2^{-(1-10^{-5})\delta k}\right).
\end{aligned}
\]
Here we used the normalization factors
\(2^{(1-2\delta)k}\) and \(2^{-2\delta k}\), respectively.
Together with \eqref{init:pressure-defect-input}, this proves
\begin{equation}
\begin{aligned}
&2^{-2(1+\delta)k}|\partial_{\hat T,\hat Z}^\beta \mathcal P|
 +2^{(1-2\delta)k}|\partial_{\hat T,\hat Z}^\beta \mathcal J^{\theta z}|
 +2^{-2\delta k}|\partial_{\hat T,\hat Z}^\beta \mathcal J^{z\theta}|\\
&\qquad\lesssim k^{C_\beta}2^{-(1-10^{-5})\delta k}.
\end{aligned}
\label{init:moment-input-bounds}
\end{equation}

\emph{The five moment equations.}
For this construction write \(x=r/\lambda\). The input condition
\eqref{intro:background-moment-patch} gives, on \(R\in I_3\),
\begin{equation}
u_B^z=0,\qquad
u_B^\theta=\lambda^{-1-\delta}A_*(Z)x^{-1-2\mu},\qquad
A_*(Z)\ge a_*/2>0,\quad\mu>0.
\label{init:moment-background-patch}
\end{equation}
We choose all profiles below so that their entire radial hull
(the smallest closed interval containing the supports of the profiles)
lies in \(\{x>0:x^2/2\in I_3\}\).
Write \(v=u_{0,\mathrm{moment}}\). Its tangential components will
satisfy
\begin{equation}
\begin{aligned}
\int_0^\infty r^2v^\theta\,dr&=0,&
\int_0^\infty rv^z\,dr&=0,\\
\int_0^\infty\frac{2u_B^\theta}{r}v^\theta\,dr
 &=\mathcal P,\\
\int_0^\infty r^2u_B^\theta v^z\,dr
 &=-\mathcal J^{\theta z},\\
\int_0^\infty ru_B^\theta v^\theta\,dr
 &=\mathcal J^{z\theta}.
\end{aligned}
\label{init:moment-five-equations}
\end{equation}
The third equation makes the contribution of
\(-2u_B^\theta v^\theta/r\) to the change in \(\mathcal P\)
equal to \(-\mathcal P\).
The fourth gives the contribution
\(\int_0^\infty r^2u_B^\theta v^z\,dr=-\mathcal J^{\theta z}\)
to the change in \(\mathcal J^{\theta z}\).
In \(\mathcal J^{z\theta}\), the same radial force term has weight
\(r^2/2\), so its contribution is
\(-\int_0^\infty ru_B^\theta v^\theta\,dr=-\mathcal J^{z\theta}\).
This explains the positive sign in the last equation.
The right-hand sides are the full incoming defects, including their
quadratic contributions.

Choose a nonnegative, nonzero smooth bump \(\eta_0\) and a number
\(d_*>0\) such that
\[
\eta_h(x)=e^{-hd_*}\eta_0(e^{-hd_*}x),\qquad h=0,1,2,
\]
have disjoint supports whose radial hull is contained in
\(\{x>0:x^2/2\in I_3\}\). This is possible by taking \(d_*\) small and
the support of \(\eta_0\) sufficiently narrow. Define
\[
m_q=\int_0^\infty x^q\eta_0(x)\,dx>0.
\]
A change of variable gives
\[
\int_0^\infty x^q\eta_h(x)\,dx=m_qe^{hd_*q}.
\]
Seek the two tangential components of $v$ in the form
\begin{equation}
A_{0,\mathrm{moment}}=\lambda^{-1-\delta}
                         \sum_{h=0}^2 b_h(t,z)\eta_h(r/\lambda),
\qquad
B_{0,\mathrm{moment}}=\lambda^{-1-\delta}
                         \sum_{h=0}^1 c_h(t,z)\eta_h(r/\lambda).
\label{init:moment-profiles}
\end{equation}
For the angular equations define a \(3\)-by-\(3\) matrix and for
the axial equations a \(2\)-by-\(2\) matrix by
\[
\begin{aligned}
(\mathsf A_\theta)_{1h}&=m_2e^{2hd_*},\\
(\mathsf A_\theta)_{2h}&=2A_*m_{-2-2\mu}e^{-(2+2\mu)hd_*},\\
(\mathsf A_\theta)_{3h}&=A_*m_{-2\mu}e^{-2\mu hd_*},
 \qquad h=0,1,2,\\
(\mathsf A_z)_{1h}&=m_1e^{hd_*},\\
(\mathsf A_z)_{2h}&=A_*m_{1-2\mu}e^{(1-2\mu)hd_*},
 \qquad h=0,1.
\end{aligned}
\]
After dividing their rows by the nonzero prefactors, these are
Vandermonde matrices. The angular powers are
\(2,-2-2\mu,-2\mu\), and the axial powers are \(1,1-2\mu\).
Writing \(q_i\) for the powers in the respective list, they are
distinct, so the determinants are nonzero: the row-normalized
angular determinant is
\(\prod_{i<j}(e^{d_*q_j}-e^{d_*q_i})\), and for the axial matrix it
is \(e^{d_*q_2}-e^{d_*q_1}\). Thus the five coefficients are
uniquely specified by
\begin{equation}
\begin{aligned}
\mathsf A_\theta\begin{pmatrix}b_0\\b_1\\b_2\end{pmatrix}
 &=\begin{pmatrix}
 0\\ \lambda^{2+2\delta}\mathcal P\\
 \lambda^{2\delta}
 \mathcal J^{z\theta}
 \end{pmatrix},\\
\mathsf A_z\begin{pmatrix}c_0\\c_1\end{pmatrix}
 &=\begin{pmatrix}
 0\\-\lambda^{-1+2\delta}
 \mathcal J^{\theta z}
 \end{pmatrix}.
\end{aligned}
\label{init:moment-coefficient-equations}
\end{equation}
Indeed, substitution of \(r=\lambda x\) gives the physical row
factors \(\lambda^{2-\delta}\), \(\lambda^{1-\delta}\),
\(\lambda^{-2-2\delta}\), \(\lambda^{1-2\delta}\), and
\(\lambda^{-2\delta}\) in the five equations, in that order.
The right-hand sides of \eqref{init:moment-coefficient-equations}
are the normalized quantities specified above multiplied by
bounded smooth powers of \(2^k\lambda\). The matrices and their inverses have
bounded fixed derivatives on the compact profile range, since
\(A_*\ge a_*/2\) and \(\mu\) is fixed. Therefore
\begin{equation}
\sum_{h=0}^2|\partial_{\hat T,\hat Z}^\beta b_h|
 +\sum_{h=0}^1|\partial_{\hat T,\hat Z}^\beta c_h|
 \lesssim k^{C_\beta}2^{-(1-10^{-5})\delta k}.
\label{init:moment-coefficient-bounds}
\end{equation}
Although reference-band normalizations were used to estimate the
right-hand sides, the five physical target integrals determine
the same coefficients on overlapping bands.

\emph{The radial component and support.}
Set
\begin{equation}
\begin{aligned}
\Psi_{0,\mathrm{moment}}(t,r,z)
 &=\frac1r\int_0^r sB_{0,\mathrm{moment}}(t,s,z)\,ds,\\
u_{0,\mathrm{moment}}
 &=A_{0,\mathrm{moment}}e_\theta
       +\nabla\times(\Psi_{0,\mathrm{moment}}e_\theta).
\end{aligned}
\label{init:moment-velocity}
\end{equation}
The second equation in \eqref{init:moment-five-equations} makes
the full integral \(\int_0^\infty sB_{0,\mathrm{moment}}\,ds\)
zero. Hence the integral from \(0\) to \(r\) vanishes both before
the first profile and after the last. The potential and all
three components are smooth and supported in their radial hull;
in particular, their \(R\)-supports lie in \(I_3\).
On this fixed support interval, \(\omega\) is bounded below by a
positive constant. Thus the unweighted estimates for these fields
and their derivatives imply the bounds below with the factor
\(\omega d_{\mathrm{edge}}^{-C_\beta}\).
Their exact component formulas are
\[
v^r=-\partial_z\Psi_{0,\mathrm{moment}},\qquad
v^\theta=A_{0,\mathrm{moment}},\qquad
v^z=\left(\partial_r+\frac1r\right)\Psi_{0,\mathrm{moment}}
     =B_{0,\mathrm{moment}}.
\]
Since these fields are independent of \(Y_{i_0}\), their physical
radial derivative is \(\partial_r\). In particular,
\[
\left(\partial_r+\frac1r\right)v^r+\partial_zv^z
 =-\partial_z\left(\partial_r+\frac1r\right)
                    \Psi_{0,\mathrm{moment}}
   +\partial_zB_{0,\mathrm{moment}}=0.
\]
Thus \(v\) is divergence-free. The first two equations in
\eqref{init:moment-five-equations} preserve the two linear
velocity moments exactly.

The fixed profiles and \eqref{init:moment-coefficient-bounds}
give, at every fixed slow derivative order,
\begin{equation}
\begin{aligned}
|\partial_{\hat R,\hat Z,\hat T}^\beta v^{\mathrm{tan}}|
 &\lesssim k^{C_\beta}2^k2^{10^{-5}\delta k}
                  \omega d_{\mathrm{edge}}^{-C_\beta},\\
|\partial_{\hat R,\hat Z,\hat T}^\beta\Psi_{0,\mathrm{moment}}|
 &\lesssim k^{C_\beta}2^{10^{-5}\delta k}
                  \omega d_{\mathrm{edge}}^{-C_\beta},\\
|\partial_{\hat R,\hat Z,\hat T}^\beta v^r|
 &\lesssim k^{C_\beta}2^{(1-\delta)k}2^{10^{-5}\delta k}
                  \omega d_{\mathrm{edge}}^{-C_\beta}.
\end{aligned}
\label{init:moment-physical-bounds}
\end{equation}
For the second line, \(r^{-1}\int_0^r s(\cdot)\,ds\) contributes
\(2^{-k}\); for the third use
\(\partial_z=2^{(1-\delta)k}\partial_{\hat Z}\). Consequently,
\begin{equation}
\begin{aligned}
2^{-(1+\delta)k}|\partial_{\hat R,\hat Z,\hat T}^\beta v^{\mathrm{tan}}|
 &\lesssim k^{C_\beta}2^{-0.99\delta k}
                  \omega d_{\mathrm{edge}}^{-C_\beta},\\
2^{-(1+\delta)k}|\partial_{\hat R,\hat Z,\hat T}^\beta v^r|
 &\lesssim k^{C_\beta}2^{-1.99\delta k}
                  \omega d_{\mathrm{edge}}^{-C_\beta}.
\end{aligned}
\label{init:moment-velocity-bounds}
\end{equation}
Moreover, on the support of \(v\), the lower bound in
\eqref{init:moment-background-patch} implies the pointwise comparison
\[
|v^\theta|\lesssim k^C2^{-(1-10^{-5})\delta k}u_B^\theta.
\]

\emph{The recomputed defects.}
First keep the pressure and flat error fixed, and put
\[
u_*=u+v,\qquad p_*=p,\qquad
f_*=\mathsf{NS}(u_*,p)-E.
\]
A star on a moment or defect denotes the quantity computed from this state.
Because \(v\) is independent of the averaging variables and
\(\langle u\rangle_{\theta,Y_{i_0}}=u_B\), the radial equation gives
\begin{equation}
\langle f_*^r-f^r\rangle_{\theta,Y_{i_0}}
 =-\frac{2u_B^\theta v^\theta}{r}+\mathsf Q^r,
\label{init:moment-radial-change}
\end{equation}
where the remaining terms are explicitly
\begin{equation}
\begin{aligned}
\mathsf Q^r={}&\partial_tv^r
 +\left(\partial_r+\frac1r\right)
                         (2u_B^rv^r+(v^r)^2)\\
 &+\partial_z(u_B^rv^z+u_B^zv^r+v^rv^z)
 -\frac{(v^\theta)^2}{r}
 -\left(\partial_r^2+\frac1r\partial_r
                       +\partial_z^2-\frac1{r^2}\right)v^r.
\end{aligned}
\label{init:moment-radial-remainder}
\end{equation}
In particular, every mixed product of \(v\) and the old velocity
correction has zero double average. Integrating
\eqref{init:moment-radial-change} and using the third moment
equation yields
\[
\begin{aligned}
\int_0^\infty\langle f_*^r\rangle_{\theta,Y_{i_0}}\,dr
 &=\underbrace{\int_0^\infty\langle f^r\rangle_{\theta,Y_{i_0}}\,dr
       -\int_0^\infty\frac{2u_B^\theta}{r}v^\theta\,dr}_{=0}
       +\int_0^\infty\mathsf Q^r\,dr.
\end{aligned}
\]
The changes in the two physical moment defects are
\[
\begin{aligned}
\mathcal J_*^{\theta z}-\mathcal J^{\theta z}
 &=\int_0^\infty r^2(u_B^\theta v^z+u_B^zv^\theta+v^\theta v^z)\,dr,\\
\mathcal J_*^{z\theta}-\mathcal J^{z\theta}
 &=\int_0^\infty r\bigl(2u_B^zv^z+(v^z)^2-u_B^\theta v^\theta\bigr)\,dr
       +\frac12\int_0^\infty r^2\mathsf Q^r\,dr.
\end{aligned}
\]
Use \(u_B^z=0\) on the entire support of \(v\) and the last two
moment equations. All three new defects are therefore given
exactly by
\begin{equation}
\begin{aligned}
\mathcal P_*&=\int_0^\infty\mathsf Q^r\,dr,\\
\mathcal J_*^{\theta z}
 &=\int_0^\infty r^2v^\theta v^z\,dr,\\
\mathcal J_*^{z\theta}
 &=
       \int_0^\infty r(v^z)^2\,dr
       +\frac12\int_0^\infty r^2\mathsf Q^r\,dr.
\end{aligned}
\label{init:moment-exact-new-defects}
\end{equation}

We now estimate each term of \(\mathsf Q^r\). Since \(v\) and the
background are auxiliary-independent, all derivatives here are
slow derivatives. The background bound
\(|\partial_{\hat R,\hat Z,\hat T}^\beta u_B^r|\lesssim k^{C_\beta}2^k\) and
\eqref{init:moment-physical-bounds} give
\[
\begin{aligned}
&2^{-(3+\delta)k}\Bigl(
 |\partial_{\hat R,\hat Z,\hat T}^\beta\partial_tv^r|
 +\left|\partial_{\hat R,\hat Z,\hat T}^\beta\!\left[
          (\partial_r+r^{-1})(2u_B^rv^r)\right]\right|
 +|\partial_{\hat R,\hat Z,\hat T}^\beta\partial_z(u_B^rv^z)|\Bigr)\\
&\qquad\lesssim k^{C_\beta}2^{-(2-10^{-5})\delta k}
                  \omega d_{\mathrm{edge}}^{-C_\beta},\\
&2^{-(3+\delta)k}
 \left|\partial_{\hat R,\hat Z,\hat T}^\beta\!\left[
 (\partial_r^2+r^{-1}\partial_r+\partial_z^2-r^{-2})v^r\right]\right|
 \lesssim k^{C_\beta}2^{-(2-10^{-5})\delta k}
                  \omega d_{\mathrm{edge}}^{-C_\beta},\\
&2^{-(3+\delta)k}\Bigl(
 |\partial_{\hat R,\hat Z,\hat T}^\beta[(\partial_r+r^{-1})(v^r)^2]|
 +|\partial_{\hat R,\hat Z,\hat T}^\beta\partial_z(v^rv^z)|\Bigr)
 \lesssim k^{C_\beta}2^{-(3-2\cdot10^{-5})\delta k}
                  \omega d_{\mathrm{edge}}^{-C_\beta},\\
&2^{-(3+\delta)k}
 \left|\partial_{\hat R,\hat Z,\hat T}^\beta\frac{(v^\theta)^2}{r}\right|
 \lesssim k^{C_\beta}2^{-(1-2\cdot10^{-5})\delta k}
                  \omega d_{\mathrm{edge}}^{-C_\beta}.
\end{aligned}
\]
The term \(\partial_z(u_B^zv^r)\) vanishes because its product is
identically zero. Thus
\begin{equation}
2^{-(3+\delta)k}|\partial_{\hat R,\hat Z,\hat T}^\beta\mathsf Q^r|
 \lesssim k^{C_\beta}2^{-(1-2\cdot10^{-5})\delta k}
                  \omega d_{\mathrm{edge}}^{-C_\beta}.
\label{init:moment-radial-remainder-bound}
\end{equation}
For example, the first equation of
\eqref{init:moment-exact-new-defects} now gives
\[
2^{-2(1+\delta)k}|\partial_{\hat T,\hat Z}^\beta\mathcal P_*|
 \lesssim k^{C_\beta}2^{-2(1+\delta)k}2^{-k}
              2^{(3+\delta)k}2^{-(1-2\cdot10^{-5})\delta k}
 =k^{C_\beta}2^{-(2-2\cdot10^{-5})\delta k}.
\]
After normalization, the other two equations give the same bound: the normalized
velocity products have size
\(k^{C_\beta}2^{-2(1-10^{-5})\delta k}\), and the weighted
\(\mathsf Q^r\) integral has the same size after normalization. Consequently,
\begin{equation}
\begin{aligned}
&2^{-2(1+\delta)k}|\partial_{\hat T,\hat Z}^\beta \mathcal P_*|
 +2^{(1-2\delta)k}|\partial_{\hat T,\hat Z}^\beta \mathcal J_*^{\theta z}|
 +2^{-2\delta k}|\partial_{\hat T,\hat Z}^\beta \mathcal J_*^{z\theta}|\\
&\qquad\lesssim k^{C_\beta}2^{-(2-2\cdot10^{-5})\delta k}
 \le k^{C_\beta}2^{-1.2\delta k}.
\end{aligned}
\label{init:moment-improved-defects}
\end{equation}

\emph{The force before pressure reconstruction.}
The force changes by the exact identity
\begin{equation}
f_*-f=\mathsf L_B(v,0)+(v\cdot\nabla)v
 +(u_{0,\mathrm{aux}}\cdot\nabla)v
 +(v\cdot\nabla)u_{0,\mathrm{aux}}
 +(w_0\cdot\nabla)v+(v\cdot\nabla)w_0.
\label{init:moment-force-update}
\end{equation}
Only the last two terms have nonzero angular modes. The derivative
bounds used in Section~\ref{init:auxiliary-correction} and
\eqref{init:moment-physical-bounds} yield
\[
\begin{aligned}
2^{-(3+\delta)k}|(v\cdot\nabla)w_0|
 &\lesssim k^C2^{10^{-5}\delta k}
                  \sqrt\omega\,d_{\mathrm{edge}}^{-C},\\
2^{-(3+\delta)k}|(w_0\cdot\nabla)v|
 &\lesssim k^C2^{-(1/2-10^{-5})\delta k}
                  \sqrt\omega\,d_{\mathrm{edge}}^{-C}.
\end{aligned}
\]
Therefore the incoming nonzero-mode bound remains
\(k^C2^{0.01\delta k}\sqrt\omega\,d_{\mathrm{edge}}^{-C}\).
To check the angular mean, define its convective flux change in
this calculation by
\[
\mathsf C^{ab}
 =(u_B^a+u_{0,\mathrm{aux}}^a)v^b
 +(u_B^b+u_{0,\mathrm{aux}}^b)v^a+v^av^b,
 \qquad a,b\in\{r,\theta,z\}.
\]
With \(\Delta_0=D_r^2+r^{-1}D_r+\partial_z^2\), the two
components satisfy
\begin{equation}
\begin{aligned}
\langle f_*^\theta-f^\theta\rangle_\theta
 &=\partial_tv^\theta-(\Delta_0-r^{-2})v^\theta
       +(D_r+2/r)\mathsf C^{r\theta}
       +\partial_z\mathsf C^{z\theta},\\
\langle f_*^z-f^z\rangle_\theta
 &=\partial_tv^z-\Delta_0v^z
       +(D_r+1/r)\mathsf C^{rz}
       +\partial_z\mathsf C^{zz}.
\end{aligned}
\label{init:moment-tangential-force}
\end{equation}
The product bounds, for every fixed \(\beta\in\mathbb Z_{\ge0}^5\), are
\[
\begin{aligned}
2^{-2(1+\delta)k}|\partial_{\hat R,\hat Z,\hat T,Y_{i_0,1},Y_{i_0,2}}^\beta\mathsf C^{ra}|
 &\lesssim k^{C_\beta}2^{-(2-10^{-5})\delta k}
                  \omega d_{\mathrm{edge}}^{-C_\beta},\\
2^{-2(1+\delta)k}|\partial_{\hat R,\hat Z,\hat T,Y_{i_0,1},Y_{i_0,2}}^\beta\mathsf C^{za}|
 &\lesssim k^{C_\beta}2^{-(1-10^{-5})\delta k}
                  \omega d_{\mathrm{edge}}^{-C_\beta},
 \qquad a=\theta,z.
\end{aligned}
\]
Taking one radial derivative costs at most
\(k^C2^k2^{10^{-5}\delta k}\), and an axial derivative costs
\(2^{(1-\delta)k}\). For the increment itself,
\(D_tv=\partial_tv\) and \(D_rv=\partial_rv\), because it is
auxiliary-independent. Hence, for \(a=\theta,z\), and
\(c_\theta=2,c_z=1\),
\[
\begin{aligned}
2^{-(3+\delta)k}\bigl(|\partial_tv^a|+|\Delta_0v^a|+r^{-2}|v^a|\bigr)
 &\lesssim k^C2^{-(1-10^{-5})\delta k}
                  \omega d_{\mathrm{edge}}^{-C},\\
2^{-(3+\delta)k}|(D_r+c_a/r)\mathsf C^{ra}|
 &\lesssim k^C2^{-(1-2\cdot10^{-5})\delta k}
                  \omega d_{\mathrm{edge}}^{-C},\\
2^{-(3+\delta)k}|\partial_z\mathsf C^{za}|
 &\lesssim k^C2^{-(1-10^{-5})\delta k}
                  \omega d_{\mathrm{edge}}^{-C}.
\end{aligned}
\]
Both the auxiliary average and the auxiliary oscillation of these
changes have the same bounds. Since \(1-2\cdot10^{-5}>0.99\), the
two incoming tangential bounds in
\eqref{init:auxiliary-output-bounds} are preserved.

The radial angular mean also satisfies the following estimate at
every fixed order of differentiation in these variables:
\begin{equation}
2^{-(3+\delta)k}
 |\partial_{\hat R,\hat Z,\hat T,Y_{i_0,1},Y_{i_0,2}}^\beta\langle f_*^r\rangle_\theta|
 \lesssim k^{C_\beta}2^{10^{-5}\delta k}
                  \omega d_{\mathrm{edge}}^{-C_\beta}.
\label{init:moment-pressure-source-bound}
\end{equation}
Indeed, the radial version of
\eqref{init:moment-tangential-force} is
\[
\langle f_*^r-f^r\rangle_\theta
 =\partial_tv^r+(D_r+r^{-1})\mathsf C^{rr}
       +\partial_z\mathsf C^{zr}-\mathsf C^{\theta\theta}/r
       -(\Delta_0-r^{-2})v^r.
\]
The largest new term has the explicit bound
\[
2^{-(3+\delta)k}
 \left|\partial_{\hat R,\hat Z,\hat T,Y_{i_0,1},Y_{i_0,2}}^\beta\frac{\mathsf C^{\theta\theta}}r\right|
 \lesssim k^{C_\beta}2^{10^{-5}\delta k}
                  \omega d_{\mathrm{edge}}^{-C_\beta}.
\]
The derivative terms have strictly better bounds by the same
product estimates. Combining this with the incoming radial
bound \eqref{init:pressure-radial-source} proves
\eqref{init:moment-pressure-source-bound}.

\subsection{Pressure reconstruction and the initialized state}
\label{init:pressure-reconstruction}

The input to this step is the state \((u_*,p_*)\) constructed in
Section~\ref{init:moment-correction}, with
\[
\mathsf{NS}(u_*,p_*)=E+f_*,\qquad p_*=p_B+p_{w_0}.
\]
We keep \(u_*\) fixed and add an axisymmetric pressure
\(\widetilde p^{[0]}\) supported in \(R\in[R_a,R_h]\)
to reduce \(\langle f_*^r\rangle_\theta\) to a prescribed radial profile.

\emph{The compactly supported pressure.}
Choose a nonnegative function
\(\rho_*\in C_c^\infty((0,\infty))\) such that
\[
\int_0^\infty\rho_*(s)\,ds=1,\qquad
\{s^2:s\in\operatorname{supp}\rho_*\}\Subset I_3,
\]
and set
\begin{equation}
\begin{gathered}
\rho_r(t,r,z)=\frac1{\sqrt2\lambda(t,z)}
             \rho_*\!\left(\frac r{\sqrt2\lambda(t,z)}\right),\\
\int_0^\infty\rho_r\,dr=1,\qquad
\operatorname{supp}_r\rho_r
 =\sqrt2\lambda(t,z)\operatorname{supp}\rho_*,\qquad
\operatorname{supp}_R\rho_r\Subset I_3.
\end{gathered}
\label{init:pressure-bump}
\end{equation}
Since \(\operatorname{supp}_R\rho_r\Subset I_3\) is fixed,
\(\omega\) is bounded below there. Thus, for each fixed \(\beta\),
\[
|\partial_{\hat R,\hat Z,\hat T}^\beta(2^{-k}\rho_r)|
 \lesssim C_\beta\omega d_{\mathrm{edge}}^{-C_\beta}.
\]
Use the cutoff \(\chi_m\) and the radial primitive
\(\mathcal I_0\) from \eqref{init:auxiliary-radial-primitive}, and define
\begin{equation}
\widetilde p^{[0]}
 =\mathcal I_0\!\left(
       -\langle f_*^r\rangle_\theta
       +\rho_r\mathcal P_*\right).
\label{init:pressure-definition}
\end{equation}
By \eqref{init:pressure-moment} and \eqref{init:pressure-bump},
\[
\begin{gathered}
\mathcal P_*=\int_0^\infty
                   \langle f_*^r\rangle_{\theta,Y_{i_0}}\,dr,\\
\int_0^\infty\bigl[-\langle f_*^r\rangle_{\theta,Y_{i_0}}
                   +\rho_r\mathcal P_*\bigr]\,dr=0.
\end{gathered}
\]
The second identity holds after auxiliary averaging. Evaluating the
unaveraged source at \(Y_{i_0}=Y_{i_0}(t,s)\) before integrating need
not give zero, so the cutoff subtraction in \(\mathcal I_0\) is still
needed.
The definition of \(\mathcal P_*\) is independent of the reference
band; the common-coordinate compatibility of \(\mathcal I_0\) then
gives a single extended pressure field. After physical evaluation,
\begin{equation}
\begin{aligned}
&\widetilde p^{[0]}(t,r,z,Y_{i_0}(t,r))\\
&=\int_0^r\Bigl[-\langle f_*^r\rangle_\theta(t,s,z,Y_{i_0}(t,s))
       +\rho_r(t,s,z)\mathcal P_*(t,z)\Bigr]\,ds\\
&\quad-\chi_m(t,r,z)\int_0^\infty
 \Bigl[-\langle f_*^r\rangle_\theta(t,s,z,Y_{i_0}(t,s))
       +\rho_r(t,s,z)\mathcal P_*(t,z)\Bigr]\,ds.
\end{aligned}
\label{init:pressure-physical-formula}
\end{equation}
Below the annulus both the integral from \(0\) to \(r\) and \(\chi_m\)
vanish; above it this integral equals the full integral and
\(\chi_m=1\). Thus the pressure vanishes outside
\(R\in[R_a,R_h]\). Its smoothness at the two edges follows from
that of the source and the fixed cutoff, as for
\eqref{init:auxiliary-fixed-shift}.

Differentiating \eqref{init:pressure-definition} along the physical
radial direction gives
\begin{equation}
D_r\widetilde p^{[0]}
 =-\langle f_*^r\rangle_\theta
       +\rho_r\mathcal P_*+\mathcal R_{0,\mathrm p},
\label{init:pressure-radial-identity}
\end{equation}
where the cutoff remainder is
\begin{equation}
\begin{aligned}
\mathcal R_{0,\mathrm p}(t,r,z,Y_{i_0})
 ={}&(\partial_r\chi_m)\Bigl[
  \int_0^\infty\langle f_*^r\rangle_\theta
       \bigl(t,s,z,Y_{i_0}+\Lambda_r^{i_0}(s^{d_0}-r^{d_0})v_r\bigr)\,ds\\
 &\hspace{38mm}-\int_0^\infty
                 \langle f_*^r\rangle_{\theta,Y_{i_0}}(t,s,z)\,ds\Bigr].
\end{aligned}
\label{init:pressure-cutoff-remainder}
\end{equation}
Translation preserves the torus average, so the two integrals
cancel after auxiliary averaging:
\[
\langle\mathcal R_{0,\mathrm p}\rangle_{Y_{i_0}}=0.
\]
Averaging \eqref{init:pressure-radial-identity} and using the zero
pressure below the annulus now gives the exact formula
\begin{equation}
\begin{aligned}
\langle\widetilde p^{[0]}\rangle_{Y_{i_0}}(t,r,z)
 ={}&-\int_0^r\langle f_*^r\rangle_{\theta,Y_{i_0}}(t,s,z)\,ds\\
 &+\left(\int_0^r\rho_r(t,s,z)\,ds\right)
       \int_0^\infty\langle f_*^r\rangle_{\theta,Y_{i_0}}(t,s,z)\,ds.
\end{aligned}
\label{init:pressure-auxiliary-average}
\end{equation}

\emph{Flatness of the cutoff remainder.}
Subtracting the auxiliary mean inside
\eqref{init:pressure-cutoff-remainder} leaves only the nonzero
Fourier modes. Write
\[
\langle f_*^r\rangle_\theta
 -\langle f_*^r\rangle_{\theta,Y_{i_0}}
 =\sum_{n\in\mathbb Z^2\setminus\{0\}}
                   a_n(t,r,z)e^{2\pi\mathrm i n\cdot Y_{i_0}}.
\]
For fixed \(t,z\), the coefficients are smooth, supported in
\(r\asymp2^{-k}\), and vanish to every order at the support edges.
The \(n\)-th mode of the full radial integral is, up to a factor
of modulus one,
\[
\int_0^\infty a_n(t,s,z)
          e^{2\pi\mathrm i\Lambda_r^{i_0}(n\cdot v_r)s^{d_0}}\,ds.
\]
One integration by parts gives
\[
\begin{aligned}
&\int_0^\infty a_n(t,s,z)
          e^{2\pi\mathrm i\Lambda_r^{i_0}(n\cdot v_r)s^{d_0}}\,ds\\
&\quad=-\int_0^\infty\partial_s\!\left(
 \frac{a_n(t,s,z)}
 {2\pi\mathrm i\Lambda_r^{i_0}d_0s^{d_0-1}(n\cdot v_r)}\right)
          e^{2\pi\mathrm i\Lambda_r^{i_0}(n\cdot v_r)s^{d_0}}\,ds.
\end{aligned}
\]
There is no boundary term. By \eqref{eq:ck-size} and
\eqref{eq:diophantine-directions}, each repetition gains
\[
\frac{2^{kd_0}}{\Lambda_r^{i_0}}
 \lesssim k^C2^{-10^{-5}\delta k},\qquad
|n\cdot v_r|^{-1}\lesssim1+|n|.
\]
After \(J\) repetitions, sufficiently many auxiliary derivatives
in \eqref{init:moment-pressure-source-bound} absorb the factor
\((1+|n|)^J\) and leave the summable factor \((1+|n|)^{-3}\).
The radial integral contributes \(2^{-k}\), while
\(|\partial_r\chi_m|\lesssim2^k\). Consequently, for every
integer \(J\ge1\),
\begin{equation}
2^{-(3+\delta)k}|\mathcal R_{0,\mathrm p}|
 \lesssim_J k^{C_J}2^{\delta k}
       \bigl(2^{-10^{-5}\delta k}\bigr)^J
       \sum_{n\ne0}(1+|n|)^{-3}.
\label{init:pressure-flat-gain}
\end{equation}
Here we used \(2^{10^{-5}\delta k}\le2^{\delta k}\).
For every fixed physical derivative, differentiating the
coefficients, the translation and the cutoff costs at most
\(k^{C_\beta}2^{C_\beta k}\). Choosing \(J\) after that derivative
order and the required power proves
\begin{equation}
|D_{t,x}^\beta(\mathcal R_{0,\mathrm p}e_r)|
 \le C_{\beta,N}\lambda^N
 \quad\text{for every fixed }\beta\text{ and }N>0.
\label{init:pressure-flat-bound}
\end{equation}

\emph{Pressure estimates and the updated force.}
The source in \eqref{init:pressure-definition}, including the bump
term, satisfies the bound \eqref{init:moment-pressure-source-bound},
by the weighted bound for \(\rho_r\) above and
\eqref{init:moment-improved-defects}.
The radial primitive gives a factor \(2^{-k}\), also after differentiation
to any fixed order in these variables. To see this without differentiating the
large auxiliary translation, use
\[
U=(2^kr)^{d_0},\qquad
\sigma=(2^ks)^{d_0}-U,\qquad
 ds=\frac{2^{-k}}{d_0}(U+\sigma)^{1/d_0-1}\,d\sigma.
\]
As in \eqref{init:auxiliary-fixed-shift}, split the primitive into
the integrals over \(\sigma<0\) and \(\sigma>0\).
At fixed \(\sigma\), its auxiliary argument is
\(Y_{i_0}+2^{-kd_0}\Lambda_r^{i_0}\sigma v_r\), independent of the
normalized slow variables. The integration length, the weight in
\(ds\), and their derivatives of any fixed order in these variables are bounded on
the normalized annulus. The endpoint estimate following
\eqref{init:auxiliary-fixed-shift} preserves the same \(\omega\),
up to a power of \(d_{\mathrm{edge}}^{-1}\). Thus
\begin{equation}
\begin{aligned}
|\partial_{\hat R,\hat Z,\hat T,Y_{i_0,1},Y_{i_0,2}}^\beta\widetilde p^{[0]}|
 &\lesssim k^{C_\beta}2^{-k}2^{(3+\delta)k}2^{10^{-5}\delta k}
                  \omega d_{\mathrm{edge}}^{-C_\beta},\\
2^{-2(1+\delta)k}
 |\partial_{\hat R,\hat Z,\hat T,Y_{i_0,1},Y_{i_0,2}}^\beta\widetilde p^{[0]}|
 &\lesssim k^{C_\beta}2^{-(1-10^{-5})\delta k}
                  \omega d_{\mathrm{edge}}^{-C_\beta}\\
 &\le k^{C_\beta}2^{-0.99\delta k}
                  \omega d_{\mathrm{edge}}^{-C_\beta}.
\end{aligned}
\label{init:pressure-size}
\end{equation}
Since \(\partial_z=2^{(1-\delta)k}\partial_{\hat Z}\), this also gives
\begin{equation}
\begin{aligned}
2^{-(3+\delta)k}|\partial_z\widetilde p^{[0]}|
 &\lesssim k^C2^{-k}2^{(1-\delta)k}2^{10^{-5}\delta k}
                  \omega d_{\mathrm{edge}}^{-C}\\
 &=k^C2^{-(1-10^{-5})\delta k}\omega d_{\mathrm{edge}}^{-C}\\
 &\le k^C2^{-0.99\delta k}\omega d_{\mathrm{edge}}^{-C}.
\end{aligned}
\label{init:pressure-axial-force}
\end{equation}

Define the completed state by
\begin{equation}
\begin{aligned}
u^{[0]}&=u_B+w_0+u_{0,\mathrm{aux}}+u_{0,\mathrm{moment}},\\
p^{[0]}&=p_B+p_{w_0}+\widetilde p^{[0]},\\
E^{[0]}&=E+\mathcal R_{0,\mathrm p}e_r,\\
f^{[0]}&=f_*+\nabla\widetilde p^{[0]}
                   -\mathcal R_{0,\mathrm p}e_r.
\end{aligned}
\label{init:moment-completed-state}
\end{equation}
Then \(\mathsf{NS}(u^{[0]},p^{[0]})=E^{[0]}+f^{[0]}\), and
\begin{equation}
\langle(f^{[0]})^r\rangle_\theta
 =\rho_r\mathcal P_*,\qquad
\mathcal P^{[0]}=\mathcal P_*.
\label{init:moment-final-radial-identity}
\end{equation}
The \(\theta z\) moment is unchanged by this pressure update.
For the \(z\theta\) moment, the change of its physical integral is
exactly
\[
\begin{aligned}
&\int_0^\infty\left(
 r\langle\widetilde p^{[0]}\rangle_{Y_{i_0}}
 +\frac{r^2}{2}\partial_r
       \langle\widetilde p^{[0]}\rangle_{Y_{i_0}}
                     \right)\,dr\\
&\qquad=\left[\frac{r^2}{2}
 \langle\widetilde p^{[0]}\rangle_{Y_{i_0}}
                                     \right]_0^\infty=0.
\end{aligned}
\]
Here the flat remainder contributes zero double mean, and the
boundary term vanishes by compact support. Thus
\eqref{init:moment-improved-defects} holds for all three defects
of the completed state.

The pressure and the added flat term are axisymmetric, so
\begin{equation}
\begin{aligned}
f^{[0]}-\langle f^{[0]}\rangle_\theta
 &=f_*-\langle f_*\rangle_\theta,\\
\langle (f^{[0]})^\theta\rangle_\theta
 &=\langle f_*^\theta\rangle_\theta,\\
\langle (f^{[0]})^z\rangle_\theta
 &=\langle f_*^z\rangle_\theta+\partial_z\widetilde p^{[0]}.
\end{aligned}
\label{init:pressure-force-identities}
\end{equation}
In particular, the change in the axial double mean is
\(\partial_z\langle\widetilde p^{[0]}\rangle_{Y_{i_0}}\),
and the change in its auxiliary oscillation is
\(\partial_z(\widetilde p^{[0]}
 -\langle\widetilde p^{[0]}\rangle_{Y_{i_0}})\).
Both satisfy \eqref{init:pressure-axial-force}, up to a fixed
constant. Combining this with the bounds for \(f_*\) proved in
Section~\ref{init:moment-correction} gives
\begin{equation}
\begin{aligned}
2^{-(3+\delta)k}
 |f^{[0]}-\langle f^{[0]}\rangle_\theta|
 &\lesssim k^C2^{0.01\delta k}\sqrt\omega\,d_{\mathrm{edge}}^{-C},\\
2^{-(3+\delta)k}\Bigl(
 |\langle (f^{[0]})^{\mathrm{tan}}\rangle_\theta
       -\langle (f^{[0]})^{\mathrm{tan}}\rangle_{\theta,Y_{i_0}}|
 +|\langle (f^{[0]})^{\mathrm{tan}}\rangle_{\theta,Y_{i_0}}|\Bigr)
 &\lesssim k^C2^{-0.99\delta k}\omega d_{\mathrm{edge}}^{-C}.
\end{aligned}
\label{init:pressure-preserved-force-bounds}
\end{equation}
By the weighted bound for \(2^{-k}\rho_r\), the radial force obeys
\[
2^{-(3+\delta)k}|\langle(f^{[0]})^r\rangle_\theta|
 =(2^{-k}\rho_r)2^{\delta k}
   \bigl[2^{-2(1+\delta)k}|\mathcal P^{[0]}|\bigr]
 \lesssim k^C2^{-0.2\delta k}\omega d_{\mathrm{edge}}^{-C}.
\]
Together with \eqref{init:pressure-preserved-force-bounds}, this proves
\begin{equation}
\begin{aligned}
2^{-(3+\delta)k}
 \bigl(|f^{[0]}|+|f^{[0]}-\langle f^{[0]}\rangle_\theta|\bigr)
 &\lesssim k^C2^{0.01\delta k}\sqrt\omega\,d_{\mathrm{edge}}^{-C},\\
2^{-(3+\delta)k}\Bigl(
 |\langle(f^{[0]})^{\mathrm{tan}}\rangle_\theta
       -\langle(f^{[0]})^{\mathrm{tan}}\rangle_{\theta,Y_{i_0}}|
 +|\langle f^{[0]}\rangle_{\theta,Y_{i_0}}|\Bigr)
 &\lesssim k^C2^{-0.2\delta k}\omega d_{\mathrm{edge}}^{-C},\\
2^{-2(1+\delta)k}|\mathcal P^{[0]}|
 +2^{(1-2\delta)k}|\mathcal J^{\theta z,[0]}|
 +2^{-2\delta k}|\mathcal J^{z\theta,[0]}|
 &\lesssim k^C2^{-1.2\delta k}.
\end{aligned}
\label{init:moment-initialized-bounds}
\end{equation}
These imply \eqref{init:stage-zero-target}, since \(0.01<0.3\).
The velocity bounds follow from
\eqref{init:primary-pointwise-bounds},
\eqref{init:auxiliary-physical-bounds} and
\eqref{init:moment-physical-bounds}. The velocity remains
divergence-free, the two linear moments are preserved, and all
velocity and pressure increments are smooth and supported in
\(R\in[R_a,R_h]\). Finally,
\(\langle E^{[0]}\rangle_{\theta,Y_{i_0}}=\mathsf E_B\), and every fixed
physical derivative of \(E^{[0]}\) is flat.

For use in the first correction stage, write
\[
f^{[0]}-\langle f^{[0]}\rangle_\theta
 =\sum_{\gamma=(\ell,\alpha,\iota)}2^{(3+\delta)\ell}
   \sum_{m=\pm1,\pm2}\mathbf f_{\gamma,1,m}
                    e^{\mathrm imN_\ell\varphi_\gamma}.
\]
Here \(\mathbf f_{\gamma,1,m}\) is normalized in force units at the
native scale \(\ell\). For every fixed \(\beta\in\mathbb Z_{\ge0}^5\),
\[
|\partial_{\hat R,\hat Z,\hat T,Y_{i_0,1},Y_{i_0,2}}^\beta
       \mathbf f_{\gamma,1,m}|
 \lesssim k^{C_\beta}2^{0.01\delta k}
       \sqrt\omega\,d_{\mathrm{edge}}^{-C_\beta}
       \mathcal G_\gamma(\tau_\gamma).
\]
The primary linear and quadratic terms satisfy this bound by the
coefficient estimates in Section~\ref{init:primary-pulse}.
In \eqref{init:auxiliary-nonzero-modes} and
\eqref{init:moment-force-update}, each new nonzero harmonic term
contains a primary coefficient or a derivative of one, so Leibniz'
rule retains the same Gaussian at every fixed derivative order.
The pressure update leaves the nonzero angular force unchanged by
\eqref{init:pressure-force-identities}.
This completes the initialization.

\section{Iteration}
\label{sec:iteration-ppt}

We retain the notation of Section~\ref{sec:initialization-ppt}.
In particular, \(u^{[j]},p^{[j]}\) include the background fields,
\(\mathrm{tan}=(\theta,z)\), and
\[
f^{[j]}=\mathsf{NS}(u^{[j]},p^{[j]})-E^{[j]}.
\]
All averages are taken before physical evaluation, in the common
variable \(Y_{i_0}\), with \(i_0=i_0(k)\) as in
Section~\ref{sec:pulse-common}. Estimates below hold on \(\mathcal B_k\),
uniformly in the angular and auxiliary variables. We use the velocity,
pressure and force units
\[
2^{(1+\delta)k},\qquad 2^{2(1+\delta)k},\qquad
2^{(3+\delta)k},
\]
respectively. The physical defects \(\mathcal P^{[j]}\),
\(\mathcal J^{\theta z,[j]}\) and \(\mathcal J^{z\theta,[j]}\)
are those defined in Section~\ref{sec:initialization-ppt}, evaluated
at the current state; their normalization factors are included in
the estimates below. We use the same radial weights
\(\omega,d_{\mathrm{edge}}\) as in
\eqref{intro:background-edge-weight} throughout the iteration.

\emph{The induction requirements.}
For \(j\ge1\), we construct the increments in
\eqref{init:later-stage-fields}, with the sizes in
\eqref{init:later-increment-sizes}, so that
\begin{equation}
\begin{aligned}
2^{-(3+\delta)k}
 \bigl(|f^{[j]}|+|f^{[j]}-\langle f^{[j]}\rangle_\theta|\bigr)
 &\lesssim k^{C_j}2^{(0.3-0.1j)\delta k}
                    \sqrt\omega\,d_{\mathrm{edge}}^{-C_j},\\
2^{-(3+\delta)k}\Bigl(
 |\langle(f^{[j]})^{\mathrm{tan}}\rangle_\theta
       -\langle(f^{[j]})^{\mathrm{tan}}\rangle_{\theta,Y_{i_0}}|
 +|\langle f^{[j]}\rangle_{\theta,Y_{i_0}}|\Bigr)
 &\lesssim k^{C_j}2^{-(0.2+0.1j)\delta k}
                    \omega d_{\mathrm{edge}}^{-C_j},\\
2^{-2(1+\delta)k}|\mathcal P^{[j]}|
 +2^{(1-2\delta)k}|\mathcal J^{\theta z,[j]}|
 +2^{-2\delta k}|\mathcal J^{z\theta,[j]}|
 &\lesssim k^{C_j}2^{-(1.2+0.1j)\delta k}.
\end{aligned}
\label{iter:targets}
\end{equation}
Moreover, we require the derivative estimates corresponding to
\eqref{init:initial-increment-sizes}, \eqref{init:later-increment-sizes}
and \eqref{iter:targets} at every fixed order, with \(k\) and \(i_0(k)\)
held fixed. The constants and loss exponents may depend on the stage
and derivative order and are uniform in \(k\). 

For each pulse
\(w\in\{w_0\}\cup\{w_{q,m\ne0},w_{q,\mathrm{stress}}\mid 1\le q\le j\}\),
write
\[
\begin{aligned}
w
 &=\sum_{\gamma=(\ell,\alpha,\iota)}2^{(1+\delta)\ell}
   \sum_{0<|m|\le 2^j}\mathbf c_{\gamma,m}(w)e^{\mathrm imN_\ell\varphi_\gamma},\\
f^{[j]}-\langle f^{[j]}\rangle_\theta
 &=\sum_{\gamma=(\ell,\alpha,\iota)}2^{(3+\delta)\ell}
   \sum_{0<|m|\le 2^{j+1}}\mathbf f_{\gamma,j+1,m}e^{\mathrm imN_\ell\varphi_\gamma}.
\end{aligned}
\]
Throughout this section, \emph{labelled} means associated with a
specified pulse label \(\gamma=(\ell,\alpha,\iota)\), which records
the scale band, the slow box, and the pulse family. For each fixed
label \(\gamma\), \(m\) denotes the harmonic order relative to the
phase \(N_\ell\varphi_\gamma\).
For each labelled coefficient \(A\),
\[
\begin{aligned}
&\left|\partial_{\hat R,\hat Z,\hat T,Y_{i_0,1},Y_{i_0,2}}^\beta A\right|\\
&\quad\lesssim k^{C_{j,\beta}}\sqrt\omega\,
 d_{\mathrm{edge}}^{-C_{j,\beta}}\mathcal G_\gamma(\tau_\gamma)
 \begin{cases}
 2^{-0.5\delta k}, & A=\mathbf c_{\gamma,m}(w_0),\\
 2^{-(0.6+0.1q)\delta k}, & A=\mathbf c_{\gamma,m}(w_{q,m\ne0}),\\
 2^{-(0.6+0.1q-10^{-5})\delta k}, & A=\mathbf c_{\gamma,m}(w_{q,\mathrm{stress}}),\\
 2^{(0.3-0.1j)\delta k}, & A=\mathbf f_{\gamma,j+1,m}.
 \end{cases}
\end{aligned}
\]
Here \(\mathcal G_\gamma\) is the Gaussian envelope in \eqref{s6n:gaussian-envelope}.

For \(a\in\{\mathrm{aux},\mathrm{moment}\}\), use the variables
\[
X_{\mathrm{aux}}=(\hat R,\hat Z,\hat T,Y_{i_0,1},Y_{i_0,2}),
\qquad X_{\mathrm{moment}}=(\hat R,\hat Z,\hat T).
\]
The axisymmetric velocity increments satisfy
\[
\begin{aligned}
&2^{-(1+\delta)k}
 \bigl(|\partial_{X_a}^{\beta}u_{q,a}^{\mathrm{tan}}|
       +2^{\delta k}|\partial_{X_a}^{\beta}u_{q,a}^r|\bigr)\\
&\quad\lesssim k^{C_{j,\beta}}\omega d_{\mathrm{edge}}^{-C_{j,\beta}}
 \begin{cases}
 2^{-(1-10^{-5})\delta k}, & q=0,\\
 2^{-(1.1+0.1q-2\cdot10^{-5})\delta k}, & 1\le q\le j.
 \end{cases}
\end{aligned}
\]

The averaged force and the three physical defects satisfy
\[
\begin{aligned}
2^{-(3+\delta)k}\Bigl(
 &\bigl|\partial_{\hat R,\hat Z,\hat T,Y_{i_0,1},Y_{i_0,2}}^{\beta}
 \bigl(\langle(f^{[j]})^{\mathrm{tan}}\rangle_\theta
      -\langle(f^{[j]})^{\mathrm{tan}}\rangle_{\theta,Y_{i_0}}\bigr)\bigr|\\
 &+\bigl|\partial_{\hat R,\hat Z,\hat T}^{\eta}
              \langle f^{[j]}\rangle_{\theta,Y_{i_0}}\bigr|\Bigr)\\
 &\lesssim k^{C_{j,\beta,\eta}}2^{-(0.2+0.1j)\delta k}
             \omega d_{\mathrm{edge}}^{-C_{j,\beta,\eta}},
\end{aligned}
\]
\[
\begin{aligned}
 &2^{-2(1+\delta)k}\bigl|\partial_{\hat T,\hat Z}^{\mu}\mathcal P^{[j]}\bigr|
 +2^{(1-2\delta)k}\bigl|\partial_{\hat T,\hat Z}^{\mu}\mathcal J^{\theta z,[j]}\bigr|
 +2^{-2\delta k}\bigl|\partial_{\hat T,\hat Z}^{\mu}\mathcal J^{z\theta,[j]}\bigr|\\
 &\quad\lesssim k^{C_{j,\mu}}2^{-(1.2+0.1j)\delta k}.
\end{aligned}
\]

We preserve incompressibility, the support condition
\eqref{init:stage-support}, and the two linear moments
\eqref{init:linear-moments}. For each \(0\le q\le j\), the auxiliary
correction satisfies \(\langle u_{q,\mathrm{aux}}\rangle_{Y_{i_0}}=0\),
and the moment correction \(u_{q,\mathrm{moment}}\) is independent of
\(Y_{i_0}\). At every completed stage we also have
\[
\langle(f^{[j]})^r\rangle_\theta=\rho_r\mathcal P^{[j]},
\]
with \(\rho_r\) from \eqref{init:pressure-bump}. The flat error
continues to satisfy
\[
\langle E^{[j]}\rangle_{\theta,Y_{i_0}}=\mathsf E_B,\qquad
|D_{t,x}^{\beta}E^{[j]}|\le C_{j,\beta,N}\lambda^N
\quad(N>0),
\]
and \(E^{[j]}-\mathsf E_B\) is supported in \(R\in[R_a,R_h]\).
Incompressibility, support, the two linear moments and these
flat-error properties are also preserved at intermediate steps.
As in Section~\ref{init:pressure-reconstruction}, the axisymmetric
pressure is reconstructed once, after all four velocity corrections.

\emph{The starting state.}
Section~\ref{init:pressure-reconstruction} supplies
\[
u^{[0]}=u_B+w_0+u_{0,\mathrm{aux}}+u_{0,\mathrm{moment}},\qquad
p^{[0]}=p_B+p_{w_0}+\widetilde p^{[0]}.
\]
By \eqref{init:moment-initialized-bounds}, this state satisfies
\eqref{iter:targets} with \(j=0\); in fact,
\begin{equation}
2^{-(3+\delta)k}
 \bigl(|f^{[0]}|+|f^{[0]}-\langle f^{[0]}\rangle_\theta|\bigr)
 \lesssim k^C2^{0.01\delta k}\sqrt\omega\,d_{\mathrm{edge}}^{-C}.
\label{iter:starting-force}
\end{equation}
Since \(u_B\), \(u_{0,\mathrm{aux}}\), and \(u_{0,\mathrm{moment}}\)
are axisymmetric, all nonzero angular modes of \(u^{[0]}\) come from
the primary pulse \(w_0\). For each label \(\gamma\), the sine and
cosine in \eqref{init:primary-fields} give the two harmonics
\(m=\pm1\). The coefficient estimates in
Section~\ref{sec:initialization-ppt} give the derivative bounds
required above. We now give the nonzero-mode correction at a general stage
\(j\ge1\).

\subsection{Correction of the nonzero angular modes}
\label{iter:nonzero-modes}

Fix \(j\ge1\) and suppose stage \(j-1\) has been completed.
The construction below is the same for every such \(j\).
For a first reading, one may set \(j=1\) throughout.
The correction cancels the principal part of
\(f^{[j-1]}-\langle f^{[j-1]}\rangle_\theta\). We construct
\(w_{j,m\ne0},p_{j,m\ne0}\) and a new flat error \(E_{\mathrm{new}}\),
and set
\begin{equation}
\begin{aligned}
u_{\mathrm{new}}&=u^{[j-1]}+w_{j,m\ne0},&
p_{\mathrm{new}}&=p^{[j-1]}+p_{j,m\ne0},\\
f_{\mathrm{new}}&=\mathsf{NS}(u_{\mathrm{new}},p_{\mathrm{new}})
                          -E_{\mathrm{new}}.
\end{aligned}
\label{iter:nonzero-updated-state}
\end{equation}
We seek the nonzero-mode bound in \eqref{iter:targets} at stage \(j\).
The axisymmetric pressure \(\widetilde p^{[j-1]}\) is unchanged in
this step.

\emph{The source and its harmonics.}
In this subsection, \(w'\) ranges over the previously constructed
pulses \(w_0\) and \(w_{q,m\ne0},w_{q,\mathrm{stress}}\),
\(1\le q<j\), and \(p_{w'}\) denotes the corresponding pulse pressure.
All sums over \(w'\) are finite. Using the estimates for the axisymmetric velocity increments
and their derivatives given just after \eqref{iter:targets}, we obtain
\begin{equation}
\begin{aligned}
2^{-(1+\delta)k}\sum_{q=0}^{j-1}
 \sum_{a=\mathrm{aux},\mathrm{moment}}
 \bigl(|u_{q,a}^{\mathrm{tan}}|+2^{\delta k}|u_{q,a}^r|\bigr)
 &\lesssim k^{C_j}2^{-(1-10^{-5})\delta k}
                    \omega d_{\mathrm{edge}}^{-C_j},\\
2^{-(2+\delta)k}\sum_{q=0}^{j-1}
 \sum_{a=\mathrm{aux},\mathrm{moment}}|\nabla u_{q,a}|
 &\lesssim k^{C_j}2^{-(1-2\cdot10^{-5})\delta k}
                    \omega d_{\mathrm{edge}}^{-C_j}.
\end{aligned}
\label{iter:old-axisymmetric-bounds}
\end{equation}
The gradient includes physical differentiation of the auxiliary
variable and of the cylindrical frame. Both bounds in \eqref{iter:old-axisymmetric-bounds} hold with
\(\partial_{X_a}^{\beta}\) applied to the components of \(u_{q,a}\)
in the first line and to the components of \(\nabla u_{q,a}\) in the
second, for every fixed multi-index \(\beta\). Here \(X_a\) is defined
just after \eqref{iter:targets}.

For each \(j\ge1\), set
\[
\mathcal H_j=\{m\in\mathbb Z:0<|m|\le 2^j\}.
\]
This set contains the nonzero harmonics of the old velocity and force
at every label. Set the coefficients of absent harmonics to zero. Each old pulse has the form
\begin{equation}
\begin{aligned}
w'&=\sum_{\gamma=(\ell,\alpha,\iota)}2^{(1+\delta)\ell}
       \sum_{m\in\mathcal H_j}\mathbf c_{\gamma,m}(w')
                          e^{\mathrm imN_\ell\varphi_\gamma},\\
|\partial_{\hat R,\hat Z,\hat T,Y_{i_0,1},Y_{i_0,2}}^\beta\mathbf c_{\gamma,m}(w')|
 &\lesssim k^{C_{j,\beta}}2^{-0.5\delta k}
                          \sqrt\omega\,d_{\mathrm{edge}}^{-C_{j,\beta}}.
\end{aligned}
\label{iter:old-pulse-coefficients}
\end{equation}
Every displayed harmonic is divergence-free, and its coefficients
retain the assigned labelled support and smooth zero extension.
The construction formulas use the native scale \(\ell\); estimates
use the reference scale \(k\). Since \(|\ell-k|\le1\), changing the slow
normalization and pulling back to the common auxiliary variable
preserves the bounds at every fixed order of differentiation.
These bounds follow from the individual increment estimates by
finite summation. At \(j=1\) we use
\eqref{init:primary-coefficient-bounds},
\eqref{init:auxiliary-physical-bounds} and
\eqref{init:moment-physical-bounds}.

Subtracting the angular mean removes all terms involving only the
background and the axisymmetric corrections, including
\(\nabla\widetilde p^{[j-1]}\). With \(\mathsf L_B\) defined
after \eqref{intro:linear-quadratic}, the exact source identity is
\begin{equation}
\begin{aligned}
f^{[j-1]}-\langle f^{[j-1]}\rangle_\theta
={}&\sum_{w'}\mathsf L_B(w',p_{w'})
                -(E^{[j-1]}-\langle E^{[j-1]}\rangle_\theta)\\
&+\sum_{w',w''}\bigl((w'\cdot\nabla)w''
                   -\langle(w'\cdot\nabla)w''\rangle_\theta\bigr)\\
&+\sum_{w'}\sum_{q=0}^{j-1}\sum_{a=\mathrm{aux},\mathrm{moment}}
 \bigl((u_{q,a}\cdot\nabla)w'+(w'\cdot\nabla)u_{q,a}\bigr).
\end{aligned}
\label{iter:source-identity}
\end{equation}
Products from distinct labels vanish by
Proposition~\ref{s6n:separation}. Thus we write
\begin{equation}
f^{[j-1]}-\langle f^{[j-1]}\rangle_\theta
 =\sum_{\gamma=(\ell,\alpha,\iota)}2^{(3+\delta)\ell}
       \sum_{m\in\mathcal H_j}
       \mathbf f_{\gamma,j,m}e^{\mathrm imN_\ell\varphi_\gamma}.
\label{iter:source-harmonics}
\end{equation}
The label range and phases are those in \eqref{init:primary-fields};
\(\mathbf f_{\gamma,j,m}\) is normalized in force units at scale
\(\ell\). The coefficient form of the incoming nonzero-mode bound is
\begin{equation}
|\partial_{\hat R,\hat Z,\hat T,Y_{i_0,1},Y_{i_0,2}}^\beta\mathbf f_{\gamma,j,m}|
 \lesssim k^{C_{j,\beta}}2^{(0.4-0.1j)\delta k}
                    \sqrt\omega\,d_{\mathrm{edge}}^{-C_{j,\beta}}
                    \mathcal G_\gamma(\tau_\gamma).
\label{iter:source-bound}
\end{equation}
The source has the normalized and transverse supports of the
\(\gamma\)-pulse and smooth zero extension across their boundaries
and the radial annulus boundaries. It also extends smoothly by zero
at the fast-time endpoints, since the previously constructed pulses
and their derivatives have support inside the assigned time interval.

The Gaussian factor in \eqref{iter:source-bound} is part of the
induction requirements. This is the source estimate required by
Proposition~\ref{s7n:forced-inverse}. We verify it directly at \(j=1\)
below and retain it in the output of the present correction.

At \(j=1\) the only old pulse is \(w_0\), and
\[
E^{[0]}-\langle E^{[0]}\rangle_\theta
 =E_{\mathrm{primary}}-\mathsf E_B,
\]
where \(E_{\mathrm{primary}}\) is given by
\eqref{init:primary-flat-error}. The source has only \(m=\pm1,\pm2\).
Each coefficient contains a primary-pulse amplitude or a derivative
of one. Retaining its time envelope in the estimates of
Section~\ref{init:primary-pulse}, the four terms in
\eqref{iter:source-identity}, divided by \(2^{(3+\delta)\ell}\), have
the following bounds, with \(a\in\{\mathrm{aux},\mathrm{moment}\}\):
\[
\begin{array}{c|c}
\text{term}&\text{bound, apart from }k^{C_\beta}
 \sqrt\omega\,d_{\mathrm{edge}}^{-C_\beta}\mathcal G_\gamma(\tau_\gamma)\\ \hline
\mathsf L_B(w_0,p_{w_0})-(E_{\mathrm{primary}}-\mathsf E_B)
     &2^{3\cdot10^{-5}\delta k}\\
(w_0\cdot\nabla)w_0-\langle(w_0\cdot\nabla)w_0\rangle_\theta
     &2^{10^{-5}\delta k}\\
\sum_a(u_{0,a}\cdot\nabla)w_0&2^{10^{-5}\delta k}\\
\sum_a(w_0\cdot\nabla)u_{0,a}&2^{-(0.5-2\cdot10^{-5})\delta k}.
\end{array}
\]
The first two rows were proved in Section~\ref{init:primary-pulse}.
For the third row, a carrier derivative contributes \(2^\ell N_\ell\),
so its normalized size is at most
\[
k^{C_\beta}2^{\delta k}2^{-(1-10^{-5})\delta k}
       2^{\delta k/2}2^{-\delta k/2}
       \sqrt\omega\,d_{\mathrm{edge}}^{-C_\beta}
       \mathcal G_\gamma(\tau_\gamma).
\]
The fourth follows from \eqref{iter:old-axisymmetric-bounds}.
The same bounds hold after differentiation to any fixed order in the
indicated variables, by Leibniz' rule,
using \(\omega\le1\) and
\(\mathcal G_\gamma^2\lesssim\mathcal G_\gamma\). Hence
\[
\begin{aligned}
|\partial_{\hat R,\hat Z,\hat T,Y_{i_0,1},Y_{i_0,2}}^\beta\mathbf f_{\gamma,1,m}|
 &\lesssim k^{C_\beta}2^{0.01\delta k}
       \sqrt\omega\,d_{\mathrm{edge}}^{-C_\beta}
       \mathcal G_\gamma(\tau_\gamma)\\
 &\le k^{C_\beta}2^{0.3\delta k}
       \sqrt\omega\,d_{\mathrm{edge}}^{-C_\beta}
       \mathcal G_\gamma(\tau_\gamma),
\end{aligned}
\]
which verifies \eqref{iter:source-bound}, including its local envelope,
at the starting stage.

\emph{The forced amplitude and the velocity correction.}
Use the native-scale quantities \(\mathbf n_\gamma\),
\(\mathsf K_{B,\ell}\) and \(d_\ell\) from
\eqref{init:primary-coefficients} and
\eqref{init:primary-principal-cancellation}. For
\(m\in\mathcal H_j\), solve
\begin{equation}
\begin{aligned}
(\partial_{\tau_\gamma}+\mathsf K_{B,\ell}+m^2d_\ell)
                      \mathbf h_{\gamma,j,m}
 +\mathrm imN_\ell\mathbf n_\gamma\pi_{\gamma,j,m}
 &=-2^{-\delta\ell}\mathbf f_{\gamma,j,m},\\
\mathbf n_\gamma\cdot\mathbf h_{\gamma,j,m}&=0,
\qquad \mathbf h_{\gamma,j,m}(0)=0.
\end{aligned}
\label{iter:forced-equation}
\end{equation}
We solve this equation on each connected component of
the full preimage \(E^{-(i(\ell)-i_0)}(\mathcal R_\gamma)\)
in the common auxiliary torus \(Y_{i_0}\). On each component, we use
the label's local coordinates \(Y_{i(\ell)}\) to apply
Proposition~\ref{s7n:forced-inverse}. Representations of the same source
are compatible on overlaps, so uniqueness gives compatible solutions.
The factor \(2^{-\delta\ell}\) converts force units to the
fast-time units of Section~\ref{s7n:forced-pulse}:
\[
2^{(1+\delta)\ell}2^{(2+\delta)\ell}
                  2^{-\delta\ell}=2^{(3+\delta)\ell}.
\]
In particular, \(\mathbf h_{\gamma,j,m}\) uses the same velocity
normalization as \(\mathbf a_\gamma\) in Section~\ref{init:primary-pulse}.
The zero-data estimate \eqref{s7n:forced-duhamel} gives
\[
\begin{aligned}
|\mathbf h_{\gamma,j,m}(\tau)|
 &\lesssim\int_0^\tau
       \frac{\mathcal G_\gamma(\tau)}{\mathcal G_\gamma(s)}
       2^{-\delta\ell}|\mathbf f_{\gamma,j,m}(s)|\,ds\\
 &\lesssim k^{C_j+2}2^{-(0.6+0.1j)\delta k}
       \sqrt\omega\,d_{\mathrm{edge}}^{-C_j}\mathcal G_\gamma(\tau),
\end{aligned}
\]
because \(\tau_{\max}(\ell)\asymp k^2\) and the radial weights are
constant in the fast-time integration. Differentiating the
projected equation as in Proposition~\ref{s7n:forced-inverse} and
using its pressure formula \eqref{s7n:forced-projection} yields
\begin{equation}
\begin{aligned}
|\partial_{\hat R,\hat Z,\hat T,Y_{i_0,1},Y_{i_0,2}}^\beta\mathbf h_{\gamma,j,m}|
 &\lesssim k^{C_{j,\beta}}2^{-(0.6+0.1j)\delta k}
       \sqrt\omega\,d_{\mathrm{edge}}^{-C_{j,\beta}}\mathcal G_\gamma(\tau_\gamma),\\
|\partial_{\hat R,\hat Z,\hat T,Y_{i_0,1},Y_{i_0,2}}^\beta\pi_{\gamma,j,m}|
 &\lesssim k^{C_{j,\beta}}2^{-(1.1+0.1j)\delta k}
       \sqrt\omega\,d_{\mathrm{edge}}^{-C_{j,\beta}}\mathcal G_\gamma(\tau_\gamma).
\end{aligned}
\label{iter:forced-amplitude-bounds}
\end{equation}
The inverse preserves the normalized and transverse supports and
smooth zero extension in those variables and in the radial variable,
as in Proposition~\ref{s7n:forced-inverse}.

Set \(\psi_\gamma=\Phi_t(\tau_\gamma/\tau_{\max}(\ell))\), and define
\begin{equation}
\begin{aligned}
\mathbf a_{\gamma,j,m}&=\psi_\gamma\mathbf h_{\gamma,j,m},&
\varpi_{\gamma,j,m}&=\psi_\gamma\pi_{\gamma,j,m},\\
\mathbf b_{\gamma,j,m}
 &=\frac{\mathrm i2^{-\ell}}{mN_\ell}\nabla\times
       \left(\frac{\mathbf n_\gamma\times\mathbf a_{\gamma,j,m}}
                         {|\mathbf n_\gamma|^2}\right),&
\mathbf c_{\gamma,j,m}&=\mathbf a_{\gamma,j,m}+\mathbf b_{\gamma,j,m}.
\end{aligned}
\label{iter:forced-coefficients}
\end{equation}
The curl here acts on the coefficient, including the cylindrical
frame. The actual fields are
\begin{equation}
\begin{aligned}
w_{j,m\ne0}
 &=\sum_\gamma\sum_{m\in\mathcal H_j}\nabla\times\left[
   \frac{\mathrm i2^{\delta\ell}}{mN_\ell}
   \frac{\mathbf n_\gamma\times\mathbf a_{\gamma,j,m}}
                         {|\mathbf n_\gamma|^2}
                   e^{\mathrm imN_\ell\varphi_\gamma}\right]\\
 &=\sum_\gamma2^{(1+\delta)\ell}
            \sum_{m\in\mathcal H_j}\mathbf c_{\gamma,j,m}
                   e^{\mathrm imN_\ell\varphi_\gamma},\\
p_{j,m\ne0}
 &=\sum_\gamma2^{2(1+\delta)\ell}
            \sum_{m\in\mathcal H_j}\varpi_{\gamma,j,m}
                   e^{\mathrm imN_\ell\varphi_\gamma}.
\end{aligned}
\label{iter:forced-fields}
\end{equation}
The second equality uses \(\mathbf n_\gamma\cdot\mathbf a_{\gamma,j,m}=0\).
The curl estimate in the proof of
\eqref{init:primary-coefficient-bounds} now gives
\begin{equation}
\begin{aligned}
|\partial_{\hat R,\hat Z,\hat T,Y_{i_0,1},Y_{i_0,2}}^\beta\mathbf b_{\gamma,j,m}|
 &\lesssim k^{C_{j,\beta}}N_\ell^{-1}2^{10^{-5}\delta k}
               2^{-(0.6+0.1j)\delta k}
               \sqrt\omega\,d_{\mathrm{edge}}^{-C_{j,\beta}}\\
 &\lesssim k^{C_{j,\beta}}2^{-(1.1+0.1j-10^{-5})\delta k}
               \sqrt\omega\,d_{\mathrm{edge}}^{-C_{j,\beta}}.
\end{aligned}
\label{iter:forced-curl-bound}
\end{equation}
The same calculation retains \(\mathcal G_\gamma(\tau_\gamma)\)
on both right-hand sides. Thus the coefficients of the new pulse
retain the local envelope as well as the radial weight.
Conjugate symmetry in Proposition~\ref{s7n:forced-inverse} makes the
fields real. They are smooth, supported in \(R\in[R_a,R_h]\), and
\[
\nabla\cdot w_{j,m\ne0}=0,\qquad
\langle w_{j,m\ne0}\rangle_\theta
 =\langle p_{j,m\ne0}\rangle_\theta=0.
\]
Their labelled supports remain separated: the inverse preserves the
transverse supports, and \(\psi_\gamma\) restricts fast time to the
same assigned interval as for the primary pulse.

\begin{prop}[The force after the nonzero-mode correction]
\label{iter:nonzero-estimates}
For the fields in \eqref{iter:forced-fields}, one can choose
\(E_{\mathrm{new}}\) in \eqref{iter:nonzero-updated-state} so that
\begin{equation}
\begin{aligned}
2^{-(1+\delta)k}|w_{j,m\ne0}|
 &\lesssim k^{C_j}2^{-(0.6+0.1j)\delta k}
                    \sqrt\omega\,d_{\mathrm{edge}}^{-C_j},\\
2^{-2(1+\delta)k}|p_{j,m\ne0}|
 &\lesssim k^{C_j}2^{-(1.1+0.1j)\delta k}
                    \sqrt\omega\,d_{\mathrm{edge}}^{-C_j},\\
2^{-(3+\delta)k}|f_{\mathrm{new}}
                  -\langle f_{\mathrm{new}}\rangle_\theta|
 &\lesssim k^{C_j}2^{-0.1j\delta k}
                    \sqrt\omega\,d_{\mathrm{edge}}^{-C_j}.
\end{aligned}
\label{iter:nonzero-output}
\end{equation}
The new force is supported in \(R\in[R_a,R_h]\).
The two linear velocity moments are unchanged, and
\[
\langle E_{\mathrm{new}}-E^{[j-1]}\rangle_\theta=0,\qquad
|D_{t,x}^\beta(E_{\mathrm{new}}-E^{[j-1]})|
 \le C_{j,\beta,N}\lambda^N\quad(N>0).
\]
Moreover, the change of the angular mean obeys
\begin{equation}
2^{-(3+\delta)k}
 |\langle f_{\mathrm{new}}-f^{[j-1]}\rangle_\theta|
 \lesssim k^{C_j}2^{-(0.1+0.1j-10^{-5})\delta k}
                    \omega d_{\mathrm{edge}}^{-C_j}.
\label{iter:nonzero-mean-change}
\end{equation}
For every fixed \(\beta\), the estimates in \eqref{iter:nonzero-output}
also hold after applying
\(\partial_{\hat R,\hat Z,\hat T,Y_{i_0,1},Y_{i_0,2}}^{\beta}\)
to the labelled harmonic coefficients, with \(C_j\) replaced by
\(C_{j,\beta}\). In \eqref{iter:nonzero-mean-change}, the same operator acts directly
on the angular mean \(\langle f_{\mathrm{new}}-f^{[j-1]}\rangle_\theta\),
again with \(C_{j,\beta}\) in place of \(C_j\).
\end{prop}

\begin{proof}
The velocity and pressure estimates follow from
\eqref{iter:forced-amplitude-bounds}--\eqref{iter:forced-curl-bound}
and the finite overlap of the native bands. We prove the force
estimates by separating the cutoff terms, the remaining linear
terms, and the velocity interactions.

\emph{The cutoff terms.}
Multiplying \eqref{iter:forced-equation} by \(\psi_\gamma\) gives
\[
\begin{aligned}
&(\partial_{\tau_\gamma}+\mathsf K_{B,\ell}+m^2d_\ell)
                 \mathbf a_{\gamma,j,m}
  +\mathrm imN_\ell\mathbf n_\gamma\varpi_{\gamma,j,m}
  +2^{-\delta\ell}\mathbf f_{\gamma,j,m}\\
&\qquad=2^{-\delta\ell}(1-\psi_\gamma)\mathbf f_{\gamma,j,m}
        +(\partial_{\tau_\gamma}\psi_\gamma)\mathbf h_{\gamma,j,m}.
\end{aligned}
\]
The source extends smoothly by zero at the fast-time endpoints,
and \(\psi_\gamma\in C_c^\infty((0,\tau_{\max}(\ell)))\), so this
identity remains valid after zero extension to the common auxiliary
torus.

Accordingly, set
\begin{equation}
\begin{aligned}
E_{\mathrm{new}}-E^{[j-1]}
 =\sum_\gamma2^{(3+\delta)\ell}\sum_{m\in\mathcal H_j}
 \Bigl[&(1-\psi_\gamma)\mathbf f_{\gamma,j,m}\\
 &+2^{\delta\ell}(\partial_{\tau_\gamma}\psi_\gamma)
                         \mathbf h_{\gamma,j,m}\Bigr]
                         e^{\mathrm imN_\ell\varphi_\gamma}.
\end{aligned}
\label{iter:nonzero-flat-error}
\end{equation}
On the support of either term in brackets,
\(\mathcal G_\gamma\le e^{-ck^2}\). Thus
\eqref{iter:source-bound}, \eqref{iter:forced-amplitude-bounds}
and the physical chain rule give
\[
|D_{t,x}^\beta(E_{\mathrm{new}}-E^{[j-1]})|
 \lesssim k^{C_{j,\beta}}2^{C_{j,\beta}k}e^{-ck^2}
 \le C_{j,\beta,N}\lambda^N.
\]
All harmonics are nonzero, so this change has zero angular mean.

\emph{The remaining linear force.}
Repeating the product-rule calculation leading to
\eqref{init:primary-linear-coefficients}, now with
\(\mathbf c_{\gamma,j,m},\varpi_{\gamma,j,m}\) and the curl
coefficient \(\mathbf b_{\gamma,j,m}\), gives the same groups of
remainder terms. The principal terms
are cancelled by \eqref{iter:forced-equation}, apart from
\eqref{iter:nonzero-flat-error}. For example, the radial viscosity
cross term is
\[
-2^{(1+\delta)\ell+1}
       (D_r\mathbf c_{\gamma,j,m})
           D_re^{\mathrm imN_\ell\varphi_\gamma}.
\]
Since \(D_re^{\mathrm imN_\ell\varphi_\gamma}
=\mathrm imN_\ell2^\ell n_\gamma^r
e^{\mathrm imN_\ell\varphi_\gamma}\), its normalized size is at most
\[
\begin{aligned}
&k^{C_j}N_\ell2^{10^{-5}\delta k}2^{-(0.6+0.1j)\delta k}
               \sqrt\omega\,d_{\mathrm{edge}}^{-C_j}\\
&\qquad\lesssim k^{C_j}2^{-(0.1+0.1j-10^{-5})\delta k}
               \sqrt\omega\,d_{\mathrm{edge}}^{-C_j}.
\end{aligned}
\]
The other terms in that expansion have the following bounds in
force units, apart from \(k^{C_j}\sqrt\omega\,d_{\mathrm{edge}}^{-C_j}\):
\[
\renewcommand{\arraystretch}{1.2}
\begin{array}{l|c}
\text{term}&\text{upper bound}\\ \hline
2^{\delta\ell}(\partial_{\tau_\gamma}+\mathsf K_{B,\ell}+m^2d_\ell)
                  \mathbf b_{\gamma,j,m}
 &2^{-(0.1+0.1j-10^{-5})\delta k}\\
\text{slow transport and lower-order background derivatives}
 &2^{-(0.6+0.1j-10^{-5})\delta k}\\
\text{phase transport defect}
 &2^{-(0.1+0.1j)\delta k}\\
\text{derivatives of pressure coefficients}
 &2^{-(0.1+0.1j-10^{-5})\delta k}\\
\text{viscosity acting twice on coefficients}
 &2^{-(0.6+0.1j-2\cdot10^{-5})\delta k}\\
\text{viscosity cross terms and frame terms}
 &2^{-(0.1+0.1j-10^{-5})\delta k}.
\end{array}
\]
For the pressure row, for instance, the bound is
\(k^{C_j}2^{\delta k}2^{10^{-5}\delta k}2^{-(1.1+0.1j)\delta k}
\sqrt\omega\,d_{\mathrm{edge}}^{-C_j}\).
The same coefficient differentiation as in
Section~\ref{init:primary-pulse} proves these bounds. In particular,
allowing the margin used there,
\begin{equation}
\begin{aligned}
2^{-(3+\delta)k}\bigl|&\mathsf L_B(w_{j,m\ne0},p_{j,m\ne0})
              +f^{[j-1]}-\langle f^{[j-1]}\rangle_\theta\\
&\hspace{18mm}-(E_{\mathrm{new}}-E^{[j-1]})\bigr|
 \lesssim k^{C_j}2^{-(0.1+0.1j-3\cdot10^{-5})\delta k}
                     \sqrt\omega\,d_{\mathrm{edge}}^{-C_j}.
\end{aligned}
\label{iter:nonzero-linear-bound}
\end{equation}

\emph{Interactions with the axisymmetric correction.}
For the rest of the proof write \(w=w_{j,m\ne0}\).
The exact update is
\begin{equation}
\begin{aligned}
f_{\mathrm{new}}
={}&f^{[j-1]}+\mathsf L_B(w,p_{j,m\ne0})
                         -(E_{\mathrm{new}}-E^{[j-1]})\\
 &+\sum_{q=0}^{j-1}\sum_{a=\mathrm{aux},\mathrm{moment}}
       \bigl((u_{q,a}\cdot\nabla)w+(w\cdot\nabla)u_{q,a}\bigr)\\
 &+\sum_{w'}\bigl((w'\cdot\nabla)w+(w\cdot\nabla)w'\bigr)
       +(w\cdot\nabla)w.
\end{aligned}
\label{iter:nonzero-force-update}
\end{equation}
The carrier gradient and \eqref{iter:forced-amplitude-bounds} give
\[
2^{-(2+\delta)k}|\nabla w|
 \lesssim k^{C_j}2^{\delta k/2}2^{-(0.6+0.1j)\delta k}
                    \sqrt\omega\,d_{\mathrm{edge}}^{-C_j}.
\]
For any old axisymmetric increment \(v=u_{q,a}\), we may weaken
the first estimate in \eqref{iter:old-axisymmetric-bounds} to
\(2^{-(1+\delta)k}|v|\lesssim k^{C_j}2^{-0.9\delta k}
\omega d_{\mathrm{edge}}^{-C_j}\). Using
\(\omega^{3/2}\le\sqrt\omega\), we obtain
\begin{equation}
\begin{aligned}
2^{-(3+\delta)k}|(v\cdot\nabla)w|
 &\lesssim k^{C_j}2^{\delta k}2^{-0.9\delta k}
                    2^{\delta k/2}2^{-(0.6+0.1j)\delta k}
                    \sqrt\omega\,d_{\mathrm{edge}}^{-C_j}\\
 &=k^{C_j}2^{-0.1j\delta k}\sqrt\omega\,d_{\mathrm{edge}}^{-C_j},\\
2^{-(3+\delta)k}|(w\cdot\nabla)v|
 &\lesssim k^{C_j}2^{\delta k}2^{-(0.6+0.1j)\delta k}
                  2^{-(1-2\cdot10^{-5})\delta k}
                    \sqrt\omega\,d_{\mathrm{edge}}^{-C_j}\\
 &=k^{C_j}2^{-(0.6+0.1j-2\cdot10^{-5})\delta k}
                    \sqrt\omega\,d_{\mathrm{edge}}^{-C_j}.
\end{aligned}
\label{iter:nonzero-axisymmetric-interactions}
\end{equation}
Finite summation gives the same bounds for the axisymmetric
interactions in \eqref{iter:nonzero-force-update}.

\emph{Interactions between pulses.}
For a fixed label, let \(\mathbf c e^{\mathrm imN_\ell\varphi_\gamma}\)
be a full divergence-free harmonic of an old pulse \(w'\) or of \(w\), with its
constant physical prefactor omitted. Then
\[
\mathrm imN_\ell2^\ell\mathbf n_\gamma\cdot\mathbf c
 =-\nabla\cdot\mathbf c
 =-(D_r+r^{-1})c^r-\partial_zc^z.
\]
For a second coefficient \(\mathbf d\) at the same label,
the product rule therefore gives the exact identity
\begin{equation}
\begin{aligned}
&(\mathbf c e^{\mathrm imN_\ell\varphi_\gamma}\cdot\nabla)
       (\mathbf d e^{\mathrm inN_\ell\varphi_\gamma})\\
&\qquad=e^{\mathrm i(m+n)N_\ell\varphi_\gamma}
 \left[(\mathbf c\cdot\nabla)\mathbf d
                   -\frac nm(\nabla\cdot\mathbf c)\mathbf d\right].
\end{aligned}
\label{iter:nonzero-carrier-cancellation}
\end{equation}
The derivatives on the right act on coefficients and on the
cylindrical frame. Here \(m,n\in\mathcal H_j\), so
\(|n/m|\le 2^j\). This factor is absorbed in the constant
depending on \(j\). Thus a spatial derivative costs at most
\(k^{C_j}2^k2^{10^{-5}\delta k}\), with no factor \(N_\ell\),
and a possible increase in the power of \(d_{\mathrm{edge}}^{-1}\).
By \eqref{iter:old-pulse-coefficients}, the old and new coefficients have sizes
\(k^{C_j}2^{-0.5\delta k}\) and \(k^{C_j}2^{-(0.6+0.1j)\delta k}\),
respectively, each with the factor
\(\sqrt\omega\,d_{\mathrm{edge}}^{-C_j}\).
Support separation removes all products of different labels. It follows that
\begin{equation}
\begin{aligned}
2^{-(3+\delta)k}
 \left|\sum_{w'}\bigl((w'\cdot\nabla)w+(w\cdot\nabla)w'\bigr)\right|
 &\lesssim k^{C_j}2^{\delta k}2^{-0.5\delta k}2^{-(0.6+0.1j)\delta k}
                         2^{10^{-5}\delta k}\omega d_{\mathrm{edge}}^{-C_j}\\
 &=k^{C_j}2^{-(0.1+0.1j-10^{-5})\delta k}\omega d_{\mathrm{edge}}^{-C_j},\\
2^{-(3+\delta)k}|(w\cdot\nabla)w|
 &\lesssim k^{C_j}2^{\delta k}2^{-(1.2+0.2j)\delta k}
               2^{10^{-5}\delta k}\omega d_{\mathrm{edge}}^{-C_j}\\
 &=k^{C_j}2^{-(0.2+0.2j-10^{-5})\delta k}\omega d_{\mathrm{edge}}^{-C_j}.
\end{aligned}
\label{iter:nonzero-pulse-interactions}
\end{equation}
Subtract the angular mean in \eqref{iter:nonzero-force-update}.
Combining \eqref{iter:nonzero-linear-bound},
\eqref{iter:nonzero-axisymmetric-interactions} and
\eqref{iter:nonzero-pulse-interactions}, and using
\(0.1+0.1j-3\cdot10^{-5}>0.1j\), proves
\[
2^{-(3+\delta)k}|f_{\mathrm{new}}
                     -\langle f_{\mathrm{new}}\rangle_\theta|
 \lesssim k^{C_j}2^{-0.1j\delta k}\sqrt\omega\,d_{\mathrm{edge}}^{-C_j}
 \le k^{C_j}2^{(0.3-0.1j)\delta k}\sqrt\omega\,d_{\mathrm{edge}}^{-C_j}.
\]
The new pulse \(w\) uses only \(m\in\mathcal H_j\), while the products
in \eqref{iter:nonzero-force-update} have nonzero harmonics
\(0<|m|\le 2^{j+1}\), still a finite set depending only on \(j\).
All the preceding residual estimates hold after differentiation
to any fixed order in these variables. In their local form, the linear terms retain
\(\mathcal G_\gamma\), and every nonlinear term contains a pulse
coefficient with that envelope; products of two envelopes are
bounded by \(C\mathcal G_\gamma\). Thus each new nonzero-mode
coefficient, in force units, has the bound
\[
k^{C_{j,\beta}}2^{-0.1j\delta k}
 \sqrt\omega\,d_{\mathrm{edge}}^{-C_{j,\beta}}
 \mathcal G_\gamma(\tau_\gamma)
\]
after \(\partial_{\hat R,\hat Z,\hat T,Y_{i_0,1},Y_{i_0,2}}^{\beta}\) derivatives. This proves the
coefficient assertion and retains the local source envelope with
the new stage gain. The assigned supports and smooth zero extensions
are preserved as well.

\emph{The angular mean and the preserved moments.}
The linear term, the flat-error change, and all interactions with
the axisymmetric increments in \eqref{iter:nonzero-force-update}
have zero angular mean.
Define, for \(a,b\in\{r,\theta,z\}\),
\[
\mathsf Q^{ab}
 =\sum_{w'}\langle (w')^aw^b+w^a(w')^b\rangle_\theta
       +\langle w^aw^b\rangle_\theta.
\]
The cylindrical divergence formulas used in
\eqref{init:primary-angular-force} give
\begin{equation}
\begin{aligned}
\langle f_{\mathrm{new}}^r-(f^{[j-1]})^r\rangle_\theta
 &=(D_r+r^{-1})\mathsf Q^{rr}+\partial_z\mathsf Q^{zr}
                                      -r^{-1}\mathsf Q^{\theta\theta},\\
\langle f_{\mathrm{new}}^\theta-(f^{[j-1]})^\theta\rangle_\theta
 &=(D_r+2/r)\mathsf Q^{r\theta}+\partial_z\mathsf Q^{z\theta},\\
\langle f_{\mathrm{new}}^z-(f^{[j-1]})^z\rangle_\theta
 &=(D_r+r^{-1})\mathsf Q^{rz}+\partial_z\mathsf Q^{zz}.
\end{aligned}
\label{iter:nonzero-mean-identities}
\end{equation}
Only opposite harmonics at the same label survive averaging.
Their carrier factors multiply to one, and
\[
2^{-2(1+\delta)k}|\partial_{\hat R,\hat Z,\hat T,Y_{i_0,1},Y_{i_0,2}}^\beta\mathsf Q^{ab}|
 \lesssim k^{C_{j,\beta}}
          (2^{-(1.1+0.1j)\delta k}+2^{-(1.2+0.2j)\delta k})
                    \omega d_{\mathrm{edge}}^{-C_{j,\beta}}.
\]
One radial derivative of a coefficient costs at most
\(k^{C_j}2^k2^{10^{-5}\delta k}\), while an axial derivative costs
\(k^{C_j}2^{(1-\delta)k}\). Substitution in
\eqref{iter:nonzero-mean-identities} proves
\eqref{iter:nonzero-mean-change}. In particular, the incoming
averaged bounds and this estimate imply
\[
\begin{aligned}
2^{-(3+\delta)k}\Bigl(
 &|\langle f_{\mathrm{new}}^{\mathrm{tan}}\rangle_\theta
       -\langle f_{\mathrm{new}}^{\mathrm{tan}}\rangle_{\theta,Y_{i_0}}|\\
 &+|\langle f_{\mathrm{new}}\rangle_{\theta,Y_{i_0}}|\Bigr)
 \lesssim k^{C_j}2^{-(0.1+0.1j-10^{-5})\delta k}
                    \omega d_{\mathrm{edge}}^{-C_j}.
\end{aligned}
\]
Finally, \(\langle w\rangle_\theta=0\) preserves both integrals in
\eqref{init:linear-moments}. The new pulse, its pressure and the new
force are smooth, with radial support in \(R\in[R_a,R_h]\).
The flat-error change has zero angular mean. This proves the
remaining assertions.
\end{proof}

The nonzero angular force has already reached its stage-\(j\)
target in \eqref{iter:targets}; the averaged terms above still
require the following corrections to reach the second target.

\subsection{The stress correction}
\label{iter:stress-correction}

Write \((u,p,E,f)\) for the state obtained in
Section~\ref{iter:nonzero-modes}. We next cancel the tangential
double average by matching its radial primitive against the fixed
primary pulse. The primitive must have compact radial support; its
two compensating terms are controlled by the moment defects at the
preceding stage.

\emph{The stress to be matched.}
Use the cutoff \(\chi_m\) from
\eqref{init:auxiliary-radial-primitive}, and put
\[
c_\theta=2,\qquad c_z=1,\qquad
\rho_a=r^{-c_a}\partial_r\chi_m\quad(a=\theta,z).
\]
These profiles are supported in a fixed interior subannulus and
\(\int_0^\infty r^{c_a}\rho_a\,dr=1\). Define the physical stress pair
\begin{equation}
\widetilde{\mathsf T}_j^a
 =-r^{-c_a}\left[
   \int_0^r s^{c_a}\langle f^a\rangle_{\theta,Y_{i_0}}(t,s,z)\,ds
   -\chi_m(t,r,z)\int_0^\infty
        s^{c_a}\langle f^a\rangle_{\theta,Y_{i_0}}(t,s,z)\,ds\right].
\label{iter:stress-primitive}
\end{equation}
It is independent of the auxiliary variable and satisfies
\begin{equation}
\left(\partial_r+\frac{c_a}{r}\right)\widetilde{\mathsf T}_j^a
 =-\langle f^a\rangle_{\theta,Y_{i_0}}
   +\rho_a\int_0^\infty s^{c_a}
                    \langle f^a\rangle_{\theta,Y_{i_0}}\,ds.
\label{iter:stress-primitive-identity}
\end{equation}
Below the annulus the integral from \(0\) to \(r\) vanishes;
above it the two integrals agree. Near the outer edge one uses the corresponding
tail integral. Thus the endpoint estimate in the proof of
\eqref{init:auxiliary-fixed-shift} applies without an auxiliary
translation. The averaged-force bound at the end of
Section~\ref{iter:nonzero-modes} gives, for every fixed \(\beta\),
\begin{equation}
\left|\partial_{\hat R,\hat Z,\hat T}^{\beta}
       \bigl(2^{-(2+\delta)k}\widetilde{\mathsf T}_j\bigr)\right|
 \lesssim k^{C_{j,\beta}}
       2^{-(0.1+0.1j-10^{-5})\delta k}
       \omega d_{\mathrm{edge}}^{-C_{j,\beta}}.
\label{iter:stress-input}
\end{equation}
In particular, the stress has smooth zero extension with the same
radial weight.

To estimate the compensating terms, let \(\mathsf Q^{ab}\) be the
angular-mean tensor in \eqref{iter:nonzero-mean-identities}.
The linear moments are preserved and the pulse pressure has zero
angular mean. In particular, the definition of the flux defect gives
\[
\begin{aligned}
\mathcal J^{z\theta}-\mathcal J^{z\theta,[j-1]}
 ={}&\int_0^\infty r\langle\mathsf Q^{zz}\rangle_{Y_{i_0}}\,dr\\
 &+\frac12\int_0^\infty r^2
       \bigl\langle f^r-(f^{[j-1]})^r\bigr\rangle_{\theta,Y_{i_0}}\,dr.
\end{aligned}
\]
In \eqref{init:weighted-force-moments}, the added radial-force
integral therefore cancels, leaving the radial integral from stage
\(j-1\). Applying that identity before and after the nonzero-mode
correction gives
\begin{equation}
\begin{aligned}
\int_0^\infty r^2\langle f^\theta\rangle_{\theta,Y_{i_0}}\,dr
 &=\partial_z\left[\mathcal J^{\theta z,[j-1]}
        +\int_0^\infty r^2\langle\mathsf Q^{\theta z}\rangle_{Y_{i_0}}\,dr
                   \right],\\
\int_0^\infty r\langle f^z\rangle_{\theta,Y_{i_0}}\,dr
 &=\partial_z\left[\mathcal J^{z\theta,[j-1]}
       -\frac{\mathcal P^{[j-1]}}2\int_0^\infty r^2\rho_r\,dr
       +\int_0^\infty r\langle\mathsf Q^{zz}\rangle_{Y_{i_0}}\,dr
                   \right].
\end{aligned}
\label{iter:stress-bump-moments}
\end{equation}
The radial term in the second line uses the exact identity at the
end of stage \(j-1\). In the respective moment units, the two
brackets are bounded by \(k^{C_j}2^{-(1.1+0.1j)\delta k}\).
Indeed, this follows from \eqref{iter:targets} at stage \(j-1\)
and the bound for \(\mathsf Q\) following
\eqref{iter:nonzero-mean-identities}. Since
\(\partial_z\) costs \(k^{C_j}2^{(1-\delta)k}\), we obtain
\begin{equation}
2^{-(3+\delta)k}\left|
 \rho_a\int_0^\infty r^{c_a}
                  \langle f^a\rangle_{\theta,Y_{i_0}}\,dr\right|
 \lesssim k^{C_j}2^{-(1.1+0.1j)\delta k}
                 \omega d_{\mathrm{edge}}^{-C_j}.
\label{iter:stress-bump-bound}
\end{equation}
The same statements hold after fixed slow derivatives. The radial
factor in this estimate is available because \(\rho_a\) has fixed
interior support.

\emph{The pulse.}
Apply \eqref{eq:assembled-signed-inverse} to
\(\widetilde{\mathsf T}_j\). For a native label
\(\gamma=(\ell,\alpha,\iota)\), denote the resulting real matching
coefficient by \(\delta g_{\gamma,j}\); thus
\eqref{eq:signed-covariance-inverse} is used with
\(\Sigma=2^{-(2+\delta)\ell}\widetilde{\mathsf T}_j\).
Set
\begin{equation}
\begin{aligned}
\mathbf a_{\gamma,j}^{\mathrm s}
 &=2^{-\delta\ell/2}\chi_\gamma
                  \delta g_{\gamma,j}\mathbf h_\gamma,&
\pi_{\gamma,j}^{\mathrm s}
 &=2^{-\delta\ell/2}\chi_\gamma\delta g_{\gamma,j}\pi_\gamma,\\
\mathbf b_{\gamma,j}^{\mathrm s}
 &=-\frac{2^{-\ell}}{N_\ell}\nabla\times
       \left(\frac{\mathbf n_\gamma\times\mathbf a_{\gamma,j}^{\mathrm s}}
                          {|\mathbf n_\gamma|^2}\right).
\end{aligned}
\label{iter:stress-coefficients}
\end{equation}
Define \(w_{j,\mathrm{stress}},p_{j,\mathrm{stress}}\) by
\eqref{init:primary-fields}, replacing
\((\mathbf a_\gamma,\mathbf b_\gamma,\pi_{0,\gamma})\) with
\((\mathbf a_{\gamma,j}^{\mathrm s},
  \mathbf b_{\gamma,j}^{\mathrm s},\pi_{\gamma,j}^{\mathrm s})\).
In particular, the velocity is the full curl in that formula.
Proposition~\ref{prop:pulse-signed-covariance} and the curl estimate
from Section~\ref{init:primary-pulse} give
\begin{equation}
\begin{aligned}
|\partial_{\hat R,\hat Z,\hat T,Y_{i_0,1},Y_{i_0,2}}^{\beta}\mathbf a_{\gamma,j}^{\mathrm s}|
 &\lesssim k^{C_{j,\beta}}2^{-(0.6+0.1j-10^{-5})\delta k}
                   \sqrt\omega\,d_{\mathrm{edge}}^{-C_{j,\beta}} \mathcal G_\gamma(\tau_\gamma),\\
|\partial_{\hat R,\hat Z,\hat T,Y_{i_0,1},Y_{i_0,2}}^{\beta}\mathbf b_{\gamma,j}^{\mathrm s}|
 &\lesssim k^{C_{j,\beta}}2^{-(1.1+0.1j-2\cdot10^{-5})\delta k}
                   \sqrt\omega\,d_{\mathrm{edge}}^{-C_{j,\beta}} \mathcal G_\gamma(\tau_\gamma),\\
|\partial_{\hat R,\hat Z,\hat T,Y_{i_0,1},Y_{i_0,2}}^{\beta}\pi_{\gamma,j}^{\mathrm s}|
 &\lesssim k^{C_{j,\beta}}2^{-(1.1+0.1j-10^{-5})\delta k}
                   \sqrt\omega\,d_{\mathrm{edge}}^{-C_{j,\beta}} \mathcal G_\gamma(\tau_\gamma).
\end{aligned}
\label{iter:stress-coefficient-bounds}
\end{equation}
The quotient by the fixed
\(g_{\mathrm{match},\gamma}\) has smooth zero extension by
Proposition~\ref{prop:pulse-signed-covariance}; no new cone condition on the stress
is needed.

Put \(w=w_{j,\mathrm{stress}}\) and set
\begin{equation}
\begin{aligned}
u_{\mathrm{new}}&=u+w,&p_{\mathrm{new}}&=p+p_{j,\mathrm{stress}},\\
E_{\mathrm{new}}&=E+\sum_\gamma
  2^{(3+2\delta)\ell-\delta\ell/2}
  (\partial_{\tau_\gamma}\chi_\gamma)\delta g_{\gamma,j}
                  \mathbf h_\gamma\cos(N_\ell\varphi_\gamma),&
f_{\mathrm{new}}&=\mathsf{NS}(u_{\mathrm{new}},p_{\mathrm{new}})
                          -E_{\mathrm{new}}.
\end{aligned}
\label{iter:stress-updated-state}
\end{equation}
As in \eqref{init:primary-flat-error}, the added error is flat with
all physical derivatives and has zero angular mean.

\begin{prop}[The force after the stress correction]
\label{iter:stress-estimates}
The new pulse satisfies the bound in
\eqref{init:later-increment-sizes}, and
\begin{equation}
\begin{aligned}
2^{-2(1+\delta)k}|p_{j,\mathrm{stress}}|
 &\lesssim k^{C_j}2^{-(1.1+0.1j-10^{-5})\delta k}
                      \sqrt\omega\,d_{\mathrm{edge}}^{-C_j},\\
2^{-(3+\delta)k}|f_{\mathrm{new}}-\langle f_{\mathrm{new}}\rangle_\theta|
 &\lesssim k^{C_j}2^{-0.1j\delta k}
                      \sqrt\omega\,d_{\mathrm{edge}}^{-C_j},\\
2^{-(3+\delta)k}
       |\langle f_{\mathrm{new}}^{\mathrm{tan}}\rangle_{\theta,Y_{i_0}}|
 &\lesssim k^{C_j}2^{-(0.2+0.1j)\delta k}
                      \omega d_{\mathrm{edge}}^{-C_j},\\
2^{-(3+\delta)k}\Bigl(
 |\langle f_{\mathrm{new}}^r\rangle_\theta|
 +|\langle f_{\mathrm{new}}^{\mathrm{tan}}\rangle_\theta
       -\langle f_{\mathrm{new}}^{\mathrm{tan}}\rangle_{\theta,Y_{i_0}}|\Bigr)
 &\lesssim k^{C_j}2^{-(0.1+0.1j-2\cdot10^{-5})\delta k}
                      \omega d_{\mathrm{edge}}^{-C_j}.
\end{aligned}
\label{iter:stress-output}
\end{equation}
These estimates hold at every fixed derivative order in the variables
indicated just after \eqref{iter:targets}.
The support, incompressibility, linear moments and flat-error
properties in the induction are preserved.
\end{prop}

\begin{proof}
The velocity and pressure bounds follow from
\eqref{iter:stress-coefficient-bounds}. For the nonzero angular
modes, use the expansion \eqref{init:primary-linear-coefficients},
the old axisymmetric bounds \eqref{iter:old-axisymmetric-bounds},
and the carrier cancellation \eqref{iter:nonzero-carrier-cancellation}.
The matching coefficient \(\delta g_{\gamma,j}\) depends only on
the slow variables, so \(\partial_{\tau_\gamma}\delta g_{\gamma,j}=0\).
Thus the principal cancellation \eqref{init:primary-principal-cancellation}
holds with the current matching coefficient and pulse amplitudes.
The proof of Proposition~\ref{prop:pulse-signed-covariance} controls every fixed
slow derivative of \(\delta g_{\gamma,j}\). Substituting
\eqref{iter:stress-coefficient-bounds} into the linear remainder
expansion \eqref{init:primary-linear-coefficients} therefore gives
the additional factor \(2^{-(0.1+0.1j-10^{-5})\delta k}\), including
terms in which derivatives fall on the matching coefficient.
When an old axisymmetric increment differentiates the new stress pulse,
the derivative can fall on the carrier. For an interaction with a pulse
already present, including \(w_{j,m\ne0}\), the carrier derivative
cancels by \eqref{iter:nonzero-carrier-cancellation}. The powers of two
in the respective force-unit bounds are
\[
\begin{aligned}
 2^{\delta k}\,2^{-(1-10^{-5})\delta k}\,
 2^{\delta k/2}\,2^{-(0.6+0.1j-10^{-5})\delta k}
 &=2^{-(0.1+0.1j-2\cdot10^{-5})\delta k},\\
 2^{\delta k}\,2^{-\delta k/2}\,
 2^{10^{-5}\delta k}\,2^{-(0.6+0.1j-10^{-5})\delta k}
 &=2^{-(0.1+0.1j-2\cdot10^{-5})\delta k}.
\end{aligned}
\]
In force units, the remaining
linear terms, the interactions with the old axisymmetric
increments, and those with the old pulses are all bounded by
\[
k^{C_j}2^{-(0.1+0.1j-4\cdot10^{-5})\delta k}
                         \sqrt\omega\,d_{\mathrm{edge}}^{-C_j}.
\]
The self-interaction is smaller. Combining this with
\eqref{iter:nonzero-output} proves the second line of
\eqref{iter:stress-output}. The same calculation retains the local
envelope, with this gain, in every nonzero harmonic coefficient.

For the mean force, linear terms and interactions with the old
axisymmetric increments vanish under angular averaging. The mean
change therefore has the divergence form
\eqref{iter:nonzero-mean-identities}, with the tensor formed from
the old pulses and \(w\). Separate the interaction with \(w_0\).
The pre-pulse identity \eqref{eq:assembled-signed-covariance} gives
exactly \(\widetilde{\mathsf T}_j^a\) for its double-mean radial
tangential flux. Since the cosine--sine cross terms have zero
angular average, the difference for the full curls is
\[
\frac12\sum_\gamma2^{2(1+\delta)\ell}
 \left\langle b_\gamma^r(b_{\gamma,j}^{\mathrm s})^a
       +(b_{\gamma,j}^{\mathrm s})^r b_\gamma^a\right\rangle_{Y_{i_0}}.
\]
Its force contribution is bounded, in force units, by
\(k^{C_j}2^{-(1.1+0.1j-3\cdot10^{-5})\delta k}
\omega d_{\mathrm{edge}}^{-C_j}\).
The axial flux of the primary cross term has the same bound:
the axial derivative costs \(2^{(1-\delta)k}\), in place of the
radial cost \(2^k\).

Every other old pulse, including \(w_{j,m\ne0}\), has coefficient
size at most
\[
k^{C_j}2^{-(0.7-10^{-5})\delta k}
\sqrt\omega\,d_{\mathrm{edge}}^{-C_j}.
\]
After double averaging, a radial derivative is a slow derivative.
Their cross terms with \(w\), and the self-interaction of \(w\),
are consequently bounded in force units by
\[
k^{C_j}\left(2^{-(0.3+0.1j-2\cdot10^{-5})\delta k}
             +2^{-(0.2+0.2j-2\cdot10^{-5})\delta k}\right)
                         \omega d_{\mathrm{edge}}^{-C_j}.
\]
For \(j\ge1\) these are smaller than the required
\(2^{-(0.2+0.1j)\delta k}\) bound. The matched term cancels the
incoming tangential double mean by
\eqref{iter:stress-primitive-identity}, leaving only the bumps
in \eqref{iter:stress-bump-bound}. This proves the third line.

Without auxiliary averaging, the primary cross term has normalized
quadratic size
\[
k^{C_j}2^{-(1.1+0.1j-10^{-5})\delta k}
                  \omega d_{\mathrm{edge}}^{-C_j}.
\]
One radial derivative of a coefficient
costs at most \(k^{C_j}2^k2^{10^{-5}\delta k}\). Thus its force
contribution is bounded by the last line of
\eqref{iter:stress-output}; all other new mean terms are smaller.
The incoming radial mean has the same bound, by
\eqref{iter:nonzero-mean-change} and the stage \(j-1\) radial
identity. This proves the remaining force estimate.
The derivative assertions follow from the same coefficient
calculations and the endpoint estimate used in
\eqref{iter:stress-input}. Finally,
\(w_{j,\mathrm{stress}}\) and \(p_{j,\mathrm{stress}}\) have zero
angular mean, so the linear moments are unchanged.
\end{proof}

We record also the size of the defects before the auxiliary and
moment corrections. At the present state,
\begin{equation}
\begin{aligned}
&2^{-2(1+\delta)k}|\partial_{\hat T,\hat Z}^{\beta}\mathcal P|
 +2^{(1-2\delta)k}|\partial_{\hat T,\hat Z}^{\beta}\mathcal J^{\theta z}|
 +2^{-2\delta k}|\partial_{\hat T,\hat Z}^{\beta}\mathcal J^{z\theta}|\\
&\qquad\lesssim k^{C_{j,\beta}}
                 2^{-(1.1+0.1j-2\cdot10^{-5})\delta k}.
\end{aligned}
\label{iter:stress-defects}
\end{equation}
For \(\mathcal P\), integrate the radial bound in
\eqref{iter:stress-output}; normalization supplies an additional
\(2^{-\delta k}\). For the two flux defects, subtract their
stage \(j-1\) values. The pulse pressures have zero angular mean,
and the velocity products gain at least
\(2^{-(1.1+0.1j-10^{-5})\delta k}\) in their respective moment
units. The additional radial-force integral in
\(\mathcal J^{z\theta}\) is bounded as for \(\mathcal P\).
These are the same integration factors as in the proof of
\eqref{init:moment-input-bounds}. Fixed derivatives in
\((\hat T,\hat Z)\) commute with these compactly supported integrals.

\subsection{The auxiliary velocity correction}
\label{iter:auxiliary-correction}

Let \((u,p,E,f)\) now be the state after the stress correction.
We cancel the remaining auxiliary oscillation in the tangential
angular mean, using the construction of
Section~\ref{init:auxiliary-correction}. Define
\begin{equation}
\begin{pmatrix}A_{j,\mathrm{aux}}\\ B_{j,\mathrm{aux}}\end{pmatrix}
 =-\Lambda_t^{-i_0}(v_t\cdot\nabla_{Y_{i_0}})^{-1}
    \bigl(\langle f^{\mathrm{tan}}\rangle_\theta
                 -\langle f^{\mathrm{tan}}\rangle_{\theta,Y_{i_0}}\bigr),
\qquad
\Psi_{j,\mathrm{aux}}=\frac1r\mathcal I_0(rB_{j,\mathrm{aux}}),
\label{iter:auxiliary-definition}
\end{equation}
where \(\mathrm{tan}=(\theta,z)\). Put
\[
v=u_{j,\mathrm{aux}}
 =A_{j,\mathrm{aux}}e_\theta
                  +\nabla\times(\Psi_{j,\mathrm{aux}}e_\theta),
\qquad
\mathcal R_{j,\mathrm{aux}}=v^z-B_{j,\mathrm{aux}}.
\]
The remainder is given by \eqref{init:auxiliary-components}, with
\(B_{0,\mathrm{aux}}\) replaced by \(B_{j,\mathrm{aux}}\).
All these fields are compatible on overlapping reference bands.

The inverse estimate \eqref{init:auxiliary-inverse-bound}, the
fixed-shift formula \eqref{init:auxiliary-fixed-shift}, and
\eqref{iter:stress-output} give, for every fixed \(\beta\),
\begin{equation}
\begin{aligned}
&2^{-(1+\delta)k}
 \bigl(|\partial_{\hat R,\hat Z,\hat T,Y_{i_0,1},Y_{i_0,2}}^{\beta}v^{\mathrm{tan}}|
                 +2^{\delta k}|\partial_{\hat R,\hat Z,\hat T,Y_{i_0,1},Y_{i_0,2}}^{\beta}v^r|\bigr)
 +2^{-\delta k}|\partial_{\hat R,\hat Z,\hat T,Y_{i_0,1},Y_{i_0,2}}^{\beta}\Psi_{j,\mathrm{aux}}|\\
&\qquad\lesssim k^{C_{j,\beta}}
       2^{-(1.1+0.1j-2\cdot10^{-5})\delta k}
                     \omega d_{\mathrm{edge}}^{-C_{j,\beta}}.
\end{aligned}
\label{iter:auxiliary-size}
\end{equation}
Indeed, \(\Lambda_t^{-i_0}\asymp k^2 2^{-(2+\delta)k}\)
supplies one additional \(2^{-\delta k}\) when converting force
units to velocity units. The radial primitive supplies \(2^{-k}\),
and \(v^r=-\partial_z\Psi_{j,\mathrm{aux}}\) then supplies the
extra radial gain. The endpoint argument in
Section~\ref{init:auxiliary-correction} preserves the same
\(\omega\) after differentiation to any fixed order in the indicated variables.
The repeated integration by parts in
\eqref{init:auxiliary-tail-bound} applies to the zero-mean source
here as well. Thus \(\mathcal R_{j,\mathrm{aux}}\) and its fast
time derivative are flat with all physical derivatives.

As in \eqref{init:auxiliary-new-fields}, set
\begin{equation}
\begin{aligned}
u_{\mathrm{new}}&=u+v,&p_{\mathrm{new}}&=p,\\
E_{\mathrm{new}}&=E+\Lambda_t^{i_0}(v_t\cdot\nabla_{Y_{i_0}})
                              \mathcal R_{j,\mathrm{aux}}e_z,&
f_{\mathrm{new}}&=\mathsf{NS}(u_{\mathrm{new}},p)-E_{\mathrm{new}}.
\end{aligned}
\label{iter:auxiliary-updated-state}
\end{equation}
The increment \(v\) is divergence-free, axisymmetric, has zero
auxiliary mean, and has smooth zero extension off the annulus.
The added flat error has zero auxiliary mean. In particular, all
the support, linear-moment and flat-error requirements are preserved.

\begin{prop}[The force after the auxiliary correction]
\label{iter:auxiliary-estimates}
The fields in \eqref{iter:auxiliary-updated-state} satisfy
\begin{equation}
\begin{aligned}
2^{-(3+\delta)k}|f_{\mathrm{new}}-\langle f_{\mathrm{new}}\rangle_\theta|
 &\lesssim k^{C_j}2^{-0.1j\delta k}
                      \sqrt\omega\,d_{\mathrm{edge}}^{-C_j},\\
2^{-(3+\delta)k}\Bigl(
 |\langle f_{\mathrm{new}}^{\mathrm{tan}}\rangle_\theta
          -\langle f_{\mathrm{new}}^{\mathrm{tan}}\rangle_{\theta,Y_{i_0}}|
 +|\langle f_{\mathrm{new}}^{\mathrm{tan}}\rangle_{\theta,Y_{i_0}}|\Bigr)
 &\lesssim k^{C_j}2^{-(0.2+0.1j)\delta k}
                      \omega d_{\mathrm{edge}}^{-C_j},\\
2^{-(3+\delta)k}|\langle f_{\mathrm{new}}^r\rangle_\theta|
 &\lesssim k^{C_j}2^{-(0.1+0.1j-2\cdot10^{-5})\delta k}
                      \omega d_{\mathrm{edge}}^{-C_j}.
\end{aligned}
\label{iter:auxiliary-output}
\end{equation}
The three defects still satisfy \eqref{iter:stress-defects}.
All estimates hold at every fixed derivative order in the variables
indicated just after \eqref{iter:targets}.
\end{prop}

\begin{proof}
In the angular-mean force formulas
\eqref{init:auxiliary-tangential-force} and
\eqref{init:auxiliary-radial-change}, replace the tensor there by
\[
\mathsf A^{ab}
 =u_B^av^b+v^au_B^b+v^av^b
  +\sum_{q=0}^{j-1}\sum_{c=\mathrm{aux},\mathrm{moment}}
       \bigl(u_{q,c}^av^b+v^au_{q,c}^b\bigr).
\]
These are exact identities: all interactions with the pulses have
zero angular mean. The fast tangential time derivative cancels
the incoming auxiliary oscillation by
\eqref{iter:auxiliary-definition}. The remainder in the axial
component is precisely the term assigned to the flat error in
\eqref{iter:auxiliary-updated-state}.

The slow time, viscous and background terms are estimated exactly
as after \eqref{init:auxiliary-tangential-force}. By
\eqref{iter:auxiliary-size}, their sum in force units is bounded by
\[
k^{C_j}2^{-(1.1+0.1j-4\cdot10^{-5})\delta k}
                          \omega d_{\mathrm{edge}}^{-C_j}.
\]
Here two radial derivatives of the coefficients account for the possible
loss \(2^{2\cdot10^{-5}\delta k}\) in the viscous term.
For the old axisymmetric increments, use
\eqref{iter:old-axisymmetric-bounds}. A radial flux contains one
radial component, with its additional \(2^{-\delta k}\), while
an axial flux is differentiated at scale \(2^{(1-\delta)k}\).
Their tangential force contributions, and the self-interaction,
are smaller than the displayed bound. This proves the required
bound for the auxiliary oscillation of the tangential force.

To check the double average, the background and all old moment
corrections have zero cross average with \(v\). The old auxiliary
corrections must be retained. More explicitly, put
\[
\mathsf C^{ab}
 =\sum_{q=0}^{j-1}
       \langle u_{q,\mathrm{aux}}^av^b
                     +v^au_{q,\mathrm{aux}}^b\rangle_{Y_{i_0}}
       +\langle v^av^b\rangle_{Y_{i_0}}.
\]
The exact double-mean changes are
\begin{equation}
\begin{aligned}
\langle f_{\mathrm{new}}^a-f^a\rangle_{\theta,Y_{i_0}}
 &=\left(\partial_r+\frac{c_a}{r}\right)\mathsf C^{ra}
                              +\partial_z\mathsf C^{za},
                                      &&a=\theta,z,\\
\langle f_{\mathrm{new}}^r-f^r\rangle_{\theta,Y_{i_0}}
 &=\left(\partial_r+\frac1r\right)\mathsf C^{rr}
                  +\partial_z\mathsf C^{zr}-\frac{\mathsf C^{\theta\theta}}r.
\end{aligned}
\label{iter:auxiliary-double-mean}
\end{equation}
The tangential changes, in force units, are bounded by
\[
k^{C_j}2^{-(2.1+0.1j-3\cdot10^{-5})\delta k}
                           \omega d_{\mathrm{edge}}^{-C_j}.
\]
This follows by adding the old decay exponent \(1-10^{-5}\) to the new one
\(1.1+0.1j-2\cdot10^{-5}\) in \eqref{iter:auxiliary-size}.
The radial gain offsets the radial derivative. Thus the tangential double-mean bound
from \eqref{iter:stress-output} is preserved.

For the nonzero modes the only new terms are
\[
\sum_{w'}\bigl((w'\cdot\nabla)v+(v\cdot\nabla)w'\bigr),
\]
where \(w'\) runs over all pulses already present, including the
two pulses added in the current stage. The calculation in
\eqref{init:auxiliary-nonzero-modes}, with
\eqref{iter:auxiliary-size} in place of
\eqref{init:auxiliary-physical-bounds}, bounds this sum in force
units by
\[
k^{C_j}2^{-(0.1+0.1j-2\cdot10^{-5})\delta k}
                         \sqrt\omega\,d_{\mathrm{edge}}^{-C_j}.
\]
The worst term differentiates the carrier of the primary pulse.
Every term retains its old labelled support and local envelope,
because it contains a pulse coefficient. This proves the first
line of \eqref{iter:auxiliary-output}, including its local
coefficient version.

In the radial angular mean, the largest new term is
\(-2u_B^\theta v^\theta/r\), whose force-unit bound is
\(2^{\delta k}\) times the velocity-unit bound in
\eqref{iter:auxiliary-size}. The radial fast time derivative is
smaller because of the extra radial gain. The remaining terms
are bounded as in \eqref{init:auxiliary-radial-change}, with
\eqref{iter:old-axisymmetric-bounds} for the old increments.
This proves the last line of \eqref{iter:auxiliary-output}.

Finally, the double-mean radial change in
\eqref{iter:auxiliary-double-mean} contains no background cross
term. Its largest contribution is
\(-\mathsf C^{\theta\theta}/r\). After radial integration and
normalization it has exponent \(2.1+0.1j-3\cdot10^{-5}\).
The changes in the two flux defects have exactly the same bound:
use \(\mathsf C^{\theta z}\), \(\mathsf C^{zz}\), and the weighted
radial integral in the definition of \(\mathcal J^{z\theta}\).
The pressure is unchanged. Consequently
\eqref{iter:stress-defects} is preserved. The same estimates hold
after applying the specified derivatives of any fixed order.
\end{proof}

\subsection{The moment correction}
\label{iter:moment-correction}

We now use the state after Section~\ref{iter:auxiliary-correction}.
As in Section~\ref{init:moment-correction}, this step reduces the
three physical defects before reconstructing the pressure. The
pressure defect controls the radial force left by that
reconstruction, while the two flux defects control the compensating
terms in the next stress correction.

\emph{The velocity.}
Define \(v=u_{j,\mathrm{moment}}\) by the five conditions
\eqref{init:moment-five-equations}, with the full current defects
\(\mathcal P,\mathcal J^{\theta z},\mathcal J^{z\theta}\) on the
right-hand sides. Use exactly the profiles and matrices in
\eqref{init:moment-profiles} and
\eqref{init:moment-coefficient-equations}, replacing the index
\(0\) by \(j\). Thus
\begin{equation}
\Psi_{j,\mathrm{moment}}
 =\frac1r\int_0^r sB_{j,\mathrm{moment}}(t,s,z)\,ds,
\qquad
v=A_{j,\mathrm{moment}}e_\theta
                  +\nabla\times(\Psi_{j,\mathrm{moment}}e_\theta).
\label{iter:moment-velocity}
\end{equation}
The matrices depend only on the fixed background patch and
profiles. Their invertibility is unchanged. The radial moment condition
\[
\int_0^\infty sB_{j,\mathrm{moment}}(t,s,z)\,ds=0
\]
makes the potential compactly supported in the same fixed interior
radial hull. The field is independent of the
auxiliary variable, is divergence-free, and has zero linear
moments. By \eqref{iter:stress-defects} and its preservation in
Section~\ref{iter:auxiliary-correction},
\begin{equation}
\begin{aligned}
&2^{-(1+\delta)k}\bigl(
 |\partial_{\hat R,\hat Z,\hat T}^{\beta}v^{\mathrm{tan}}|
       +2^{\delta k}|\partial_{\hat R,\hat Z,\hat T}^{\beta}v^r|\bigr)
 +2^{-\delta k}
       |\partial_{\hat R,\hat Z,\hat T}^{\beta}\Psi_{j,\mathrm{moment}}|\\
&\qquad\lesssim k^{C_{j,\beta}}
       2^{-(1.1+0.1j-2\cdot10^{-5})\delta k}
                         \omega d_{\mathrm{edge}}^{-C_{j,\beta}}.
\end{aligned}
\label{iter:moment-size}
\end{equation}
Here the radial weight is available on the fixed interior
support, as in \eqref{init:moment-physical-bounds}.
Keep the pressure and flat error unchanged and set
\[
u_*=u+v,\qquad p_*=p,\qquad
f_*=\mathsf{NS}(u_*,p)-E.
\]

\emph{The new defects.}
At this stage the double-mean velocity is
\[
\langle u\rangle_{\theta,Y_{i_0}}
       =u_B+\sum_{q=0}^{j-1}u_{q,\mathrm{moment}}.
\]
Thus the calculation in \eqref{init:moment-radial-change} gives
\begin{equation}
\langle f_*^r-f^r\rangle_{\theta,Y_{i_0}}
       =-\frac{2u_B^\theta v^\theta}{r}+\mathsf Q_j^r,
\label{iter:moment-radial-change}
\end{equation}
where
\begin{equation}
\begin{aligned}
\mathsf Q_j^r={}&\partial_tv^r
 +\left(\partial_r+\frac1r\right)(2u_B^rv^r+(v^r)^2)
 +\partial_z(u_B^rv^z+u_B^zv^r+v^rv^z)\\
&-\frac{(v^\theta)^2}{r}
 -\left(\partial_r^2+\frac1r\partial_r
                         +\partial_z^2-\frac1{r^2}\right)v^r\\
&+\sum_{q=0}^{j-1}\left[
  \left(\partial_r+\frac1r\right)(2u_{q,\mathrm{moment}}^rv^r)
  +\partial_z(u_{q,\mathrm{moment}}^rv^z
                         +u_{q,\mathrm{moment}}^zv^r)
  -\frac{2u_{q,\mathrm{moment}}^\theta v^\theta}{r}\right].
\end{aligned}
\label{iter:moment-radial-remainder}
\end{equation}
The terms before the sum over \(q\) are the remainder
\eqref{init:moment-radial-remainder}; the sum records the interactions
with the old moment corrections. All auxiliary corrections have
zero cross average with \(v\).

Use the five moment conditions and \(u_B^z=0\) on the support of
\(v\). The exact counterpart of
\eqref{init:moment-exact-new-defects} is
\begin{equation}
\begin{aligned}
\mathcal P_*&=\int_0^\infty\mathsf Q_j^r\,dr,\\
\mathcal J_*^{\theta z}
 &=\int_0^\infty r^2\left[v^\theta v^z
   +\sum_{q=0}^{j-1}(u_{q,\mathrm{moment}}^\theta v^z
                          +u_{q,\mathrm{moment}}^zv^\theta)\right]dr,\\
\mathcal J_*^{z\theta}
 &=\int_0^\infty r\left[(v^z)^2
                  +2\sum_{q=0}^{j-1}u_{q,\mathrm{moment}}^zv^z\right]dr
                    +\frac12\int_0^\infty r^2\mathsf Q_j^r\,dr.
\end{aligned}
\label{iter:moment-new-defects}
\end{equation}
In particular, the linear equations still cancel the complete
incoming defects; the old moment cross terms remain in the new
defects.

For completeness, all terms in \eqref{iter:moment-radial-remainder}
can be bounded using slow derivatives. The time and viscous terms,
and the terms involving the background meridional velocity, have
force-unit bound
\[
k^{C_j}2^{-(2.1+0.1j-2\cdot10^{-5})\delta k}
                              \omega d_{\mathrm{edge}}^{-C_j}.
\]
The meridional products with old moment corrections are smaller.
The largest possible term is the old swirl cross term. From
\eqref{iter:old-axisymmetric-bounds} and \eqref{iter:moment-size},
\begin{equation}
2^{-(3+\delta)k}
 |\partial_{\hat R,\hat Z,\hat T}^{\beta}\mathsf Q_j^r|
 \lesssim k^{C_{j,\beta}}
       2^{-(1.1+0.1j-3\cdot10^{-5})\delta k}
                          \omega d_{\mathrm{edge}}^{-C_{j,\beta}}.
\label{iter:moment-remainder-bound}
\end{equation}
Indeed, the centrifugal factor gives \(2^{\delta k}\) in force
units, while the two velocities have exponents
\(1-10^{-5}\) and \(1.1+0.1j-2\cdot10^{-5}\).
The new swirl square has a strictly better exponent for \(j\ge1\).
Radial integration supplies one more \(2^{-\delta k}\) after
normalization of the defects. The other products in
\eqref{iter:moment-new-defects} obey the same resulting bound.
Consequently,
\begin{equation}
\begin{aligned}
&2^{-2(1+\delta)k}|\partial_{\hat T,\hat Z}^{\beta}\mathcal P_*|
 +2^{(1-2\delta)k}|\partial_{\hat T,\hat Z}^{\beta}\mathcal J_*^{\theta z}|
 +2^{-2\delta k}|\partial_{\hat T,\hat Z}^{\beta}\mathcal J_*^{z\theta}|\\
&\qquad\lesssim k^{C_{j,\beta}}
                   2^{-(2.1+0.1j-3\cdot10^{-5})\delta k}.
\end{aligned}
\label{iter:moment-improved-defects}
\end{equation}
This is stronger than the defect requirement in \eqref{iter:targets}.

\emph{The force before pressure reconstruction.}
The bounds \eqref{iter:auxiliary-output} remain valid with
\(f_{\mathrm{new}}\) replaced by \(f_*\). To see this, use
\eqref{init:moment-force-update}, including the interactions with
all old axisymmetric increments and pulses. The new nonzero
angular terms are the interactions of \(v\) with the pulses.
The estimate in \eqref{init:auxiliary-nonzero-modes} bounds them,
in force units, by
\[
k^{C_j}2^{-(0.1+0.1j-2\cdot10^{-5})\delta k}
                          \sqrt\omega\,d_{\mathrm{edge}}^{-C_j}.
\]
Each term retains the labelled support and envelope of its pulse
factor. The tangential double-mean change is bounded by
\eqref{init:moment-tangential-force} and
\eqref{iter:moment-size}, with exponent
\(1.1+0.1j-2\cdot10^{-5}\). The additional old moment products
are smaller. The only new auxiliary oscillations in the angular
mean come from products with the auxiliary corrections already present.
Let \(U=u_{q,\mathrm{aux}}\), \(0\le q\le j\). For \(a\in\{\theta,z\}\),
the radial flux \(U^rv^a+v^rU^a\) contains one radial velocity component
in each product. For example, \eqref{iter:old-axisymmetric-bounds} and
\eqref{iter:moment-size} give, using the stronger bound
\eqref{iter:auxiliary-size} for \(q=j\),
\[
\begin{aligned}
 &2^{-(3+\delta)k}
   \left|\left(D_r+\frac{c_a}{r}\right)(U^rv^a)\right|\\
 &\quad\lesssim k^{C_j}\omega d_{\mathrm{edge}}^{-C_j}\,
     2^{\delta k}2^{10^{-5}\delta k}
     2^{-(2-10^{-5})\delta k}2^{-(1.1+0.1j-2\cdot10^{-5})\delta k}\\
 &\quad= k^{C_j}\omega d_{\mathrm{edge}}^{-C_j}\,
     2^{-(2.1+0.1j-4\cdot10^{-5})\delta k}.
\end{aligned}
\]
The same bound holds for \(v^rU^a\), with the radial gain now supplied
by \(v^r\), while the axial derivative incurs no auxiliary radial loss.
These estimates preserve the
second line of \eqref{iter:auxiliary-output}.

Finally, before auxiliary averaging the main new radial term is
\(-2u_B^\theta v^\theta/r\), of force-unit exponent
\(0.1+0.1j-2\cdot10^{-5}\). All the remaining radial terms,
including the old auxiliary cross terms, are smaller. Thus
\begin{equation}
2^{-(3+\delta)k}|\partial_{\hat R,\hat Z,\hat T,Y_{i_0,1},Y_{i_0,2}}^{\beta}\langle f_*^r\rangle_\theta|
 \lesssim k^{C_{j,\beta}}
       2^{-(0.1+0.1j-2\cdot10^{-5})\delta k}
                         \omega d_{\mathrm{edge}}^{-C_{j,\beta}}.
\label{iter:pressure-source}
\end{equation}
The support and linear moments are preserved, and the flat error
has not changed. This completes the velocity corrections at
stage \(j\).

\subsection{Pressure reconstruction and completion of the stage}
\label{iter:pressure-reconstruction}

Retain the starred state from Section~\ref{iter:moment-correction}.
It already contains the old axisymmetric pressure
\(\widetilde p^{[j-1]}\). Define the pressure increment by
\begin{equation}
q_j=\mathcal I_0\!\left(-\langle f_*^r\rangle_\theta
                                      +\rho_r\mathcal P_*\right),
\qquad
\widetilde p^{[j]}=\widetilde p^{[j-1]}+q_j.
\label{iter:pressure-increment}
\end{equation}
Here \(\rho_r\) and \(\mathcal I_0\) are the fixed profiles and
operator used in Section~\ref{init:pressure-reconstruction}.
In particular, this is a single pressure reconstruction after the
four velocity corrections.

The radial identity \eqref{init:pressure-radial-identity} becomes
\begin{equation}
D_rq_j=-\langle f_*^r\rangle_\theta
                       +\rho_r\mathcal P_*+\mathcal R_{j,\mathrm p}.
\label{iter:pressure-radial-identity}
\end{equation}
The remainder has the formula \eqref{init:pressure-cutoff-remainder}
with the present \(f_*\). It has zero auxiliary mean and is flat
with every physical derivative. Indeed, the proof of
\eqref{init:pressure-flat-bound} applies using
\eqref{iter:pressure-source}; the number of integrations by parts
is chosen after the desired derivative order and decay power.
The increment \(q_j\) has smooth zero extension off the annulus,
by the same argument using the integral from \(0\) to \(r\) and
the tail integral as in Section~\ref{init:pressure-reconstruction}.

The normalized bump term is smaller than the source in
\eqref{iter:pressure-source}, by
\eqref{iter:moment-improved-defects}. Applying the fixed-shift
estimate from the proof of \eqref{init:pressure-size} gives
\begin{equation}
\begin{aligned}
&2^{-2(1+\delta)k}|\partial_{\hat R,\hat Z,\hat T,Y_{i_0,1},Y_{i_0,2}}^{\beta}q_j|
 +2^{-(3+\delta)k}|\partial_{\hat R,\hat Z,\hat T,Y_{i_0,1},Y_{i_0,2}}^{\beta}\partial_zq_j|\\
 &\qquad\lesssim k^{C_{j,\beta}}
        2^{-(1.1+0.1j-2\cdot10^{-5})\delta k}
                         \omega d_{\mathrm{edge}}^{-C_{j,\beta}}.
\end{aligned}
\label{iter:pressure-size}
\end{equation}
This estimate is for \(q_j\), rather than the accumulated
axisymmetric pressure. The radial integral gives \(2^{-k}\),
which becomes an additional \(2^{-\delta k}\) in pressure units;
the axial derivative then converts pressure units to force units.

Define the completed stage by
\begin{equation}
\begin{aligned}
u^{[j]}&=u_*,&p^{[j]}&=p_*+q_j,\\
E^{[j]}&=E+\mathcal R_{j,\mathrm p}e_r,&
f^{[j]}&=f_*+\nabla q_j-\mathcal R_{j,\mathrm p}e_r.
\end{aligned}
\label{iter:completed-stage}
\end{equation}
These definitions agree with \eqref{init:later-stage-fields}.
They give \(\mathsf{NS}(u^{[j]},p^{[j]})=E^{[j]}+f^{[j]}\)
and, exactly,
\begin{equation}
\langle(f^{[j]})^r\rangle_\theta=\rho_r\mathcal P_*,
\qquad
\mathcal P^{[j]}=\mathcal P_*,\qquad
\mathcal J^{a,[j]}=\mathcal J_*^a\quad(a=\theta z,z\theta).
\label{iter:completed-defects}
\end{equation}
For the pressure defect, integrate the first identity and use
\(\int\rho_r\,dr=1\). The \(\theta z\) defect contains no pressure.
For the \(z\theta\) defect, the change is the integral of
\[
r\langle q_j\rangle_{Y_{i_0}}
       +\frac{r^2}{2}\partial_r\langle q_j\rangle_{Y_{i_0}}
 =\partial_r\left(\frac{r^2}{2}\langle q_j\rangle_{Y_{i_0}}\right),
\]
whose boundary values vanish. This is the pressure invariance
proved in Section~\ref{init:pressure-reconstruction}.

The nonzero angular force and the angular component are unchanged.
The axial angular mean changes by \(\partial_zq_j\), so
\eqref{iter:pressure-size} preserves the tangential estimates in
\eqref{iter:auxiliary-output}. For the radial mean,
\eqref{iter:moment-improved-defects} and the weighted bound for
\(2^{-k}\rho_r\) give
\begin{equation}
2^{-(3+\delta)k}|\partial_{\hat R,\hat Z,\hat T,Y_{i_0,1},Y_{i_0,2}}^{\beta}
                         \langle(f^{[j]})^r\rangle_\theta|
 \lesssim k^{C_{j,\beta}}
        2^{-(1.1+0.1j-3\cdot10^{-5})\delta k}
                         \omega d_{\mathrm{edge}}^{-C_{j,\beta}}.
\label{iter:completed-radial-bound}
\end{equation}
It is independent of the auxiliary variable and is stronger than
the required double-mean bound. Combining the preceding estimates
yields
\begin{equation}
\begin{aligned}
2^{-(3+\delta)k}
 \bigl(|f^{[j]}|+|f^{[j]}-\langle f^{[j]}\rangle_\theta|\bigr)
 &\lesssim k^{C_j}2^{-0.1j\delta k}
                         \sqrt\omega\,d_{\mathrm{edge}}^{-C_j},\\
2^{-(3+\delta)k}\Bigl(
 |\langle(f^{[j]})^{\mathrm{tan}}\rangle_\theta
      -\langle(f^{[j]})^{\mathrm{tan}}\rangle_{\theta,Y_{i_0}}|
 +|\langle f^{[j]}\rangle_{\theta,Y_{i_0}}|\Bigr)
 &\lesssim k^{C_j}2^{-(0.2+0.1j)\delta k}
                         \omega d_{\mathrm{edge}}^{-C_j}.
\end{aligned}
\label{iter:completed-force}
\end{equation}
The first bound is stronger than \eqref{iter:targets}. The defect
bound follows from \eqref{iter:moment-improved-defects} and
\eqref{iter:completed-defects}. All assertions also hold after applying the derivatives
specified at the start of this section.

The fields \(u^{[j]}-u_B\), \(p^{[j]}-p_B\),
\(f^{[j]}\), and \(E^{[j]}-\mathsf E_B\) remain smooth,
with radial support in \(R\in[R_a,R_h]\).
The velocity \(u^{[j]}\) remains divergence-free, and both linear
moments are preserved. Every added flat term has zero double mean; hence
\(\langle E^{[j]}\rangle_{\theta,Y_{i_0}}=\mathsf E_B\).
Each velocity harmonic of the two new pulses is a curl, hence
divergence-free. The pulse and nonzero-mode force coefficients
retain their assigned labelled supports, smooth zero extensions
and Gaussian bounds, as verified at each velocity step.
All new nonzero harmonics therefore satisfy \(0<|m|\le 2^{j+1}\).
Write the nonzero angular part of \(f^{[j]}\) as in
\eqref{iter:source-harmonics}, with \(j+1\) in place of \(j\).
The coefficient estimates from the velocity corrections, together with the
fact that the pressure correction leaves the nonzero angular force unchanged,
give, for each label \(\gamma\), each \(m\in\mathcal H_{j+1}\), and every
fixed \(\beta\in\mathbb Z_{\ge0}^5\),
\[
 \left|\partial_{\hat R,\hat Z,\hat T,Y_{i_0,1},Y_{i_0,2}}^\beta
                  \mathbf f_{\gamma,j+1,m}\right|
 \lesssim k^{C_{j,\beta}}2^{-0.1j\delta k}
          \sqrt\omega\,d_{\mathrm{edge}}^{-C_{j,\beta}}
          \mathcal G_\gamma(\tau_\gamma).
\]
Since \(-0.1j\le 0.4-0.1(j+1)\), this gives \eqref{iter:source-bound}
at stage \(j+1\). We have therefore established the complete
induction, including the increment bounds
\eqref{init:later-increment-sizes}.

\subsection{Smooth realization of the iteration}
\label{iter:smooth-realization}

We construct \(u^\infty,p^\infty\) so that on each successive
range of \(\lambda\) they agree with a completed stage plus a
single cutoff increment. To prove that the resulting force is flat,
it will suffice to estimate the effect of this one increment on the
equation and then choose the scale at which it is applied.

All finite stages are defined on the same domain
\(0<\lambda\le 2^{-k_0}\), \(|Z|\le1\), with \(k_0\) fixed before
the iteration. Indeed, the forced equations with the fixed background
coefficients are invertible for each nonzero harmonic, and each stage
uses only finitely many harmonics. The stress inverse uses the fixed
positive primary amplitudes and accepts arbitrary signed targets.
The torus inverse, radial primitive and moment matrices are fixed
as well. Thus increasing the stage changes the constants in the
estimates without imposing a new restriction on \(\lambda\).

\emph{The fields between consecutive radii.}
After physical evaluation, the velocity increment has the form
\begin{equation}
 u^{[j]}-u^{[j-1]}=\nabla\times A_j+b_j e_\theta,\qquad
 \partial_\theta b_j=0.
\label{iter:increment-potentials}
\end{equation}
Here \(A_j\) is the sum of the two pulse potentials and
\((\Psi_{j,\mathrm{aux}}+\Psi_{j,\mathrm{moment}})e_\theta\), and
\(b_j=A_{j,\mathrm{aux}}+A_{j,\mathrm{moment}}\).
These potentials, \(b_j e_\theta\) and the pressure increment
\(p^{[j]}-p^{[j-1]}=p_{j,m\ne0}+p_{j,\mathrm{stress}}+q_j\)
extend smoothly by zero off \(R\in[R_a,R_h]\). For the moment
potential this uses the radial moment condition after
\eqref{iter:moment-velocity}.

Fix a smooth cutoff \(\chi\) on \([0,\infty)\), equal to one on
\([0,1/2]\) and zero on \([1,\infty)\). For now let
\(0<\rho_1\le 2^{-k_0}\) and \(0<\rho_{j+1}\le\rho_j/2\).
Put \(\chi_j(t,z)=\chi(\lambda(t,z)/\rho_j)\) and set
\begin{equation}
\begin{aligned}
 u^\infty&=u^{[0]}+\sum_{j\ge1}
        \bigl[\nabla\times(\chi_j A_j)+\chi_j b_j e_\theta\bigr],\\
 p^\infty&=p^{[0]}+\sum_{j\ge1}\chi_j(p^{[j]}-p^{[j-1]}).
\end{aligned}
\label{iter:realized-fields}
\end{equation}
Both sums are locally finite for \(\lambda>0\).
The velocity is divergence-free since the cutoff is applied to
the potential before taking the curl and \(\chi_j b_j\) is
independent of \(\theta\). The corrections to \(u_B,p_B\)
retain their smooth zero extensions off the annulus.

On \(\rho_{j+1}\le\lambda\le\rho_j\), the cutoffs with index
less than \(j\) are all one, and those with index greater than
\(j\) are all zero. The first \(j-1\) increments telescope, giving
\begin{equation}
\begin{aligned}
 u^\infty&=u^{[j-1]}+\nabla\times(\chi_j A_j)+\chi_j b_j e_\theta,\\
 p^\infty&=p^{[j-1]}+\chi_j(p^{[j]}-p^{[j-1]}).
\end{aligned}
\label{iter:fields-between-radii}
\end{equation}

\emph{The force between consecutive radii.}
We first record the physical derivative
estimates needed in \eqref{iter:fields-between-radii}. Throughout,
\(B_n\) denotes a nonnegative loss depending only on the derivative
order \(n\) and the fixed construction. It may increase from one
estimate to the next. The constants \(C_{j,n}\) may
also depend on the stage. For \(j\ge1\),
\begin{equation}
 \max_{|\nu|\le n}\Bigl(
 |D_{t,x}^\nu A_j|+|D_{t,x}^\nu(b_j e_\theta)|
 +|D_{t,x}^\nu(p^{[j]}-p^{[j-1]})|\Bigr)
 \le C_{j,n}\lambda^{0.1\delta j-B_n}.
\label{iter:physical-increments}
\end{equation}
Indeed, \(2^{-k}\asymp\lambda\), and the increment estimates
\eqref{init:later-increment-sizes} give the gain \(2^{-0.1j\delta k}\)
up to a fixed power of \(2^k\). The pulse potentials are given by
their explicit formulas, and the other potentials obey
\eqref{iter:auxiliary-size} and \eqref{iter:moment-size}.
At fixed derivative order the powers of \(2^k\) from the physical
chain rule, including the frequencies and the cylindrical frame,
are independent of \(j\). Harmonic factors involving \(2^j\)
enter \(C_{j,n}\). The edge weight bounds any fixed power of
\(d_{\mathrm{edge}}^{-1}\) occurring with \(\sqrt\omega\) or
\(\omega\), and \(k^{C_{j,n}}\) can be absorbed into one power
of \(\lambda^{-1}\) with a stage-dependent constant.

For the completed stage in \eqref{iter:fields-between-radii},
the corresponding bounds are
\begin{equation}
\begin{aligned}
 \max_{|\nu|\le n}|D_{t,x}^\nu u^{[j-1]}|
   &\le C_{j,n}\lambda^{-B_n}
      \qquad (R_a\le R\le R_h),\\
 \max_{|\nu|\le n}|D_{t,x}^\nu\mathsf{NS}(u^{[j-1]},p^{[j-1]})|
   &\le C_{j,n}\lambda^{0.1\delta(j-1)-B_n}.
\end{aligned}
\label{iter:physical-stages}
\end{equation}
The velocity bound follows from the background estimate
\eqref{intro:background-velocity-bounds} and the finite sum of
velocity corrections. For the force, apply the same physical
chain rule to \eqref{iter:completed-force} and use the flatness
of \(E^{[j-1]}\). At \(j=1\), use the initialized force bound
\eqref{init:moment-initialized-bounds}, whose factor
\(2^{0.01\delta k}\) costs only a fixed power of \(\lambda^{-1}\).
Outside the annulus the residual is the flat background error in
\eqref{intro:flat-background-remainder}, and is zero for \(R\ge R_h\).

To compare the two fields on this band, write \(v,q\) for the
velocity and pressure increments on the right of
\eqref{iter:fields-between-radii}. Their contribution to the
equation is
\begin{equation}
\begin{aligned}
 \mathsf{NS}(u^\infty,p^\infty)-\mathsf{NS}(u^{[j-1]},p^{[j-1]})
 ={}&\partial_t v-\Delta v+\nabla q+(u^{[j-1]}\cdot\nabla)v\\
    &+(v\cdot\nabla)u^{[j-1]}+(v\cdot\nabla)v .
\end{aligned}
\label{iter:band-force-difference}
\end{equation}
The constants in estimating this difference must be independent
of the radii. By \eqref{intro:background-variables},
\(\lambda^2-z^2\lambda^{2\delta}=1-t\). Implicit differentiation,
using \(1-\delta Z^2\ge1-\delta\), gives
\(|D_{t,x}^\nu\lambda|\lesssim_\nu\lambda^{1-2|\nu|}\).
On the support of a derivative of \(\chi_j\) we have
\(1/2\le\lambda/\rho_j\le1\), so the chain rule gives
\(|D_{t,x}^\nu\chi_j|\lesssim_\nu\lambda^{-2|\nu|}\)
uniformly in \(\rho_j\).

After \(n\) derivatives, the linear terms in
\eqref{iter:band-force-difference} use \(v\) through order \(n+2\)
and \(q\) through order \(n+1\). The curl in \(v\) thus uses the
cutoff and \(A_j\) through order \(n+3\). By the product rule and
\eqref{iter:physical-increments} each linear term retains a factor
\(\lambda^{0.1\delta j}\) up to a loss depending only on \(n\).
The mixed terms use \(u^{[j-1]}\) through order \(n+1\) only on the
annulus supporting \(v\), where \eqref{iter:physical-stages} gives
the same type of loss. The last term contains two factors of the
increment and hence gains \(\lambda^{0.2\delta j}\) before these losses.
Adding the residual of the previous stage gives
\begin{equation}
 \max_{|\nu|\le n}|D_{t,x}^\nu\mathsf{NS}(u^\infty,p^\infty)|
 \le C_{j,n}\lambda^{0.1\delta(j-1)-B_n}
 \qquad(\rho_{j+1}\le\lambda\le\rho_j).
\label{iter:band-force}
\end{equation}
Here \(B_n\) is independent of \(j\), and \(C_{j,n}\) is independent
of every radius. In particular these constants can be fixed
before choosing \(\rho_j\).

\emph{Choice of radii and flatness.}
Now choose \(\rho_1,\rho_2,\ldots\) successively. In addition to
\(\rho_1\le2^{-k_0}\) and \(\rho_{j+1}\le\rho_j/2\), require
\begin{equation}
 \rho_j\max_{0\le n\le j} C_{j,n}\le1.
\label{iter:cutoff-choice}
\end{equation}
At stage \(j\) there are only finitely many constants to absorb.
Since \(\lambda\le\rho_j\) on this band, \(C_{j,n}\lambda\le1\)
for \(n\le j\). Absorbing the constant at the cost of one power of \(\lambda\) in
\eqref{iter:band-force} therefore bounds its left-hand side by
\(\lambda^{0.1\delta(j-1)-B_n-1}\).

With the radii fixed, prescribe a derivative order \(n\) and a
power \(N>0\). Choose \(j_*\ge n\), \(j_*\ge1\), such that
\(0.1\delta(j_*-1)-B_n-1\ge N\). Every \(0<\lambda\le\rho_{j_*}\)
belongs to a band with \(j\ge j_*\), so the bound there is at
most \(\lambda^N\). On the rest of a smaller fixed working
neighborhood, \(\lambda\) is bounded below and only finitely many
stages and dyadic scales occur. Increasing the constant on this
part gives
\begin{equation}
 f^\infty=\mathsf{NS}(u^\infty,p^\infty),\qquad
 |D_{t,x}^\nu f^\infty|\le C_{\nu,N}\lambda^N
 \qquad(N>0).
\label{iter:realized-flat-force}
\end{equation}

The fields have smooth one-sided extensions to \(t=1\) away
from the origin, as verified in Section~\ref{iter:localization}.
Near \((t,x)=(1,0)\), the coordinate relation gives
\[
\lambda\lesssim (1-t)^{1/2}+|z|^{1/(1-\delta)}.
\]
Thus \eqref{iter:realized-flat-force}, for every physical
derivative and every power, extends \(f^\infty\) smoothly to
that point with zero jet. This proves the required local smooth
forcing statement. Since all corrections vanish for \(R>R_h\),
the realized velocity still agrees with the heat exterior in
\eqref{intro:background-heat}.

\subsection{Localization and the blowup conclusion}
\label{iter:localization}

We finish the proof of Theorem~\ref{thm:part-ii-main} by localizing
the velocity and pressure in physical space and time.

\emph{Regularity up to \(t=1\) away from the origin.}
Fix \(x_0\ne0\) in our spatial neighborhood of the origin, with
cylindrical radius \(r_0\) and height \(z_0\). The scale \(\lambda(t,z)\)
depends only on \(t,z\), so even at positive radius it tends to zero
as \((t,z)\to(1,0)\).

If \(z_0\ne0\), then \(|Z|\le1\) in \eqref{intro:background-variables}
gives \(\lambda\ge |z|^{1/(1-\delta)}\), a positive lower bound in
a neighborhood of the point. To see that the change of variables
is smooth there up to \(t=1\), put \(s=\lambda^2\). Its defining
equation and the derivative with respect to \(s\) are
\[
 s-z^2s^\delta=1-t,\qquad
 \partial_s(s-z^2s^\delta)=1-\delta Z^2\ge1-\delta.
\]
Thus \(s\), and hence \(\lambda,R,Z\), are smooth functions of the
physical variables there, including at \(Z=\pm1\).
Only finitely many bands \(k\) meet this neighborhood, and only
finitely many stage cutoffs \(\chi_j\) in \eqref{iter:realized-fields}
are nonzero there. The background is smooth by
Section~\ref{intro:flux-cone}. For the corrections, the stage-zero
estimates in Section~\ref{sec:initialization-ppt} and the physical
derivative bounds \eqref{iter:physical-increments} give locally
uniform control of every space--time derivative for \(t<1\). Applying
the same bounds with one more time derivative gives locally uniform
one-sided limits of all these derivatives as \(t\uparrow1\).
This proves smoothness of \(u^\infty,p^\infty\) up to \(t=1\).

If \(z_0=0\), then \(r_0>0\). The bound
\(\lambda\lesssim(1-t)^{1/2}+|z|^{1/(1-\delta)}\) from
Section~\ref{iter:smooth-realization} gives a neighborhood in which
\(r>r_0/2\) and \(R=r^2/(2\lambda^2)>R_h\) for \(t<1\).
Throughout this neighborhood all annular corrections vanish.
Thus \(u^\infty,p^\infty\) agree with the heat exterior
\eqref{intro:background-heat}, which is smooth up to \(t=1\)
at positive radius.

\emph{A potential on the cutoff shell.}
Following the potential cutoff method of \cite[Section~10.1]{OpenAI},
choose \(0<a<b\) small and \(t_0\in(0,1)\) close enough to one that
\(\overline B_b(0)\times[t_0,1)\) lies in the working region.
We need a vector potential only on the shell \(a<|x|<b\). For
\(a<\rho<b\) and \(t_0\le t<1\), the divergence theorem in \(B_\rho(0)\)
gives
\[
 \int_{|x|=\rho}u^\infty(t,x)\cdot n\,dS=0.
\]
This also holds at \(t=1\) by smooth convergence in time on the sphere.
Fix one intermediate sphere \(S\). The zero flux allows us to solve
\(\Delta_S\phi=u^\infty\cdot n\) on \(S\) with zero mean. Setting
\(A=n\times\nabla_S\phi\) there gives the prescribed normal curl.
In the gauge \(A\cdot n=0\), the tangential curl equations determine
\(A\) on the shell by radial integration. The difference
\(\nabla\times A-u^\infty\) is then radial and divergence-free,
so its vanishing on \(S\) implies its vanishing throughout the shell.
The fixed spherical inverse and radial integration commute with time
derivatives, hence \(A\) is smooth up to \(t=1\).

Choose a smooth spatial cutoff \(h\) equal to one on \(B_a(0)\)
and zero outside \(B_b(0)\), with \(\nabla h\) compactly supported
in the open shell. Choose \(0<\varepsilon<a\) with
\(t_0<1-\varepsilon^2\), and a smooth time cutoff \(\tau\) equal
to zero for \(t\le t_0\) and to one for \(t\ge1-\varepsilon^2\).
Set
\begin{equation}
 u=\tau\bigl(h u^\infty+\nabla h\times A\bigr),\qquad
 p=\tau h p^\infty,\qquad f=\mathsf{NS}(u,p),
\label{iter:localized-fields}
\end{equation}
with the cutoff terms extended by zero off their supports.
On the shell the expression in parentheses is
\(\nabla\times(hA)\). Inside \(B_a(0)\) it is \(u^\infty\).
Thus \(u\) is divergence-free, and \(u,p\) vanish for \(t\le t_0\)
and have spatial support in the fixed compact set \(\overline B_b(0)\).

\emph{The force, energy and blowup.}
On \(Q_\varepsilon\) the cutoffs equal one, so \(u,p\) retain the
annular correction support of \eqref{iter:realized-fields}, and
\(f=f^\infty\) satisfies \eqref{iter:realized-flat-force}.
In particular all derivatives of \(f\) extend with value zero at
\((1,0)\). At every other point of the terminal slice, smoothness of
the localized fields gives smoothness of \(f\) up to \(t=1\).
Extend \(f\) smoothly across \(t=1\) and multiply the extension
by spatial and later-time cutoffs equal to one on
\(\overline B_b(0)\times[t_0,1]\). This gives
\(f\in C_c^\infty(\mathbb R^3\times(0,\infty);\mathbb R^3)\).

For each \(t<1\), spatial integration of the equation is
justified by smoothness and compact support. The energy identity
and Cauchy's inequality give
\[
\frac{d}{dt}\|u(t)\|_2^2+2\|\nabla u(t)\|_2^2
 \le\|u(t)\|_2^2+\|f(t)\|_2^2.
\]
Since \(u=0\) initially and the extended force is smooth and
compactly supported, Gronwall's inequality yields
\(\sup_{0\le t<1}\|u(t)\|_2<\infty\). Integrating the same
inequality gives
\[
2\int_0^T\|\nabla u(t)\|_2^2\,dt
 \le\int_0^T\bigl(\|u(t)\|_2^2+\|f(t)\|_2^2\bigr)\,dt
 \le C,\qquad T<1,
\]
with \(C\) independent of \(T\). Letting \(T\uparrow1\) proves
\(\nabla u\in L^2(\mathbb R^3\times[0,1])\).

Finally fix \(R_*>R_h\). For \(t\) sufficiently close to one,
the curve \(z=0,\ r=\sqrt{2R_*(1-t)}\) lies in
\(Q_\varepsilon\), where the cutoff equals one and all annular
corrections vanish on that curve. Consequently,
\eqref{intro:background-heat} gives
\[
\|u(t)\|_\infty
 \ge c_\infty R_*^{-(1+\delta)/2}H_\delta(2/R_*)
                                  (1-t)^{-(1+\delta)/2}.
\]
The coefficient is positive. This proves
\eqref{intro:main-blowup} and completes
Theorem~\ref{thm:part-ii-main}.

\section*{Acknowledgement}
ZL was supported in part by  NSFC (No.12431007), NSFC Excellent Research Group Project
(No.62588101), NSFC (No.12494544) and
New Cornerstone Science Foundation through the XPLORER PRIZE. XR was supported in part by NSFC (No.12601444) and the National Key R\&D Program of China (No. 2023YFA1010700).  We used AI tools to help prepare and polish the article.

\end{CJK*}
\end{document}